\documentclass[12pt]{amsart}
\usepackage{longtable}
\usepackage{amsfonts,amssymb,amscd,amsmath,epsfig}
\usepackage{latexsym}
\usepackage{mathrsfs}
\usepackage{tocvsec2}

\usepackage[]{youngtab}

\usepackage[
hyperindex=true,pagebackref=true,bookmarks=true,
colorlinks=true,linkcolor=blue,citecolor=red]{hyperref}

\let\lx\l

\begin{document}

\renewcommand{\tilde}{\widetilde}
\renewcommand{\hat}{\widehat}

\newcommand{\BR}{{\mathbb R}}
\newcommand{\BQ}{{\mathbb Q}}
\newcommand{\BC}{{\mathbb C}}
\newcommand{\BP}{{\mathbb P}}
\newcommand{\BZ}{{\mathbb Z}}
\newcommand{\BN}{{\mathbb N}}
\newcommand{\BS}{{\mathbb S}}

\newcommand{\cH}{{\mathcal H}}
\newcommand{\cA}{{\mathcal A}}
\newcommand{\cB}{{\mathcal B}}
\newcommand{\ccF}{{\mathfrak F}}
\newcommand{\cD}{{\mathcal D}}
\newcommand{\cL}{{\mathcal L}}
\newcommand{\cF}{{\mathcal F}}
\newcommand{\cP}{{\mathcal P}}
\newcommand{\cX}{{\mathcal X}}
\newcommand{\cY}{{\mathcal Y}}
\newcommand{\cS}{{\mathcal S}}
\newcommand{\cSol}{\hbox{$\mathcal Sol$}}
\newcommand{\cT}{\hbox{$\mathcal T$}}

\newcommand{\Z}{{\mathbb Z}}
\newcommand{\Q}{{\mathbb Q}}
\newcommand{\N}{{\mathbb N}}
\newcommand{\C}{{\mathbb C}}
\newcommand{\R}{{\mathbb R}}
\newcommand{\X}{{\mathbb X}}
\newcommand{\Y}{{\mathbb Y}}

\newcommand{\CH}{{\mathcal H}}
\newcommand{\CA}{{\mathcal A}}

\def\HH{\mbox{${\mathcal H}$\kern-5.2pt${\mathcal H}$}}

\newcommand{\binomial}[2]{\genfrac{(}{)}{0pt}{}{ #1 }{ #2 }}
\newcommand{\qbinomial}[2]{\genfrac{[}{]}{0pt}{}{ #1 }{ #2 }_q }
\newcommand{\qbinom}[3]{\genfrac{[}{]}{0pt}{}{ #1 }{ #2 }_{ #3 } }


\def\der{\partial}
\def\tensor{\otimes}
\def\gam{\gamma} \def\Gam{\Gamma}
\def\del{\delta} \def\Del{\Delta}
\def\kap{\kappa}
\def\lam{\lambda} \def\Lam{\Lambda}
\def\Comp{{\mathbb C}}
\def\sM{{\mathcal M}}

\newtheorem{theorem}{Theorem}[section]
\newtheorem{maintheorem}[theorem]{Main Theorem}
\newtheorem{proposition}[theorem]{Proposition}
\newtheorem{definition}[theorem]{Definition}
\newtheorem{lemma}[theorem]{Lemma}
\newtheorem{corollary}[theorem]{Corollary}
\newtheorem{notation}[theorem]{Notation}
\newtheorem{remark}[theorem]{Remark}
\newtheorem{example}[theorem]{Example}

\newtheorem{theorem }{Theorem}[section]
\newtheorem{maintheorem }[theorem]{Main Theorem}
\newtheorem{proposition }[theorem]{Proposition}
\newtheorem{definition }[theorem]{Definition}
\newtheorem{lemma }[theorem]{Lemma}
\newtheorem{corollary }[theorem]{Corollary}
\newtheorem{notation }[theorem]{Notation}
\newtheorem{remark }[theorem]{Remark}
\newtheorem{example }[theorem]{Example}

\newtheorem{ maintheorem }[theorem]{Main Theorem}
\newtheorem{ theorem}{Theorem}[section]
\newtheorem{ proposition}[theorem]{Proposition}
\newtheorem{ definition}[theorem]{Definition}
\newtheorem{ lemma}[theorem]{Lemma}
\newtheorem{ corollary}[theorem]{Corollary}
\newtheorem{ notation}[theorem]{Notation}
\newtheorem{ remark}[theorem]{Remark}
\newtheorem{ example}[theorem]{Example}

\newtheorem{thm}{Theorem}[section]
\newtheorem{prop}[thm]{Proposition}
\newtheorem{lem}[thm]{Lemma}
\newtheorem{cor}[thm]{Corollary}
\newtheorem{conj}[thm]{Conjecture}
\newtheorem{con}[thm]{Conjecture}
\newtheorem{dfn}[thm]{Definition}
\newtheorem{df}[thm]{Definition}
 \newcommand{\rem}{{\bf Comment.\ }}
 \newcommand{\rmk}{{\bf Comment.\ }}
 \newcommand{\exmp}{{\bf Example.\ }}
 \newcommand{\ex}{{\bf Example.\ }}
 \newcommand{\prob}{{\bf Problem.\ }}

\newtheorem{note}{Note} 
\renewcommand{\thenote}{}
\newtheorem*{acka}{Acknowledgments}
\newtheorem{ack}{Acknowledgments}
\renewcommand{\theack}{}
\renewcommand{\appendixname}{\bf Appendix}
\renewcommand{\proof}{{\em Proof.\ }}

\hyphenation{
ap-pen-dix as-ymp-tot-ic at-trib-uted at-trib-ut-able
Bry-li-n-sky com-mu-ta-tion de-ge-ne-rate
de-riv-a-tive dis-trib-ute equi-vari-ant ex-tra-or-di-nary  
geo-met-ric griev-ance griev-ous grad-ed ho-lo-no-my ho-mo-thetic
in-fin-ite-ly in-fin-i-tes-i-mal Ha-rish Cha-n-dra mul-ti-plic-able 
non-euclid-ean non-iso-mor-phic non-smooth par-a-digm 
par-a-bol-ic pa-rab-o-loid pa-ram-e-trize phe-nom-e-non 
post-script pseu-do-dif-fer-en-tial pseu-do-fi-nite 
qua-drat-ics quad-ra-ture Han-kel rec-tan-gle semi-def-i-nite 
set-up wide-spread Euler-ian Feb-ru-ary Gauss-ian Grothen-dieck 
Hamil-ton-ian Her-mi-t-ian her-mi-t-ian Jan-u-ary 
Japan-ese Ka-shi-wa-ra Kor-te-weg Le-gendre No-vem-ber Rie-mann-ian 
Sep-tem-ber Za-mo-lo-d-chi-kov Kni-zh-nik quan-tum Op-dam
Mac-do-nald Ca-lo-ge-ro Su-ther-land Mo-ser 
Ol-sha-net-sky  Pe-re-lo-mov in-de-pen-dent ope-ra-tors 
cy-clo-to-mic ra-tio-nal de-gen-er-a-tion 
in-ter-est-ing de-for-ma-tions de-for-ma-tion pro-ce-dure 
fol-lows ope-ra-tors  pre-serve suf-fices ap-proach 
for-mu-las con-sider its com-ple-tion cor-re-spond-ing 
au-to-mor-phism be-cause pro-por-tional fi-nal-ly let-ting 
equi-v-a-lence ge-n-er-al-ized Mac-do-nald iden-ti-ties 
cor-re-s-pond sub-dia-grams par-ti-tion na-t-u-ral-ly 
or-dered stan-dard de-for-ma-tion ar-gu-ment com-bined 
sphe-r-i-cal rep-re-sen-ta-tions tri-go-no-me-t-ric
ge-n-er-al-ly speak-ing pri-m-it-ive ir-re-du-cible 
sum-ma-tion  rep-re-sen-ta-tives pro-por-ti-o-na-li-ty
ultra-sphe-ri-cal Ro-gers}

\def\ffor{\quad\hbox{ for }\quad}
\def\wwhen{\quad\hbox{ when }\quad}
\def\wwhere{\quad\hbox{ where }\quad}
\def\aand{\quad\hbox{ and }\quad}
\def\for{\  \hbox{ for } \ }
\def\iif{ \ \hbox{ if } \ }
\def\when{ \ \hbox{ when } \ }
\def\where{\  \hbox{ where } \ }
\def\and{\  \hbox{ and } \ }
\def\and{\  \hbox{ and } \ }
\def\oor{\  \hbox{ or } \ }
\def\proof{{\em Proof. \  }}

\def\equal{\stackrel{\,\mathbf{def}}{= \kern-3pt =}}

\def\la{\lambda}
\def\La{\Lambda}
\def\om{\omega}
\def\Om{\Omega}
\def\Th{\Theta}
\def\th{\theta}
\def\al{\alpha}
\def\be{\beta}
\def\ga{\gamma}
\def\ep{\epsilon}
\def\up{\upsilon}
\def\Up{\Upsilon}
\def\de{\delta}
\def\De{\Delta}
\def\ka{\kappa}
\def\kapp{\hbox{\bf \ae}}
\def\si{\sigma}
\def\Si{\Sigma}
\def\Ga{\Gamma}
\def\ze{\zeta}
\def\io{\iota}
\def\bio{b^\iota}
\def\aio{a^\iota}
\def\twio{\tilde{w}^\iota}
\def\hwio{\hat{w}^\iota}
\def\gio{\g^\iota}
\def\Bio{B^\iota}

\def\del{\delta}
\def\pa{\partial}
\def\vp{\varphi}
\def\ve{\varepsilon}
\def\inf{\infty}

\def\vph{\varphi}
\def\vps{\varpsi}
\def\vPh{\varPhi}
\def\vep{\varepsilon}
\def\vpi{{\varpi}}
\def\vth{{\vartheta}}
\def\vsi{{\varsigma}}
\def\vrh{{\varrho}}

\def\bph{\bar{\phi}}
\def\bsi{\bar{\si}}
\def\bvp{\bar{\varphi}}

\newcommand{\bS}{{\mathbf S}}
\newcommand{\bH}{{\mathbf H}}
\newcommand{\bF}{{\mathbf F}}
\newcommand{\bE}{{\mathbf E}}

\def\tal{\tilde{\alpha}}
\def\tbe{\tilde{\beta}}
\def\tde{\tilde{\delta}}
\def\tpi{\tilde{\pi}}
\def\txi{\tilde{\xi}}
\def\tPi{\tilde{\Pi}}
\def\tPhi{\tilde{\Phi}}
\def\tV{\tilde{V}}
\def\tJ{\tilde{J}}
\def\tla{\tilde{\lambda}}
\def\tga{\tilde{\gamma}}
\def\tGa{\tilde{\Gamma}}
\def\tvs{\tilde{{\varsigma}}}
\def\tu{\tilde{u}}
\def\tU{\tilde{U}}
\def\tw{\widetilde w}
\def\tW{\widetilde W}
\def\tB{\tilde B}
\def\tv{\tilde v}
\def\tV{\tilde V}
\def\tz{\tilde z}
\def\tb{\tilde b}
\def\ta{\tilde a}
\def\tih{\tilde h}
\def\trh{\tilde {\rho}}
\def\tx{\tilde x}
\def\tf{\tilde f}
\def\tg{\tilde g}
\def\tG{\tilde G}
\def\tk{\tilde k}
\def\tl{\tilde l}
\def\tL{\tilde L}
\def\tD{\tilde D}
\def\tR{\tilde R}
\def\tP{\tilde P}
\def\tH{\tilde H}
\def\tp{\tilde p}

\def\hH{\hat{H}}
\def\hh{\hat{h}}
\def\hR{\hat{R}}
\def\hY{\hat{Y}}
\def\hX{\hat{X}}
\def\hP{\hat{P}}
\def\hT{\hat{T}}
\def\hV{\hat{V}}
\def\hG{\hat{G}}
\def\hF{\hat{F}}
\def\hw{\widehat{w}}
\def\hW{\widehat{W}}
\def\hu{\hat{u}}
\def\hs{\hat{s}}
\def\hv{\hat{v}}
\def\hb{\hat{b}}
\def\hB{\widehat{B}}
\def\hze{\hat{\zeta}}
\def\hsi{\hat{\sigma}}
\def\hrh{\hat{\rho}}
\def\hth{\hat{\theta}}
\def\hy{\hat{y}}
\def\hx{\hat{x}}
\def\hz{\hat{z}}
\def\hg{\hat{g}}
\def\he{\hat{e}}
\def\hE{\widehat{E}}

\def\B{\mathbf{B}}
\def\I{\mathbf{I}}
\def\P{\mathbf{P}}
\def\G{\mathbf{G}}
\def\S{\mathbf{S}}
\def\F{\mathbf{F}}
\def\one{\mathbf{1}}
\def\Sn{\mathbf{S}_n}
\def\0{\mathbf{0}}
\def\H{\mathbf{H}}
\def\V{\mathbf{V}}

\def\f{\mathcal{F}}
\def\çF{\mathcal{F}}
\def\o{\mathcal{O}}
\def\t{\mathcal{T}}
\def\r{\mathcal{R}}
\def\l{\mathcal{L}}
\def\m{\mathcal{M}}
\def\k{\mathcal{K}}
\def\n{\mathcal{N}}
\def\d{\mathcal{D}}
\def\p{\mathcal{P}}
\def\cP{\mathcal{P}}
\def\a{\mathcal{A}}
\def\h{\mathcal{H}}
\def\c{\mathcal{C}}
\def\y{\mathcal{Y}}
\def\e{\mathcal{E}}
\def\v{\mathcal{V}}
\def\z{\mathcal{Z}}
\def\x{\mathcal{X}}
\def\s{\mathcal{S}}
\def\g{\mathcal{G}}
\def\u{\mathcal{U}}
\def\w{\mathcal{W}}
\def\i{\mathcal{I}}
\def\j{\mathcal{J}}
\def\b{\mathcal{B}}

\def\lan{\langle}
\def\llb{(\!(}
\def\ran{\rangle}
\def\rrb{)\!)}
 \def\dim{{\hbox{\rm dim}}_{\mathbb C}\,}
\def\lng{\hbox{\rm{\tiny lng}}}
\def\sht{\hbox{\rm{\tiny sht}}}
\def\sph{\hbox{\rm{\tiny sph}}}
\def\inv{\hbox{\rm{\tiny inv}}}

\def\br#1{\langle #1 \rangle}

\def\rank{\hbox{rank}}
\def\gl{\mathfrak{gl}_N}

\newcommand{\Aut}{\operatorname{Aut}}
\newcommand{\Hom}{\operatorname{Hom}}
\newcommand{\End}{\operatorname{End}}
\newcommand{\Ind}{\operatorname{Ind}}
\newcommand{\ad}{\operatorname{ad}}
\newcommand{\pr}{\operatorname{pr}}
\newcommand{\aweyl}{\tilde{\mathbb S}_n}
\newcommand{\hec}{{\mathcal H}^t_n}
\newcommand{\Func}{{\mathcal F}({\mathbb C}^n,{\mathcal H}^t_n)}
\newcommand{\tr}{\operatorname{tr}}
\newcommand{\Out}{\operatorname{Out}}
\newcommand{\Rad}{\operatorname{Rad}}
\newcommand{\Spec}{\operatorname{Spec}}
\newcommand{\id}{\operatorname{id}}
\newcommand{\Int}{\operatorname{Int}}
\newcommand{\ct} {\operatorname{ct}}

\newcommand{\rat}{{\mathbb Q}}
\newcommand{\real}{{\mathbb R}}
\newcommand{\cplx}{{\mathbb C}}
\newcommand{\zint}{{\mathbb Z}}

\newcommand{\sq}{\phantom{1}\hfill$\qed$}
\newcommand{\Rea}{\Re}
\newcommand{\Ima}{\Im}

\newcommand{\st}{\bowtie}
\newcommand{\modd}{\mbox{\,mod\,}}
\newcommand{\lr}{\langle}
\newcommand{\rr}{\rangle}
\newcommand{\eps}{\varepsilon}
\newcommand{\phk}{\phi^{(k)}}
\newcommand{\psk}{\psi^{(k)}}
\newcommand{\Res}{\mbox{Res}\;}
\newcommand{\sgn}{\mbox{sgn}}
\newcommand{\mn} {\left\{ \begin{array}{c}m\\
n\end{array}\right\}}

\def\sX{\mathscr{X}}
\def\sH{\mathscr{H}}
\def\sY{\mathscr{Y}}
\def\TT{\mathfrak{T}}
\def\JJ{\mathfrak{J}}
\def\HH{\mathfrak{H}}
\def\FF{\mathfrak{F}}
\def\GG{\mathfrak{G}}
\def\CC{\mathfrak{C}}
\def\LL{\mathfrak{L}}

\def\BB{\mathfrak{B}}
\def\AA{\mathfrak{A}}
\def\ZZ{\mathfrak{Z}}
\def\HH{\hbox{${\mathcal H}$\kern-5.2pt${\mathcal H}$}}
\def\HHH{\hbox{${\mathbb H}$\kern-4.2pt${\mathbb H}$}}
\def\tHH{\widetilde{\HH\ }}

\font\smm=msbm10 at 12pt 
\def\symbol#1{\hbox{\smm #1}}
\def\lsmash{{\symbol n}}
\def\rsmash{{\symbol o}}
\def\#{\sharp}

\font\tenbf=cmbx10
\font\tenrm=cmr10
\font\tenit=cmti10
\font\ninebf=cmbx9
\font\ninerm=cmr9
\font\nineit=cmti9
\font\eightbf=cmbx8
\font\eightrm=cmr8
\font\eightit=cmti8
\font\sevenrm=cmr7
\font\sevenbf=cmbx7



\title [Superpolynomials of algebraic links]
{Superpolynomials of algebraic links}
\author[Ivan Cherednik]{Ivan Cherednik $^\dag$}

{
\centering 
{\sf\em Dedicated to the memory of Ian Macdonald
\medskip\par}
\vskip 0.3cm}

\begin{abstract}
Theory of motivic superpolynomials is developed, including
its extension to algebraic links colored by rows, relations
to $L$-functions of plane curve singularities, the justification
of the motivic versions of Weak Riemann Hypothesis, and 
recurrences for iterated torus links. The key theme is 
the conjectural coincidence of motivic 
superpolynomials with the DAHA ones, which can be 
interpreted as a far-reaching generalization of the Shuffle 
Conjecture. Applications include 
affine Springer fibers of type $A_n$ and compactified Jacobians
in the most general case
(for arbitrary characteristic polynomials) and extended
rho-invariants of algebraic knots. 
The 2nd connection conjecture relates the 
superpolynomials to the Galkin-St\"ohr $L$-functions,
which is some counterpart of the ORS conjecture. 
The corresponding theory of plane curve singularities
is systematically exposed and developed, which 
can be seen in the case of Hopf links as
a generalized version of Schubert Calculus. 
\end{abstract}

\thanks{$^\dag$ \today.
\ \ \ Partially supported by NSF grant
DMS--1901796}

\address[I. Cherednik]{Department of Mathematics, UNC
Chapel Hill, North Carolina 27599, USA\\
chered@email.unc.edu}

\newcommand{\hga}{\hat{\ga}}
 \def\sht{\raisebox{0.4ex}{\hbox{\rm{\tiny sht}}}}
 \def\bysame{{\bf --- }}
\let\oldt\~   
 \def\~{{\bf --}}
 \def\rr{{\mathsf r}}
 \def\ss{{\mathsf s}}
 \def\tt{{\mathsf t}}
 \def\mm{{\mathsf m}}
 \def\pp{{\mathsf p}}
 \def\ll{{\mathsf l}}
 \def\aa{{\mathsf a}}
 \def\bb{{\mathsf b}}
 \def\cc{{\mathsf c}}
 \def\NS{\hbox{\tiny\sf ns}}
 \def\ssum{\hbox{\small$\sum$}}
\newcommand{\comment}[1]{}
\renewcommand{\tilde}{\widetilde}
\renewcommand{\hat}{\widehat}
\renewcommand{\V}{\mathbb{V}}
\renewcommand{\S}{\mathbb{S}}
\renewcommand{\F}{\mathbb{F}}
\newcommand{\q}{\mathcal{Q}}
\newcommand{\dagx}{\hbox{\tiny\mathversion{bold}$\dag$}}
\newcommand{\ddagx}{\hbox{\tiny\mathversion{bold}$\ddag$}}
\newtheorem{conjecture}[theorem]{Conjecture}
\newcommand*\toeq{
\raisebox{-0.15 em}{\,\ensuremath{
\xrightarrow{\raisebox{-0.3 em}{\ensuremath{\sim}}}}\,}
}
\newcommand{\unknot}{\hbox{\tiny\!\raisebox{0.2 em}
{$\bigcirc$}}\!}
\newcommand{\mmu}{\hbox{\mathversion{bold}$\mu$}}
\newcommand{\lla}{\hbox{\mathversion{bold}$\lambda$}}
\newcommand{\dde}{\hbox{\mathversion{bold}$\delta$}}

\newcommand\rightthreearrow{\hbox{\tiny
        $\mathrel{\vcenter{\mathsurround0pt
         \ialign{##\crcr
         \noalign{\nointerlineskip}$\rightarrow$\crcr
         \noalign{\nointerlineskip}$\rightarrow$\crcr
         \noalign{\nointerlineskip}$\rightarrow$\crcr
                }}}$ }}
\newcommand\rightfourarrow{\hbox{\tiny
        $\mathrel{\vcenter{\mathsurround0pt
         \ialign{##\crcr
         \noalign{\nointerlineskip}$\rightarrow$\crcr
         \noalign{\nointerlineskip}$\rightarrow$\crcr
         \noalign{\nointerlineskip}$\rightarrow$\crcr
         \noalign{\nointerlineskip}$\rightarrow$\crcr
                }}}$ }}

\newcommand\rightdotsarrow{\hbox{\small
        $\mathrel{\vcenter{\mathsurround0pt
         \ialign{##\crcr
         \noalign{\nointerlineskip}$\,\rightarrow$\crcr
         \noalign{\nointerlineskip}$\cdots$\crcr
         \noalign{\nointerlineskip}$\,\rightarrow$\crcr
                }}}$ }}
\newcommand\rightdotsarrowtiny{\hbox{\tiny
        $\mathrel{\vcenter{\mathsurround0pt
         \ialign{##\crcr
         \noalign{\nointerlineskip}$\,\rightarrow$\crcr
         \noalign{\nointerlineskip}$\cdots$\crcr
         \noalign{\nointerlineskip}$\,\rightarrow$\crcr
                }}}$ }}

\newcommand*{\vect}[1]{\overrightarrow{\mkern0mu#1}}
\newcommand{\twoone}
{\hbox{\rm
$\circ\!$\raisebox{-2.6pt}{$\rightarrow$}
\raisebox{-2.6pt}{$\!\!\circ\!\!\rightarrow$}
\kern-33pt\raisebox{+2.6pt}{$\rightarrow$}
\kern+14pt}}
\newcommand{\twoonetiny}
{\hbox{\tiny
$\circ\!$\raisebox{-2.pt}{$\rightarrow$}
\raisebox{-2.pt}{$\!\!\circ\!\!\rightarrow$}
\kern-23pt\raisebox{+2.pt}{$\rightarrow$}
\kern+10pt}}

\newcommand{\twotwo}
{\hbox{\rm $\circ\!\rightrightarrows\!
\raisebox{2.5pt}{\hbox{\small $\circ$}}
\kern-5.5pt\raisebox{-2.5pt}{\hbox{\small $\circ$}}
\!\rightrightarrows\,$}}
\newcommand{\twotwotiny}
{\hbox{\tiny $\circ\!\rightrightarrows\!
\raisebox{2.pt}{\hbox{\tiny $\circ$}}
\kern-4.2pt\raisebox{-2.pt}{\hbox{\tiny $\circ$}}
\!\rightrightarrows$}}
\newcommand{\tax}{\hbox{\sf[r,s]}}
\newcommand{\lxi}{\raisebox{0.5pt}{${}^\xi$}\!}
\vskip -0.0cm

\newcommand{\qbin}[2]{\begin{bmatrix}{#1}\\ {#2}\end{bmatrix}_q}

\maketitle
\smallskip

{\footnotesize
\centerline{{\bf MSC} (2010): 14H50, 17B22, 17B45, 20C08,
20F36, 33D52, 30F10, 55N10, 57M25}
}
\smallskip

\vskip -0.5cm
\renewcommand{\baselinestretch}{1.2}
{\vbadness=10000\textmd
{
\tableofcontents}
}
\renewcommand{\baselinestretch}{1.2}
\vfill\eject

\renewcommand{\natural}{\wr}

\setcounter{section}{0}
\setcounter{equation}{0}
\section{\sc Introduction}
\vbadness=3000
\hbadness=3000

\subsection{\bf Overview} The main line of the paper is the conjectural
coincidence of the {\em DAHA superpolynomials} of algebraic links
colored by rows 
with the corresponding {\em motivic superpolynomials}. The definition
and study of the latter in such a generality is new.
This is part
of the program outlined in \cite{ChW,ChS}. For instance, we
consider applications to $L$-functions of
plane curve singularities, prove {\em Weak Riemann Hypothesis}
in the unibranch case (for $a=0$), and prove the motivic
counterpart of the DAHA recurrence relation connecting 
DAHA superpolynomials of algebraic cables with those of
colored torus links. Many examples are provided in this
paper; the main focus is on the machinery of counting
{\em standard modules}.  Trefoil, double trefoil,
colored trefoil, Hopf links and certain families of torus links and
cables are considered in detail. A lot of interesting ``linear algebra"
is involved even in relatively simple cases. This theory
can be seen as unification of 
{\em Schubert Calculus} with the theory of {\em plane curve 
singularities}, classical and new.
 
\vfil

The DAHA superpolynomials were defined in \cite{CJ,CJJ} for 
torus knots colored by any Young diagrams;
the main construction was for any reduced root systems and even
for DAHA for  $C^\vee C_1$, the root system with four $t$
(in \cite{CJJ}). These papers were
extended in \cite{ChD1} to colored {\em iterated} torus knots 
and  to iterated torus {\em links} in \cite{ChD2}. Algebraic links
are iterated torus links, but the latter class is significantly
wider. The
{\em Khovanov-Rozansky polynomials}, the origin of 
superpolynomials, can be generally defined for any links, not
only iterated torus ones, but this theory remains challenging for 
links and with (arbitrary) colors, especially in the reduced
setting. 
\vfil

The {\em motivic} superpolynomials are known by now only for 
{\em algebraic} links and when the colors are ``rows". They
were defined for {\em irreducible} plane curve singularities
in \cite{ChP1}, which was extended to arbitrary ranks 
in \cite{ChP2}; see also \cite{ChW}. We extend \cite{ChP2}
in this paper  to any plane curve singularities. This is
a significant development, corresponding to  
the consideration of
{\em arbitrary} affine Springer fibers of type $A$, the ones
for any characteristic polynomials, possibly reducible and
not square-free. The $p$-adic {\em orbital integrals} of type $A$
are motivic superpolynomials when  $t=1,a=0$ and in the absence
of colors; see  \cite{Cha, Yun} for some generalizations.
A challenge is to extend their theory 
and the corresponding Fundamental Lemma 
to general $t,a$ and arbitrary, possibly reducible and
non-reduced, {\em spectral curves}.  

\subsection{\bf Main developments and connections}
There are various developments and improvements in this
paper. For instance, the usage of flags of {\em standard modules} 
in the definition of motivic superpolynomials from 
\cite{ChP1,ChP2,ChW} is complemented
by the {\em decomposition formula}
 where only standard modules $M$ and their $q$-ranks 
$r\!k_q(M)$ occur. 

The notion of standard modules is the key for us. 
Given a ring of singularity $\r$ over a base field $\F$
with $\kappa$ irreducible components,  
they are $\r$-invariant submodules of 
$\Om=\oplus_{i=1}^\kappa \o_i^{\oplus c_i}$ such
that $M\otimes \o=\Om$, where 
$\o_i=\F[[z_i]]$ are normalizations of 
the rings $\r_i$ for irreducible components. 
The numbers $c_i\in \N$, called $\o$-ranks in this paper,
correspond to the weights $c_i\om_1$ (colors)
in the DAHA theory. 
\vfil

The main decomposition formulas are 
in terms of $\prod_{i=c}^{r\!k_q(M)-1} (1-q^i a)$, where $c=\max\{c_i\}$. 
The original
definition from \cite{ChP1,ChP2} was in terms of certain
standard flags with $a^\ell$ associated with their lengths $\ell$. 
For colored algebraic links, the formulas directly in terms of $r\!k_q$
appeared more convenient,  though the usage of flags is possible too.

There are various applications. 
The proof of {\em Weak Riemann Hypothesis}
from \cite{ChW} in this paper (for uncolored knots and $a=0$) is one
on them. Among other applications, let us mention
motivic formulas for the {\em DAHA-Jones polynomials} (before the
$a$-stabilization). For instance,
DAHA-Jones polynomials for  $A_1$ 
become (conjecturally) certain sums over standard modules 
of $r\!k_q\le 2$.

For algebraic knots, another important
decomposition is obtained, which is with respect to  
$\prod_{i}(1-q^i a/t)$.
There is one more, in terms of $\prod_i (1-(q/t)^i a)$, 
assuming the conjectural  relations of the 
motivic superpolynomials
to generalized Galkin-St\"ohr $L$-functions and those
from \cite{ORS}. There are other decompositions 
if the superduality is used:  $q\leftrightarrow 1/t$.
All of them  are meaningful geometrically.
They are related
to {\em differentials} in the theory
of HOMFLY-PT homology \cite{KhR1,KhR2,Rou,Ras,DGR}, but there
are deviations; see \cite{CJ}. 

\vfil
 
Theory of torsion-free modules of any ranks 
over singular curves is generally a 
difficult one, even for small ranks
and relatively simple curve singularities. 
We use that this theory
is not that involved for {\em pure}
singularities. The corresponding {\em Quot schemes} 
result in {\em compactified Jacobians}, projective
schemes, but we only need their stratification in
terms of {\em standard} modules for superpolynomials. 
We note that {\em all} our motivic constructions can be 
transferred to $p$-adic setting (to $char=0$), which is based
on the usage of {\em Witt vectors}. 
This is parallel to {\em Fundamental Lemma}. The 
$p$-adic superpolynomials obtained so far are no different
from those over $\F_q$, but we do not conjecture this.
\vskip 0.2cm

The definition
of generalized (flagged) Galkin-St\"ohr $L$-functions and
their conjectural relation to motivic superpolynomials is
another line of this paper. They
are Dedekind sums over ideals of $\r$ of finite colength
enhanced by adding $r\!k_q$. This is very general, but 
the connection with motivic superpolynomials is expected
only for plane curve singularities; for instance, it
does not hold for Gorenstein $\r$. The
corresponding  coincidence
conjecture is stated only in the uncolored case in
this paper. The $L$-functions
of arbitrary $\o$-ranks require further 
study; see \cite{ChS} for preliminary definitions
and some examples.
\vfil

Such $L$-functions are related to the {\em nested Hilbert
schemes of ideals} of plane curve singularities
and the {\em ORS Conjecture} from \cite{ORS, GORS}.
Our motivic superpolynomials are associated with generalized
{\em compactified Jacobians}. This is a different direction,
though with ample connections with  Hilbert schemes.
\vskip 0.2cm
\vfil

The construction of {\em DAHA superpolynomials},
the first part of this paper,
goes in two steps:
First, we define  the {\em DAHA-Jones
polynomials}, which is for arbitrary reduced 
root systems, ``colored" by arbitrary dominant weights and depending
on the DAHA-parameters $q,t$. This theory was triggered by 
paper \cite{AS}. 

Second, the {\em $a$-stabilization theorem} states that there exist 
certain {\em superpolynomials}
in terms of $a,q,t$ such that the DAHA-Jones
polynomials for $A_n$  are their specializations at $a=-t^{n+1}$. 
This is expected to hold for classical
series of root systems and even for the $E$-series.
In contrast to DAHA superpolynomials, where the
$a$-stabilization requires effort, 
the $a$-parameter is fully transparent in
the motivic construction; however the polynomial
dependence on $q$ is an open problem for them.
\vfil

The definition of the colored DAHA-Jones polynomials 
requires the theory of Macdonald polynomials; they
incorporate colors, which are Young diagrams in the $A$-case.
The combinatorics of these polynomials is complicated unless for
minuscule weights, when these polynomials 
 become symmetric monomial functions.
All fundamental weights are minuscule  for $A_n$.
For such weights, the DAHA construction becomes entirely combinatorial;
only the action of the 
projective $SL(2,\Z)$ on the generators of DAHA and the
definition of {\em coinvariant} are needed. Due to the superduality,
this is when the {\em First Coincidence Conjecture} with
motivic superpolynomials is stated.

This conjecture is
a far-reaching generalization of the so-called
{\em Shuffle Conjecture} proved in \cite{CaM}. Its ``combinatorial
half" is a special case of the  DAHA construction
(upon the usage of $E$-polynomials). The 
formulas for motivic superpolynomials generalize
(a great deal) its  ``formula-half". Explicit
formulas for the motivic superpolynomials
can be generally obtained for any given {\em family}
of plane curve singularities, but they can be involved. This
is mostly based on 
the {\em method of syzygies}, similar to the machinery
from \cite {Pi}. In DAHA theory, many explicit formulas 
for superpolynomials are known,
including quite involved colored cables and links; see
\cite{ChD2,ChW} and references there. 

\subsection{\bf Motivic construction}
Plane curve singularities are those  for local complete
rings  $\r$ over the base field $\F$
with $1$ and two generators $x,y$ 
(as complete rings). They will be considered
as subrings in the normalization ring $\o^\kap=
\oplus_{i=1}^\kap \o_i$, where 
$\o_i\equal\F[[z_i]]$ and $\kap$ is the number of
irreducible components. Namely, $F(x,y)=0$ for
$F(x,y)=\prod_{i=1}^\kap F_i(x,y)$ for $F_i$ over
$\F$  assumed absolutely irreducible in this paper
and non-proportional. 
The localization of $\r$ must be $\oplus_{i=1}^\kap\F((z_i))$.

If $\F\subset \C$,  $\kap$ is the number of
connected components of the corresponding algebraic link, which is
$\{F(x,y)\!=\!0, (x,y)\!\in\! \C^2\}\cap S^3_\ep$ for a small 
$3$-sphere $S^3_\ep$ centered at $(0,0)$. 

The units $e_i\in \o_i$ are naturally 
idempotents in $\r$:\, 
$1=\sum_{i=1}^\kap e_i$, $e_i e_j=
\de_{ij}e_i$. We set  $z_i=z e_i$ for
$z\equal z_1+\cdots+ z_\kap$, and always assume that 
$x_i=x e_i$ and $y_i=y e_i$ belong to $(z_i)=z_i\o=z\o_i=z_i\o_i$
for all $i$. The maximal ideal $\mathfrak{m}_{\r}$ of $\r$ 
is the one generated by the (formal) series
in terms of $x$ and $y$ with zero constant term.

The corresponding {\em standard}
$\r$-submodules
$M\in \o$, the key objects for us, are such that 
$M\otimes_{\r}\o_i=\o_i$ for all $i$, where these and similar
tensor products are naturally 
considered as $\o$-submodules of $\o$ and 
its projections onto $\o_i$.   
Explicitly, this definition means that 
the projections $M_i$
of $M$ onto $\o_i$ must {\em not} be  inside $(z_i)$ for
any $1\le i\le \kap$.

We mostly take  $\F=\r/\mathfrak{m}_{\r}=\F_q$ as the base
field (for $q=p^m$ elements for prime $p$), assuming
that $\r_i/\mathfrak{m}_i=\F_q$ for $1\le i\le \kap$. Generally,
the latter can be extensions of $\F_q$ as in \cite{Sto},
but we do not know how to interpret the
corresponding superpolynomials topologically. 

Such  $M$ are of
$\o$-rank one, but can have several generators 
over $\r$. The (minimal) number of generators  equals 
$r\!k_q(M)\equal \text{\,dim\,}_{\F_{\tilde{q}}}M/\mathfrak{m}_{\r}M$,
  the {\em $q$-rank} of $M$. The {\em invertibles} are $M$
with $r\!k_q(M)=1$. 

The next step is to go to 
$\r$-modules $M$ of any $\o$-ranks.
For $\si=(c_1\ge c_2\ge \cdots \ge c_\kap>0)$, 
{\em standard modules} are defined then as $\r$-submodules 
$M\subset \Om\equal
\o_1^{c_1}\oplus \o_2^{c_2}\oplus\cdots\oplus \o_\kap^{c_{\kap}}$ 
such that
$M\otimes_{\r}\o_i=\o_i^{c_i}$ for any $1\le i\le \kap$. Accordingly,
$r\!k_q(M)\equal \text{\,dim\,}_{\F_{\tilde{q}}}\,
M/\mathfrak{m}M\ge c_1$. This is the way to address 
non-reduced $\r$, which are for not square-free equations. The
modules are
called {\em invertible} if $rk_q(M)=c_1$, though this does not mean the
usual invertibility if $c_1>1$.

The 
{\em motivic superpolynomial} of $\r,\si$ is defined as follows:
$$
\h^{mot}_{\si}=\sum_{M} t^{\hbox{\tiny dim}(\Om/M)}
\prod_{j=r\!k_{min}}^{r\!k_q(M)-1}(1+aq^j)
\text{ summed over standard } M.
$$
The first product $\Pi$ here (for $r\!k_q(M)=c_1$) is $1$.
The conjecture is that $\h^{mot}_\si$  depend polynomially on $q,t,a$
and coincide with the DAHA superpolynomials
$\hat{\h}^{{\la}}(q,t,a)$ for the corresponding links colored
by the sequences $\{\la\}=\{c_1\om_1,c_2\om_1,\ldots,
c_\kappa \om_1\}$ (only pure rows). It implies that
$\h^{mot}_{\si}$  are topological
invariants because so are the DAHA superpolynomials (a relatively
simple corollary of basic DAHA facts).

The
standard modules form quasi-projective
varieties when \,dim\,$\Om/M$ is fixed.   
For  $\kap=1$, these varieties 
 can be arranged
as some strata of a {\em single} quasi-projective variety following 
\cite{ChP2}, which can be extended to any $\kap$. For $\kap=1$
and $c_1=1$, i.e. in the uncolored unibranch case, we 
arrive at the
{\em compactified Jacobian} of $\r$, which is an irreducible
projective variety due to  \cite{PS,Rego}. 

This variety and its generalizations to arbitrary $\kap$
and $\si$ are of obvious interest, but the strata above
 are sufficient for motivic superpolynomials 
(for any motivic constructions). 
\vskip 0.2cm

The {\em Second Coincidence Conjecture} states
that $\h^{mot}(q\mapsto qt, t, a)$ coincide
with the generalized {\em Galkin-St\"ohr
$L$-functions}. See Conjecture \ref{conj:HLlinks}
and \cite{KT}.
This conjecture is entirely in terms
of plane curve singularities; however, it includes
a nontrivial passage from  $q$ and $q t$
and does not hold for more general (Gorenstein)
curve singularities.
\vskip 0.2cm

The coincidence conjectures have 
various aspects. Some features can be
justified in one theory but remain conjectures
in other theories.  For instance, 
the {\em Weak Riemann Hypothesis} stated for DAHA
superpolynomials in 
\cite{ChW} is proven motivically in this paper 
(for $a=0$). Also, the dependence on $a$
is quite transparent for motivic superpolynomials;
by contrast with the $a$-stabilization in the DAHA approach.  
On the other hand, the topological invariance of
superpolynomials is a DAHA theorem, which remains a conjecture in the
motivic setting. 

The superduality of superpolynomials is of conceptual importance. 
It is proven in full generality for DAHA. It is, generally, 
the invariance (up to some multipliers) with respect to  
 $q\leftrightarrow t^{-1}$, when $a$ remains unchanged
and the corresponding Young diagrams are transposed.   
However, the superduality remains a 
conjecture for motivic uncolored superpolynomials (when $c_i=1$).
Importantly, it becomes
the {\em functional equation} for Galkin-St\"ohr $L$-functions.
St\"ohr found 
a relatively simple ``combinatorial" proof, which can be
extended to our $a$-generalizations.

As another implication of the coincidence conjectures,
the passage to the
generalized  
$L$-functions results in 
the definition of the
{\em quasi-rho invariants of algebraic knots}. The latter are 
a certain integer variant
of the classical von Neumann {\em rho-invariants} of knots.
Our motivic superpolynomials provide their ``triple
deformations", which is a combination
of certain algebraic manipulations and the
usage of the superduality.

\subsection{\bf Geometric and other aspects}
The origin of the theory of $L$-functions for
curve singularities is the following classical formula:
$$\ze(X,t)\equal
\exp(\,\sum_{n=1}^\infty\,t^n\,|X(\mathbb F_{q^n})|\,/\,n\,)=  
\frac{P_1(t)}{(1-t)(1-qt)},
$$
where $X$ is a smooth projective curve over
$\F_q$ and $P_1$
is the usual polynomial
of degree $2g$ for the genus $g$ of $X$. The smoothness of $X$ 
is mostly needed for the Riemann Hypothesis; the rationality
holds without this assumption (Dwork and others).  

The Galkin-St\"ohr $L$-functions are 
counterparts of $P_1(t)$. They satisfy the same
Hasse-Weil functional
equation, which is for $t\mapsto 1/(qt)$. This is for 
Gorenstein curve singularities; plane curve singularities are 
such. The corresponding
zeta-function is $Z(q,t)=L(q,t)/(1-t)^{\kappa}$; it does not
satisfy the functional equation.

\vskip 0.2cm

{\em Kapranov's zeta} is the most general motivic 
version of $Z(q,t)$. The justification 
of the functional equation 
from Proposition 3 from \cite{ORS} is
parallel
to that proof of the {\em motivic} functional equation 
in \cite{Kap}; see also Section 6 in \cite{Gal} and 
Section 3 in \cite{Sto}. 

Adding $a$ and the colors is beyond the Hasse-Weil-Deligne theory 
\cite{Del1,Del2}, as well as the following feature.
All coefficients of DAHA superpolynomials are presumably 
positive integers for rectangle Young diagrams and
algebraic knots. This is the last 
unproven {\em intrinsic} conjecture from the list
in \cite{CJ}, extended to any algebraic knots 
in  \cite{ChD1}.  For links and non-rectangle 
Young diagrams, the positivity conjecturally holds only 
upon the division by
some powers of $(1-t), (1-q)$. For instance, the division by
$(1-t)^{\kap-1}$ is presumably sufficient
for uncolored algebraic links with $\kap$ components.
\vskip 0.2cm

In the DAHA theory of the  superpolynomials, the rationality in
terms of $q,t$ is manifest. The polynomiality for arbitrary weights
is more subtle because the
coefficients of Macdonald polynomials have $q,t$-denominators.
This is  Theorem 1.2, \cite{CJJ},  extended later to 
arbitrary colored iterated torus links. The $a$-stabilization
for torus knots was announced in \cite{CJ,CJJ} as a corollary of
some lemma in \cite{SV}); its complete proof
was provided in \cite{GoN}. See \cite{ChD2} for the most 
general version:\, iterated torus links colored by any Young
diagrams.
\vskip 0.1cm

The DAHA {\em superduality} was conjectured for torus knots in \cite{CJ};
it was proved in \cite{GoN} on the basis of the
$q\leftrightarrow t$\~duality of the
modified Macdonald polynomials.
An alternative approach to the proof via roots of unity and the
generalized {\em level-rank
duality\,} was presented in \cite{CJJ}; it is expected to work 
for classical root systems. Thus,
our connection conjectures $\h^{daha}(q,t,a)=\h^{mot}(q,t,a)=
L(q/t,t,a)$
link the DAHA superduality 
to  the classical {\em Hasse-Weil's
functional equation}: $\ze(X,q^{-1}t^{-1})=
q^{1-g}t^{2-2g}\ze(X,t)$.

Importantly, the superduality
of physical superpolynomials (with colors) was conjectured
in \cite{GS}. This links the super-duality and
the functional equation for $\h^{daha},\h^{mot},L$ to the
$S$-duality in {\em SCFT}, super-conformal field theory.
In fact, the DAHA theory is closely connected with 
{\em Knizhnik-Zamolodchikov equations} and, as such, to {\em SCFT}.
Accordingly, $Z(q,t,a)$ can be expected
a partition function for motivic topological {\em LGSM}, Landau-Ginzburg
sigma model.  See \cite{ChS}
for further discussion and other relations to physics. The
correspondence {\em SCFT} $\leftrightsquigarrow$ {\em LGSM}
 is a fundamental
feature of {\em string theory}. 

\subsection{\bf Perspectives}
The superpolynomials have surprisingly many incarnations
in mathematics and physics. {\em Khovanov-Rozansky polynomials\,}
are their origin,  
Poincar\'e polynomials of the HOMFLY-PT
triply-graded link homology \cite{Kh,KhR1,Rou}. They depend on
$q,t,a$ and are infinite series in
the {\em unreduced\,} case. The reduced setting is needed to make
them polynomials. Generally, they are difficult
to calculate and there are
unsettled problems for links and
arbitrary colors. The best tool by now is the usage of {\em Soergel modules}; 
see also \cite{Wed}. 
They serve all links.
For uncolored unreduced algebraic links, they were conjectured
in \cite{ORS} to coincide with 
certain generating series
for {\em nested Hilbert schemes\,} of the corresponding
{\em plane curve singularities}. 

Our motivic and DAHA
superpolynomials  are conjectured to give
{\em reduced stable} Khovanov-Rozansky polynomials of algebraic links. 
However, these are different theories.
For instance, DAHA are much 
more suitable to incorporate colors and deal with links and
arbitrary root systems, but can be used (by now) only
for iterated torus links. The {\em superduality} is a theorem for
DAHA, but, generally, an open question in topology.
Only relatively recently, the {\em Soergel modules} allowed to 
check it for some simple colors. The topological invariance of
DAHA superpolynomials is a theorem, but it is still a conjecture
for motivic superpolynomials. However, the latter were used
to justify {\em Weak Riemann Hypothesis}, and so on.

We note here that 
Corollary 3.4 in  \cite{Mel2},
proves Conjecture 2.7 $(ii)$ from \cite{CJ}
(uncolored torus knots) on the connection
with Khovanov-Rozansky theory. 
The proof is essentially a direct combinatorial identification of
Gorsky's combinatorial formulas with those in \cite{CJ} using
Soergel modules. The usage of {\em nonsymmetric} Macdonald polynomials
in \cite{CJ} gives the identification of what this approach with
our conjecture. See also \cite{Kh,EH,Hog,HM}.

The motivic superpolynomials 
conjecturally  coincide with {\em geometric superpolynomials\,}
introduced in \cite{ChP1,ChP2}. They were defined there
for any algebraic knots colored by rows, developing \cite{ChD1}.
See \cite{GoN,GM,GMO,GMV} for various geometric-combinatorial
aspects. 
Their definition is similar to that of  motivic superpolynomials, 
when {\em Borel-Moore homology} is used instead of 
counting the points over $\F_q$. An alternative approach is
the usage of the Serre-Deligne weight filtration as in
\cite{ORS,GORS,KT}, which is better connected with our motivic
constructions.  Generally, we have here 
different choices of  {\em motivic measures}.
See also prior
works, especially \cite{Kap,GSh}. For physics origins, see 
\cite{DGR,AS,GS,DMMSS,FGS}; paper \cite{AS} triggered \cite{CJ}.

\setcounter{equation}{0}
\section{\sc Double Hecke algebras}
\subsection{\bf Affine root systems}
Until the definition of superpolynomials in 
Section \ref{sec:sup-def}, which
are in type $A$, 
any reduced irreducible root systems 
$R=\{\al\}\subset \R^n$ can be taken; $\al_i$ for $1\le i\le n$ 
are the simple roots.

The root lattice is $Q=\oplus_{i=1}^n \Z\al_i$.
Replacing $\Z$ by $\Z_{\pm}=\{m\in\Z, \pm m\ge 0\}$, we obtain
$P_\pm,Q_\pm$. See  e.g., \cite{Bo} and \cite{C101}.

The vectors $\ \tal=[\al,j]\in \R^{n+1}$
for $\al \in R, j \in \Z $ form the
{\em affine root system\,}
$\tR \supset R$, where $\al\in R$ are identified with
$[\al,0]$. Generally, $\tal=[\al,\nu_{\al} j]$, where
$\nu_\al=1$ for short roots, and $\nu_\al=2,3$ for long 
ones ($3$ is for $G_2$). 
The root $\al_0 \equal [-\th,1]$ is added to the simple
roots, where $\th$ is the  {\em maximal short root\,}.
The corresponding set
$\tR_+$ of positive roots is 
$R_+\cup \{[\al,j],\ \al\in R, \ j > 0\}$.

In the case of $A_n$, 
$R\!=\!\{\al\!=\!\vep_i\!-\!\vep_j, i\!\neq\! j\}$
for the basis $\{\vep_i, 1\le i\le n\!+\!1\}\in \R^{n+1}$,
orthonormal with respect to the usual euclidean form
$(\cdot,\cdot)$. The Weyl group is $W=\S_{n+1}$; it is
generated by the reflections (transpositions)
$s_\al$ for the set of positive  roots
$R_+=\{\vep_i-\vep_j, i<j\}$; $R_-=-R_+$.
The simple roots are $\al_i=\vep_i\!-\!\vep_{i+1},$
$\th=\vep_{1}-\vep_{n+1}$. Generally, 
the weight lattice is
$P=\oplus^n_{i=1}\Z \om_i$,
where $\{\om_i\}$ are fundamental weights, satisfying
$ (\om_i,\al_j)=\de_{ij}$. For $A_n$:\,
\begin{align}
&\om_i=\vep_1+\cdots+\vep_i-\frac{i}{n+1}(\vep_1+\cdots+\vep_{n+1})
\for i=1,\ldots,n,\\
&\rho=\om_1\!+\cdots+\!\om_n=\frac{1}{2}\bigl(
(n\!-\!1)\vep_1+(n\!-\!3)\vep_2+\cdots+(1\!-\!n)
\vep_n\bigr).\notag
\end{align}


{\sf Affine Weyl groups.} 
Given $\tal=[\al,j]\in \tR,  \ b \in P$, let
\begin{align}
&s_{\tal}(\tz)\ =\  \tz-(z,\al^\vee)\tal,\
\ b'(\tz)\ =\ [z,\ze-(z,b)]
\label{ondon}
\end{align}
for $\tz=[z,\ze]\in \R^{n+2}$.
The
{\em affine Weyl group\,}
$\tW=\lan s_{\tal}, \tal\in \tR_+\ran$
is the semidirect product $W\lsmash Q$ of
its subgroups $W=$ $\lan s_\al,
\al \in R_+\ran$ and $Q$, where $\al$ in the latter are identified with
\begin{align*}
& s_{\al}s_{[\al,\,1]}=\
s_{[-\al,\,1]}s_{\al}\text{ for any }
\al\in R.
\end{align*}

The {\em extended Weyl group\,} $ \hW$ is $W\lsmash P$, where
the action is
\begin{align}
&(wb)([z,\ze])\ =\ [w(z),\ze-(z,b)] \for w\in W, b\in P.
\label{ondthr}
\end{align}
It is isomorphic to $\tW\rsmash \Pi$ for $\Pi\equal P/Q$.
The latter group consists of $\pi_0=$\,id\, and the
images $\pi_r$ of $\om_r$ in $P/Q$.
Note that
$\pi_r^{-1}$ is $\pi_{r^\iota}$,  where $\iota$ is
the standard involution of the {\em non-affine\,}
Dynkin diagram of $R$,
induced by $\al_i\mapsto \al_{n+1-i}$ for $A_n$.
Generally, we set
$\iota(b)=-w_0(b)=b^{\,\iota}\,$, where $w_0$ is the
longest element
in $W$. For $A_n$:\, the element  $w_0$
 sends $ \{1,2,\ldots,n+1\}$
to $\{n+1,\ldots,2,1\}$, and $\Pi=\Z_{n+1}$.

The group $\Pi$
is naturally identified with the subgroup of $\hW$ of the
elements of the length zero; the {\em length\,} is defined as
follows:
\begin{align*}
&l(\hw)=|\La(\hw)| \for \La(\hw)\equal\tR_+\cap \hw^{-1}(-\tR_+).
\end{align*}
One has $\om_r=\pi_r u_r$ for $1\le r\le n$,
where $u_r$ is the
element $u\in W$ of {\em minimal\,} length such that
$u(\om_r)\in P_-$,
equivalently, $w=w_0u$ is of {\em maximal\,} length such that
$w(\om_r)\in P_+$. The elements $u_r$ are very explicit.
Let $w^r_0$ be the longest element
in the subgroup $W_0^{r}\subset W$ of the elements
preserving $\om_r$; this subgroup is generated by simple
reflections. One has:
\begin{align}\label{ururstar}
u_r = w_0w^r_0 \hbox{\,\, and\,\, } (u_r)^{-1}=w^r_0 w_0=
u_{r^\iota} \for 1\le r\le n.
\end{align}
\smallskip

Setting $\hw = \pi_r\tw \in \hW$ for $\pi_r\in \Pi,\, \tw\in \tW,$
\,$l(\hw)$ coincides with the length of any reduced decomposition
of $\tw$ in terms of the simple reflections
$s_i,\, 0\le i\le n.$ Thus, indeed, $\Pi$ is a subgroup of
$\hW$ of the elements of length $0$. 


\subsection{\bf Definition of DAHA}
We follow \cite{CJJ,CJ,C101}.
The {\em double affine Hecke algebra, DAHA\,}, of
type $A_n$ will be defined
over the ring
$\Z_{q,t}\equal\Z[q^{\pm 1/(n+1)},t^{\pm 1/2}]$
formed by
polynomials in terms of $q^{\pm 1/(n+1)}$ and
$t^{\pm1/2}.$ We note that the coefficients of the
Macdonald polynomials will belong to
$
\Q(q,t).
$
Generally, we have $t_{sht}$ and $t_{lng}$ (but one $q$)
and the fractional
power of $q$ is adjusted to $(v)$ below and (\ref{tauplus}).

It is convenient to use the following notation:
\begin{align*}
&t=q^k,\ \rho_k\equal \frac{k}{2}\,\sum_{\al>0}\al=
k\sum_{i=1}^n  \om_i.
\end{align*}
Generally, $\rho_k$ is defined in terms of  
$k_{sht}$ and $k_{lng}$.

For pairwise commutative $X_1,\ldots,X_n,$
\begin{align}
& X_{\tb}\ \equal\ \prod_{i=1}^nX_i^{l_i} q^{ j}
\iif \tb=[b,j],\ \hw(X_{\tb})\ =\ X_{\hw(\tb)},
\label{Xdex}\\
&\hbox{where\ } b=\sum_{i=1}^n l_i \om_i\in P,\ j \in
\frac{1}{n+1}\Z,\ \hw\in \hW.
\notag \end{align}
For instance, $X_0\equal X_{\al_0}=qX_\th^{-1}$ for maximal
{\em short} root $\th$.
\medskip

\begin{definition}
The double affine Hecke algebra $\HH\ $
is generated over $\Z_{q,t}$ by
the elements $\{ T_i,\ 0\le i\le n\}$,
pairwise commutative $\{X_b, \ b\in P\}$ satisfying
(\ref{Xdex})
and the group $\Pi,$ where the following relations are imposed:

(o)\ \  $ (T_i-t^{1/2})(T_i+t^{-1/2})\ =\
0,\ 0\ \le\ i\ \le\ n$;

(i)\ \ \ $ T_iT_jT_i...\ =\ T_jT_iT_j...,\ m_{ij}$
factors on each side;

(ii)\ \   $ \pi_rT_i\pi_r^{-1}\ =\ T_j \iif
\pi_r(\al_i)=\al_j$;

(iii)\  $T_iX_b \ =\ X_b X_{\al_i}^{-1} T_i^{-1} \iif\,
(b,\al_i)=1,\
0 \le i\le  n$;

(iv)\ $T_iX_b\ =\ X_b T_i\, $ when $\, (b,\al_i)=0
\for 0 \le i\le  n$;

(v)\ \ $\pi_rX_b \pi_r^{-1}= X_{\pi_r(b)}\, =\,
X_{ u^{-1}_r(b)}
 q^{(\om_{\iota(r)},b)},\ \, 1\le r\le n$.
\label{double}\end{definition}

Given $\tw \in \tW, 1\le r\le n,\ $ the product
\begin{align}
&T_{\pi_r\tw}\equal \pi_r T_{i_l}\cdots T_{i_1},\where
\tw=s_{i_l}\cdots s_{i_1} \for l=l(\tw),
\label{Twx}
\end{align}
does not depend on the choice of the reduced decomposition
of $\tw$.
Moreover,
\begin{align}
&T_{\hv}T_{\hw}\ =\ T_{\hv\hw}\  \hbox{ whenever\,}\
 l(\hv\hw)=l(\hv)+l(\hw) \for
\hv,\hw \in \hW. \label{TTx}
\end{align}
In particular, we arrive at the pairwise
commutative elements
\begin{align}
& Y_{b}\equal
\prod_{i=1}^nY_i^{l_i} \iif
b=\sum_{i=1}^n l_i\om_i\in P,\
Y_i\equal T_{\om_i},b\in P.
\label{Ybx}
\end{align}
When acting in the polynomial representation
(see below), they are called {\em difference
Dunkl operators}.

\subsection{\bf The automorphisms}\label{sect:Aut}
The following maps can be (uniquely) extended to
automorphisms of $\HH\,$, where
$q^{1/(2(n+1))}$ must be added to $\Z_{q,t}$
(see \cite{C101}, (3.2.10)--(3.2.15))\,:
\begin{align}\label{tauplus}
& \tau_+:\  X_b \mapsto X_b, \ T_i\mapsto T_i\, (i>0),\
\ Y_r \mapsto X_rY_r q^{-\frac{(\om_r,\om_r)}{2}}\,,
\\
& \tau_+:\ T_0\mapsto  q^{-1}\,X_\th T_0^{-1},\
\pi_r \mapsto q^{-\frac{(\om_r,\om_r)}{2}}X_r\pi_r\
(1\le r\le n),\notag\\
& \label{taumin}
\tau_-:\ Y_b \mapsto \,Y_b, \ T_i\mapsto T_i\, (i\ge 0),\
\ X_r \mapsto Y_r X_r q^\frac{(\om_r,\om_r)}{ 2},\\
&\tau_-(X_{\th})=
q T_0 X_\th^{-1} T_{s_{\th}}^{-1};\ \
\si\equal \tau_+\tau_-^{-1}\tau_+\, =\,
\tau_-^{-1}\tau_+\tau_-^{-1},\notag\\
&\si(X_b)=Y_b^{-1},\   \si(Y_b)=
T_{w_0}^{-1}X_{b^{\,\iota}}^{-1}T_{w_0},\ \si(T_i)=T_i (i>0).
\label{taux}
\end{align}
These automorphisms fix $\ t,\ q$
and their fractional powers, as well as the
following {\em anti-involution\,}:
\begin{align}
&\vph:\
X_b\mapsto Y_b^{-1},\, Y_b\mapsto X_b^{-1},\
T_i\mapsto T_i\ (1\le i\le n).\label{starphi}
\end{align}

\comment{
The following anti-involution results directly from
the group nature of the DAHA relations:
\begin{align}\label{star-conj}
H^\star= H^{-1} \for H\in \{T_{\hw},X_b, Y_b, q, t\}.
\end{align}
To be exact, it is naturally extended to the fractional
powers of $q,t$:
$$
\star:\ t^{\frac{1}{2}} \mapsto t^{-\frac{1}{2}},\
q^{\frac{1}{2(n+1)}}\mapsto  q^{-\frac{1}{2(n+1)}}.
$$
This anti-involution serves the inner product in the theory
of the DAHA polynomial representation.
}

Let us list the matrices corresponding to the automorphisms and
anti-automorphisms above upon the natural projection
onto $SL_2(\Z)$, corresponding to
$\,t^{\frac{1}{2}}=1=q^{\frac{1}{2(n+1)}}$.
The matrix {\tiny
$\begin{pmatrix} \al & \be \\ \ga & \de\\ \end{pmatrix}$}
will then represent the map $X_b\mapsto X_b^\al Y_b^\ga,
Y_b\mapsto X_b^\be Y_b^\de$ for $b\in P$. One has:
\smallskip

\centerline{
$\!\!\!\tau_+\rightsquigarrow$
{\tiny
$\begin{pmatrix}1 & 1 \\0 & 1 \\ \end{pmatrix}$},\
$\tau_-\rightsquigarrow$
{\tiny
$\begin{pmatrix}1 & 0 \\1 & 1 \\ \end{pmatrix}$},\
$\si\rightsquigarrow$
{\tiny
$\begin{pmatrix}0 & 1 \\-1 & 0 \\ \end{pmatrix}$},\
$\vph\rightsquigarrow$
{\tiny
$\begin{pmatrix}0 & -1 \\-1 & 0 \\ \end{pmatrix}$}.\
}
\smallskip

The {\em projective\,} $PSL_2(\Z)$ (due to Steinberg)
is the group generated by $\tau_{\pm}$ subject to the
relation $\tau_+\tau_-^{-1}\tau_+=
\tau_-^{-1}\tau_+\tau_-^{-1}.$ The notation will
be $PSL_{\,2}^{\wedge}(\Z)$; it is isomorphic to the braid
group $B_3$.

{\sf The coinvariant.}
The projective $PSL_2(\Z)$ and the {\em coinvariant\,},
to be defined now, are the main ingredients of our approach.

Any $H\in \HH$ can be uniquely represented in the form
$$
H=\sum_{a,w,b} c_{a,w,b}\, X_a T_{w} Y_b \for w\in W,
a,b\in P
$$
(the DAHA-PBW Theorem, see \cite{C101}). Using this
presentation, the
{\em coinvariant\,} is a unique functional $\HH\to \C\,$
defined as follows:
\begin{align}\label{evfunct}
\{\,\}_{ev}:\ X_a \ \mapsto\  q^{-(\rho_k,a)},\
Y_b \ \mapsto\  q^{(\rho_k,b)},\
T_i \ \mapsto\  t^{1/2}.
\end{align}
The main symmetry of the coinvariant is
\begin{align}\label{evsym}
&\{\,\vph(H)\,\}_{ev}\,=\,\{\,H\,\}_{ev} \hbox{\, for\, }
H\in \HH.
\end{align}
Also,
$\,\{\,\iota(H)\,\}_{ev}=\{\,H\,\}_{ev}$, where 
$\iota$ is naturally extended to $\HH$:
\begin{align}\label{iotaXY}
\iota(X_b)\!=\!X_{\iota(b)}, \
\iota(Y_b)\!=\!Y_{\iota(b)}, \
T_i^\iota\!=\!T_{\iota(i)},\ 1\le i\le n.
\end{align}

\subsection{\bf Polynomial representation}
It is isomorphic to $\Z_{q,t}[X_b]$
as a vector space with the action of $T_i(0\le i\le n)$
given by
the {\em Demazure-Lusztig operators\,}:
\begin{align}
&T_i\  = \  t^{1/2} s_i\ +\
(t^{1/2}-t^{-1/2})(X_{\al_i}-1)^{-1}(s_i-1),
\ 0\le i\le n.
\label{Demazx}
\end{align}
The elements $X_b$ become the multiplication operators
and  $\pi_r (1\le r\le n)$ act via the general formula
$\hw(X_b)=X_{\hw(b)}$ for $\hw\in \hW$. 

The {\em coinvariant} is closely related to the 
polynomial representation.
For any $H\in \HH$,
one has: $\{H T_w Yb\}_{ev}=\{H\}_{ev}\, \chi (T_w Y_b)$,
where $\chi$ is the standard character (one-dimensional
representation) of the affine Hecke algebra
$\h_Y$, generated by $T_w, Y_b$ for
$w\in W, b\in P$;\, $\chi$\,  sends\,
$Y_b\mapsto q^{(\rho_k,b)}$ and $ T_i\mapsto t^{1/2}$.
Thus, $\{\ldots\}_{ev}$ acts via the projection
$H\!\mapsto\! H\!\!\Downarrow \,\equal H(1)$ of $\HH\,$
onto the {\em polynomial representation \,}. Indeed, the
latter is
isomorphic to the $\HH$\~module induced from $\chi$
due to \cite{C101,CJ,CJJ}.

Note that $\tau_-$
naturally acts in the polynomial representation.
See formula (1.37) from \cite{CJJ}, which is based on the
identity
\begin{align}\label{taumin-poly}
\tau_-(H\!\Downarrow)\,=\,\tau_-(H)\!\Downarrow\ =\,
\Bigl(\tau_-\bigl(H\bigr)\Bigr)(1) \for H\in \HH.
\end{align}
\smallskip

{\sf Symmetric Macdonald polynomials.}
The standard notation is  $P_b(X)$ for $b\in P_+$;
see \cite{Mac1,C101}. Let me mention Kevin Kadell's key
contributions for the classical root systems and 
Rogers polynomials for $A_1$. 
The usual definition is via the orthogonality
conditions. DAHA provides the following alternative 
approach. They are such that 
\begin{align}\label{macdopers}
L_f(P_b)=f(q^{-\rho_k-b})P_b,\  L_f\equal f(X_a\mapsto Y_a)
\end{align}
for any symmetric ($W$\~invariant) polynomial
$f\in \C[X_a,a\in P]^W$, where the coefficient
of $X_b$ in $P_b$ is assumed $1$. These conditions fix them
uniquely (for generic $q,t$). The coefficients of
$P_b$ belong to the field $\Q(q,t)$, and 
$P_b(X^{-1})\,=\,P_{b^{\iota}}(X)$. The polynomials 
$P_b^\circ\equal P_b/P_b(q^{-\rho_k})$ for $b\in P_+$ 
are called 
{\em spherical Macdonald polynomials\,}.
\smallskip

\comment{
Let $c_+$ be such that $c_+\in W(c)\cap P_+$
(it is unique);
recall that $Q_+=\oplus_{i=1}^n \Z_+\al_i$.
 For $b\in P_+$, the following are the defining relations:
\begin{align*}
&P_b\! -\!\!\!\!\sum_{a\in W(b)}\!\!\! X_{a}
\in\, \oplus^{c_+\neq b}_{c_+\in b-Q_+}\,\Q(q,t) X_c
\hbox{\, and\, }
\lan \,P_b \,X_{c^{\,\iota}}\,\mu(X;q,t)\,\ran\!=\!0 \for
\\
&\hbox{all $c$ in $\oplus$ above;\ \ }
\mu(X;q,t)\!\equal\!\prod_{\al \in R_+}
\prod_{j=0}^\infty \frac{(1\!-\!X_\al q^{j})
(1\!-\!X_\al^{-1}q^{j+1})
}{
(1\!-\!X_\al t q^{j})
(1\!-\!X_\al^{-1}t^{}q^{j+1})}\,.
\end{align*}
Here and further $\lan f\ran$ is the {\em constant term\,}
of a Laurent series or polynomial $f(X)$;
$\mu$ is considered
a Laurent series of $X_b$ with
the coefficients expanded in terms of
positive powers of $q$. 
}

One has (see (3.3.23) from \cite{C101}):
\begin{align}\label{macdeval}
&P_{b}(q^{-\rho_k})=
P_{b}(q^{\rho_k})=
&q^{-(\rho_k,b)}
\prod_{\al>0}^{(\al,b)>0}\,\prod_{j=0}^{(\al,b)-1}
\Bigl(
\frac{
1- q^{j}\,t\, X_\al(q^{\rho_k})
 }{
1- q^{j}X_\al(q^{\rho_k})
}
\Bigr).
\end{align}

\subsection{\bf \texorpdfstring{{\mathversion{bold}$J$}}
{J}-polynomials}
They are necessary for managing torus iterated
{\em links\,}; spherical polynomials $P^\circ_b$
are sufficient for such {\em knots}.

For $\,b=\sum_{i=1}^n b_i \om_i\in P_+$, i.e.
for a {\em dominant\,} weight with
$\,b_i\ge 0$ for all $\,i$, the corresponding
{\em Young diagram\,} is  as
follows:
\begin{align}\label{omviavep}
&\la\!=\!\la(b)\!=\!\{\la_1\!=\!b_1\!+\!\cdots\!+\!b_n,\,
\la_2\!=\!b_2\!+\!\cdots\!+\!b_n,\,\ldots,\, \la_n\!=\!b_n\},\\
&b=\sum_{i=1}^n \la_i \vep_i-
\frac{|\la|}{n+1}\,\bigl(\vep_1+\cdots+\vep_{n+1}\bigr)
\for |\la|\equal\sum_{i=1}^n \la_i.\notag
\end{align}
One has: $(b,\vep_{i}\!-\!\vep_j)=b_i\!+\!\cdots\!+\!b_{j\!-\!1}=
\la_i\!-\!\la_j$. Also, $(b,\rho)=n|\la|/2-
{\sum_{i=1}^{n-1}i\la_{i+1}}$, and
$b^2\equal(b,b)=\sum_{i=1}^n \la_i^2-|\la|^2/(n\!+\!1)$.


Let us provide the  set of all $[\al,j]$ in the
product from (\ref{macdeval}); it is
\begin{align*}
&\{\,[\al,j],\  \al=\vep_l\!-\!\vep_{m}\in R_+,\, j>0\, \mid\,
b_l+\cdots+b_{m-1}>j> 0\}.
\end{align*}

The {\em $J$\~polynomials\,}
are as follows:
\begin{align}\label{P-arms-legs}
J_\la\equal h_\la P_b \for \la=\la(b),\,
h_\la=\prod_{\Box\in\la}
(1-q^{arm(\Box)}t^{leg(\Box)+1});
\end{align}
they are $q,t$\~integral.
Here $arm(\Box)$ is the {\em arm number\,}, which is
the number of boxes in the same row as $\Box$
strictly after it; $leg(\Box)$ is the {\em leg number\,},
which is the number of boxes in the column of $\Box$
strictly below it. Here $\la=(\la_i),$\, 
$\la_1\!\ge \la_2\!\ge\!\ldots\!\ge\!\la_{n-1}\!\ge \la_n\,$,
which are the numbers of
boxes in the corresponding rows; the $i${\tiny th} row
is assumed above the  $(i+1)${\tiny th}. 

The following formula is equivalent to (\ref{P-arms-legs}):
\begin{align}\label{j-polynom}
&J_\la= t^{-(\rho,b)}
\prod_{p=1}^n \prod_{j=0}^{\la_{p^*\,}-1}
\Bigl( 1\!-\!q^{j} t^{\,p+1}\Bigr) P_b^\circ, \ \,
p^*= n\!-\!p\!+\!1,\, b\in P_+.\
\end{align}
See, for instance, \cite{ChD2} and Theorem 2.1 from \cite{GoN}.
Note that the {\em arms and\, legs\,} 
are not needed in this presentation (in terms of
$P_b^\circ$). In this form, 
counterparts of $J$\~polynomials can be defined for any 
reduced root systems; see \cite{ChD2}, Section 2.6.
Note that 
$P(t^{-\rho})=P(t^{\rho})$ 
for any symmetric polynomials $P$ due to $w_0(\rho)=-\rho$.

The following the simplest instance of the 
{\em $a$-stabilization\,}; $J_\la(t^{-\rho})$ for $A_n$ is the value
of the formula in (\ref{stabevalx}) for  $a\!=\!-t^{n+1}$:\,
\begin{align}\label{stabevalx}
J_\la(t^{-\rho})\!=\!
(a^2)^{-\frac{|\la|}{4}}\,t^
{\sum_{i=0}^{n-1}(i+1/2)\la_{i+1}}\,
&\prod_{p=1}^n \prod_{j=0}^{\la_{p\,}-1}
\Bigl( 1\!+\!q^{j}\, a\, t^{-p+1}\Bigr).
\end{align}

The general $a$-stabilization  in the definition 
of superpolynomials below can be
reduced to the following proposition from 
\cite{ChD2}. 

\begin{proposition}\label{PROP:stab-values}
Given two Young diagrams $\la$ and $\mu$,\,
the values $P_\la(q^{\mu+\rho_k})$ are $a$\~stable.
Namely,  
there exists a universal rational function in terms
of $\,q,t,a\,$ times a certain  
power of $a^{1/2}$ and power of $t^{1/2}$ 
such that 
 $P_b(q^{c+\rho_k})$ for $A_n$ is its value 
at  $a=-t^{n+1}$.
Here $n\ge 0,\, \la=\la(b),\mu=\la(c)$ are the Young
diagrams for  $b,c\in P_+$ for $A_n$,
$n+1$ is no smaller 
than the number of rows in $\la$ and in $\mu$,
 and $P_{m\om_{n+1}}=1$ in $A_{n}$ by definition. 
Also, the $a$-stabilization holds for 
$P_\la^\circ(q^{\mu+\rho_k})$. \sq
\end{proposition}

\section{\sc DAHA-Jones theory}
\subsection{\bf Iterated torus knots}\label{sec:ITER-KNOTS}
The torus knots $T(\rr,\ss)$ are defined for any integers
assuming that \,gcd$(\rr,\ss)=1$. In the standard torus $T^2$
it is a path with $\rr$ horizontal turns (along the torus)
and $\ss$ vertical ones.  
There are obvious symmetries:
 $\,T(\rr,\ss)=T(\ss,\rr)=T(-\rr,-\ss)$, where 
we use ``$=$" for the {\em ambient isotopy equivalence}.
Also, $\,T(\rr,\ss)=\unknot\,$\, if $\rr=1$ or $\ss= 1$, since
this torus knot is 
{\em unknot\,} $\unknot\,$. We mention that
changing $\rr\mapsto -\rr$
results in the {\em mirroring} of $T(\rr,\ss)$. The DAHA
construction is for any torus knots, but the motivic theory is
only for {\em algebraic torus
knots}, which are for positive $\rr,\ss$. Here and below see
\cite{RJ,EN,ChD1}.
\vskip 0.2cm

 Following  \cite{ChD1}, 
the DAHA construction of invariants of iterated torus knots
is based one their {\em $\tax$-presentations}. Algebraically, such
a presentation is an arbitrary pair of sequences of
integers:
\begin{align}\label{iterrss}
\vec{\rr}=\{\rr_1,\ldots \rr_\ell\}, \
\vec{\ss}=\{\ss_1,\ldots \ss_\ell\} \hbox{\, such that\,
gcd}(\rr_i,\ss_i)=1,
\end{align}
where $\ell$ will be called the {\em length.}
For algebraic knots, $[\rr_i,\ss_i]$ will 
be interpreted below as 
{\em characteristic\,} or {\em Newton pairs\,}
for the corresponding plane curve singularities. They  
are not topological invariants. The following sequences are.

\vskip 0.2cm
{\sf Cabling parameters.} 
The {\em cabling\,} $C\!ab(\aa,\bb)(K)$ of 
any oriented
knot $K$ in (oriented) $S^3$ is defined as follows;
see e.g., \cite{EN} and references therein.
We consider a small $2$\~dimensional torus
around $K$ and put there the torus knot $T(\bb,\aa)$
in the direction of $K$,
which is $C\!ab(\aa,\bb)(K)$. The parameter
$\aa$ is the number of
turns around $K$, $\bb$ is that along $K$.
For instance, $C\!ab(a,0)K=\unknot$ and $Cab(a,1)K=K$.

This procedure depends on
the order of numbers $\aa$ and $\bb$, the orientation of $K$,
and its {\em framing}. We take the natural
zero-framing  when  the parallel copy of the knot has
zero linking number with a given knot.

Using this operation, iterated torus knots (also called
satellite knots)
are as follows.
Given two sequences
$\vec\aa=(\aa_1,\ldots,\aa_\ell),\vec\rr=(\rr_1,\ldots,\rr_{\ell})\}$, 
called a {\em cable presentation}, 
we begin with  $T(\rr_1,\aa_1)=C\!ab(\aa_1,\rr_1)(\unknot)$
(note that we transpose $\aa$ and $\rr$ here), and then set:
\begin{align}\label{Knotsiter} 
\t(\vec\rr,\vec\ss)=C\!ab(\vec{\aa},\!\vec{\rr})(\unknot)=
\Bigl(C\!ab(\aa_\ell,\rr_\ell)\cdots \bigl(C\!ab(\aa_2,\rr_2)
(T(\rr_1,\!\aa_1))\bigr)\Bigr),
\end{align}
where the connection with the $\tax$-presentation is as follows:
\begin{align}\label{Newtonpair}
\aa_1=\ss_1,\,\aa_{i}=\aa_{i-1}\rr_{i-1}\rr_{i}+\ss_{i}\,\
(i=2,\ldots,\ell).
\end{align}
See  e.g., \cite{EN}. The $\tax$-presentation will be the
one used later. However, to prove the topological invariance
of the DAHA construction (a theorem) the passage to the cable presentation
is necessary. 

\Yboxdim7pt
Recall that knots and links are considered up to
{\em ambient isotopy\,}. 
Generally, they are colored by dominant weights, which  
will be Young diagrams in what will follow; uncolored links
are those for $\yng(1)\,$. 

\subsection{\bf From knots to links}\label{sec:knots-links}
Generally, {\em iterated torus links\,} are 
disjoint unions of {\em colored} iterated torus knots described by
the corresponding {\em splice diagrams}. {\em Cabling} 
can be naturally extended to torus iterated links: we
begin then with an oriented {\em unlink} with the standard framing. 
Torus knots $T(\rr,\ss)$ can be now for \, gcd\,$(\rr,\ss)>1$.
For algebraic links, the main subject of this paper, 
the canonical orientations of their connected components 
is when all pairwise
linking numbers are positive. 

The data now will be 
unions of the sequences of \tax-parameters
for  iterated torus knots subject to
the following {\em incidence conditions}:
\begin{align}\label{tau-link}
&\l_{(\vec\rr^{\,j},\vec\ss^{\,j})}^{\,\Up,\, (b^j)}=
\Bigl(\{\t(\vec\rr^{j},\vec\ss^{j}),\, b^j\in P_+\},
j=1,\ldots,\kap\Bigr)\,\&\,
\Up=(\up_{j,k})\\ 
&\text{such that }
0\le \up_{j,k} \le \min\{\ell^j,\ell^{k}\},\, \up_{j,k}=\up_{k,j},
1\le j,k \le \kap, \notag\\
&\text{and we assume that}\ \, [\rr_{i}^j,\ss_{i}^j]=
[\rr_{i}^k,\ss_{i}^k]\hbox{\,\, for all\,\, }
1\le i \le \up_{j,k}.
\notag
\end{align}
Here $\ell^j$ is the length of $\vec\rr^j=\{\rr_i^j\}$ and
$\vec\ss^j=\{\ss_i^j\}$ for $1\le j\le \kap$; we
naturally set $\up_{j,j}=\ell^j$. 

Equivalently, the incidence matrix $\Up$ determines a 
\tax-labeled 
graph with the {\em vertices}
$v=(i,j)$ identified as in (\ref{tau-link}) and supplied 
with the \tax-labels.
The {\em paths\,} are sequences
of vertices with fixed $j$ and increasing consecutive $i$.
The {\em edges} are from $i$ to $i+1$ in
any paths; so this graph is oriented naturally. We see that 
$\Up$ is a disjoint
union of its {\em maximal subtrees}. Each of them contains a unique 
{\em initial vertex\,}, which is $v_1^j=(i=1,j)$ for any $j$ in
this subtree.  The  $i$\~index of $v$ is 
then $1$ plus the distance to $v$ from the initial vertex
in the corresponding maximal subtree. We conclude that the graph
$\Up$ is a disjoint union of 
 \tax-labeled trees that are bunches of paths from the same 
initial vertex (which may ramify).

\comment{ 
For $i\le \ell^j$,
the pairs $[\rr_i^j,\ss_i^j]$ become labels, called
{\em \tax-labels\,} of the 
{\em vertices\,} $(i,j)$ of $\l$; the 
square brackets will be used for them.
}
\vskip 0.2cm

Accordingly, $\l$ will be the graph $\Up$ where 
we add one or several {\em arrowheads\,}
at the end of {\em every\,} path. The ends  are the vertices 
$(i=\ell^j,j)$. The colors $b^j$ are assigned to
the arrowheads. 
Topologically, the $j${\tiny th} {\em maximal path} is  interpreted as
the knot $\t(\vec\rr^{j},\vec\ss^{j})$
colored by $b^j\in P_+$ at its arrowhead.
Later, we will switch to 
Young diagrams $\la^j=\la(b^j)$ from dominant weights $b^j$.

The maximal paths can totally coincide in $\Up$, 
but have
different arrows; then their $j$\~indices will be different and
they will be treated as different paths.
If $\Up$ contains no vertices, then $\l$  will
be a collection of arrowheads; the  corresponding link is 
the {\em unlink} colored by $\{b^j\}$.


\vfil
The $\aa$-parameters are calculated
along the maximal paths exactly as we did for the knots, 
starting with $i=1, \aa_{1}=\ss_{1}$ for all $j$. 
The corresponding 
$\aa_i^j$ depend only on the corresponding vertices $v$ in $\Up$.
The pairs
$\{\aa_i^j,\rr_i^j\}$
will be called the {\em cab-labels} of the vertices.

Only the \tax-labels and the colors are needed in the DAHA
construction below. The topological invariance of our
DAHA construction includes that it actually depends only
the {\em cab-labels}.
We mostly omit  topological aspects
in this paper. See \cite{ChD2,ChW} for details and (many) examples.
For instance, the {\em pairs of graphs} below can be generally 
reduced (isotopically) to one graph, but this can be ``unnatural".
\comment{
The torus knot colored by $b\in P_+$ (or by the corresponding
$\la$) is denoted by $T_{\rr,\ss}^b$;
$C\!ab_{\aa,\rr}^b(\mathbf L)$,
equivalently $C\!ab_{0,1}^b Cab_{\aa,\rr}(\mathbf L)$, is
the cable $C\!ab(\aa,\rr)(\mathbf L)$
of a link $\mathbf L$ colored by $b$.
The color is attached to the last $C\!ab$ in
the sequence of cables.
In the absence of vertices, the notation is $\unknot^{\,b}$
(the unknot colored by $b\in P_+$) or $C\!ab(0,1)^b$. {\em
We will use the same notation $\l$ for the graph and the
corresponding link $\mathbf{L}$.}
}


\comment{
\begin{align}\label{Knotsiterx}
\l\bigl(\vec\rr^{j},\vec\ss^{j},\, 1\le j\le \kap\bigr)
\rightsquigarrow
\Bigl(\coprod_{j=1}^\kap C\!ab(\vec{\aa}^j,\!\vec{\rr}^j)
\Bigr)(\unknot),
\end{align}
where the composition and coproduct of cables is with respect to
the graph structure and $C\!ab(\vec{\aa}^j,\!\vec{\rr}^j)=
\cdots C\!ab(\aa^j_2,\rr^j_2)T(\rr^j_1,\ss^j_1)$ is as in
(\ref{Knotsiter}).
In this work, the coproduct symbol $\,\hbox{\small$\coprod$}\,
which stands for the union of cables, will be simply replaced
by comma; we set $\bigl(C\!ab(\aa,\rr),C\!ab(\aa',\rr')\bigr)$
instead of $C\!ab(\aa,\rr)\coprod C\!ab(\aa',\rr')$.
The color is attached to the last $C\!ab$ in
the sequence of cables.

The $\aa$\~parameters and the corresponding cables
are calculated as above along the
corresponding paths.
See \cite{ChD2,ChW} for details and (many) examples.
We do not need much the topological aspects
in this work; the graphs are sufficient for the DAHA
construction.
}

\vfil
{\sf Pairs of graphs.}
Let $\{\mathcal{L},\,'\!\mathcal{L}\}$ be a pair of labeled
graphs with colored arrows as above.
The cabling construction provides a canonical
embedding of the iterated torus links into
the solid torus. 
Their {\em twisted union\,} is as follows:
we put the link for $\l$ into the horizontal solid torus
and $'\!\mathcal{L}$ 
into the complementary vertical one. 

Since this paper is mostly about {\em algebraic\,} links, we will
always change the natural orientation
of the second link (if not empty) by the opposite one.
Without this switch, the resulting union
is never algebraic. The corresponding twisted union will be
denoted by $\{\mathcal{L},\,'\!\mathcal{L}\}$. 

This would be 
$\{\mathcal{L},\,'\!\mathcal{L}^\vee\}$ in the
notation from \cite{ChD2}; only such pairs 
will be considered below 
(we omit $\vee$ in this paper). Note that  
$\{\mathcal{L},\,'\!\mathcal{L}\}$ and 
$\{'\!\mathcal{L},\,\mathcal{L}\}$ result in isotopic links,
which corresponds to formula (4.24) from \cite{ChD2} in the 
DAHA setting.

For instance,
\Yboxdim7pt
$\bigl\{\mathcal{L}=\{\circ\rightarrow \yng(1)\},\,
'\!\mathcal{L}=\{\circ\rightarrow\yng(1)\}\bigr\}$
represents in this paper the uncolored Hopf
$2$\~links with the linking number $lk=+1$, which
was denoted by $\{\circ\rightarrow\yng(1),
\circ\rightarrow\yng(1)^\vee\}$
in \cite{ChD2}. It is an algebraic link corresponding 
to the singularity $xy=0$ (see below). 
\vskip 0.2cm

The theory in \cite{EN} is without colors, as well as
that in \cite{ObS,ORS} (for algebraic links).
Generally, attaching colors to HOMFLY-PT polynomials 
can be incorporated using proper
framed links; {\em the skein} is used.
This becomes complicated for the triply-graded HOMFLY-PT
homology. We need the {\em $q,t$-skein} of the torus $T^2$, which
is essentially spherical DAHA. 
{\em Colors} will be arbitrary in the DAHA construction.

Note that the transposition of $\rr_1$ and
$\ss_1$ (only for the first pair!) does not change
the isotopy type of the corresponding {\em knot}, but 
may influence the resulting {\em twisted union}. 
To give an example, let $'\!\l$ be one pure arrow (no vertices)
colored by $\la$ and $\l$ an arbitrary graph as above. Then 
the twisted union for the pair  $\{\l,\,'\!\l\}$ 
corresponds to adding {\em the meridian} colored by $\la$
to $\l$. The linking number of the meridian with $\l$
may change when $\rr_1,\ss_1$
are transposed.
\vskip 0.2cm

{\sf Algebraic links.}
We will provide here only some basic facts; 
see \cite{EN} for details and
references, especially Theorem 9.4 there.

Arbitrary algebraic links can be described using
the cabling construction above for {\em positive pairs of trees\,}
$\{\l,\,'\!\l\}$;  only trees occur.
The positivity means the positivity of all
$\rr,\ss$ and the inequalities
$\,'\ss_{1} \ss_{1}>\,'\rr_1 \rr_1$ for the labels $[\rr_1,\ss_1]$
and  $['\rr_1,'\ss_1]$ of the {\em initial vertices} of the trees 
$\l$ and $'\!\l\}$. This claim is upon the 
usage of the symmetries
of the corresponding {\em splice diagram\,}. 
See \cite{EN,New} and \cite{ChD2} for details.

Such links are disjoint unions of algebraic iterated torus {\em knots} 
associated with 
the maximal paths in $\l$ and $'\!\l$ for the indices
$1\le j\le \kap$ and
$1\le j\le\, '\!\kap$. These components and the corresponding 
pairwise linking numbers uniquely determine the corresponding
link due to the {\em Reeve Theorem}. The linking number between the 
knots for the paths with the indices $1\le j<k\le \kap$ in $\l$ 
(or those in $'\!\l$) is
\begin{align}\label{linkjk}
lk(j,k)=\aa^{j}_{i_o}\rr_{i_o}^j\,
\Bigl(\,\prod_{i=i_o+1}^{\ell^j} \rr^j_i\Bigr)
\Bigl(\,\prod_{i=i_o+1}^{\ell^k} \rr^k_i\Bigr),
\where i_\circ=\up(j,k).
\end{align}
The linking number between knot $j$ from $\l$ and knot $j'$ in
$'\!\l$ are with full
products of $\rr_i^j$ and $'\rr_i^{j'}$ over $i$ in this formula.
\smallskip

Algebraically, isolated plane curve singularity are considered,
presented by polynomial equations $F(x,y)=0$ in a neighborhood 
of $0=(x=0,y=0)$. Their intersections
with a small $3$-dimensional sphere in
$\C^2$ around $0$ are the corresponding {\em algebraic link\,}.

Given a \tax-labeled positive tree
 $\l=\l^\Up_{(\vec\rr^j,\vec\ss^j)}$, the 
irreducible components of the singularity are
associated with its maximal paths 
(numbered by $j$) and are given by the equations:
\begin{align}\label{yxcurve}
y = c^j_1\,x^{\ss^j_1/\rr^j_1}
(1+c^j_2\,x^{\ss^j_2/(\rr^j_1\rr^j_2)}
\bigl(1+c^j_3\,
x^{\ss^j_3/(\rr^j_1\rr^j_2\rr^j_3)}
\Bigl(\ldots\Bigr)\bigr)) \hbox{\, at\, } 0.
\end{align}
Here the numbers $\rr,\ss$ are from the corresponding
labels; the parameters $c_i^j\in \C$ must be
sufficiently general.

The simplest example is the equation
$y^{\kap\rr}= x^{\kap\ss}$ for $\rr,\ss>0$ such that
 gcd$(\rr,\ss)=1$. Its link is $T(\kap\rr,\kap\ss)$:\,
$\kap$ copies of  the torus knot $T(\rr,\ss)$ with the 
coinciding pairwise linking numbers, which are $\rr\,\ss$.

The unibranch components and their
(pairwise) linking numbers uniquely determine the topological
type of the corresponding
germ due to the Reeve Theorem; see e.g., \cite{EN}.
{\em All linking numbers are strictly positive for
algebraic links.} If $\l$ or $'\l$ is a disjoint union
of {\em several} trees, then the corresponding 
iterated torus link is {\em non-algebraic}; they must be trees (not 
their disjoint unions)
for algebraic links.

\subsection{\bf DAHA-Jones invariants}
They, the $J\!D$-polynomials,  are uniformly defined for any (reduced,
irreducible) root systems $R$. To be more exact, this is for 
twisted affine root systems  $\tilde{R}$. We will 
need only $R=A_{n}$ to construct {\em DAHA superpolynomials}.

The {\em combinatorial data} will be the \tax-labeled 
colored graphs
$\l_{(\vec\rr^{j},\vec\ss^{j})}^{\,\Up, (b^j)}$
from (\ref{tau-link}) and pairs $\{\l,\,'\!\l\}$ of such graphs
(the second can be empty).  Recall that
$
1\le j\le \kap,\,\, \vec\rr^{\,j}=\{\rr_i^j\},\,\,
\vec\ss^{\,j}=\{\ss_i^j\},\,\, 1\le i\le \ell^j,
$
for $\l$ (similarly for $'\!\l$),
and the {\em arrowheads\,} are added at the ends of {\em paths\,},
which are colored by $b^j\in P_+$; see (\ref{tau-link}).
Generally, these graphs are  
disjoint unions
of oriented labeled subtrees with the set of unique initial vertices
(one for a subtree)
and ``colored" arrows at their ends. We will do in this section any iterated
torus links, not only algebraic.


For the latter reference,
let $\la^j=\la(b^j)$ for dominant $b^j$. We set:
\begin{align}\label{P-arms-legs-LCM}
&ev(b^1,\ldots,b^\kap)=ev(\la^1,\ldots,\la^\kap)=
LC\!M\bigl(J_{\la^1}(t^{\rho}),\ldots,
J_{\la^\kap}(t^{\rho})\bigr),
\end{align}
where $LC\!M$ is normalized by the condition that it is a
$q,t$\~polynomial with the constant term $1$.

One has the following combinatorially transparent formula:
\begin{align}\label{P-arms-legs-union}
&ev(\la^1,\ldots,\la^\kap)\ =\
ev(\la^1\!\vee\cdots\vee\la^\kap)\,, \hbox{\,\, where}\\
&\text{the diagram\,\,} \la^1\!\vee\cdots\vee\la^\kap
\hbox{\, is the  union of } \{\la^j\}.\notag
\end{align}

We note that the $J$\~polynomials in the $A_n$\~case
are not {\em minimal\,} integral 
forms of the $P$-polynomials,
those proportional to $P_\la$ but without $q,t$\~denominators.
They are not minimal even  
for $t=q$, when the $P$-polynomials become monomial symmetric
polynomials, obviously ``minimal". The switch from $P$ to $J$ 
is a must for invariants of iterated torus {\em links}. 
Algebraically, they
are needed for the $a$-stabilization of Macdonald polynomials.
See \cite{GoN},\cite{ChD2}; the latter reference contains
their theory for any reduced root systems $R$.
\smallskip

{\sf Pre-polynomials.}
Recall that $H\!\!\Downarrow\, \equal H(1)$, where the
action of $H\in \HH$ in the polynomial representation
 is used.
We represent torus knots $T(\rr,\ss)$
by the matrices $\ga[\rr,\ss]=\ga_{\rr,\ss}\in PSL_{\,2}(\Z)$
with the
first column $(\rr,\ss)^{tr}$ ($tr$ is the transposition)
for $\,\rr,\ss\in \N$ such that \,gcd$(\rr,\ss)=1$.
Let $\hat{\ga}_{\rr,\ss}\in PSL_{\,2}^{\wedge}(\Z)$ be
{\em any\,} pullback of $\ga_{\rr,\ss}$ to 
$PSL_{\,2}(\Z)$; see Section \ref{sect:Aut}.

{\em Pre-polynomials}  are defined
for any pair $\{\l,\,'\!\l\}$ of labeled graphs with colored
arrowheads from  
(\ref{tau-link}). The {\em positivity\,} of $\l,\,'\!\l$
is not needed for their definition and the $J\!D$-polynomials 
below. 
\comment{
\begin{align}\label{links-x-y}
\l=\l_{(\vec\rr^j,\vec\ss^j)}^{\Up,(b^j)},\
'\!\l=\, '\!\l_{('\vec\rr^j,'\vec\ss^j)}^{\,'\Up,('\!b^j)}
\where b^j,\,'\!b^j\in P_+,
\end{align}
$$
1\le j\le \kap,\,'\!\kap\for \l,\, '\!\l,\
\vec\rr^j= (\rr_i^j \mid 1\le i\le \ell^j),\,
'\vec\rr^j=('\!\rr_i^j \mid 1\le i\le\, '\!\ell^j).
$$
}
\smallskip

For any labels of $\l,'\!\l$, we will lift 
$(\rr,\ss)^{tr}$ and 
 $('\!\rr,'\!\ss)^{tr}$\, to
the corresponding 
$\hat{\ga}$ and $'\hat{\ga}\in PSL_{\,2}^{\wedge}(\Z)$.
Let us begin with the definition of
{\em prepolynomials} for  $\l$. 
We set $\tilde{\p}^j=J_{b^j}$ for $i=\ell^j$
for the corresponding colored arrowheads, treated below as
additional edges of $\Up$.

Given a vertex $v$ with the label  $[\rr,\ss]$,
$\tilde{\p}_v=
\biggl(\hga\,\bigl(\prod_{v\to v'}
\tilde{\p}_{v'}\bigr)\biggr)\!\!\Downarrow$, by induction, 
for $\hga$ associated  with $[\rr,\ss]$, and
$v\to v'$ is when either $i'=i+1$ in the corresponding 
path, or $v'$ is  an arrowhead from $v$.
The product is over all such $v'$ (including the
 arrowheads from it if any).

The pre-polynomials $\tilde{\p}_{s}$ for
maximal subtrees $\Up^s$ of $\Up$ are those for $v_s=v_1^j$, the
corresponding initial vertices (serving all $j$ in this subtree). Finally,  
$\tilde{\p}_{\l}\equal\prod_s \tilde{\p}_{s}$. The definition for
pairs $\{\l,'\!\l\}$ is as follows. 

Let $\tilde{\q}$ be the pre-polynomial for $'\!\l$, the one 
corresponding to the
graph $'\Up$ and the sequence $\{'\la^j=\la('b^j)\}$ for
its arrowheads.  
Then the pre-polynomial of $\{\l,\,'\!\l\}$ is 
$\tilde{\p}_{\l,'\!\l}\equal
\bigl(\tilde{\q}(Y) \tilde{\p}(X)\bigr)\!\!\Downarrow$.
We substitute $X\mapsto Y$ in $\tilde{\q}$; 
using the DAHA automorphisms, 
$\tilde{\q}(Y)=\vph\hbox{\tiny $\circ$\,}
\iota\, (\tilde{\q})$.
\smallskip

{\sf DAHA-Jones construction.}
Let $\mathbf{b}=(b^j),\,\mathbf{'b}=('b^j)$;
(\ref{P-arms-legs-union}) will be used.
The {\em DAHA-Jones polynomial\,} 
are as follows:
\begin{align}\label{jones-bar}
&J\!D_{\l,\,'\!\l}=
J\!D_{\,(\,\vec\rr^j,\,\vec\ss^j\,)
\, , (\,'\vec\rr^j,\,'\vec\ss^j\,)}
(\mathbf{b},\mathbf{'b}\, ; \,q,t)\,\equal\,
\frac{\bigl\{\tilde{\p}_{\l,\,'\!\l}\bigr\}_{ev} }
{ev(\mathbf{b},\mathbf{'b})}.
\end{align}
Recall that $\{\p\}_{ev}=\p(t^{-\rho})=\p(t^{\rho})$ 
for symmetric $\p$.

In the case of iterated torus {\em knots\,},
when $\l$ is one path with one $b$ (and there is no
$'\!\l$), 
we arrive at formula (2.12) from \cite{ChD1}:
\begin{align}\label{jones-ditx}
& J\!D_{\vec\rr,\vec\ss}
(b;q,t)\! =\!
\Bigl\{\hat{\ga_{1}}\Bigr(
\cdots\Bigl(\hat{\ga}_{\ell-1}
\Bigl(\bigl(\hat{\ga}_\ell(P_b)/
P_b(t^{-\rho})\bigr)\!\Downarrow
\Bigr)\!\Downarrow\Bigr) \cdots\Bigr)\Bigr\}_{ev}.
\end{align}
Only $P_b^\circ=P_b/P_b(t^{-\rho})$ are needed in this case
(no $J$-polynomials). 

The simplest link is for the union 
of any number of arrowheads
colored by $\la^1,\ldots, \la^\kap$. Then
$J\!D=\prod_{j=1}^\kap J_{\la^j}(t^{\rho})/
J_{\la^1\vee\cdots\vee \la^\kap}(t^\rho)$.

\smallskip

{\sf Polynomiality.}
The following theorem and related statements are
mostly from \cite{ChD1} and \cite{ChD2}.

\begin{theorem}\label{THM-integr-Jones}
The invariants  $J\!D_{\l,\,'\!\l}\,$ 
are polynomials in terms of $q,t$
up to a factor $q^\bullet t^\bullet$,
where the powers $\bullet$ can be rational.
Modulo such factors,
it does not depend
on the particular choice of the lifts \,$\ga\in PSL_2(\Z)$
and $\hat{\ga}\in PSL_{\,2}^{\wedge}(\Z)$\ for the vertices
of $\Up,'\Up$. 

Up to the $q^\bullet t^\bullet$\~equivalence
and a possible factor from $\Q$, the following
{\em\sf hat-normalization} $\hat{J\!D}_{\l,\,'\!\l}\,$ exists.
It is a $q,t$\~polynomial not divisible by
$\,q\,$ and by $\,t\,$ with integral
coefficients of its $q,t$\~monomials such that 
their total $GCD$ equals $1$, and  
the coefficient of the minimal pure
power of $t$ must be positive. For algebraic links colored
by rows $m\om_1$, this normalization becomes
$\hat{J\!D}_{\l,\,'\!\l}\,(q\!=\!0,t\!=\!0)=1$. \sq
\end{theorem}


{\sf Topological symmetries.}
The polynomial $\hat{J\!D}_{\l,\,'\!\l}$ defined in
Theorem \ref{THM-integr-Jones}\,
depends
only on the topological type of
the link corresponding to the pair of graphs
$\{\l,\,'\!\l\}$.
For instance, the
vertices with $\rr=1$ can be omitted in 
$\l$ and $'\!\l$ and the transposition
$[\rr_1^j,\ss_1^j]\mapsto $ $[\ss_1^j,\rr_1^j]$  (only
for $i=1$ and in the absence of $'\!\l$) does not influence
$\hat{J\!D}_{\l}$. Also, the pairs
$\{\l,'\!\l\}$ and $\{'\!\l,\l\}$ result in coinciding 
polynomials.
\vskip 0.2cm

The justification of such symmetries 
for torus knots is in  Theorem 1.2 from
\cite{CJJ}. It suffices to check that 
DAHA-Jones polynomials coincide for $T(\rr,\ss)$,  
$T(\ss,\rr)$ and $T(-\rr,-\ss)$, and they
are trivial for $T(\rr,1)$. 

Concerning torus iterated {\em knots}, it
suffices to check that the $J\!D$\~polynomials for
$C\!ab(m\rr+\ss,\rr)T(m,1)$ and $T(m\rr+\ss,\rr)$ 
must be the same since $T(m,1)$. The corresponding
DAHA fact is the commutativity of  $\tau_-^m$
with $\Downarrow$, which simply means that $\tau_-$ acts
in the polynomial representation. 
Formula (\ref{tau-minus}) below
is sufficient to extend the topological invariance  to iterated
torus {\em links}. See \cite{ChD2} for detail. 

We note that 
 $\hat{\ga}\bigl(\tilde{\p}\bigr)\in \HH$,
``knot operators", 
and their products 
are invariants of the corresponding iterated torus 
links considered in $T^2$, in
$T^2\times [0,\ep]$ for small $\ep>0$ to be exact. One level down,
the pre-polynomials themselves 
are invariants of such links considered 
in the horizontal solid torus. Also, the substitution $X\mapsto Y$
is the passage to the vertical solid torus. Combining them
and applying the coinvariant $\{\cdots\}_{ev}$, 
the $J\!D$-polynomials 
$\bigl\{\tilde{\q}(Y)\tilde{\p}(X)\bigr\}_{ev}$  
become invariants of the corresponding
twisted products in $\S^3$. The {\em coinvariant} 
plays here the role of ``trace" in parallel 
topological considerations. 

\smallskip

{\sf Specialization $q\!=\!1$.}
In the formula for $J\!D_{\l,\,'\!\l}(q\!=\!1)$,
let the number of maximal paths in $'\!\l$ be
$'\!\kap=\nu$ (for the readability); $ev(b^1,\ldots,b^\kap)$ from
(\ref{P-arms-legs-LCM}) will be used.  Then for generic $t$,
\begin{align}\label{q-1-prod}
&\frac{ev(b^1,\ldots,b^\kap,\,'b^1,\ldots,
'\!b^\nu)}
{ev(b^1)\cdots ev(b^\kap)
\,ev('b^1)\cdots
ev('b^{\,\nu})}(q\!=\!1)\,\,
J\!D_{\l,\,'\!\l}(q\!=\!1)\\
=\ \,&
\prod_{j=1}^{\kap} J\!D_{\,\vec\rr^j,\,\vec\ss^j}
\,\bigl(\,b^j\,;\,q\!=\!1,\,t)\
\prod_{j=1}^{\nu} J\!D_{\,'\vec\rr^j,\,'\vec\ss^j}
\,\bigl(\,'b^j\,;\,q\!=\!1,\,t),\notag\\
&\hbox{where\,\,}
J\!D_{\,\vec\rr,\,\vec\ss}
\,\bigl(b;\,q\!=\!1,t
\bigr)\!=\!
\hbox{\small$\prod$}_{p=1}^n J\!D_{\,\vec\rr,\,\vec\ss}\,
(\om_p;\,q\!=\!1,t)^{\be_p}\notag
\end{align}
for $b=\hbox{\small$\sum$}_{p=1}^n \be_p \om_p\in P_+$.
The $J\!D$\~polynomials (without the hat-normalization) are used;
see (\ref{jones-ditx}) above and formula (2.18)
in \cite{ChD1}.

Let us comment on the factor in terms of $ev(\cdots)$ 
in the left-hand side. If $J_\la$ are
replaced in the definition of  $J\!D_{\l,\,'\!\l}$ by
the spherical polynomials $P_b^\circ$ and no division by
$LC\!M$ is performed, then this factor becomes $1$ and 
$J\!D_{\l,\,'\!\l}(q\!=\!1)$ becomes a (pure)  product 
from the theorem. In particular, this factor is $1$
for iterated torus  {\em knots}.

The same formula at $q\!=\!1$ holds for superpolynomials
$\hat{\h}$ defined below. For instance, 
$\hat{\h}_K(b;\,q\!=\!1,t,a)\!=\!
\hbox{\small$\prod$}_p \hat{\h}_K
(\om_p;\,q\!=\!1,t,a)^{\be_p}$ for iterated torus knots $K$, where
$\be_p$ is the number of $p$-columns in $\la(b)$.

\subsection{\bf DAHA superpolynomials}\label{sec:sup-def}
Following \cite{CJ,GoN,CJJ,ChD1,ChD2}, the construction from
Theorem \ref{THM-integr-Jones} and other statements above can
be extended to the {\em DAHA-superpolynomials\,},
the result of the stabilization  of
$\hat{J\!D}^{A_n}_{\l,\,'\!\l}$, which are
$J\!D$-polynomials for $A_n$ under the {\em hat-normalization.}

The $a$\~stabilization for torus knots
was announced and outlined in \cite{CJ}; its proof was published in
\cite{GoN}.
Both approaches use \cite{SV}.
The superduality conjecture was proposed in \cite{CJ}
(let us also mention \cite{GS})
and proven in \cite{GoN} for torus knots; see also
\cite{CJJ} for an alternative approach based on the
generalized level-rank duality. The justifications
of the $a$\~stabilization and
the superduality was extended to
arbitrary iterated torus knots in \cite{ChD1},
and to links in \cite{ChD2}.

In contrast to knots, the 
polynomiality of the superpolynomials for links
requires the usage of  $J_\la$\~polynomials.  Without the latter,
the superpolynomials for  (colored) {\em links\,}, 
may have non-trivial $t$\~denominators if $P$\~polynomials
are used. Also, the usage of  $J_\la$\~polynomials is
important (for knots too) to justify 
the stabilization and superduality; see \cite{GoN, ChD1,ChD2}. 
However, the definition of the 
$J\!D$\~polynomials and $\hat{\h}$-polynomials for knots
requires only 
spherical $P^\circ_\la$. Moreover, one can use
spherical {\em nonsymmetric} Macdonald polynomials  $E_\la$
for iterated torus knots, which have many advantages; see 
Section \ref{sec:exa} below.

\smallskip

{\sf Stabilization.}
The sequences $\vec\rr^j,\,\vec\ss^j$ of length $\ell^j$
for the graph $\l$ and $'\vec\rr^j,'\!\vec\ss^j$ of
length $'\!\ell^j$ for the graph $'\!\l$ will be as above.
We 
use the DAHA-Jones polynomials $\hat{J\!D}$, which
are under the {\em hat-normalization}.
Recall that $\la=\la(b)$ is the Young diagram
representing $b\in P_+$.
\smallskip

We consider now $P_+\ni b=$ $\sum_{i=1}^n \be_i \om_i$
for $A_n $ as
(dominant) weights for {\em any\,} $A_m$ 
with $m\ge n-1$, setting $\om_{n}=0$ for $A_{n-1}$ (when $m=n-1$).
See \cite{CJ,GoN,CJJ,ChD1} concerning the prior versions
of the following general theorem for torus knots and
iterated torus knots.


\begin{theorem}\label{STABILIZ}\cite[Theorem 2.3]{ChD2}
Given a pair $\{\l,\, '\!\l\}$
colored by $\mathbf{b}=(b^j),\, '\mathbf{b}=('\!b^j)$,
there exists a
unique polynomial from $\Z[q,t^{\pm 1},a]$
\begin{align}\label{h-polynoms-hat}
&\hat{\h}_{\l,\,'\!\l}\ =\
\hat{\h}_{(\vec\rr^j,\,\vec\ss^j),
('\vec\rr^j,\,'\vec\ss^j)}(\mathbf{b},\,
'\mathbf{b};\,q,t,a)
\end{align}
such that for any  $m\!\ge\! n\!-\!1$ and proper powers
of $q,t$ (possibly rational)\,:
\begin{align}\label{jones-sup-hat}
&\hat{\h}_{\l,\,'\!\l}(q,t,a\!=\!-t^{m+1})\,=\,
\pm\, q^\bullet t^\bullet\,
\hat{J\!D}^{A_m}_{\l,\,'\!\l}(q,t).
\end{align}
Here $\hat{\h}(a\!=\!0)$
is automatically 
hat-normalized. We note that this normalization becomes
$\hat{\h}(q\!=\!0,t\!=\!0,a\!=\!0)=1$ in Theorem
\ref{con:motknot} and further  coincidence
conjectures below.
 The correction factors 
$\pm\, q^\bullet t^\bullet$ are not needed for sufficiently
large $\,m$ for links and for any  $m\!\ge\! n\!-\!1$ in the case of 
iterated torus knots. Also,  the consideration of
one sufficiently large
$\,m$\, is sufficient to fix $\,\hat{\h}\,$
uniquely. \sq
\end{theorem}

The polynomials $\hat{\h}$ depend only on the isotopy
class of the corresponding iterated torus links.
All symmetries for the $J\!D$\~polynomials
hold for $\hat{\h}$, including the product
formula at $q\!=\!1$
from (\ref{q-1-prod}).
\vskip 0.2cm

Let us discuss the $a$\~degree of $\hat{\h}$\~polynomials.
\comment{
We will assume
that $\rr_i^j,\,'\!\rr_i^j\neq 0$ for $i\!>\! 1$.
Then \,
deg${}_a \hat{\h}^{min}_{\l,\,'\!\l}
(\,\lla,\,'\!\lla\,;\,q,t,a)\,$  is no greater than
\vskip -0.3cm
\begin{align}\label{deg-a-j}
\sum_{j=1}^\kap
\max\{1,|\,\rr^j_1|\}|\rr^j_2\cdots\rr^j_{\ell^j}|\,|\la^j|+
\sum_{j=1}^{'\!\kap}
\max\{1,|\,'\rr^j_1|\}|\,'\rr^j_2\cdots\,'\rr^j_{\,'\!\ell^j}|
\,|'\!\la^j| -\De,&\notag\\
\for \De\!=\!|\la^1\!\vee\!\ldots\!\vee\!\la^\kap\!\vee\,
'\!\la^1\!\vee\!\ldots\!\vee\, '\!\la^{'\!\kap}|
 \hbox{\, for \,} \hat{\h}^{min},\
\De\!=\!|\la^{j_o}| \hbox{\, for \,} \hat{\h}^{j_o},&
\end{align}
where $|\la|$ is the number of boxes in $\la$.
\smallskip
}
We conjectured in \cite{ChD2} for algebraic links that
\begin{align}\label{deg-a-jj}
\hbox{deg}_a\hat{\h}_{\l,\,'\!\l}
=&\,\hbox{\small $\sum$}_{j=1}^\kap\
\min\{\,\rr^j_1,\ss^j_1\}\, \rr^j_2\cdots
\rr^j_{\ell^j}\,|\la^j|\!\\
&+
\hbox{\small $\sum$}_{j=1}^{\,'\!\kap}
\ \min\{\,'\rr^j_1,\,'\ss^j_1\}\,
\,'\rr^j_2\cdots
\,'\rr^j_{\,'\!\ell^j} \,|\,'\!\la^j|-\De\notag\\
\hbox{for\ } \De&=\!|\la^1\!\vee\!
\ldots\!\vee\!\la^\kap\!\vee\,
'\!\la^1\!\vee\!\ldots\!\vee\, '\!\la^{'\!\kap}|
\hbox{\, from \,(\ref{P-arms-legs-union})},
\notag
\end{align}
where $|\la|$ is the number of boxes in $\la$.
A somewhat weaker statement can be justified. This is
compatible with the specialization $q\!=\!1$, which gives that
the $a$\~degree is no smaller than that in (\ref{deg-a-jj}).

The right-hand side of this
formula is the multiplicity
of the corresponding singularity generalized to the colored
case. For instance, $deg_a=|\la|(\ss-1)$ for torus knots
$T(\rr,\ss)$ colored by $\la$, where $\ss<\rr$. 
\vskip 0.2cm

{\sf Superduality.}
As above, we switch from $\mathbf{b},\,'\mathbf{b}$ to the
corresponding sets of Young diagrams
$\{\la^j\}$,\, $\{'\!\la^j\}$; their transpositions
will be denoted by $\{\cdot\}^{tr}$.
Up to powers of $q$ and $t$, one has:
\begin{align}\label{iter-duality}
\hat{\h}_{\l,'\!\l}(\{\la^j\},\{'\!\la^j\}; q,t,a)=
q^{\bullet}t^{\bullet}
\hat{\h}_{\l,'\!\l}(\{\la^j\}^{tr},
\{'\!\la^j\}^{tr};t^{-1},q^{-1},a).
\end{align}
This is {\em superduality}, one of the key properties
of superpolynomials. It was conjectured in \cite{GS,CJ} and
proven in the context of DAHA superpolynomials in 
\cite{GoN,ChD1,ChD2}.

\subsection{\bf Some examples}\label{sec:exa}
An important simplification of the construction of the $J\!D$-polynomials
and $\hat{\h}$-polynomials for iterated torus {\em knots}, is due to
the usage of nonsymmetric Macdonald polynomials $E_\la$ for dominant $\la$
instead of $P_\la$. One simply replaces $P_\la^\circ
=P_\la/P_\la(t^{-\rho})$ by  $E_\la^\circ
=E_\la/E_\la(t^{-\rho})$; no other changes are
necessary. Recall, that $E_\la$ are eigenvectors of $Y$-operators.
Also, one can do the calculations for $GL_{n+1}$, which
is somewhat simpler than those for $SL_{n+1}$ (for $A_n$).  
\vskip 0.2cm

{\sf The case of trefoil}.
Let us calculate $\hat{\h}_{3,2}$ for uncolored trefoil.
We begin with the case of $A_1$. 

In this case,
$\HH$ is generated by $X^{\pm1},Y^{\pm1},T$ subject to  
group relations $TXTX=1=TY^{-1}TY^{-1}, Y^{-1}X^{-1}YXT^2=q^{-1/2}$
and the quadratic one 
$(T-t^{1/2})(T+t^{-1/2})=1$. The action of $\tau_{\pm}$ is:
\begin{align*}
&\tau_+\!:\!
Y\!\mapsto\! q^{-\frac{1}{4}}XY,\ X\!\mapsto\!X,\ T\!\mapsto\!T,\ 
\tau_-\!:\! X\!\mapsto\! q^{\frac{1}{4}}YX,\  
Y\!\mapsto\!Y,\  T\!\mapsto\!T.
\end{align*}
The polynomial representation is given by the formulas
$T\mapsto t^{1/2}s\!+\! 
\frac{t^{1/2}\!- t^{-1/2}}{X^2-1}(s\!-\!1)$, 
$X\mapsto X,\ \, Y\mapsto spT$, where $ s(X)=X^{-1},\,
\, p(X)=q^{1/2}X.$
\vskip 0.2cm

By $\,\sim$\,, 
we mean ``\,up to $q^{\bullet}t^{\bullet}$\,"
in the following calculation: 
\begin{align*}
&J\!D_{3,2}\!=\!
\{\tau_+\tau_-^2(X)\}_{ev}\!\sim\!\{(XY)(XY)X(1)\}_{ev}\!\sim\!
\{Y(X^2)\}_{ev}\\
&=t^{-\frac{1}{2}}q^{-1}X^2-
t^{\frac{1}{2}}+t^{-\frac{1}{2}}|_{X^2\mapsto t^{-1}}\sim
1+qt-qt^2,
\end{align*}
where we replaced  $P_1$ 
by the {\em nonsymmetric} Macdonald polynomial $E_1$.
One has:  $Y(X)\!=\!(qt)^{-\frac{1}{2}}X$ and $E_1=X$.
Recall that 
$\{H\}_{ev}\!\equal\! H(1)(X\!\mapsto\! t^{-\rho})$.
The {\em Jones hat-polynomial} is
$\hat{J\!D}_{3,2}(q\mapsto t)= 1\!+\!t^2\!-\!t^3.$

\vskip 0.2cm

For $GL_{n+1}$, the corresponding
$\HH$ is generated by $X_i^{\pm 1}, Y_j^{\pm 1}, T_k$,  where 
$1\le i,j\le n+1$ and
$1\le k \le n$ for pairwise commutative $\{X_i\}$ and 
$\{Y_j\}$. One has:   
 $\tau_+(Y_1)\!=\!q^{-1/2} 
X_1 Y_1$,\  $\tau_-(X_1)\!=\!q^{+1/2}
Y_1 X_1$ and so on. The action of $Y_1$ in the polynomial
representation is via
the formula 
$Y_1\!=\!\pi T_{n}\cdots T_1$, where 
$\pi: X_1\!\mapsto\! X_2, X_2\!\mapsto\! X_3,
\ldots, X_n\!\mapsto\! q^{-1} X_1$. Here $T_k$ are
those for $A_n$, and $X_{\al}=X_i X_j^{-1}$ for
$\al=\ep_i-\ep_j$. 
These formulas are almost equally simple 
for any $Y_i$, but we need only $X_1, Y_1$.

The calculation of $J\!D$
remains practically the same as it was for $A_1$.
The polynomial $E_1=X_1$ can be used instead of $P_1$.
We obtain that  $\hat{J\!D}^{A_n}_{3,2}
=1+qt-qt^{n+1}$ and 
$\hat{\h}_{3,2}=1+qt+aq$. Indeed, the values of the latter 
at $a=-t^{n+1}$ are $\hat{J\!D}^{A_n}_{3,2}$.

Note that the relations
$\hat{\h}_{3,2}(a\!\mapsto\! -t)=1$ and    
$\hat{\h}_{3,2}(a\!\mapsto\! -t^{2})=1+qt-qt^2$ (only two)
are sufficient to fix $\hat{\h}_{3,2}$ uniquely using
that $deg_a=1$. Generally, 
$deg_a\hat{\h}_{\rr,\ss}(\la)\!=\!
|\la|\bigl(\text{Min}(\rr,\ss)\!-\!1\bigr)$. 

The DAHA construction 
in the case of {\em uncolored}  $T(2n+1,2)$ is very much similar. 
The calculations are simple,
but it is not clear from the DAHA approach why
$\hat{\h}_{3,2}$ and $\hat{\h}_{2n+1,2}$ are so simple.
The motivic approach gives these formulas right away.

\vfil
{\sf Hopf $2$-links}. The $J\!D$-invariants
of such links colored by $\la,\mu\in P_+$ are essentially
$\pi^{\pm}_{\la,\mu}=\{\tau^{\pm 1}_-(P_\la P_{\mu})\}_{ev}$.
The sign is $+$ for $2$-plus-links (algebraic),
those with the
linking number $+1$. 
To be exact, we need to use here $J_\la$ and $J_\mu$, and divide
by the corresponding $LC\!M$. 

Alternatively, the Hopf $2$-link  
can be obtained by  adding 
the {\em meridian} to $\unknot$\,. 
This is the simplest
case of the usage of pairs of trees, $\{\l,\,'\!\l\}$.
 The topological invariance of
$\hat{\h}$-polynomials results in an alternative formula
for $\pi^{\pm}_{\la,\mu}$: it must be  
$q^{\bullet}t^{\bullet}
\Bigl\{P_{\mu}(Y^{\pm 1})\bigl(P_{\la}(X)\bigr)\Bigr\}_{ev}$ . 
The coincidence is equivalent
to following fundamental  DAHA identity:
\begin{align}\label{tau-minus}
\{\tau^{-1}_-(fg)\}_{ev}=
\{\dot{\tau}_-^{-1}(f),\dot{\tau}_-^{-1} (g)\}_{ev}, 
\end{align}
where $\{f,g\}_{ev}\equal\bigl(f(Y^{-1})(g(X))\bigr)(t^{-\rho})$
for any Laurent polynomials $f,g$. 
See Theorem 3.9 from \cite{ChD2}. 

For $2$-plus-links, the identity
$
\{\dot{\tau}_-(fg)\}_{ev}=
\{\dot{\tau}_-(f)(Y)\dot{\tau}_-(g)\}_{ev}$ is used ({\em ibid.}).
We obtain that the 
corresponding $\hat{\h}$ is essentially 
$P_{\la}(q^{\iota(\mu)+\rho_k})$ for $\iota=-w_0$.  
Combining this formula 
with the approach  based on the product
formulas for Macdonald polynomial, 
we obtain that such evaluations are $a$-stable. This is claimed
in Proposition \ref{PROP:stab-values}, the key for the 
$a$-stabilization of  $J\!D$-polynomials in the $A$-case
and the existence of the DAHA superpolynomials. 

\vskip 0.2cm
The simplest superpolynomials for 
Hopf $2$-plus-links (algebraic) are when $\mu=\om_i$
(generally, a minuscule weight). 
Up to the passage from $P_\la$ to $J_\la$ and
the division by the $LC\!M$, the $J\!D$-polynomial 
becomes essentially 
$\{\tau_-(P_\la P_{\om_i})\}_{ev}$.
The calculation of such coinvariants  can be
readily reduced to the {\em Pieri rules}; expansions
of $P_\la P_{\om_i}$ in terms of Macdonald polynomials
are very explicit. Pieri rules are an important part of
the theory of Macdonald polynomials. Their justification
was a demonstration
of the power of DAHA:\, generally, the DAHA-Fourier transform was used.
See \cite{C101}.
Then we use that $\tau_-$ acts in the
polynomial representation:
$\dot{\tau}_-(P_\la)=q^{-(\la,\la)/2-(\la,\rho_k)}P_\la$ 
for $\la\in P_+$,  where
 $\dot{\tau}_-$ is the corresponding restriction. 

As a demonstration,  $\hat{\h}\!=\!
(1\!+\!q^ma)\!+\!(q^m\!-\!1)t$\, for the Hopf $2$-plus-link
colored by $\om_1$ and $m\om_1$ (the $m$-row).  We use
here the {\em Pieri rules}. We will obtain
the same formula in the motivic theory below.
\vskip 0.2cm

{\sf HOMFLY-PT polynomials.}
The definition of $H\!O\!M(t,\aa;\la)$
is especially simple in the uncolored case, which is for 
$\la=\yng(1)=\om_1$.
The following {\em skein  relation} is sufficient to define
them (the reduced ones):
{\small
$$
\aa^{1/2}
H\!O\!M(\nwarrow\kern-10pt\nearrow\kern-10.5pt
\nearrow\kern-11pt\nearrow
\kern-11.5pt\nearrow)
\!-\!\aa^{-1/2}
H\!O\!M(\nearrow\kern-10.2pt\nwarrow\kern-10.7pt
\nwarrow\kern-11.2pt\nwarrow
\kern-11.7pt\nwarrow)\!=\!
(t^{1/2}\!-\!t^{-1/2})
H\!O\!M(\uparrow\uparrow),\ H\!O\!M(\bigcirc)=1.
$$
}

Basically, 
$H\!O\!M(t,\aa;\la)\,\,\widehat{}\!=\!\hat{\h}(q=t,t,\aa=-a; \la)$
for {\em iterated torus knots}, 
where $\,\widehat{}\,$ is the {\em hat-normalization} 
of the reduced $H\!O\!M$: must be $1$ for the unknot.
We note that the parameter $\sqrt{\aa}$ is mostly used in
topology.
The coincidence is due to the author 
for torus knots, Morton-Samuelson (iterated torus knots), and
Cherednik-Danilenko (iterated torus links); see \cite{CJ, MoS,ChD2}.

The following  adjustments are necessary for {\em links}.
The Macdonald polynomials $P_\la$ 
must be used instead
of  the $J$-polynomials $J_\la$ in the DAHA construction. Also, 
the division by $LC\!M$ in 
our $\hat{\h}(\{\la\},\{\mu\})$ must be replaced with
the division by {\em one} of
$P_{\la^i}(t^{-\rho})$ (or one of those for $\mu^j$) 
to match the corresponding {\em reduced} HOMFLY-PT
polynomial.  This is because 
the latter are defined 
in topology with respect to one ``distinguished"  component of
a link. The division by $LC\!M$, which gives the 
$a,t$-polynomiality, is, generally, ``non-topological". This
normalization gives the same as the topological normalization above  
if all non-distinguished components are uncolored
or if all colors coincide.   

\vskip 0.2cm
For the uncolored trefoil, i.e. for $T(3,2)$ and when 
$\la=\om_1=\yng(1)\, $:\,
$H\!O\!M\!=\!\aa(t+t^{-1}-\aa),\, 
H\!O\!M\,\,\widehat{}\!=\! 1+t^2-t\aa$; 
recall that  
$\h\!=\!1+qt+qa$. 

The Alexander polynomials $Alex(t)$ are, generally,
the following specialization: 
$H\!O\!M(t,\aa\!=\!1)(1-t)^{1-\kappa}$, which is  for links
with $\kappa$ components. In particular,
$Alex\!=\!t^{-1}-1+t,\, Alex\,\widehat{}\!=\!1-t+t^2$ for trefoil.
 
The simplest {\em link} is the Hopf $2$-plus-link, $2$
uncolored unknots with the linking number $+1$. Then:
$H\!O\!M=\aa^{1/2}$ 
{\large $\,\frac{1+\aa-t-t^{-1}}{t^{1/2}-t^{-1/2}}$}
(they can be rational in $t$), 
$\hat{\h}=1+(q-1)t+aq$, 
and $Alex\,\widehat{}=1.$
\vskip 0.2cm

According to (\ref{P-arms-legs}), $J_{\Box}=(1-t)P_{\Box}$,
which generally gives that $H\!O\!M\,\,\widehat{}=
\hat{\h}(t,t,-\aa)/(1-t)^{\kappa-1}$ in the uncolored case.
The latter becomes {\large $\frac{1+t^2-t-\aa t}{1-t}=
\frac{1+\aa -t -t^{-1}}{1-t^{-1}}$} for trefoil, which
is $H\!O\!M\,\,\widehat{}$.

The {\em superduality} 
becomes $t^{\frac{1}{2}}\to -t^{-\frac{1}{2}}, 
\aa^{\frac{1}{2}}\to \aa^{-\frac{1}{2}}$ for  $H\!O\!M$;
it is obviously compatible with 
the skein relation above (in the uncolored case).
Also, the Young diagram  $\la$ must go 
to its transpose.

The symmetry $t^{\frac{1}{2}}\to -t^{-\frac{1}{2}}$
holds for $Alex(t)$ too. However, it does not hold for
the Quantum Group $A_n$-invariants.
 The latter are basically
$H\!O\!M(t,\aa=t^{n+1};\la)$ and the substitution $\aa=t^{n+1}$
is obviously incompatible with the superduality unless $q=t$.
Recall that $\aa=-a$.

We note that {\em mirroring\,} of iterated torus links results in  
changing $\aa_i\mapsto -\aa_i$ in (\ref{Newtonpair}), and 
in $q\mapsto q^{-1}, t\mapsto t^{-1}, a\mapsto a^{-1}$
in the superpolynomials (followed by the hat-normalization).
Also, one must replace $Y$ by $Y^{-1}$ for the second graph
if it is present. See \cite{ChD2} for details.

\vskip 0.2cm

\section{\sc Motivic superpolynomials}
We switch now to motivic theory. The main objective is
to define and study {\em motivic superpolynomials} 
in full generality:\,
for multibranch plane curve singularities and for
arbitrary ranks. They conjecturally coincide
with the DAHA superpolynomials
for the corresponding algebraic links colored by ``rows". 
From the viewpoint of the theory 
of Springer fibers of type $A$, this is the case of the most
general characteristic polynomials. This section and the next
one are for {\em uncolored  unibranch} singularities;
the general case will be addressed later.


\subsection{\bf Unibranch singularities}\label{sec:basic}
For any base field $\F$, plane curve singularities are
those with the {\em singularity rings} 
$\r=\F[[x,y]]/(F(x,y))$ for 
$F(x,y)\in \F[[x,y]]$ 
such that $F(0,0)=0$ and $(F(x,y))=F(x,y)\F[[x,y]]$.
 The {\em unibranch} ones
are for irreducible $F(x,y)$, assumed absolutely  
irreducible (over the algebraic closure  $\overline{\F}$ of $\F$).

Equivalently, unibranch singularities are
for 
$\r=\F[[x,y]]\subset \o\equal \F[[z]]$, where  $z$ is
the {\em uniformizing parameter}, 
$x,y$ are in the maximal ideal $\mathfrak{m}_\o=(z)=z\o$ 
of $\o$, the localization of
$\r$ is the whole $\F((z))$, and $\r/\mathfrak{m}=\F$ for
the maximal ideal $\mathfrak{m}=\r \cap (z)$ of $\r$.
 We use this 
definition in this paper;
the equation $F(x,y)=0$ will not be needed.
 
Quite a few considerations below will be for any Gorenstein
singularities; plane curve singularities are such.
However, the main properties of superpolynomials, including the
coincidence conjecture and superduality,  hold only
for (multibranch) plane curve singularities.  

\vskip 0.2cm
The simplest invariants of a singularity are the
{\em multiplicity}, which is  $mult=$
dim$\,\F[[z]]/\F[[z]]\mathfrak{m}$, and the arithmetic
genus \,$\de=$\,dim$\,\F[[z]]/\r$, the {\em Serre number}.


\vskip 0.2cm
{\em Valuation semigroup.}
It is one of the key in the theory of unibranch singularities.
The definition of
this semigroup is as follows:
$\Ga\equal\bigl\{\,\nu_z(f), 0\neq f\in \r\subset \o=\F[[z]]
\,\bigr\}$, where
$\nu_z$ is the valuation, the order of $z$. Then
$\de=|G|$ for the {\em set of gaps}
$G\equal\Z_+\setminus \Ga$, and the
multiplicity of the corresponding singularity
is $mult=\min(\Ga\setminus \{0\})$.
\vskip 0.2cm

The {\em conductor} of $\r$ is 
the greatest
ideal $\mathfrak{c}=(z^\cc)=z^\cc\o$ in $\o$ that belongs to $\r$.
Generally, 
$\de+1\le \cc\le 2\de$, and $\cc=2\de$ if and only if $\r$ is 
{\em Gorenstein}, which includes plane curve singularities. 
Equivalently, $\r$ is such if and only if the map 
$(\Ga\setminus \{\cc+\Z_+\}) \ni \ga\mapsto \cc-1-\ga
\in G=\Z_+\setminus \Ga$ is an isomorphism.  Its image is in $G$,
which is of size $\de$,
since  $(\cc-1)$ is a gap. This gives the inequality above and that
this map is an isomorphism if and only if $\cc=2\de$. 
\vskip 0.2cm

{\sf Algebraic links.} Let $\F=\C$. The link corresponding to
$\r$ is the 
intersections of $\{(x,y)\in \C^2 \mid  F(x,y)=0\}$
considered near $(0,0)\in \C^2$
with the $3$-sphere 
$S_\ep^3\subset \C^2$ centered at $(0,0)$ of small
radius $\ep$.
They are {\em knots} (connected) for unibranch 
singularities, which are for irreducible $F(x,y)$.
The invariants $\de$ and $mult$ are topological: they depend
only on the isotopy class of the corresponding link.

Importantly, $\Ga$ totally determines
the topological type  
of the corresponding algebraic knot for $\F=\C$ 
(considered up to isotopy), 
which fundamental fact is due to Zariski and others. Thus, 
topological invariants
of rings $\r$ are exactly those expressed in terms of $\Ga$.

To give an example, the hat-normalized {\em Alexander polynomial} 
from Section \ref{sec:exa}
is $Alex\,\widehat{}=$
$(1-t)\sum_{\nu\in \Ga}t^\nu$ in terms of $\Ga$.
 For instance, 
it becomes  $(1-t)(\frac{1}{1-t}-t)=1-t+t^2$ for trefoil, $T(3,2)$.
\vfil

The theory of topological equivalence  
of algebraic {\em links} is significantly
more ramified. The coincidence of
semigroups for the components and the corresponding 
pairwise linking  numbers is sufficient
due to the Reeve Theorem.
The linking numbers are the indices of
the corresponding branches and can be algebraically calculated
via the ring of singularity. All must be positive
for algebraic links. Generally, {\em splice diagrams} 
provide full topological classification.
\vfil

{\em Cables.} Torus knots $T(\rr,\ss)$ are those for the
singularities $x^\ss\!=\!y^\rr$, where $\rr,ss>0$ and
gcd$(\rr,\ss)\!=\!1$. The corresponding rings are
$\r=\C[[x\!=\!z^\rr,y\!=\!z^\ss]]$; they are called
{\em unibranch
quasi-homogeneous} due to their invariance with respect to
the action $x\mapsto \varsigma^\ss, y\mapsto \varsigma^\rr$
for $\varsigma\neq 0$. For them, 
the multiplicity is $\min(\rr,\ss)$ and  
$\de =$ {\large$\frac{(\rr-1)(\ss-1)}{2}$.} The latter formula is actually 
due to Sylvester (from Frobenius coin problem). 

The simplest ``non-torus" family 
is $\r\!=\!\C[[z^4,z^6\!+\!z^{7+2m}]]$ for $m\!\in\! \Z_+$,
which are 
of multiplicity $4$ and with $\de_m\!=\!8\!+\!m$. The corresponding
cables are $C\!ab(13+2m,2)C\!ab(2,3)$, where 
$C\!ab(2,3)=T(3,2)$. One can replace here 
$T(3,2)$ by $T(2,3)$ because it does not change
the isotopy type of the cable; the cable parameters
are topological invariants. 
\vskip 0.2cm

The simplest ``triple" cable is 
$C\!ab(53,2)C\!ab(13,2)C\!ab(3,2)$.
 The corresponding
ring is $\r\!=\!$
$
\!\C[[x\!\!=\!z^8,y\!=\!z^{12}\!+\!z^{14}\!+\!z^{15}]].
$ 

let us provide the corresponding {\em Newton's pairs.}
Generally, they are
$\{\rr_1,\ss_1\},\{\rr_2,\ss_2\}, \ldots\,$ such that
the Puiseux-type equation of the singularity is
{  $y=x^{\frac{\ss_1}{\rr_1}}\Bigl(1+ c_1 
x^{\frac{\ss_2}{\rr_1 \rr_2}}
\bigl(1+c_2 x^{\frac{\ss_3}{\rr_1 \rr_2 \rr_3}}(\cdots)\bigr)\Bigr)$}
for generic $c_i$. 
The Newton's pairs are $\{(2,3),(2,1),(2,1)\}$ in this example.

The valuation semigroup is $\Gamma=\lan 0,8,12,26,53\ran$, the
linear span of these numbers over  $\Z_+=\{\Z\ni n\ge 0\}$.
The Serre number is $\de\!=\!42$  and $mult=8$.
Generally, $\Ga=\lan 0,\, \rr_1 \rr_2 \rr_3,\,
\aa_1 \rr_2 \rr_3,\, \aa_2 \rr_3,\, \aa_3\ran$ for the cable parameters 
$(\aa_i,\rr_i)$ above (here $\rr_1<\ss_1$ is used). Recall that
$\aa_1\!=\!\ss_1$, $a_2\!=\!\rr_1 \ss_1 \rr_2+\ss_2,$\, 
$\aa_3\!=\!\aa_2 \rr_2 \rr_3+\ss_3\,$, and so on for any number 
of $\{\rr_i,\ss_i\}$. 
\vskip 0.2cm

As another example, let us calculate $\Ga$ and $\de$ for 
the ring 
$\r=\F[[z^{\upsilon \rr}, z^{\upsilon \ss} +z^{\upsilon \ss+p}]]$,
where, gcd$(\rr,\ss)=1$, $\rr<\ss$,
$\upsilon>1$ and gcd$(\upsilon,p)=1$ for $p\ge 1$. The
corresponding cable is $C\!ab(\upsilon \rr \ss\!+\!p,\,\upsilon)\,
C\!ab(\ss,\rr)$. One has:
\begin{align}\label{Ga-cab}
&\Ga=\lan \upsilon \rr,\,
\upsilon \ss,\,\upsilon \rr \ss\!+\!p\ran,\ 
2\de=\upsilon^2 \rr \ss\!-\!\upsilon (\rr\!+\!\ss)+
(\upsilon\!-\!1)p+1.
\end{align}

We will mostly assume below that the first $C\!ab(\rr_1,\ss_1)$ in any 
cables satisfies the inequality $\rr_1\le \ss_1$, as it was in these
examples. 
\vskip 0.2cm

{\sf Good reduction.}
The motivic superpolynomials will be associated with
plane curve singularities over finite 
fields $\F=\F_q$ for $q=p^k$ and prime $p$.
The passage from the base field $\C$ to $\F_q$ will be
necessary below to relate DAHA and motivic superpolynomials.
It is as follows. We begin with 
$\r$ over $\C$,  define it over $\Z$, which is always 
doable {\em within a given isotopy type}, and then
consider $\r\otimes _{\Z} \F_q$ for $q=p^k$.  
Then $p$ is called a {\em prime (place) of good
reduction} if
$F(x,y)$ remains irreducible over $\F_q$ 
upon this procedure and $\Gamma$ remains unchanged.  
\vskip 0.2cm

All primes $p$ are good in this sense for the
rings $\C[[x\!=\!z^\rr,y\!=\!t^\ss]]$ considered above.
More generally, it is expected that for any algebraic knot,
there are no primes $p$  of bad reduction.
This means that  given any $p$, there exists
$\r$ within a given topological type
(representing this knot) where this $p$ is good.
\vskip 0.2cm

To give an example, let us consider 
$\Z[[x\!=\!t^4,y\!=\!t^6\!+\!t^7]]$.
Then  $\Ga\!=\!$ {\small $\{0,4,6,8,10,13,14,16,17,18,\ldots\}$} 
and $\de=8$ for $\F=\C$.  This ring has bad reduction only at $p=2$. 
Indeed, 
$\nu_z(y^2-x^3)=14$ over $\F_2$,  which is $13$ for $p\neq 2$. However,
this singularity is equivalent over $\C$
(analytically, not only topologically) to the one for
$\Z[[t^4+t^5,t^6]]$, where $p=3$ is the
only place of bad reduction. Thus, the corresponding cable
has no primes of bad reduction. 

\subsection{\bf Standard flags}\label{sec:stand}
Let $\r\subset \o\equal\F[[z]]$ be a unibranch 
plane curve singularity over a field $\F$, which will be 
assumed irreducible over $\overline{\F}$. 
An $\r$-invariant submodule $M\subset \o$ is called a {\em standard
module} if $\o\supset M\ni \phi=1+z(\cdot)$.
Equivalently, $\De(M)\ni 0$, where 
$\De(M)\equal \{\nu_z(v)\mid 0\!\neq\!v\!\in\! M\}$.
The sets $\De(M)$ are $\Ga$-modules for any $\r$-modules
$M$, which means that $\Ga+\De\subset \De$.
{\em Standard}  sets $\De$ are those in $\Z_+$
containing $0$; equivalently,  $\Ga\subset\De$ and $\Ga+\De=\De$.
Obviously, a module $M$ is
standard if and only if $\De(M)$ is standard.

Given $M\subset \o$, we will constantly use
$deg(M)\equal$ 
dim$_\F(\o/M)=|Z_+\setminus \De(M)|$, and the {\em deviation} of $M$ 
from $\r$, which is 
$dev(M)\equal \de-deg(M)=|\De(M)\setminus \Ga|$. Also,
the {\em $q$-rank} of $M$ is 
$r\!k_q(M)\equal M/\mathfrak{m}M$, which will be used when $\F=\F_q$.
To provide an example,
$mult(\r)=r\!k_q(\o)=\min(\Ga\setminus \{0\})$.
The $q$-rank equals the minimal number of generators of $M$ over $\r$,
which is no greater than 
the number of generators of $\De(M)$ over the action
of $\Ga$. Indeed, the $\r$-span $M'$  of (any) 
elements in $M$ with the valuations generating $\De(M)$  
has the same $\De$ and therefore must coincide with $M$. 

For standard modules: 
$dev(M)\ge 0$
and it is $0$ if and only if $M=\phi\r$
for some $\phi\in 1+(z)\in \o$. The latter modules are called
{\em invertible}. 
They form the {\em generalized
Jacobian variety} $Jac$ of $\r$, which is 
the algebraic group $U^\bullet_\o\equal \{\phi\}$ modulo 
$U_\r^\bullet=U^\bullet_\o \cap \r$ canonically isomorphic to 
$z\o/(z\o\cap \r)\cong \F^\de$ (the 
group scheme $\mathbb{G}_a^\de$, to be exact).
\vskip 0.2cm

For $\ell\ge 0$, the  
{\em flagged Jacobian}
 $\j_\ell$, considered as a set of $\F$-points by now, 
is formed by {\em standard flags} 
$\vect{M}=
M_0\!\!\subset\!\! M_1\!\!\subset\!\! \cdots\!\!\subset\!\!
 M_\ell\!\subset\!
\o\!=\!\F[[z]]$ of 
$\r$-submodules $M_i$ of $\o$ such that\  \,
(a)\, $M_0\ni \phi=1\!+\!z(\cdot)$ (where $(\cdot)\in \o$), 
\, (b)\,
dim$\,M_i/M_{i-1}\!=\!1 \hbox{ and } M_i=M_{i-1}\oplus 
 \C\, z^{g_i}(1+z(\cdot))$, \,  and (c)\, (important)\,
$g_i<g_{i+1},\, \text{ where } i\ge 1$. We will call them
$\ell$-flags. Standard modules $M$ are then $0$-flags.
\vfil

\vskip 0.2cm
{\em Piontkowski cells.}
An abstract increasing  sequence
of $\Ga$-modules 
$\vect{\De}\!=\!
\{\De_0\!\subset\!\cdots \!\subset\!\De_\ell\subset \Z_+\}
$ is called {\em standard} if $\De_0$ contains $\Ga$, 
$\De_i=\De_{i-1}\cup\{g_i\}$,
and $g_i<g_{i+1}$ for $1\le i\le \ell$.
We set $\De(\vect{M})=\{\De(M_i)\}$. The flags $\vect{M}$ are 
standard if and only if the corresponding $\vect{\De}$ are standard.

Given a standard 
$\vect{\De}$, the corresponding  {\em Piontkowski cell}
is defined as 
$\j_\ell(\vect{\De})
\equal\bigl\{ \vect{M}\in \j_\ell \mid \De(\vect{M})=
\vect{\De}\bigr\}$. One has:\,
$\j_\ell=\cup\, \j_\ell(\vect{\De})$, 
where the union is disjoint. Some cells can be empty.

For quasi-homogeneous
singularities $\r=\F[[x=z^\rr,y=z^\ss]]$, where 
gcd$(\rr,\ss)=1,\, \rr,\ss>1$,
all standard $\Ga$-modules $\De$ 
come from some standard $M$, i.e. there are no
empty Piontkowski cells.  Their
number is the rational Catalan number 
$\frac{1}{\rr+\ss}\binom{\rr+\ss}{\rr}$. We note
that standard $\De$ are described
by {\em Dyck paths} in rectangles ``$\rr\times \ss$" due
to \cite{GM}. 

Piontkowski in \cite{Pi} employed the
method of syzygies to justify  that 
the corresponding
cells in $\j_{\ell=0}$
are nonempty affine spaces ${\mathbb A}^m$ and find formulas for 
their dimensions $m$. The generalization to flags is in \cite{ChP1}. 
Presumably, the $\C^*$-action in $\r$ 
and  $\j_0$ can be used to justify that these
cells are affine spaces due to a proper singular version  
of {\em Bia{\lx}yncki-Birula Theorem} as in \cite{Dr};
alternatively, Gr\"obner cells can be used as in \cite{ChG}.
\vskip 0.2cm

Generally, it seems that empty $\j_0(D)$
for some standard $D$ are always present
unless for quasi-homogeneous singularities.
For instance, let $\r=\F[[z^4,z^6+z^7]]$, the simplest non-torus
example. Then two from $25$ standard $\De$
are not in the form $\De(M)$ for standard $M$, 
which phenomenon is due to Piontkowski.   
\vskip 0.2cm

If all (non-empty) Piontkowski cells are affine
spaces, then
dim$\,\j_\ell(\vect{\De})$ and the list of empty cells
are sufficient to know. This happens for relatively ``small"
families, including $\r\!=\!\C[[z^4,z^6\!+\!z^{7+2m}]]$ for 
$m\!\in\! \Z_+$. See \cite{Pi,ChP1}. 
The method of syzygies always gives some
combinatorial dimension formulas for such families,
but they can be involved. 
\vskip 0.2cm

The varieties $\j_\ell(\vect{\De})$   are
conjectured to be {\em configurations of affine
spaces}, i.e. 
unions and differences of affine spaces $\mathbb{A}^m$ 
in a bigger $\mathbb{A}^N$, not always equidimensional and
connected. This conjecture  includes            
multibranch generalizations and any $\o$-ranks considered below. 
\vskip 0.2cm

Let us mention that the connection between the dimensions
of these cells in $\j_0$ and the {\sf deviations}
in the case of $\F_q[[z^\rr,z^\ss]]$ was observed in \cite{LS},
which  is a special case of motivic
{\em superduality}.
This is of obvious importance:\, the deviations
are immediate in terms of $\De(M)$, the dimensions are
generally difficult to calculate. 


\subsection{\bf Compactified Jacobian and {\em ASF}}\label{sec:asf-nil}
Let us supply $\j_0$ with a structure of a {\em projective}
variety. 
By construction, standard modules constitute  
a {\em disjoint} unions of
quasi-projective varieties, those for 
different values of the deviations of $M$. Importantly, they
can be combined in one irreducible projective variety. 
The main steps are as follows. 
\vskip 0.2cm

Any standard $M$ contains the {\em conductor} of $\r$.
 Recall that it
is the greatest
ideal $(z^\cc)=z^\cc\o$ in $\o$ that belongs to $\r$. We use that 
$\cc=2\de$, which holds if and only if $\r$ is 
Gorenstein. Equivalently, 
the map $$
(\Ga\setminus \cc+\Z_+) \ni \ga\mapsto \cc-1-\ga
\in G=\Z_+\setminus \Ga=\text{``the set of gaps" } 
$$ 
is an isomorphism if and only if $\r$ is Gorenstein. 

Thus,
 $\phi\in 1+(z)\in M$ (for standard $M$) implies that 
 $\phi\cdot (z^{\cc})=
(z^{\cc})\subset M$, which is $(z^{2\de})$ for Gorenstein $\r$.

The key is the fact that 
$z^{dev(M)}M\supset 
(z^{2\de})$ for standard $M$ due to 
Pfister-Steenbrink \cite{PS}. Equivalently, 
$dev(M)+\De(M)\supset 2\de+\Z_+$. 
\vskip 0.2cm

Using this inequality,  let $M\mapsto M'\equal z^{dev(M)}M$.
Then $dev(M')=dev(M)-dev(M)=0$. This map establishes an 
identification of standard $M$ with $\r$-modules  
$(z^{2\de})\subset M'\subset \o$ such that  $dev(M')=0$. 
The inverse map is 
$M'\mapsto z^{-d} M'$ for $d=\text{Min}\{
\nu_z(m)\!\mid\! 0\!\neq\!m\in M'\}$. 

Finally, 
 $\j_0$ becomes $\{M' \mid dev(M')=0\}$ as a set, 
a subset of 
 the Grassmannian Gr$(2\de,\de)$ of the
subspaces of the middle dimension in
 $\o/(z^{2\de})$. These subspaces must be invariant with
respect to the action of $\r$ (for $x$ and $y$ in the case
of plane curve singularities). This supplies $\j_0$
with a natural structure of a projective scheme over $\F$, 
a closed subscheme of Gr$(2\de,\de)$. This is 
the {\em compactified Jacobian}.
The same notation $\j_0$ will be used for it. We will mostly
need its description via standard $M$ in this paper.

\vfil
The irreducibility of $\j_0$
for plane curve singularities is due to Rego.
 The flagged Jacobians  $\j_\ell$ naturally become some
fiber spaces over $\j_0$; see Theorem \ref{thm:a-formula}.
The projectivizations of   $\j_\ell$
for any $\ell$ are of importance, but we do not
discuss them in this paper.

We note that this construction can be applied to any irreducible curve
singularities; the conductor must be changed correspondingly. 
The compactified Jacobians are reducible for non-planar singularities,
but they are always
connected due the action of
the generalized Jacobian there \cite{PS}; so are $\j_\ell$ for
any $\ell\ge 0$.

\vfil
{\sf Affine Springer fibers.} 
Their definition requires the knowledge of the
equation $F(x,y)=0$ for the
generators $x,y$ of $\r$.
Adjusting the $z$-series for $x,y$ in
our approach, 
this equation can be assumed a polynomial 
in terms of $x,y$ over $\F$, irreducible
since we do the unibranch case.

Let $n$ and $m$ be the top $x$-degree and 
$y$-degree of this equation. Then
our compactified Jacobian $\j_0$ can be interpreted as a (parahoric)
affine Springer fiber $\x_\ga$ from \cite{KL}
either for $GL_n$ or  for $GL_m$;
the equation connecting $x$ and $y$ becomes
the corresponding {\em characteristic equation}. The abbreviation
$ASF$ will be used.



Generally,  $ASF$ from  \cite{KL} is for  
a reductive Lie algebra $\mathfrak{g}$, either 
$GL_n$ or  $GL_m$ in our case.
For any field $\F$, let $\mathfrak{g}[[x]]=
\mathfrak{g}\otimes_\F \F[[x]]$ and 
 $\mathfrak{g}((x))=
\mathfrak{g}\otimes_\F \F((x))$. Accordingly,
we define $G[[x]]$ and $G((x))$
for simply-connected $G$ with Lie$(G)=\mathfrak{g}$.

Given $\ga\in  \mathfrak{g}[[x]]$, let 
$\x_\ga^\circ \equal \{g\in  
G((x))/G[[x]]
\mid g^{-1}\ga g\in  \mathfrak{g}[[x]]\},$
 where we assume that
the centralizer of $\ga$ in $G((x))$ is {\em anisotropic}
(the so-called nil-elliptic case).  For $GL_n$,
$$\x_\ga^\circ\cong 
\{M\subset \k=\F((z))
\mid \r M=M,\, \k M=\k\},$$ where
$F(x,y)=F_\ga(x,y)=$\,det\,$(\mathbf 1 y-\ga)$, which is  
of degree $n$ with respect to $y$.
Finally, 
$\x_\ga=\x_\ga^\circ/\Z\cong \j_0$, where the action of $v\in \Z$
is $M\mapsto z^v M$, a particular case of the action  
 due to \cite{KL}.

Our definition of $\j_0$ is directly in terms of $\r$ and
does not require $F(x,y)$. The fact 
that the transposition $x\leftrightarrow y$ and $n\leftrightarrow m$,
where $m$ is the degree of $F(x,y)$ with respect to $x$, does not
change the corresponding $ASF$ is not immediate 
from the definition of $\x_\ga$. 
 
The corresponding 
{\em orbital integral}
is  $\h^{mot}(q,t\!=\!1,a\!=\!0)$ for
the motivic superpolynomials defined below, where $\F=\F_q$. 
Conjecturally, $\h^{mot}(q\!=\!1,t\!=\!1,a\!=\!0)$ 
gives the Euler characteristic $e(\x_\ga)$;
this is the case of  ``the field with $1$ element".

\comment{
Also, let us mention here the Oblomkov-Shende-Rasmussen conjecture
(2012), 
which is about {\em nested Hilbert schemes} of pairs of
ideals in $\r$; see below.
 It is an extension of the prior Oblomkov-Shende 
conjecture proved by Maulik for HOMFLY-PT polynomials; see below.
}
 
\subsection{\bf Unibranch superpolynomials}\label{sec:unibr}
The rings $\r\subset \o$ will be now over 
 $\F=\F_q$. This section is mostly for any irreducible curve
singularities.
 However the conjectures about the
structure of Piontkowski cells, the topological invariance
of motivic superpolynomials and their connection to the
DAHA ones will be only for plane curve singularities. 

Following \cite{ChP1}, the
{\em motivic superpolynomial} of $\r$  for $\ell\ge 0$ is 
$$\h^{mot}\!\equal\!
\sum_{\{M_0\subset\cdots\subset
 M_\ell\}\in\j_\ell(\F)}t^{deg(M_\ell)}a^\ell,
$$
where it suffices
to assume that $M_0$ is standard. The number of such flags is finite
because $M_0$ (any standard $M$) contains $z^{2\de} \F_q[[z]]$.
\vskip 0.2cm

{\sf Reduction formula.}
This definition can be reduced to the consideration
of standard $M$ containing $\r$.
Following Lemma 3.1 from \cite{ChG}, we consider the pairs
$\tilde{M}=\{\vect{M},I\subset M_{\ell}\}$, where $\vect{M}$
are standard flags and $I$ are invertible modules, and define 
$\tilde{\h}^{mot}=
\sum_{\tilde{M}}\, t^{deg(M_\ell)}a^\ell.$
The number of $I\subset M_{\ell}$ equals $q^m$ for
$m=dev(M_{\ell})=
\de-deg(M_{\ell})$, which results in the formula
$$\tilde{\h}^{mot}=q^{\de}\h^{mot}(q,t\mapsto t/q,a)=
q^{\de}\sum_{\r\subset \vect{M}} 
t^{deg(M_\ell)}a^\ell.$$
Note that all flags here are standard because 
they contain $1$.

This provides 
an expression for motivic superpolynomials in
terms of certain {\em ideals} in
$\r$; i.e. the usage of $\o$
can be bypassed.
The {\em reciprocity map} is used:\, $M\mapsto 
M^\ast\equal\r:M=\{f\in \F((z)) \mid fM\subset \r\}$, where
$N:M\equal \{f\in \F((z)) \mid fM\subset N\}\simeq M^\ast:N^\ast$.
For instance, $\r^\ast=\r$ and $\o^\ast=(z^{2\de})$.
We claim that $\De(M^\ast)=(\De(M))^\ast$, where 
\begin{align}\label{Dedual}
\De^\ast\equal \{p\in \Z\mid p+\De\subset \Ga\}=
(2\de-1)-(\Z\setminus \De). 
\end{align}
See e.g.,  Lemma 2.7 in \cite{GM}; 
$v\pm X=\{v\pm x \mid x\in X\}$ for $X\subset \Z \ni v$.
The modules here are $z^m M$ for standard $M$ and  $m\in \Z$. If they
are standard, then $\Ga\subset\De\subset \Z_+$ and  
$\De^\ast=(2\de-1-(\Z_+\setminus \De))\cup (2\de+\Z_+)$. 
The ring $\r$ is assumed Gorenstein.
Using (\ref{Dedual}), $\De(M^\ast)\subset \De^\ast$ and 
dim\,$ N/M=$ dim\,$M^\ast/N^\ast$, which  results in the coincidence
$\De(M^\ast)=\De^\ast$.
See e.g., (2.5),(2.6) from \cite{Sto}. In particular,  $\r\subset M\subset \o$ is equivalent to
$\r\supset M^\ast \supset (z^{2\de})$:\, the passage
to ideals mentions above.
\vskip 0.2cm

We will state now the uncolored {\em Coincidence Conjecture} 
for unibranch plane curve singularities, to be extended in further
sections.  Claims $(ii)-(iv)$ follow from $(i)$, the key,
and some DAHA facts, though they can be
considered as independent conjectures. 

\begin{conjecture}\label{con:motknot}
Given $\r_\Z$ defined over $\Z$, let $K$ be the corresponding
algebraic knot for $\r_\C=\r\otimes_{\Z} \C$. Assume that  
$\r_\F=\r_\Z\otimes \F$ for $\F=\F_q$ is a good reduction of $\r_\Z$,
which means (as defined above) that $\r_\C\subset \o_\C=\C[[z]]$ and 
$\r_\F\subset \o_F=\F[[z]]$ have the same semigroup $\Ga$.

(i) Then the uncolored DAHA superpolynomial $\hat{\h}(q,t,a)$ 
constructed for $K$,  which is 
a polynomial in terms of $q,t,a$ with integer coefficients,
coincides with $\h^{mot}(q,t,a)$ defined for $\F=\F_q$.

(ii) Conjecturally 
$q^\de t^\de \h^{mot}(t^{-1}, q^{-1}, a)\!=\! 
\h^{mot}(q, t, a)$ (the superduality);
this is a theorem for  $\hat{\h}(q,t,a)$.
It conjecturally becomes the functional equation for the corresponding
$L$-function introduce below. 

(iii) The superpolynomial $\h^{mot}$ is a polynomial in terms
of $q,t,a$ with integer coefficients; recall that it is a polynomial
in terms of 
$t,a$ by construction. Since $\hat{\h}$ is a polynomial,
$(i)$ implies $(iii)$.

(iv) Furthermore, $\h^{mot}$ is a topological invariant of
the singularity, i.e. depends only on $\Ga$. This is known
for the DAHA superpolynomials; this claim follows 
from $(iii)$ upon minor assumptions.
\end{conjecture}

It will be extended later to arbitrary plane curve
singularities colored by rows, which
is for arbitrary $F(x,y)=0$, possibly reducible and
not square-free,
corresponding to
the most general $ASF$ of type $A$.

\section{\sc Unibranch theory continued}
The main topic of this section will be the justification
of the motivic Weak Riemann Hypothesis for $a=0$ from \cite{ChW}.
It will be based on the {\em decomposition formulas} below
 for $\h^{mot}$, which
are of clear independent importance. It is not known how to
obtain them in the DAHA setting. We continue to 
consider unibranch 
$\r\subset \o=\F_q[[z]]$.

\subsection{\bf The {\bf\em q}-rank decomposition}
The following theorem will be generalized 
below for any plane curve singularities $F(x,y)=0$, 
not only irreducible and square-free.

For an arbitrary module $M$, let $M^{\{i\}}=M\cap
(z^i)$, an $\r$\~module, and 
$\overline{M}^{\{i\}}$ be
the image of $M^{\{i\}}$ in $\overline{M}=M/\mathfrak{m} M$,
a vector space over $\F_q$.  Then
dim\,$M^{\{g\}}/M^{\{g\!+\!1\}}\!=\!1$ for
$g\!\in\! \De(M)$\, and $0$ otherwise. Recall that $\mathfrak{m}=
\r\cap \mathfrak{m}_\o=\r\cap (z)$ is the maximal ideal in $\r$,
and $\r/\mathfrak{m}=\F_q$.

\begin{theorem}\label{thm:a-formula}
For any  $\r$, not only Gorenstein, and
standard $M$:
\begin{align}\label{hmotrk}
\h^{mot}(q,t,a)=\sum_{M} t^{deg(M)}\prod_{i=1}^{r\!k_q(M)-1}
(1+q^i a), \text{\, where\, } \prod_{i=1}^0=1.
\end{align}
\end{theorem}
{\it Proof.} We follow Proposition 2.4 from \cite{ChP2}.
The $\r$\~modules $M_i$ from a standard $\ell$-flag
are uniquely determined by 
their images
$\overline{M}_i\equal M_i/\mathfrak{m}M_{\ell}$ in
$\overline{M}_{\ell}$, which follows
from the {\em Nakayama Lemma}.
Given $\vect{\De}=\{\De_0\!\subset\! \De_1\!\subset
\cdots\subset\!\De_\ell\}$ and a 
standard module $M_{\ell}$, 
standard $\ell$\~flags $\vect{M}'$ such that $M'_\ell=M_\ell$
and $\De(\vect{M}')=\vect{\De}$
exist if and only if
\begin{align}\label{mocondd}
\hbox{dim\,}_{\F_q}\,\!\overline{M}_{\ell}^{\{g_i\}}/
\overline{M}_{\ell}^{\{g_i\!+\!1\}}=1
\for 1\le i\le \ell.
\end{align}
These flags can be identified via the Nakayama Lemma with full
$\ell$\~flags of $\F_q$\~subspaces
$\overline{M}_{0}\!\subset\cdots\!\subset
\overline{M}_{i}\!\subset\cdots
\!\subset\overline{M}_{\ell}$
such that
\begin{align}\label{moconddd}
\overline{M}_i\!+\!\overline{M}_{\ell}^{\{g_i\!+\!1\}}=
\!\overline{M}_{\ell} \for 1\le i\le \ell-1.
\end{align}
Assuming(\ref{mocondd}),
this space of such flags is biregular to $\mathbb{A}^{\ell(\ell-1)/2}$
if  $M_0$ is fixed. Thus, it suffices to count the
number of all possible $M_0$ for a given $M_{\ell}$,
which will be $M$ in (\ref{hmotrk}). The module $M_0\subset M_{\ell}$
must be standard, which gives the condition $\overline{M}_0\not\subset
z\F_q[[z]]$, where the latter is of co-dimension $1$ in $\F_q[[z]]$.
This is sufficient for the following calculation.

The contribution 
of one $M_{\ell}$ to $\h^{mot}$ with $n=r\!k_q(M)$ is:
\begin{align}\label{rk-formula1}
&h^{mot}(q,t,a)=1+\sum_{k=1}^n q^{\frac{k(k-1)}{2}}\Bigl(\,
\text{\footnotesize $\qbin{n}{k}-\qbin{n-1}{k-1}$}\ \Bigr) a^k\\
=&\,1+q^{\frac{k(k+1)}{2}}
\text{\footnotesize $\qbin{n-1}{k}$}\ a^k
=(1+aq)\cdots (1+q^{n-1}a),\notag
\end{align}
where the formula {\footnotesize $\qbin{n}{k}-\qbin{n-1}{k-1}$}
$=q^k\text{\footnotesize $\qbin{n-1}{k}$}$ for $q$-binomial coefficients
was used and the $q$-binomial theorem. Finally, we consider
the summation over all standard $M=M_\ell$.
\sq

\vskip 0.2cm

\vskip 0.2cm
{\sf Recovering $J\!D$-polynomials.}
The following Corollary is of conceptual importance: it
provides a {\em motivic} (conjectural) interpretation
of $J\!D$-polynomials, {\em i.e. before the  $a$-stabilization.}
This can  make the verification of the first (main!) claim 
of Conjecture \ref{con:motknot} more direct.

\begin{corollary} 
(i) Assuming Conjecture \ref{con:motknot},
the connection of
uncolored $\hat{J\!D}_{A_n}(q,t)$ for $n\ge 0$
from formulas (\ref{jones-bar}) and  (\ref{jones-sup-hat})
(under the hat-normalization) with the motivic superpolynomials
is as follows:
\begin{align}
&\h^{mot}(q,t,a=-1/q^{n+1})=(qt)^\de \hat{J\!D}_{A_n}(t^{-1},q^{-1})\\
=&\sum_{M} t^{deg(M)}\prod_{i=1}^{r\!k_q(M)-1}
(1-q^{i-n-1}), \text{ where } r\!k_q(M)\le n+1. \notag
\end{align}

(ii) In the case of $A_1$, i.e. for uncolored refined (hat-normalized)
Jones polynomials,  the following motivic
formula (conditionally) holds:
\begin{align*}
&\hat{J\!D}_{A_1}(t^{-1},q^{-1})
\!=\!1+(1-q^{-1})q^{-\de}\sum_{M} t^{-dev(M)}, \text{where\ }
r\!k_q(M)\!=\!2.\text{\sq} 
\end{align*}
\end{corollary}

The summations here are over standard $M$.
Recall that generally 
$deg(M)=|\Z_+\setminus \De|=\de-dev(M)$. 

Concerning $(ii)$,
the standard modules of $r\!k_q(M)=2$ are ``mostly"
with $\De(M)=\Ga\cup \{g+\Ga\}$ for $g\in G$,
those of $dev=1$, but not always.
 Such modules can emerge from the
cells with the number of $\Ga$-generators of $\De$, the $\Ga$-rank
of $\De$, greater than $2$. This makes the problem of counting
the standard modules in $(ii)$ not entirely in terms of (relatively
simple) Piontkowski cells of $dev=1$. 
\vskip 0.2cm

Note that negative $q,t$-monomials can occur in $(i,ii)$.
This is a standard for Jones
polynomials and their quantum group generalizations.
The specializations of
(uncolored, reduced)  Khovanov-Rozansky polynomials to $A_n$ 
are positive $q,t$-polynomials by construction; they is different from
for our $J\!D$-polynomials.
 
Generally, the coefficients
of the products $\prod_{i=1}^{r-1}(1+q^i a)$, which are 
$\sum_{r\!k_q(M)=r} t^{deg(M)}$, can have
negative $q,t$-monomials. We note here that 
the positivity of the
$q,t,a$-coefficients for DAHA superpolynomials
is the last unproven {\em intrinsic} conjecture 
from the initial list in \cite{CJ}. This is for
any algebraic knots colored by rectangle diagrams. 

\vskip 0.2cm
 
\subsection{\bf Examples}
Among torus knots, $T(9,4)$ seems the simplest
examples with negative terms in the coefficients of
the decomposition
from Theorem \ref{thm:a-formula}.
 Then $\r=\F_q[[x=z^4,y=z^9]]$ and 
$\h^{mot}_{9,4}=$
\renewcommand{\baselinestretch}{0.5} 
{\small
\(
q^{12} t^{12}+(1+a q) q^7 t^4 (q\mathbf{-t}+q^2 t+2 q^2 t^2+q^3 t^2
\mathbf{-q t^3}
+q^2 t^3+q^3 t^3+2 q^3 t^4+q^3 t^5+q^4 t^5+q^3 t^6+q^4 t^6+q^4 t^7)
+ (1+a q)(1+a q^2)q^4 t^2 (1+q+q^2+q t+2 q^2 t+q^3 t
+2 q^2 t^2+2 q^3 t^2+q^4 t^2+q^2 t^3+3 q^3 t^3+
2 q^4 t^3+q^3 t^4+2 q^4 t^4+q^3 t^5+2 q^4 t^5+q^5 t^5+q^4 t^6
+q^5 t^6+q^5 t^7)
+(1+a q) (1+a q^2) (1+a q^3) (1+q t+q^2 t+q^3 t+q^2 t^2
+q^3 t^2+q^4 t^2+q^3 t^3+q^4 t^3+q^5 t^3+q^4 t^4+q^5 t^4+q^5 t^5
+q^6 t^6).
\)
}
\renewcommand{\baselinestretch}{1.2} 

\vskip 0.2cm
{\sf The case of $T(7,3)$.}
One has: 
$\r=\F_q[[z^3,z^7]]$, 
$\de=6$, $\Ga\setminus \{2\de+\Z_+\}=\{0,3,7,6,9,10\}$, and
$G=\{1,2,4,5,8,11\}$ in this case.
Recall that all Piontkowski cells are affine spaces for
quasi-homogeneous singularities. Generally, the number of 
standard modules for  $\F[[z^r,z^s]]$
is the rational slope Catalan number $\frac{1}{r+s}\binom{r+s}{r}$
(Beauville). Table \ref{Table7-3} provides the contributions to the
$(1+a q^i)$\~decomposition in this example.

The sets of gaps $D=\De\setminus \Ga$ for standard $\De$ are shown.
Recall that $|D|=dev(\De)$ and $deg(D)=\de-|D|$, which gives the
power of $t$. The
brackets $[\ldots]$ will be used below for the $D$-sets.

The cell with $r\!k_q=1$ is unique:\, the generalized 
Jacobian,
which contributes $q^6 t^6$. The simplest one with $r\!k_q=2$
is for $D=[11]$; it contributes $q^5 t^5$ to $\h^{mot}$ (the
first column). The module $M=\o$, which contributes $1$, is
of $q$-rank $3$; it is not unique such.

{\footnotesize
\begin{table*}[ht!]
\[
\centering
\begin{tabular}{|l|l|}
\hline                   
$D$-sets, $r\!k_q=2$ & $q,t$-terms \\
\hline 
 & \\   
$\widetilde{11}$ & $q^5t^5$ \\  
$\widetilde{4}, 11$ & $q^5 t^4$\\ 
$\widetilde{8},11$ & $q^4 t^4$\\
$\widetilde{5},8,11$ & $q^3 t^3$ \\
$\widetilde{2},5,8,11$ & $q^3 t^2$\\ 
$\widetilde{1},4,8,11$ & $q^4 t^2$\\      
$\widetilde{4},8,11$ & $(q^4-q^3) t^3$\\   
   
\hline
\end{tabular}
\hspace{0.1cm}
\begin{tabular}{|l|l|}
\hline    
$D$-sets, $r\!k_q=1$ & $q,t$-terms  \\
\hline 
$\emptyset$ & $q^6 t^6$ \\
\hline 
$D$-sets, $r\!k_q=3$ & $q,t$-terms  \\
\hline  
& \\ 
$\widetilde{4},\widetilde{5},8,11$ & $q^2 t^2$\\
$\widetilde{1},4,\widetilde{5},8,11$ & $q^2 t$ \\
$\widetilde{2},\widetilde{4},5,8,11$ & $q t$ \\
$\widetilde{1},\widetilde{2},4,5,8,11$ & $q^0 t^0$\\    
$\widetilde{4},\widetilde{8},11$ & $q^3 t^3$\\
\hline   
\end{tabular}
\]
\caption{Ranks and $q,t$-terms for $T(7,3)$}
\label{Table7-3}
\end{table*}
}

\comment{
{\footnotesize
\begin{table*}[ht!]
\[
\centering
\begin{tabular}{|l|l|}
\hline      
$D$-sets, $r\!k_q=1$ & terms  \\
\hline  
$\emptyset$ & $q^6 t^6$ \\
\hline                   
$D$-sets, $r\!k_q=2$ & terms \\
\hline 
 & \\
$\widetilde{4}, 11$ & $q^5 t^4$\\ 
$\widetilde{8},11$ & $q^4 t^4$\\
$\widetilde{5},8,11$ & $q^3 t^3$ \\
$\widetilde{2},5,8,11$ & $q^3 t^2$\\ 
$\widetilde{1},4,8,11$ & $q^4 t^2$\\      
$\widetilde{4},8,11$ & $(q^4-q^3) t^3$\\   
   
\hline
\end{tabular}
\hspace{0.1cm}
\begin{tabular}{|l|l|}
\hline    
$D$-sets, $r\!k_q=2$ & terms  \\
\hline 
 & \\  
$\widetilde{11}$ & $q^5t^5$ \\  
\hline 
$D$-sets, $r\!k_q=3$ & terms  \\
\hline  
& \\
$\widetilde{4},\widetilde{5},8,11$ & $q^2 t^2$\\
$\widetilde{1},4,\widetilde{5},8,11$ & $q^2 t$ \\
$\widetilde{2},\widetilde{4},5,8,11$ & $q t$ \\
$\widetilde{1},\widetilde{2},4,5,8,11$ & $q^0 t^0$\\    
$\widetilde{4},\widetilde{8},11$ & $q^3 t^3$\\
\hline   
\end{tabular}
\]
\caption{Ranks and $q,t$-terms for $T(7,3)$}
\label{Table7-3}
\end{table*}
}
}

One has:
$\h^{mot}_{7,3}=$
{\small
$q^6 t^6+(1+a q) (1+a q^2) (1+q t+q^2 t+q^2 t^2+q^3 t^3)+
(1+a q) (q^3 t^2+q^4 t^2+q^4 t^3+q^4 t^4+q^5 t^4+q^5 t^5)$}.
 The coefficients are all
positive in this decomposition, but this is a special feature
of knots $T(3n\pm 1,3)$.

The valuations of the generators of the corresponding $M$
are $0$ (not shown) and those marked by {\em tilde} in this table.
Their number is smaller than the number 
of $\Ga$-generators of $\De=\Ga\cup D$ for some $M$
if $D=[4,8,11]$; the generators of the corresponding $\De$ 
are $0,4,8$. Namely, 
$\j_0(\De)\simeq \mathbb{A}^4$ is a union of
$\mathbb{A}^4\setminus \mathbb{A}^3$ formed by the modules of
$q$-rank $2$ (contributing $(q^4-q^3)t^3$\,),  
and those of $r\!k_q=3$, the remaining ones, which 
form $\mathbb{A}^3$ and contribute $q^3 t^3$. 
\vskip 0.2cm

{\sf The case of $\F[[z^4,z^6+z^7]]$.}
This is $\r$ with the simplest non-torus algebraic 
cable, which is  $K=C\!ab(13,2)C\!ab(2,3)$.
All cells are affine spaces in this case;
so dim$\,=$dim\,$J_0(\De)$ are sufficient in the table
below for the corresponding
$D\equal \De\setminus \Ga$. Recall that
dim\,$\o/M=\de-|D|$, which gives the
power of $t$.
Two standard $\De$
from $25$ have no standard $M$, namely for
$D=[2,15]$ and $D=[2,11,15]$. 

\comment{
We provide the top ranks of $M$ corresponding to $D$-sets.
{\footnotesize
\begin{table*}[ht!]
\[
\centering
\begin{tabular}{|l|l|l|}
 \hline 
\hbox{$D$-sets} & $dim$ & $rk$\\
\hline
$\varnothing$ & 8 & 1\\
15 & 7 & 2\\
11,15 & 6 & 2\\
7,11,15 & 6 & 2\\
9,15 & 7 & 2\\
9,11,15 & 5 & 3\\
7,9,11,15 & 4 & 3\\
3,7,9,11,15 & 4 & 2\\
5,9,11,15 & 5 & 2\\
5,7,9,11,15 & 3 & 3\\
3,5,7,9,11,15 & 2 & 3\\
1,5,7,9,11,15 & 4 & 2\\
\hline
\end{tabular}
\hspace{0.1cm}
\begin{tabular}{|l|l|l|}
\hline
\hbox{$D$-sets} & $dim$ & $rk$\\
\hline
1,3,5,7,9,11,15 & 2 & 3\\
2,7,11,15 & 6 & 2\\
2,9,15 & 7 & 2\\
2,9,11,15 & 6 & 3\\
2,7,9,11,15 & 5 & 3\\
2,3,7,9,11,15 & 4 & 3\\
2,5,9,11,15 & 5 & 3\\
2,5,7,9,11,15 & 3 & 4\\
2,3,5,7,9,11,15 & 1 & 4\\
1,2,5,7,9,11,15 & 3 & 3\\
1,2,3,5,7,9,11,15 & 0 & 4\\
\,2,15 and 2,9,15 & $\emptyset$ & 0\\
\hline   
\end{tabular}
\]
\caption{Dimensions and ranks$}
\label{Table4-6-13-rk}
\end{table*}
}
}

The whole (uncolored) superpolynomial is: 
$\h^{mot}(q,t,a)=$
\renewcommand{\baselinestretch}{1.2} 
{\small
\(
1 + q t + q^8 t^8 + q^2 \bigl(t + t^2\bigr) 
 + 
 q^3 \bigl(t + t^2 + t^3\bigr) + q^4 \bigl(2 t^2 + t^3 + t^4\bigr) 
+  q^5 \bigl(2 t^3 + t^4 + t^5\bigr) + q^6 \bigl(2 t^4 + t^5 
+ t^6\bigr) + 
 q^7 \bigl(t^5 + t^6 + t^7\bigr) + 
a \bigl(q + q^2 \bigl(1 + t\bigr) + q^3 \bigl(1 + 2 t + t^2\bigr) 
+ q^4 \bigl(3 t + 2 t^2 + t^3\bigr) + 
    q^5 \bigl(t + 4 t^2 + 2 t^3 + t^4\bigr) + q^6 \bigl(t^2 
+ 4 t^3 + 2 t^4 + t^5\bigr) + 
    q^7 \bigl(t^3 + 3 t^4 + 2 t^5 + t^6\bigr) + q^8 \bigl(t^5 
+ t^6 + t^7\bigr)\bigr)+  
a^2 \bigl(q^3 + q^4 \bigl(1 + t\bigr) + q^5 \bigl(1 + 2 t + t^2\bigr) 
+ q^6 \bigl(2 t + 2 t^2 + t^3\bigr) + q^7 \bigl(2 t^2 + 2 t^3 
+ t^4\bigr) + q^8 \bigl(t^3 + t^4 + t^5\bigr)\bigr) +
a^3 \bigl(q^6 + q^7 t + q^8 t^2\bigr).
\)
}
\renewcommand{\baselinestretch}{1.2} 

{\footnotesize
\begin{table*}[ht!]
\[
\centering
\begin{tabular}{|l|l|}
 \hline 
\hbox{$D$-sets} & $dim$\\
\hline
$\varnothing$ & 8\\
15 & 7\\
11,15 & 6\\
7,11,15 & 6\\
9,15 & 7\\
9,11,15 & 5\\
7,9,11,15 & 4\\
3,7,9,11,15 & 4\\
5,9,11,15 & 5\\
5,7,9,11,15 & 3\\
3,5,7,9,11,15 & 2\\
1,5,7,9,11,15 & 4\\
\hline
\end{tabular}
\hspace{0.1cm}
\begin{tabular}{|l|l|}
\hline
\hbox{$D$-sets} & $dim$\\
\hline
1,3,5,7,9,11,15 & 2\\
2,7,11,15 & 6\\
2,9,15 & 7\\
2,9,11,15 & 6\\
2,7,9,11,15 & 5\\
2,3,7,9,11,15 & 4\\
2,5,9,11,15 & 5\\
2,5,7,9,11,15 & 3\\
2,3,5,7,9,11,15 & 1\\
1,2,5,7,9,11,15 & 3\\
1,2,3,5,7,9,11,15 & 0\\
2,15 and 2,11,15 & $\emptyset$\\   
\hline   
\end{tabular}
\]
\caption{Dimensions for $\Ga=\lan 4,6,13\ran$}
\label{Table4-6-13-0}
\end{table*}
}

\vskip 0.2cm
The table below provides all $D$ and the
corresponding dimensions of the cells $\j_0(\De)$.
For instance, there are $3$ cells of dimensions $7$ in $\j_0$
(for $a=0$). Namely, those  with $D=[15],[9,15],[2,9,15]$
and  $t^7, t^6, t^5$. Generally, the number of cells of
dim\,$=\de-1$ is the multiplicity of singularity;
it equals the coefficient of $t$ for $q\!=\!1,a\!=\!0$  due to the
superduality, which is for $\Z_+\!\setminus \De=\{1\},\{2\},\{3\}$ in this
example. Only $\{1\}$ results in dim\,$=1$ and $qt$ in $\h^{mot}$. 

We note 
the {\em reciprocity involution} of standard $\De$, which sends
$\De\mapsto \De^\ast- \min\{\De^\ast\}$ for
$\De^\ast= \Ga\setminus\{(2\de-1)-D\}$. It preserves ``dim". 
For instance, it sends $[15]\mapsto [2,9,15]$, and
$[2,9]\mapsto [2,9]$.


Let us provide the $(1\!+\!a q^i)$-decomposition for this 
superpolynomial: $\h^{mot}=$
\renewcommand{\baselinestretch}{1.2} 
{\small
\(
q^8 t^8\ +\  (1\!+\!a q) (1\!+\!a q^2) (1\!+\!a q^3) 
\bigl(1+q t+q^2 t^2\bigr)\ +\ 
(1\!+\!a q) (1\!+\!a q^2)\bigl(q^2 (1+q) t+q^2 (q+q^2) t^2
+q^2 (q+2 q^2+q^3) t^3
+q^2 (q^2+q^3) t^4+q^5 t^5\bigr)\ +\ (1\!+\!a q) \bigl(q^4 t^2+
(-1+q) q^4 t^3+2 q^6 t^4+q^4 (q^2+q^3) t^5+q^4 (q^2+q^3) t^6
+q^7 t^7\bigr).
\)
}
\renewcommand{\baselinestretch}{1.2}

\subsection{\bf Weak Riemann Hypothesis}
The field $\F$ can be arbitrary in the
following proposition. 
We will assume there and below that $\de\ge 1$.

Let 
$E_M=\{f\in \o\mid fM\subset M\}$, 
$U_M=E_M\cap \o^*$, where $\o^*$ is the group of invertible
elements in $\o$. For instance, $U_\r=\r^*\subset U_M$ for any 
$\r$-module $M$.
We set 
$U_M^\bullet=U_M/U_\r$, which is isomorphic to the
group $U_M^1=\{f\in 1+z\o \mid \phi M= M\}$ modulo $U_\r^1$.
For instance,  $U_\o^\bullet\cong \mathbb{G}_a^\de(\F_q)$. It  
acts simply transitively on invertible modules $I\subset \o$,
and $U_M^\bullet$ acts in the sets 
$I(M)\equal\{I\subset M\}\cong \mathbb{A}^{dev(M)}$ 
without fixed points. We use Lemma \ref{lem:inv}, 
which provides the reduction of standard $M$ to 
modules $M\supset \r$.
In particular, dim$\, U_M^\bullet\le dev(M)$.

\begin{proposition}\label{prop:iandj}
(i) Let us assume that $I(M)$ is not one orbit of $U_M^\bullet$,
equivalently, the dimension of the latter (unipotent) group is smaller
than $dev(M)$. Then the dimension of the orbit
$U_\o^\bullet (M)$, an affine space isomorphic to
$U_\o^\bullet/U_M^\bullet$, is greater than  
$deg(M)\!=\!\de\! -\! dev(M)$.

(ii) Assume that $I(M)=U_M^\bullet(I_0)$ for some
$I_0\in I(M)$, i.e. there is only one $U_M^\bullet$-orbit there.
Let $\phi_\circ=1+z(\cdots)\in M$ and 
$M_\circ=\phi_\circ^{-1}M$, which is a ring containing $\r$. If
$\De=\De(M)$ is not in the form $\Ga\cup \{m+\Z_+\}$ for $m\ge 0$,
then $\j_0(\De)$ contains an algebraic fibration with the
fibers isomorphic to $\mathbb{A}^{deg(M)+1}$.

(iii) Let $I(M)=U_M^\bullet(I_0)$ as above and 
$\De(M)=\Ga\cup \{m+\Z_+\}$ for $m\ge 0$. Then 
$\j_0(\De)=U_\o^\bullet(M)$, which is
the affine space of dimension $deg(M)=|\{g\in G
\mid 1\le g<m\}|$. Recall that $G=\Z_+\setminus \Ga$.
\end{proposition}
{\em Proof.} The orbits $U_\o^\bullet(M)$
for any given $\r$-module $M$ in $\F((z))$ are biregular to 
affine spaces due to the Theorem of Chevalley-Rosenlicht.
This is for any unipotent groups; see Theorem 1 
in \cite{Pi}
and the discussion there. 
Since $\j_0(\De)$ is a union of the corresponding $U_\o^\bullet$-orbits,
the dimensions of these orbits (affine spaces) 
must be strictly  greater than $\de-dev(M)=deg(M)$ in $(i)$.
Recall that $deg(M)\equal$ 
dim$_\F(\o/M)=|G(M)|$, where $G(M)=|Z_+\setminus \De(M)|$. 
\vskip 0.2cm

In the case of $(ii)$, the orbit $U_\o^\bullet(\phi M)$ is isomorphic
to $\mathbb{A}^{deg(M)}$ for any invertible $\phi$.
The module $M_\circ$ 
there contains $1$ and, therefore, the whole
$\r$. The assumption is
that any ratio $(1+f)/(1+g)$ for $f,g \in  
M_\circ\cap (z)$ must belong to
$U_{M_\circ}$. This gives (formally) that $fg\in M_\circ$ for
any such $f,g$. Thus,
$M_\circ$ is a subring of $\o$ (containing $\r$). 

Let $\{e_g, g\in D=\De\setminus \Ga\}$ for $\De=\De(M_\circ)$
be any system of 
elements in $M_\circ$ such that $\nu_z(e_g)=g$. They and $1$
generate $M_\circ$ as an $\r$-module. If $e_g$ is 
in this system, then $e_{g+\ga}$ can be omitted for any $\ga\in \Ga$,
but we will take all of them. Let $m'=\min\{D\}$ and
$c'=\min\{n\in \Z_+ 
\mid n+\Z_+\subset \De\}$, the conductor of $M_\circ$.
Note that $c'-1\not\in \De$.

For any $\la\in \F$, 
consider the span  $M^\la$ of $e_1(\la)=e_1+\la z^{c-1}$ for
$e_1=e_{m'}$ and $e_g$ for all $g>m'$ (taken as such). 
Then $M^\la$
is a ring and an $\r$-module with the same $\De$.
Indeed, the ``correction terms" in any
products of $e_1(\la)$ with $e_g$ will be of valuation
greater than $c'-1$, which belong to $M\cap M^\la$.
 Similarly, $\De$ serves $M^\la$  because
any relations that may result in adding extra elements to
$\De$ may have corrections only due to the multiplication
of $e_1(\la)$ by elements of $\r$, i.e.
with valuations greater than $c'-1$, which do not influence $\De$.

Also, $M^\la$ is not $\phi M_\circ$ for any
$\phi=1+z(\cdots)$ if $\la\neq 0$. Indeed, such $\phi$ must
then belong to $M^\la$ because $\phi\cdot 1\in M^\la$.
However, the latter is a ring and contains $\phi^{-1}$,
which series must be from $M_\circ$. 

We obtain that  $\j_0(\De)$ is a disjoint union of 
copies of $\mathbb{A}^{deg(M)+1}$. The dimension $deg(M)+1$
is due to the action of $U_M$ and
the one-parametric group formed by the $\la$-transformations
of the rings $M_\circ$ considered above.  
\vskip 0.2cm

Claim $(iii)$ is a straightforward calculation of 
the parameters of the
corresponding $\j_0(\De)$. The generators $e_g$ of $M$
can be taken as follows: $e_0=1+\sum_{g\le m}{\la_g z^g}$ and
$e_g=z^g$ for $g>c'-1$, where $g\in G=\Z_+\setminus \Ga$.
The parameters $\la_g$ are arbitrary from $\F$ there. \sq
 
\vskip 0.2cm
As an application, we can verify the {\em Weak Riemann Hypothesis}
from \cite{ChW} for motivic superpolynomials modulo 
the conjecture that $|\j_0(\De)|$
are polynomials in terms of $q$ for any standard $\De$.
The next theorem can be reformulated to hold for any
field $\F$, which we will omit.
Part $(i)$ is direct from 
Proposition  \ref{prop:iandj}. Part $(ii)$ will be proven below.

\begin{theorem} \label{thm:q-rank}
We assume that $|\j_0(\De)|$ are polynomials
in terms of $q$ for any standard $\Ga$-modules $\De$, which 
holds conjecturally for any plane curve singularities.
Let $G=\Z_+\setminus \Ga =\{0<g_1<g_2<\cdots<g_{\de}\}$.
As above, $\de\ge 1$ and we set $m=\min\{\Ga\setminus\{0\}\}$.

(i) Any monomials 
$q^i t^j$ in $\h^{mot}(q,t, a=0)$ satisfy $i\ge j$, when 
$i=j$ holds only for $\De$  from part  $(iii)$ of 
Proposition  \ref{prop:iandj}. 
The sum of such extremal ($i=j$)  monomials  in
$\h^{mot}(q,t, a=0)$ is 
$\sum_{i=0}^\de q^i t^i$. 


\comment{
(ii) Let $\h^{mot}_k$ be the coefficient of $a^k$ in $\h^{mot}$,
where $0\le k\le mult(\r)-1$.
Then $q^{-k(k+1)/2}\h^{mot}_k =\sum_{i=0}^{r_k}q^i t^i$ modulo 
the monomials
$q^i t^j$ with $i>j\ge 0$, where $r_k$ is as follows. 
It is the greatest $r\ge 0$
such that the $\Ga$-module
$\De_r=\Ga\cup\{g_{r+1},g_{r+2},\cdots,g_{\de}=2\de-1\}$ 
is of $q$-rank
$k+1$. For instance, $\De_\de=\Ga, r_1=\de, \ 
\De_{\de-1}=\Ga\cup \{g_\de\}, r_2=\de-1$. 
}

(ii) Let $\h^{mot}_\dag=\h^{mot}(q,t, a=-t/q)$. It is 
a sum of monomials $\pm q^i t^j$ with $i,j\ge 0$ such that  $i+1\ge j$.
The extremal terms ($i+1=j$)  can be only 
due to any modules of $q$-rank $2$ contributing $q^i t^{i}$
or modules of $q$-rank $3$ contributing $q^i t^i$ 
to $\h^{mot}(a=0)$. Respectively, they contribute 
 $-t(q^i t^i)$ for 
$0\le i\le \de-1$and $+t^2 q (q^i t^i)$ to $\h^{mot}_\dag$, where
 for the former 
 and
$0\le i\le \de-3$ for the latter. Their total contribution
to  $\h^{mot}_\dag$ is then  
$-t$ if $\de=1$, and $-t-t^{\de}q^{\de-1}$ otherwise. \sq
\end{theorem}

{\sf Weak RH}. It  was stated in \cite{ChW} for
DAHA superpolynomials of algebraic links colored by rectangle 
Young diagrams. Let $\h^{mot}(q,t,a)=\sum_{k=0}^{mult-1}
 \h^{mot}_k(q,t)a^k$, where $mult$ is that for $\r$.
 Then the claim in the unibranch case is that 
all $t$-zeros $z\neq 0$ of $\H_k(q,t)\equal
\h^{mot}_k(qt,t)$ satisfy {\em Weak RH}:\,
$|z|=q^{-1/2}$ for $q>0$ sufficiently close to $0$.

\begin{conjecture}\label{conj:weak}
The sums of extremal monomials of
$\h^{mot}_k(q,t)$ are 
$q^{k(k+1)/2}(1+(qt)+\cdots+(qt)^{n_k})$ for certain $n_k$. 
Here $n_0=\de$ due to Part $(i)$ of the theorem.
In particular, Weak RH holds for $\H_k(q,t)$. 
\sq
\end{conjecture}

We did not provide (conjecture) this formula for the {\em extremal
polynomials} in \cite{ChW}. It
readily gives that $|z|=q^{-1/2}$ for sufficiently small $q$.
Indeed, the
complex conjugation $\overline{\ze}$ and $\ze'=1/\ze$ are 
roots of $\H_k(q,q^{-1/2}\ze)$ if $\ze=q^{1/2} t$ is such a root.
Thus, $\ze$ and $\overline{\ze}'$
converge in the limit $q\to 0_+$ to the same root $\varsigma$
of the polynomial
$1+\varsigma^2+\varsigma^4+\cdots+\varsigma^{2 n_k}$.
Then $\varsigma$ must be a root
of multiplicity at least $2$, which is impossible.

Similarly, Part $(ii)$ gives {\em Weak RH} for nonzero $z$ when
$a\mapsto -t/q$. Indeed, we arrive at 
the equation $\varsigma+\varsigma^{2\de-1}=0$
for $\de>1$. We did not consider {\em RH} in this case
in \cite{ChW}; though see \cite{ChS}.

\subsection{\bf Proof of Part (ii)}
%
As above, the $\Ga$-rank, $r\!k_\Ga(\De)$ is
the number of $\Ga$-generators of $\De$. They are
uniquely determined as follows.
The smallest nonzero element $g^1$ in $\De$
is automatically the smallest generator. Consider 
$\De\setminus \{g^1+\Ga\}$. The smallest element $g^2$ there
will be the $2${\small nd} generator. Then, consider 
$\De\setminus \bigl\{\{g^1+\Ga\}\cup \{g^2+\Ga\}\bigr\}$ 
and proceed by induction. These generators are exactly
all primitive (indecomposable)
elements in $\De$, i.e. those not in the form $a+\ga$ for
$a\in \De$ and $0<\ga\in\Ga$. Equivalently, they are $g$ such 
that 
$\De\setminus\{g\}$ are $\Ga$\~modules.

Let $\De_r\equal\Ga\cup\{g_{r+1},g_{r+2},\ldots,g_{\de}=2\de-1\}$
for $r=0,1,\cdots,\de$.  
For instance, $\De_0=\Z_+$ is  of 
maximal possible $\Ga$-rank $m$, 
$\De_{\de-1}=\Ga\cup \{g_\de\}$ is of rank $2$, and 
$\De_\de=\Ga$ is of rank $1$.

Note that 
$\De_{\de-2}=\Ga\cup \{g_\de,g_\de-m\}$ and
$r\!k_\Ga(\De_{\de-2})=2$ for Gorenstein $\r$.  Indeed,
$\De_{\de-2}=\Ga\cup \{2\de-1,2\de-1-m\}$
due to the minimality of $m$, where we use the map 
$\Ga\setminus \{2\de+\Z_+\}\ni \ga\mapsto 2\de-1-\ga\in G$. 

Certain modifications of
Parts $(ii)$ and $(iii)$ of the following Lemma are expected
to be sufficient to prove Conjecture \ref{conj:weak}.

\begin{lemma}\label{lem:Der}
(i) As above, $m=g_1=\min\{\Ga\setminus \{0\}\}$, the multiplicity
of the singularity associated with  $\r$. 
For any $0\le r\le \de$, one has:  $r\!k_q(\De_r)$ 
$=|\{g_i \mod m \mid g_i \in \De_r\}|$. For instance,
the sequence $r\!k_q(\De_r)$ 
monotonically decreases for $r=0,1,\ldots, \de$ with steps
$0,-1$.  As above,
$\De_r=\Ga\cup \{n>g_r\}$, where 
we set $g_0=0$. This claim holds for any $\r$ (not only
Gorenstein). 

(ii) Let $k$ be such that $g_{\de-k+1}-g_{\de-k}$ is not divisible 
by $m$, but
$g_{\de-i}=g_{\de}-m i$ for $1\le i\le k-1$. Then 
modules $\De_\de,\De_{\de-1},\ldots,\De_{\de-k}$ 
are the only ones of
$q$-rank $2$ among all $\De_r$, which holds  
for any Gorenstein rings. 
Let $\De=\De_{\de-k}=
\Ga\cup \{g_{\de-k+1},\ldots,g_\de\}$, 
$\De'=\De_{\de-k-1}=\Ga\cup \{g_{\de-k},\ldots,g_\de\}$. Then
$r\!k_q=2$ for $\De$ 
and $r\!k_q=3$ for $\De'$. 
The contributions of $\De$ and $\De'$ to $\h^{mot}(a=0)$ will
be $q^k t^k$ and $q^{k-1}t^{k-1}$ respectively.

(iii) One has $r\!k_\Ga(\De_{r})=r\!k_\Ga(\De_{r-1})-1$
if and only if $g_r+m\in \Ga$ for $m$ as above.
 If this is the case,
assume additionally that $g_{r+1}=g_r+1$. Let $\De'_r$=
$\Ga\cup \{g_{r},g_{r+2},\ldots,g_{\de}\}$, which is
deleting
$g_{r+1}$ from $\De_{r-1}$; its degree is  $deg=r$. Then 
$\j_0(\De'_r)\cong\mathbb{A}^{r+1}$ for the corresponding
Piontkowski cell. This cell is
a disjoint union of 
$\mathbb{A}^r$ for $M$ (in this cell) with
$r\!k_q(M)=r\!k_\Ga(\De_{r})$
and $\mathbb{A}^{r+1}\setminus \mathbb{A}^r$ for $M$
with $r\!k_q(M)=r\!k_\Ga(\De_{r})-1$.
\end{lemma}
{\it Proof.} The $\Ga$-generators $\{g^i\}$ of any  $\De$ must
have different
residues $g^i \!\mod m$. In particular, $r\!k_\Ga(\De_r)\le
|\{g \!\mod m \mid g \in \De\}|$, where the inequality is
quite possible.
Claim $(i)$ means that these ranks coincide and are equal to 
$r\!k_q$ for $\De=\De_r$. 

To justify $(i)$, we proceed by induction.
Assuming that  $r\!k_\Ga(\De_r)$ 
$=|\{g_i \!\mod m \mid g_i \in \De_r\}|$, the passage
to $\De_{r+1}$ is by removing $g_{r+1}$ from $\De_r$,
the first nonzero element.
Then the remaining generators of $\De_r$ and
$g_{r+1}+m$ generate $\De_{r+1}$ and are its generators
unless $g_{r+1}+m \in \Ga$ .
Indeed, it is impossible for $g_{r+1}+m$ to be represented
in the form $g'+\ga$ for $g'\in \De$ and $0\neq \ga\in \Ga$
because $\ga\ge m=\min\{\Ga\setminus \{0\}\}$.

The justification of
the coincidence of  $r\!k_\Ga(\De_r)$ with  $r\!k_q(\De_r)$
is  as follows. The corresponding standard $M$ have the generators
$e^0=1+\sum_{j=1}^r \la_j e_i$, where $\nu_z(e_i)=g_i$, for
$\la_i\in \F$ and $e^i$ such that  $\nu_z(e^i)=g^i$ for the set 
$\{g^i\}$ 
of $\Ga$-generators of $\De_r$, where
 $1\le i\le r\!k_\Ga(\De_r)-1$.
They  are $\r$-independent in this case (not for any $\De$).
Accordingly, the corresponding cell $\j_0(\De_r)$
is isomorphic to  $\mathbb{A}^r$.
\vskip 0.2cm

Part $(ii)$. The generators of  $\De$ are 
$\{0,g_{\de-k+1}\}$ and those of $\De'$ are
$\{0,g_{\de-k},g_{\de-k+1}\}$; 
their $\Ga$-ranks are $2$ and $3$. 
Such $\De$ are all among $\{\De_r\}$ 
of rank $2$. Indeed, the elements of $\De\setminus \Ga$ must
be an arithmetic sequence with step $m$, which is possible only
for $D$ as above among $\{\De_r\}$. 

Rank $2$ changes to rank $3$
in the sequence $\{\De_\de, \De_{\de-1},\cdots\}$ exactly at 
$\De'$, i.e. right after
$\De=\De_{\de-k}$.  The dimension
of the cell $\j_0(\De)$ is $\de-k$ by Part $(i)$, 
the degree of $\De$ is $deg=\de-dev(\De)=
\de-k$, and this cell  contributes  $-t(tq)^{\de-k}$
to $\h^{mot}_\dag$, where $k=1,2,\ldots,\de$.

The degree of $\De'$ is $\de-k-1$, where $k=2,\ldots, \de-1$.
 The coordinates of the corresponding
cell are the coefficients of $g_0=1+\sum_{i=1}^{\de-k-1}\la_i e_i$.
Its contribution to  $\h^{mot}_\dag$ is 
$q^{1+2}(t^2/q^2)(tq)^{\de-k-1}=t(tq)^{\de-k}$;
Theorem \ref{thm:q-rank} $(i)$ has been used.
There can be no other extremal terms in $\h^{mot}_\dag$. 
\vskip 0.2cm

Part $(iii)$. This is the case when the generators $\{e^i\}$
can be linearly dependent over $\r$ and
the $q$-rank is smaller than the $\Ga$-rank.
For $\De'_r$, let 
$e^0=1+\sum_{i=1}^{r-1}\la_i e_i+\la_{r+1}e_{i+1}$, 
$e^1=e_{r}+\mu e_{r+1}$, $e^2=e_{r+2}$ and so on. There are
free $r+1$ parameters. However, if $\la_1=\mu$, then
the $q$-rank becomes $r\!k_q(\De_r)$ due to the relation
$(z^{m+g_r}+\cdots)e^0-(z^m+\cdots)e^1=
(\la_1-\mu)z^{g_{r+1}}+\cdots$. We use that $g_r+m\in \Ga$ and
$g_{r+1}=g_r+1$. \sq 
\vskip 0.2cm

Let us finalize $(ii)$ of the theorem. 
Upon the cancelations of the contributions of the terms 
from Lemma \ref{lem:Der} $(ii)$ 
to $\h^{mot}_\dag$, their 
total contribution becomes
$-t-q^{\de-1}t^{\de}$, unless for  $\de=1$ when it is $-t$.
As we already mentioned, an extension
of Parts $(ii)$ and $(iii)$ of the lemma is
expected to be applicable
to Conjecture \ref{conj:weak}. 

\section{\sc Superpolynomials in any ranks}\label{sec:ranks}
This is when the characteristic polynomial in $ASF$ becomes
a power of irreducible $F(x,y)$, which corresponds to considering
the corresponding algebraic knots colored by rows
in the DAHA setting. 
  
\subsection{\bf Basic definitions}
Let us consider a liner space $\Om\equal \o^{\oplus c}=\oplus_{i=1}^c
\o\vep_i$ over $\o$, where the basis 
$\vep_1,\ldots, \vep_c$ will be fixed in this section. The
{\em valuation} in $\Om$ will be
$\nu(\sum_{i=1}^c f_i \vep_i)=\min_{i=1}^c(\nu_z(f_i)+\nu^i)$,
where we take $0\le \nu^i=\nu(\vep_i)<1$ and assume that 
$\nu^1<\nu^2<\cdots<\nu^c$. Then $\nu(\Om)=\cup_{i=1}^c 
\{\nu^i+\Z_+\}$, $\Ga^c\equal \nu(\r^{\oplus c})=
\cup_{i=1}^{c}\{\nu^i +\Ga\}$.
Actually, $\nu^i$ can be arbitrary here such
that  $\nu^i-\nu^j\not\in \Z$ for $i\neq j$ (see also below).

For any $\r$-submodule
$M\subset \Om$: 
dim$_{\F}(\Om/M)=|G^c(M)|$, where $G^c(M)$, the {\em set of gaps},
is defined as $\nu(\Om)\setminus\De^{\!c}(M)$ 
for $\De^{\!c}(M)\equal \nu(M)$. 
The latter are $\Ga$-modules. 
These definitions are fully parallel to the case of $\o$-rank one;
let $\nu^1=0$ for the compatibility with the latter. 

The {\em standard modules} $M\subset \Om$
are $\r$-invariant  such that $\o\,M=\Om$. 
Equivalently, $\De^{\!c}(M)$ must be standard, where
a $\Ga$-module  $\De^{\!c}\subset \nu(\Om)$ is called standard
if $\{\nu^i\}\subset \De^{\!c}$. The number of such $\De$
 is $(std\,_\Ga)^c$,
where $std\,_\Ga$ is the number of standard $\Ga$-modules
$\De\subset \Z_+$.

Let  $dev(M)=\de\, c-deg(M)$ for $deg(M)=\,
\text{dim}_{\F}(\Om/M)=|D^c(M)|$. Here
$D^c(M)\equal\De^{\!c}(M)\setminus \Ga^c$, which is 
the set of {\em added gaps} from $\Ga^c$ to $\De^c(M)$.   
In particular, $|G^c|=\de\, c=dev(\Om)=$\,dim$_{\F}\Om/\r^{\oplus c}$
for $G^c\equal D^c(\Om)=\cup_{i=1}^{c}\{\nu^i +G\}$,
where $G=\Z_+\setminus \Ga$ as above.

The {\em invertible modules} are by definition
 standard ones with $c$ generators
over $\r$. Equivalently, 
$\De^{\!c}(M)=\Ga^c$,
which is obviously the smallest standard $\Ga$-module.
 Equivalently, invertible modules are standard of deviation $0$
or standard such that  $deg(M)=\de\,c.$

\vskip 0.2cm
Let 
$\{e_g, g\in D^c \}$ be a system of any elements $e_g\in M$
such that $\nu(e_g)=g$, where $D^c=D^c(M)$ for a given standard $M$.
They generate $M$ as an $\r$-module.
For any {\em invertible} submodule
$I\subset M$, its  $\r$-generators can be taken as follows:
$\ep^i=e_{\nu^i}+\sum_{g\ge 1}a_g^i e_g$, 
where  $1\le i \le c$,\, $g\in D^c$ and  $a_g^i\in \F$. 
Recall that $0\le \nu^i<1$ for $1\le i\le c$.

\begin{lemma}\label{lem:inv}
Different choices of $\{a_g^i\}$, which can be arbitrary,
 result in different $I$.
Therefore, the set $\{I\subset M\}$ is naturally an 
affine space $\mathbb{A}^{dim}$ (over $\F$) for $dim=dev(M)\,c$,
which identification  
depends on the choice of  $\{e_g, g\in D^c(M)\}$.  
In particular, the set  $Jac^{\,c}$ of all invertible modules
is isomorphic to  the affine space $\mathbb{A}^{\de c^2}$. 
\end{lemma}
{\em Proof.} The coefficients $a_g^i$ for $g\ge 1$
can be arbitrary because we can
add any such $e_g\in M$ to the generators of $I$.
 If $I$ contains $2$ elements in such form 
with the same $e_{\nu^i}$, then it contains
$\ep=e_g+\sum_{h>g} b_{h}^i e_{h}\in I$
for some $g,h\in D^c$, and  $\nu(\ep)=g\not\in \Ga^c$
due to the definition of $D$-sets.
This is impossible since
$\De^{\!c}(I)=\Ga^c$ for invertible $I$. 
 \sq
\vskip 0.2cm

Alternatively, the second claim can be verified
using the natural action of $GL_c(\o)$,
the group of invertible matrices $A$ with the entries in  $\o$,
on standard modules. These matrices are considered as
$\o$-linear transformation of $\Om$;
the $i${\tiny th} column of $A$ gives the image
of $\vep_i$. Obviously, standard and invertible modules
go to standard and invertible ones. This action is
transitive on the set of invertible modules. The stabilizer
of $\r^{\oplus c}\subset \Om$ is $GL_c(\r)$ and 
$Jac^{\,c}=GL_c(\o)/GL_c(\r)\cong \F^{\,\de c^2}.$
\vskip 0.2cm

We note that the dimension $\de c^2$ from the Lemma
is a singular counterpart of the
famous dimension formula $(g-1)c^2+1$ for the space of stable 
bundles of rank $c$ over a smooth curve of genus $g\ge 2$.
We do not ``divide" by the action of $GL_c$, which explains 
the deviation.
\vskip 0.2cm

The definition of {\em standard
flags} $\vect{M}=\{M_0,\ldots,M_{\ell}\}$ is 
a direct extension of that for $c=1$. They are such that 
$\De^{\!c}_i=\De^{\!c}_{i-1}\cup\{g_i\}$, where
$g_i<g_{i+1}$ and $\De^{\!c}_i=\De^{\!c}(M_i)$ for $1\le i\le \ell$.
Now $g_i$ are from the set
$G^c=\cup_{i=1}^{c}\{\nu^i +G\}$ of cardinality $\de c$
for $G=\Z_+\setminus \Ga$ (defined above). 

Accordingly, the Piontkowski cells are the schemes 
$\j_\ell(\vect{\De}^c)$
 of standard flags $\vect{M}$ 
such that  $\De^{\!c}(\vect{M})=\vect{\De}^c$. 
The motivic superpolynomials
will not depend on the ordering of $\nu^i$, but such cells
depend on it. Generally, they can be arbitrary real numbers,
not only from  $[0,1)$,
such that $\nu^i-\nu^j\not\in \Z$ for $\i\neq j$. 
The passage from one such ordering to another 
is an example of {\em wall crossing}. Similarly to the case of
$c=1$, 
all $\j_\ell(\vect{\De}^c)$ can be combined as strata in 
a single variety; see \cite{ChP2}. Study of its algebraic closure 
(it is not projective) and the
closures of $\j_\ell(\vect{\De}^c)$ is an interesting new direction.  

\subsection{\bf First decomposition theorem}
Upon these adjustments, the definition of motivic superpolynomials
remains basically the same. Following \cite{ChP2} and
Theorem \ref{thm:a-formula}, we arrive at the following 
definition-theorem. Recall that $deg(M)\!=$dim\,$_{\F_q}(\Om/M)$
and $\prod_{i=1}^0\!=\!1$. Now $\F\!=\!\F_q$.

\begin{theorem}\label{thm:a-formula-c}
Considering the summation over all standard $\ell$-flags 
$\vect{M}=\{M_0\subset M_1\subset\cdots\subset M_\ell\}\subset \Om$
over $\F=\F_q$, let
\begin{align}\label{motc-def}
&\h^{mot}=\h^{mot}(q,t,a, c;\r)\equal\!\!\!
\sum_{\{M_0\subset\cdots\subset
 M_\ell\}}t^{deg(M_\ell)}a^\ell,
\end{align}
which sum is finite. 
Alternatively:
\begin{align}\label{hmotrk-c}
\h^{mot}(q,t,a, c;\r)\,=\,\sum_{M} t^{deg(M)}\prod_{i=c}^{r\!k_q(M)-1}
(1+q^i a),
\end{align}
where the summation is over all standard $M$, and 
$\r$ can be any singularity ring, not only Gorenstein or
with $2$ generators. 
Recall that the $q$-ranks of such $M$  are
$\ge c$. 
\end{theorem}
{\it Proof.} First of all, let us justify that the
sums in  Theorem \ref{thm:a-formula-c}
are finite. 
Any standard module contains $z^{2\de}\Om$
because $\r$ contains $(z^{2\de})$, the conductor. Indeed
$M$ contains $\ep^{\,i}\in \vep_i+z\o$ for any $1\le i\le c$
and all their $\r$-linear combinations. In particular,
it contains $(z^{2\de})\ep^i$ and, therefore, the whole
$z^{2\de}\Om$.  Thus, all modules and flags can be
identified with their images 
in the space $\Om/ z^{\de}\Om$, which is finite-dimensional.

\vfil
Let us mention that this consideration is the starting point when
supplying $\{M\}$ with a structure of (a single) 
quasi-algebraic variety; see \cite{ChP2}. However, the corresponding
construction is more involved than that for $c=1$. 
\vfil

The first step toward the justification of
formula (\ref{hmotrk-c}) is the relation 
$\mathfrak{m}_\r M_{\ell}\subset M_0$, which then must
hold for all $M_i$. Let us prove it.

Recall that $\mathfrak{m}_\r=\r\cap z\o$, and
$\De^{\!c}_i=\De^{\!c}_{i-1}\cup\{g_i\}$, where $g_i<g_{i+1}$
and $\De^{\!c}_i=\De^{\!c}(M_i)$ for $1\le i\le \ell$.
All $g_i$ are
{\em primitive}: not in the form $\ga+g_j$ for any $j$ and
$\ga\in \Ga\setminus \{0\}$. Indeed, if $g_i=\ga+g_j$ then
$j<i$ due to the inequalities we imposed. Then
$D^c_{j}$ is obtained from $D^c_{j-1}$
by adding both, $g_j$ and $g_i$, which is impossible.

We obtain that
 $\{g_i,1\le i\le \ell\}=\De^{\!c}(M_\ell)\setminus \De^{\!c}(M_0)$ 
does not intersect
$\ga+\{g_i,1\le i\le \ell\}$ for any $0\neq \ga\in \Ga$.
Thus, the  $\ga$-translation of this set, which belongs to 
$\De^{\!c}(M_\ell)$ by
construction, actually belongs to $D^c(M_0)$. The same holds
if $\{g_i,1\le i\le \ell\}$ is extended by the
generators of $\De^{\!c}(M_0)$.
Thus, $\ga+\De^{\!c}(M_\ell)\subset \De^{\!c}(M_0)$ for 
$\Ga\ni \ga\neq 0$. Therefore
$\mathfrak{m}_\r M_{\ell}\subset M_0$. 
\vfil

Next, the space of all standard flags for a fixed 
$M_0$ is biregular to $\mathbb{A}^{\ell(\ell-1)/2}$, which is
exactly the same formula as it was for $c=1$; no change here. 
Thus, we only need to interpret the conditions for $M_0$
that make it {\em standard}. This is as follows.
\vskip 0.2cm

Let $M\!=\!M_\ell, 
M^{\{1\}}\!=\!M\cap z\Om,\, M_0^{\{1\}}\!=\!M_0\cap z\Om$,
$\overline{M}^{\{1\}}\!=\!
M^{\{1\}}/\mathfrak{m}_\r M$,  and $\overline{M}_0^{\{1\}}\!=\!
M_0^{\{1\}}/\mathfrak{m}_\r M$. For $n\equal r\!k_q(M)$, one has:
dim\,$\overline{M}/\overline{M}^{\{1\}}=c$ (since $M$ is standard)
and dim\,$\overline{M}^{\{1\}}=n-c$. 
 
The condition that  $\overline{M}_0$ is standard reads: 
$\overline{M}_0+\overline{M}^{\{1\}}=\overline{M}$. Equivalently, 
dim\,$\overline{M}/\overline{M}_0=$\,
dim\,$\overline{M}^{\{1\}}/\overline{M}_0^{\{1\}}=\ell$. 
Combining the latter  with
the equalities above, we obtain that
$\overline{M}_0^{\{1\}}=n-c-\ell$. 

Let us fix the subspace
$\overline{M}_0^{\{1\}}$ inside $\overline{M}^{\{1\}}$. The number
of choice for this subspace is {\footnotesize
$\qbin{n-c}{n-c-\ell}=\qbin{n-c}{\ell}$}. For each such choice,
we need to pick $\overline{\vep}_k+\sum_{i=1}^\ell 
a^k_i\, \overline{e}_{g_i}\in
\overline{M}_0$, 
where $1\le k\le c$. Here 
 $e_{g_i}\in M$ are such that $\nu(e_{g_i})=g_i$ and
$\overline{\vep}_k, \overline{e}_{g_i}$ are their images in $\overline{M}$. 
Different choices of the coefficients $a_i^k\in \F_q$, which can be
arbitrary, result in 
different $\overline{M}_0$, provided that $\overline{M}^{\{1\}}_0$
is fixed.

We obtain that there are $q^{\,c\,\ell}${\footnotesize 
$\qbin{n-c}{\ell}$} choices for $\overline{M}_0$, and  there
are $q^{\frac{\ell(\ell-1)}{2}}$ continuations of 
$\overline{M}_0$ to standard flags with $M=M_\ell$. Thus,
$\h^{mot}$ from (\ref{motc-def}) becomes
\begin{align*}
&\sum_{\ell=0}^{n-c} 
q^{\frac{\ell(\ell-1)}{2}} q^{\,c\,\ell}\text{\footnotesize 
$\qbin{n-c}{\ell}$} a^{\ell}
=
(1+q^{c}a)(1+q^{c+1}a)&\cdots (1+q^{n-1}a), 
\end{align*}
where 
$r\!k_q(M)=n$.
As above, the last product is $1$ if $n=r\!k(M)=c$, i.e. for 
invertible $M$. Then we add the summation $\sum_{M} t^{deg(M}$ over all
standard $M$ and obtain the required formula. 
\sq
\vskip 0.2cm

The greatest standard module is $\Om$. One has:
$r\!k_q(\Om)=$\,dim\,$\Om/\mathfrak{m}_\r\Om=m\,c$
for  the maximal ideal $\mathfrak{m}_\r=
\r\cap z\o$ in $\r$ and $m=\min\{\Ga\setminus \{0\}\}$. 
The module $M=\Om$  
contributes $(1+a q^c)\cdots (1+a q^{m\,c-1})$, the longest such
product. Accordingly, $deg_a \h^{mot}=r\!k_q(\Om)-c=(m-1)c$, though
there can be several families of standard modules of maximal
$q$-rank. 

When $a=-q^{-c}$, all products $\prod_{i=c}^{r\!k_q-1}(1+a q^i)$  
vanish unless $r\!k_q=c$, and  we obtain the contribution of invertible
modules, which is $q^{\de\,c^2}t^{\de\,c}$.

\vskip 0.2cm
We note that similar calculations involving the $q$-binomial
theorem were used in Conjecture 2 and Proposition 3 of \cite{ORS}
in the context of {\em nested} Hilbert schemes of plane
curve singularities. This proposition
states that the corresponding {\em ORS series} can be made a 
Laurent polynomial
upon the multiplication by $(1-t)$ in the unibranch case (in our
notation), and the functional equation upon
$t\mapsto 1/(qt)$ holds.

 When $a=0$, their series and these facts
are direct counterparts of those for the $L$-functions from
\cite{Gal,Sto}:\, for any Gorenstein rings over $\mathbb{F}_q$ 
and the classical  motivic measure $X\mapsto |X(\mathbb{F}_q)|$.
The proof of the functional equation from \cite{Sto} can be
adjusted to incorporate the dependence of $a$, but this can be
done only for plane curve singularities. See Theorem
\ref{thm-feq} below.

\vskip 0.2cm

{\sf Colored trefoil.} Let us calculate $\h^{mot}$
for $\r=\F_q[[z^2,z^3]]$ for arbitrary $c\ge 1$; recall that
$\de=1, m=2$ and $deg_a=(m-1)c=c$ in this case. Any standard
module contains $z^{2}\Om$ because $\r\supset (z^{2\de})$
(the conductor). In the case of trefoil,
the condition  $M\supset (z^{2\de})$
is sufficient for being an $\r$-submodule. Thus,
standard $M$ are uniquely determined by
$\F_q$-linear subspaces $\tilde{M}=M/(z^2\Om)\subset
\tilde{\Om}=\Om/(z^2 \Om)$ such that $\tilde{M}/(z\Om)$ 
coincides with $\Om/(z \Om)$. We obtain that $M$ can be described by
$\tilde{M}'=\tilde{M}\cap (z\Om)$ and vectors
$\tilde{\ep}^i=\vep_i+\sum_{i=1}^c a^i_j (z\vep_j)$ 
in $\tilde{M}$, 
where the coefficients are arbitrary
and different $\{a^i_j\}\subset \F_q$ correspond to different $M$.
This gives that  the number of standard $M$  is $q^{c\, k}$ when
$\tilde{M}'$ is fixed; the number of choice for the latter
is {\tiny
$\qbin{c}{k}$}. Here
$k=$\,dim\,$_{\F_q}\tilde{M}'$  (the $q$-rank)
can be from $0$ to $c$, the $q$-rank
of $M$ is $2c-k$, and the length of the $a$-products
is $c-k$. We arrive at:
 \begin{align*}
&\h^{mot}_{3,2}(c)=\sum_{k=0}^{c}q^{c\,k} t^k 
\text{\footnotesize 
$\qbin{c}{k}$}
(1+q^{c}a)(1+q^{c+1}a)&\cdots (1+q^{2c-k-1}a). 
\end{align*}
It coincides with 
$\h^{daha}(c\om_1)$. From now on, we will 
mostly use the notation $\h^{daha}$ instead of $\hat{\h}$
for DAHA superpolynomials.

Interestingly, the formula (2.52) for $\h^{daha}$
from \cite{CJJ} was presented differently:
\begin{align}\label{trefa-t}
&\h^{daha}_{3,2}(c\om_1)=\sum_{k=0}^{c}q^{c\,k} t^k 
\text{\footnotesize 
$\qbin{c}{k}$}
(1+\frac{a}{t}) \cdots (1+\frac{a}{t}q^{k-1}). 
\end{align}
It  was calculated there via the DAHA-Jones invariants
for the root system $(C\!\vee\! C)_1$. This
root systems results in one $q$ and four $t$-parameters. 
Their certain 
specialization provides the DAHA superpolynomials
for $T(2n+1,2)$ colored by $c$ (a theorem). The  
DAHA-formula there appeared coinciding
 with the formula predicted 
in \cite{DMMSS,FGS}. 

When $a=-t^2, t=q$. i.e. in the case of quantum $sl_2$, 
formula (\ref{trefa-t})  becomes
that due to Habiro. Generally, we obtain $2$ different
decompositions for 
$\h^{daha}_{3,2}(c\om_1)=\h^{mot}_{3,2}(c)$. The corresponding
identity can be deduced from the $q$-binomial theorem,
but this is not immediate. The following general theorem 
fully clarifies  (\ref{trefa-t}). 

We note that both formulas can be naturally extended to $c\in \C$
from their values at $c\in \Z_+$. The standard formulas 
for the $q$-binomial coefficients 
$\text{\tiny $\qbin{c}{k}$}$ in terms of
$(x;q)_\infty$  can be used.

\subsection{\bf Second decomposition theorem}
Let us generalize (\ref{trefa-t}). We note that
the $1${\small st} claim of the next theorem readily 
results in the identity
$\h^{mot}(c)(q,t,a=-t)=1$. This is manifest 
for $\h^{daha}(c\om_1)$ because
$\hat{J\!D}^{A_0}=1$ for $c\om_1$ (and only for such weights).
However, the DAHA meaning of this decomposition theorem 
and the previous one remains unknown beyond some examples.

\begin{theorem}\label{thm:2ndexp}
For any standard $M$, we define an $\r$-module
 $M^\Diamond\equal\{f\in \Om \mid \mathfrak{m}_\r f
\subset M\}$ and set
 $d^\Diamond(M)=$\,dim$_{\F_q}(M^\Diamond/M)$.
For instance, $M^\Diamond=
\Om$ if $M\supset \mathfrak{m}_\r$, and $M^\Diamond=M$ if and only
if $M=\Om$. The claims (i,ii) below will be 
for arbitrary singularity rings $\r$
(not only for those with $2$ generators or Gorenstein). 

(i) Given standard  $M$,
the contribution of the flags with  $M_0=M$ to $\h^{mot}$ is
$t^{deg(M)}\,
(1+\frac{a}{t}) (1+\frac{a}{t}q)\cdots (1+\frac{a}{t}q^{d-1})$,
so it depends only on $deg(M)$ and $d=d^\Diamond(M)$. We arrive
at the following presentation:
\begin{align}\label{2ndHfomula}
\h^{mot}(c)= \sum_{M} t^{deg(M)}\,
\Bigl(1+\frac{a}{t}\Bigr) \Bigl(1+\frac{a}{t}q\Bigr)\cdots 
\Bigl(1+\frac{a}{t}q^{d^\Diamond(M)-1}\Bigr),
\end{align}
where the summation is over all standard $M\subset \Om=\o^{\oplus c}$.

(ii) Assuming that the flags are inside a standard $\r$-module $N$
and $\mathfrak{m}_\r N\subset M$, 
the contribution of such  flags with $M_0=M$ will be then 
$\sum_{d=0}^{r}\, t^{deg(N)+d} \, q^{c\, d}\text{\footnotesize 
$\qbin{r}{d}$}
(1+\frac{a}{t}) (1+\frac{a}{t}q)\cdots (1+\frac{a}{t}q^{d-1})$
for $d=$\,dim\,$_{\F_q} N/M$, $r=r\!k_q(N)\!-\!c$  and \,$deg(N)+d=deg(M)$.

(iii) In the case of $\r=\F_q[[z^2,z^3]]: M^\Diamond=\Om$ 
for any standard $M$. Thus, we can take $N=\Om$ and 
$r=c(m-1)=c$ in (ii), which gives:
\begin{align*}
&\h^{mot}(c)=\sum_{d=0}^{c} q^{c\,d}\, t^d
\text{\footnotesize 
$\qbin{c}{d}$}
(1+\frac{a}{t}) \cdots (1+\frac{a}{t}q^{d-1}). 
\end{align*}
\end{theorem}
{\it Proof}. The equivalence $\{M^\Diamond = M\} \Leftrightarrow 
\{M=\Om\}$ for standard $M$ is as follows. If $M\neq \Om$, then let
$g'=\max\{g\in \nu(\Om)\setminus \De(M)\}$ and $v'\in \Om$ 
be any element such that $\nu(v')=g'$. Then $m_\r v'\subset
M$ and $v'\in M^{\Diamond}$, which is impossible.

Claim $(i)$ is justified following the proof of
Theorem \ref{thm:a-formula-c}, where we now fix $M_0$ instead
of $M_{\ell}$ and calculate the number of the
corresponding standard flags. Concerning the number of flags
inside $\mathfrak{m}_\r N\subset N$, 
let $\hat{M}\subset \hat{N}\subset \hat{\Om}$ be the images of
$M\subset N\subset \Om$ 
modulo $\mathfrak{m}_\r N$; the latter belongs to $M$.
Also,  let $\hat{M}'\subset N^{\prime}$ be the images of
their intersections with $\Om'=z\Om$, which
contains  $\mathfrak{m}_\r N$.

Using that $N,M$ are standard, 
dim\,$_{\F_q}(\hat{N}^{\prime}/\hat{M}')=d=$ dim$_{\F_q} 
(N/M)$ and 
dim\,$\hat{M}/\hat{M}^{\prime}=c=$\,dim\,$\hat{N}/\hat{N}^{\prime}$. 
Accordingly, the number
of all possible $\hat{M}^{\prime}$ of codimension $d$ inside
$\hat{N}^{\prime}$, which is of $q$-rank $r$, is
{\footnotesize $\qbin{r}{d}$}. Then the number of lifts of
$\hat{M}^{\prime}$ to $\hat{M}\subset \hat{N}$ is $q^{c\, d}$. 

The case of $\r=\F_q[[z^2,z^3]]$ is exceptional because
all standard $M$ contain $\mathfrak{m}_\r \Om=z^2\Om$.
Thus, $M^\Diamond=\Om$, $r(M)=c$ for all of them, and 
$\hat{M}^{\Diamond\,\prime}=(z\Om)/(z^2\Om)$. This gives the
required formula. \sq
\vskip 0.2cm

\Yboxdim5pt
{\sf Trefoil of rank 2.} 
It makes sense to show explicitly how the valuation $\nu$,
$\Ga$-modules and Piontkowski cells are used in this case;
we will consider only $c=2$. 
The corresponding DAHA superpolynomial from \cite{CJ} is 
{\small 
$
\h^{daha}_{3,2}(2\om_1)=
1  + q^2 t + q^3 t + q^4 t^2 + a(q^2 + q^3 + q^4 t + q^5 t)+
a^2 q^5$.} 
\Yboxdim7pt

Its motivic counterpart is as follow; see \cite{ChP2}.
For $\r=\F_q[[z^2,z^3]]$ and $\Ga=\{0,2,3,\ldots\}$,
let $\nu^1=0, \nu^2=\frac{1}{2}$. Since any standard $M$ contains 
$z^2\Om$, there are four standard $\Ga$-modules
$\De^{\!c}$, which are:
\begin{align*}
&\De_{tot}=\{\Z_+, \frac{1}{2}+\Z_+\},\ 
\De'=\{\Z_+, \frac{1}{2}, \frac{1}{2}+2, \frac{1}{2}+3,\ldots\},\ 
\De''=\\
\{0&\,,2,3,\ldots, \frac{1}{2}+\Z_+\ldots\},\ 
\De_{inv}=\{0,2,3,\ldots, \frac{1}{2}, \frac{1}{2}+2, 
\frac{1}{2}+3,\ldots\}.
\end{align*}

The corresponding  $\tilde{M}=M/z^2\Om$ are  as follows:
$\tilde{M}_{tot}=\Om/z^2\Om$, 
$\tilde{M}'(\al,\be)=\lan \vep_1+\al z\vep_1, 
\vep_2+\be z \vep_1, z\vep_2\ran $,
$\tilde{M}''(\al,\be,\la)=\lan \vep_1+\al z\vep_2, 
\vep_2+\be z \vep_2, z \vep_1 +\la z \vep_2\ran$,
$\tilde{M}_{inv}(\al,\be,\la,\mu)=\lan \vep_1+\al z\vep_1+\be z\vep_2, 
\vep_2+ \la z \vep_1 + \mu z \vep_2\ran$. 
Here $\al,\be,\la,\mu\in \F_q$.
The span $\lan\cdot\ran$ must be generally over $\r$, but
the span over $\F_q$ modulo
$z^2\Om$ is sufficient here (only for trefoil).

Notice that the first terms in the vectors above are 
those with minimal possible $\nu$ and the other terms must have
higher $\nu$. This explains the difference
between $\tilde{M}'$ and $\tilde{M}''$; the corresponding 
dimensions are different.
We obtain that the dimensions
of $\j_0(\De)$ are $0,2,3,4$ for 
$\De=\De_{tot}, \De',\De'', \De_{inv}$.

The $1$-flags 
are $\tilde{M}'\subset \tilde{M}_{tot}$, 
$\tilde{M}''\subset \tilde{M}_{tot}$, 
and the
extensions of $\tilde{M}_{inv}(\al,\be,\la,\mu)$ to $\tilde{M}'$ or 
$\tilde{M}''$
by $z\vep_2$ or by $z\vep_1+\ga z\vep_2$. The corresponding
cells are $\mathbb{A}^d$ for $d=2,3,4,5$. 
The $2$-flags can be only of type
$\tilde{M}_{inv}\subset \tilde{M}''\subset \Om$; they constitute  
$\mathbb{A}^5$.
\vskip 0.2cm

This calculation can be extended to bottom
and top $t$-terms for arbitrary $\r$. For the sake of simplicity, 
let $a=0$ (no flags). Then the lowest $t$-term is 
$1$, which is for $M=\Om$, 
the top $t$-term of $\h^{mot}$ at $a=0$
is  $q^{c\de^2} t^{c\de}$, where $c\de^2$ 
is the dimension of the cell of
invertible modules. 

Next, the coefficient of $t^{c\de -1}$ 
counts the modules obtained from $\Ga^c$ by {\em adding} one
gap in the form 
$\nu^i+(\de-1)$ for $i=1,2,\ldots, c$.  The corresponding
contribution to $\h^{mot}$ will be then  $t^{c\de -1}
(q^{c(\de^2-1)}+q^{c(\de^2-1)+1}+\cdots+
q^{c\de^2-1})$. Similarly, the coefficient of $t$ will be 
$(q^c+q^{c+1}+\cdots+q^{c\,m-1})$, where 
$m=\min\{0\neq \ga\in\Ga\}$ is the
multiplicity of the singularity.
 The number of $q$-powers
here is $c(m-1)$, which is 
due to {\em removing}  $\nu^i+g$ (one such value) from 
$\nu(\Om)$, where $0<g<\mu$
and $1\le i\le c$.

\subsection{\bf The case of 
\texorpdfstring{{\mathversion{bold}$\,T(2n+1,2)$}}{ T(2n+1,2)}}
Let us discuss general $T(2n+1,2)$ colored by
$c\om_1$. We have:\, $\r=\F_q[[x\!=\!z^2,y\!=\!z^{2n+1}]]$,\,
 $\de=n$,\,
$r\!k_q(\Om)=m\, c=2c$,\, and 
$deg_a\h^{mot}=(m-1)\, c=c$ because $m=2$ (as for the trefoil).
The contribution of invertible modules is
$t^{\de\, c}q^{\de\,c^2}=t^{2n}q^{4n}$.

We will discuss only  $c\!=\!2$, i.e. 
the standard modules of $\o$-rank $2$. An important
(and predictable) feature of 
the formulas provided below is their stabilization with respect
to $n$, which corresponds to the inclusions
$\F_q[[x=z^2,y=z^{2n+1}]]\supset \F_q[[x=z^2,yx=z^{2n+3}]]$.
\vskip 0.2cm

Let $\h^{mot}_{2n+1,2}(c\!=\!2)=
A_0+A_1(1+q^2 a)+A_2(1+q^2 a)(1+q^3 a^2)$.
Then,
{\small
\begin{align}\label{2c2n+1} 
&A_0=t^{2n}q^{4n},\ A_1=
t^n q^{3n-1}(1+q)\bigl(1+qt+q^2t^2+\cdots+q^{n-1}t^{n-1}\bigr),\\ 
&A_2=\sum_{i=0}^{n-1} t^i q^{2i}\bigl(1+q+\cdots+q^i\bigr)+
\sum_{i=n}^{2n-2} t^i q^{2i}\bigl(1+q+\cdots+q^{2n-2-i}\bigr). \notag
\end{align}
}

We present these terms in this form to demonstrate the stabilization
with respect to $n$. For the term $A_2$, we have pure uniformization
and the embeddings of the corresponding formulas for different $n$ 
(no rescaling is needed). We will omit the calculations.

The corresponding formulas can be obtained
for any $c$. The approach from
\cite{CJJ} is based
on the recurrence relations (2.54):\, the ``stabilization"
with respect to $c$ in the decomposition in terms
of $(1+ q^ia/t)$.  The coincidence with motivic
superpolynomials was justified (theoretically) for several 
$2n+1, c$, but not in full generality
by now. 
\vskip 0.2cm

{\sf Uncolored examples.}
The simplest example beyond the trefoil is uncolored $r=\F_q[[z^2,z^5]]$.
Then $\de\!=\!2,\, \Ga\!=\!\{0,2,4,5,\ldots\},\, G\!=\!\{1,3\}$. The
Piontkowski cells are for  $D_0=[1,3], D_1=[3], D_3=[\emptyset$],
where $deg_0\!=\!0, deg_1\!=\!1, deg_2\!=\!2$, and the number of modules
in the corresponding cells are $1, q^1, q^2$. Here and
below we use the brackets $[\cdots]$ for the $D$-sets.

The corresponding
$q$-ranks are $r\!k_q=2, 1, 1$. One has $M_0^\Diamond=M_1^\Diamond=
\o$, and $M_2^\Diamond=\lan 1+\al z, z^3\ran$ for (invertible)
$M_1=\lan 1+\al z+\be z^3\ran$. Accordingly, $d^\Diamond_0=0$ and
$d^\Diamond_1=d^\Diamond_2=1$. The corresponding decompositions are:
$\h^{mot}_{5,2}=(1+qt)(1+qa)+q^2t^2=1+(qt+q^2t^2)(1+a/t)$.  

It is almost equally simple to justify  that $\h^{mot}_{2n+1,2}=$
{\small
$$(qt)^n+(1+qt+\cdots+(qt)^{n-1})(1+aq)=1+
(qt+\cdots+(qt)^n)(1+a/t).
$$
}

Generally, the coefficients of $q^it^j$ can be negative in the 
decomposition via the products
$\prod_{k=0}^{r\!k_q-1}(1+q^k a/t)$ for $0\le deg_a(\h^mot)= (m-1)+c$.
 The simplest {\em uncolored} 
torus example
we found is $T(7,6)$ for $\r=\F_q[[z^6,z^7]].$ Namely,
the coefficient
of $(1+a/t)(1+qa/t)$ contains $-q^{11}t^{10}$. All other $q^it^j$
in this decomposition are with non-negative coefficients. 
This is an indirect confirmation that
the $q$-rank of $M^\Diamond$ is not always constant in Piontkowski 
cells (even) for torus knots.

Similarly,
such negative coefficients occur in the decompositions
(the main for us) with respect to the products 
$\prod_{k}(1+q^k a)$ (always, for sufficiently large knots).
 For instance, two negative
terms occur for uncolored $T(9,4)$. We note that,
generally, the $q$-rank of $M$ or $M^\Diamond$
can increase when the limits of standard modules $M$
are considered. Thus, the method due to Bia{\lx}ynicki-Birula 
is not applicable:\, the
subvarieties of $\j_{\ell}(\vect{\De})$  with 
fixed  $r\!k(M)$ or fixed $d^\Diamond(M)$ are not always 
affine spaces (otherwise, all terms $q^it^j$ would be 
with positive coefficients). 

Let us complement the 
$\prod_{k\ge 1}(1+aq^k)$\~decomposition 
of $\h^{mot}_{7,3}$ provided above 
with that in terms of 
$\prod_{k\ge 0}(1+aq^k/t)$.
One has:
$\h^{mot}_{7,3}=$
{\small
$q^6 t^6+(1+a q) (1+a q^2) (1+q t+q^2 t+q^2 t^2+q^3 t^3)+
(1+a q) (q^3 t^2+q^4 t^2+q^4 t^3+q^4 t^4+q^5 t^4+q^5 t^5)
=1+(q t+q^2 t+q^3 t^2+q^4 t^2+q^5 t^4+q^6 t^6)(1+a/t)+
(q^2 t^2+q^3 t^3+q^4 t^3+q^4 t^4+q^5 t^5)(1+a/t)(1+qa/t)
$}.

\comment{
The coefficient of $\h^{mot}$ of the
decompositions with respect to the products 
 $t^d \prod_{i=c}^{r-1}(1+ q^ia/t)$ corresponds to counting $\ell$-flags 
$\vect{M}$ with $r\!k_q(M_0)=r-\ell$ and $deg(M_0)=d$ for any
$\ell$. This makes sense geometrically: 
 we count the number of all extension of all standard
 $M=M_0$ of degree $d$ to $M'=M_{\ell}$ or $q$-rank $r$
such that $\mathfrak{m}_\r M'\subset M$ (for all possible $\ell$). 
This interpretation can explain some remarkable
features of the decompositions for products above
with respect to $a/t$ instead of $a$. 

For instance, for uncolored trefoil, we obtain 
$(1+a/t)qt+1=1+qt+aq$. For trefoil colored by $c=2$,
\begin{align*}
&\h^{mot}_{3,2}(c\!=\!2))=\sum_{k=0}^{2}q^{2\,k} t^k 
\text{\footnotesize 
$\qbin{2}{k}$}
(1+q^2 a)\cdots (1+q^{4-k-1}a)\\
&=q^4t^2+q^2t(1+q)(1+q^2 a)+(1+q^2 a)(1+q^3 a)=
1+a^2 q^5+q^2 t+q^3 t+q^4 t^2+a (q^2+q^3+q^4 t+q^5 t).
\end{align*}
The number of all extensions of
invertible modules to $r\!k_q=4$, which can be only $\Om$,
 results in the formula for $a/t$
$t^2 q^{4}(1+a/t)(1+qa/t)$.
}

\subsection{\bf Discussion}
The usage of  ``motivic" in this paper  requires a comment.
The most abstract version of motivic superpolynomials is when 
the Piontkowski cells $\j_\ell(\vect{\De}^c)$ are considered
as their classes 
in the Grothendieck ring  $K_0(V\!ar/\mathbb{F})$ 
of varieties over $\F$.  The class 
$[\mathbb A^1]$ replaces $q$ and the coefficients of $t^i a^j$  
become from $K_0(V\!ar/\mathbb{F})$. This is 
a ``true" motivic definition. 
The rationale for this generalization is the conjecture that the
Piontkowski cells are {\em affine configuration spaces}:
some their intersections in bigger affine spaces and unions.
Compare with  Lemma 9 from \cite{ORS}.

The coefficients of such superpolynomials
contain (a lot of) information
beyond what our $\h^{mot}$ provide. 
The latter is the specialization of such abstract superpolynomials
under the {\em motivic
measure}\, $K_0(V\!ar/\mathbb{F}_q)\supset [X]\mapsto
|X(\mathbb F_q)|$. Other natural choices of motivic measures
are based on 
the usage of the {\em weight
filtration} as in \cite{ORS,GORS} 
or Borel-Moore homology as in \cite{ChD1,ChD2}. The polynomial
dependence on $q$ is granted  for the latter two 
measures.
\vfil

{\sf Polynomiality, etc.}
An important part of
the  {\em Coincidence Conjecture} is that 
$\h^{mot}$ is a polynomial in 
terms of $q,t,a$;  this is manifest for $t,a$, 
{\em but no for $q$.} Also,
its topological invariance is a conjecture.

Both claims are theorems for 
$\h^{daha}=\hat{\h}_{\l,'\!\l}$ 
in the notation from 
Theorem \ref{STABILIZ}, where the pair
$\{\l,'\!\l\}$ encodes the corresponding
algebraic link and the colors are arbitrary Young diagrams.
\vfil

The $q,t$-polynomial dependence 
of the $J\!D$-polynomials for torus knots and their topological 
invariance are from \cite{CJ,CJJ}.
These theorems are for arbitrary (reduced) root systems. 
See \cite{CJ, GoN} concerning the
 ``$a$-stabilization" of $J\!D$-polynomials of type $A$
to superpolynomials;
  \cite{SV} was used (for torus knots). Also, see \cite{CJ} for 
examples of the $a$-stabilization for the  $BCD$-series
and \cite{ChE} for the $E$-series from \cite{DG}. 
The stabilization is expected to hold for these series, 
but there is no such theory by now (only examples).  

The generalization of these properties to colored iterated 
torus links is due to \cite{ChD2}. 
We note that the topological invariance
of $J\!D$-polynomials and the corresponding theory of {\em DAHA-vertex}
for iterated torus {\em links} appeared directly related to
the generalized Cherednik-Macdonald-Mehta formula, one of the key 
tools in DAHA theory. 
\vfil

The uniform (and polynomial) dependence of 
$\h^{mot}(c=1)$ on $q$ and its
topological invariance for 
plane curve singularities was justified in \cite{ChG} 
when $\Ga$ has  $2$ generators (the case of torus knots),
or $3$ generators (with some restrictions). 
Generally, plane curve singularities 
depend on ``continuous" parameters, {\em modes} in physics
literature. 
The formal-analytic classification of such singularities
is essentially known; 
in formal series, the convergence is ignored. We used the
 elements of this classification in \cite{ChG} 
when considering $\Ga$ with $2$ and $3$ generators.

\vskip 0.2cm
The Coincidence Conjecture $\h^{daha}=\h^{mot}$ 
was checked in many examples, 
including various cases 
when not all 
Piontkowski cells are  affine spaces $\mathbb{A}^m$. 
If these cells
are affine spaces, the method of {\em syzygies}, generally,  
provides formulas
for their dimensions. All cells $\j_\ell(\vect{\De}^c)$ are
nonempty affine spaces for $T(\rr,\ss)$
colored by any $c$. Generally, the number of standard
modules is $(std)^c$,
where $std$ is that for $c=1$. We discuss  below
$std=\frac{1}{\rr+\ss}\binom{\rr+\ss}{\rr}$  for $T(\rr,\ss)$. Recall that
$\h^{daha}(q\!=\!1;c)=\h^{daha}(q\!=\!1;c\!=\!1)^c$, which is not known
for $\h^{mot}$ in full generality; see Section 4.1.3 from \cite{ChP2}.
 We note that non-affine cells
occur for $c>1$ even in the simplest non-torus case 
$\F[[z^4,z^6\!+\!z^7]]$, {\em ibid.} 

\vfil

{\sf Shuffle Conjecture and $E$-polynomials.}
The coincidence $\h^{daha}=\h^{mot}$ for $m\om_1$
is, actually, an advanced
version of the {\em Shuffle Conjecture}, proved in \cite{CaM},
where the left-hand
side is the ``combinatorial half"  and the right-hand side is the
``formula half". More exactly,  the weights $m\om_1$ in $\h^{daha}$
must be replaced by minuscule $\om_m$ due to the superduality.
The polynomials $P_{\om_m}$
are monomial symmetric functions (in type $A$).
The corresponding $E$-polynomials are simply monomials.
Switching to them,
the DAHA half of the
Coincidence Conjecture becomes entirely combinatorial.
The usage of Macdonald polynomials for arbitrary
(non-minuscule) weights makes the DAHA half involved 
combinatorially, and the motivic half is not known in such
generality. 

The passage to {\em nonsymmetric Macdonald $E$-polynomials} is
an important feature of \cite{CJ} and further works. They
are much simpler to deal with than their symmetrizations, the
$P$-polynomials. 
The $E$-polynomials are sufficient to obtain the $J\!D$-polynomials
and DAHA superpolynomials for torus iterated {\em knots} 
due to \cite{CJ}.  Namely, the invariants 
remain {\em exactly the same} 
if spherical $P^\circ_b$ are replaced by spherical $E^\circ_b$.
However, the  $P$-polynomials and $J$-polynomials
are still needed
for {\em links} (see below).  Given a link, 
only one of $J_\la$ 
can be replaced by the 
corresponding $E$-polynomial. If {\em all} $J_\la$ for the 
connected components 
are replaced by 
$E_\la$, then our construction works perfectly, but we
do not know the geometric
meaning of the resulting polynomials.
\vskip 0.2cm

{\sf Relations to $ASF$.}
The conjecture  $\h^{mot}(c)=\hat{\h}(c\om_1)$
from \cite{ChP1,ChP2} will be extended below to 
{\em non-unibranch} plane curve singularities. 
This is the case of {\em affine Springer 
fibers}, $ASF$,  of type $A$
with the most general characteristic polynomials:\, possibly reducible
and not square-free. Accordingly, this is the case of
reducible and non-reduced {\em spectral curves}
of the corresponding Hitchin fibers. 
It is interesting to extend the {\em Fundamental Lemma}
to this setting:\, with  $q,t,a$ and for arbitrary
``spectral curves". 

\vskip 0.2cm
 
{\sf Using Witt vectors.} Direct counterparts
of the motivic superpolynomials for knots
can be defined for any
$p$-adic integral domains $\o$ instead of $\F_q[[z]]$
and its {\em orders} $\r$,
subrings in $\o$ with the same localization field. 
Our definitions can be readily transferred to this setting.
The corresponding $p$-adic superpolynomials
count {\em standard} $\ell$-flags of $\r$-modules  
$M\subset \o^{\oplus c}$,
those containing a unit in $\o$,  
with the weights $t^{deg(M)} a^{\ell}$. Here $|\o/M|=q^{\deg(M)}$
for $\r/\mathfrak{m}=\F_q$, where  
$\mathfrak{m}\subset \r$ is the maximal ideal.
 The $q$-rank is as above:\,
$r\!k_q=\text{dim}_{\F_q}
M/\mathfrak{m} M$.

There will be no quasi-projective varieties and Gorenstein rings; 
{\em Witt vectors} will be used ``instead". However, both theories
all parallel. However, plane curve singularities
have $p$-adic counterparts:\, complete subrings in $\o$ with 
two generators (one of them can be $p$). 


For instance, let $\o=\Z_p[[\pi]]$ for the
$p$-adic $\Z_p$, 
$\pi^\ss=p$ and $\r=\Z_p[[x\!=\!p,y\!=\!\pi^\rr]]\subset \o$, where 
 $\rr,\ss>0$ and gcd$(\rr,\ss)=1=$gcd$(p,\ss$ 
(the tamely ramified case). 
 The corresponding
superpolynomials coincide with those  for
$\r=\F_p[[z^\ss,z^\rr]]\subset \o=\F_p[[z]]$ and for several other
families. The corrections due to Witt vectors do not influence
the output, but we do not conjecture this by now.

\vskip 0.2cm
\vfil

{\sf Geometric superpolynomials.}
We note that motivic superpolynomials 
are always {\em significantly} more convenient to
deal with and faster to calculate 
than the {\em flagged $L$-functions} defined below.
The Piontkowski cells are known in many examples; if they
are affine spaces then their dimensions 
dim$\,\j_\ell(\vect{\De})$ are sufficient. This is
always the case for quasi-homogeneous singularities.
The calculations of $\h^{mot}$
are almost immediate for  sufficiently small quasi-homogeneous
singularities.
 
Concerning $\h^{daha}$,
 practical finding DAHA superpolynomials is, generally,
more uniform than that for motivic
ones, and they exist in much greater generality. Due to the
$a$-stabilization needed for $\h^{daha}$, the calculations
can be involved. However, they
are, generally, faster than those for motivic ones,
unless all cells are affine spaces and
the exact formulas for their dimensions are used, 
 as those in \cite{Pi,ChP1,ChP2} .

\vskip 0.2cm
If  $\j_0$ is known to be a disjoint union of affine cells,  
then the coefficient of $q^i$
in  $\h^{mot}$ for $t\!=\!1,a\!=\!0$ is the {\em Betti number}
$b_{2i}\!=$\,rk\,$H_{2i}(\j_{0};\R)$ and $b_{2i+1}\!=\!0.$
In particular, $\h^{mot}(q\!=\!1,t\!=\!1,a\!=\!0)$ 
is the Euler number $e(\j_0)$. 
The simplest application is the formula for $e(\j_0)$
for $\r=\F[[z^\rr,z^\ss]]$, where gcd$(\rr,\ss)=1$.
It is the rational Catalan number 
$\frac{1}{\rr+\ss}\binom{\rr+\ss}{\rr}$ 
due to Beauville. This is 
the number of all standard $\De$ for such 
 $\r$, which are 1-1 with Dyck paths 
in the rectangles ``$r\times s$". 
\vskip 0.2cm

We conjectured in \cite{ChD1} that the relation
to Betti numbers of $\j_0$ always holds for the 
corresponding DAHA 
superpolynomials.
More generally, the conjecture is
that the {\em geometric superpolynomials} defined in terms
of {\em Borel-Moore homology} of $\j_{\ell}$ 
coincide with the DAHA superpolynomials for any algebraic knots. 
Basically, these theories correspond
to different choices of the {\em motivic
measure}. They coincide with our
motivic ones if $\j_{\ell}$ can be covered by affine spaces,
which follows from the definition of the {\em Borel-Moore homology}. 

\vskip 0.2cm

\section{\sc From knots to links} The consideration of
{\em non-unibranch} plane singularities 
colored by $c\om_1$ (pure rows) corresponds to the theory
of $ASF$ of type $A$ with arbitrary (not only irreducible or
square free) characteristic polynomials.
Also, such an extension is needed in Rosso-Jones-type 
formulas \cite{RJ}. A basic problem is
to reduce the superpolynomials
of {\em iterated} torus knots to 
torus ones. It cannot be done as such. However, 
Proposition 4.3 in \cite{ChW} expresses DAHA
superpolynomials of certain uncolored iterated torus knots 
in terms of uncolored
torus {\em links} and torus knots
colored by $c\om_1$.
In this paper, Proposition
\ref{thm:h-iter} is a motivic counterpart of 
one case of this DAHA formula. Generally,
iterated torus links colored by $c\om_1$ are
expected to be closed for such formulas.

\subsection{\bf Adding idempotents}
The ring will be now $\r\subset \o\equal 
\oplus_{i=1}^\kappa \o_i e_i$,
where $\o_i=\F[[z_i]]$ and $e_i e_j=\de_{ij} e_i$,
and $\kappa$ is the number of (absolutely)
irreducible components
of a singularity.
We set
$e\equal\sum_{i=1}^\kappa e_i$ (the unit in $\o$),
 $z\equal \sum_{i=1}^\kappa z_i e_i$ and identify $z_i$ with
$z e_i$. 
Generally, $f_i$ will be the projection $f e_i$ for any $f\in \o$. 
As above, $\r$ contains $1=e$ and have $2$ generators: 
$x=\sum_{i=1}^\kappa x_i$ and 
$y=\sum_{i=1}^\kappa y_i$ in $\mathfrak{m}_{\o}=z\o$. 
Also, the localizations of the projection 
$\r_i$ of $\r$ must be $\F((z_i))$, and $\r_i/\mathfrak{m}_i=\F$
for $\mathfrak{m}_i=\r_i\cap z_i\o_i$, the maximal ideals of $\r_i$.

By construction, 
$\prod_{i=1}^\kappa
F_i(x,y)=0$, where $F_i(u,v)$ are irreducible polynomials over $\F$,
assumed irreducible over its algebraic closure $\overline{\F}$, 
such that 
$F_i(x_i,y_i)=0$.  The polynomiality of $F_i(u,v)$ 
can be achieved by deforming $x_i$ and $y_i$ 
if necessary (the Weierstrass Preparation Theorem). 
The product $F(u,v)=
\prod_{i=1}^\kappa F_i(u,v)$ 
must be {\em square-free}, i.e. $F_i$ and $F_j$ are 
assumed non-proportional
for $i\neq j$. This is standard for curve singularities.
The general case (with multiplicities) will be addressed later via
arbitrary $\o$-{\em ranks}:\, the passage from $\o$ to $\Om$. 

Let us provide some general facts and definitions. A convenient
reference is \cite{Hil} with the following reservation:
there must be $\de=(\mu+\kappa-1)/2$
at the bottom of page 248 and throughout \cite{Hil}. 
For the sake of completeness,
the {\em Milnor number} is 
$\mu=$\,dim\,$_{\F} (\F[u,v]/(F_u,F_v)$ for the ideal $(F_u,F_v)$
in $\F[u,v]$ generated by the partial derivatives $F_u,F_v$
provided that
characteristic of $\F$ is $0$ or sufficiently
general $p>0$. Actually, we do not need it in this paper.

Let $\de=\de_\r\equal\text{dim}_{\F}(\o/\r)$. 
Then $\mu=2\de\!-\!\kappa\!+\!1$
due to Milnor (a topological proof) and Risler. To give
an example, consider {\em quasi-homogeneous} singularities,
which are for square-free 
$F(u,v)$ such that $F(\la^a u,\la^b v)=\la^N$ provided gcd$(a,b)=1$.
 Then, $\de=\frac{(N-a)(N-b)}{2ab}+\frac{\kappa-1}{2}$ and
$\mu=$\,dim\,$_{\F}\bigl(\F[u,v]/(F_u,F_v)\bigr)=\frac{(N-a)(N-b)}{ab}$.

In particular, let $F(u,v)=u^{m}-v^{n}$, where 
$n=\kappa a$, $m=\kappa b$,  gcd$(a,b)=1$, and $N=\kappa a b$; the
corresponding link is $T(n,m)$. 
\vskip 0.2cm

{\sf Intersection numbers}.
For any polynomials $f,g$ in terms of $u,v$, their 
{\em intersection number} is 
$I(f,g)\equal$\,dim\,$_{\F} \bigl(\F[u,v]/(f,g)\bigr)$. Accordingly,
 $I_{i,j}=I(F_i,F_j)$ for the components $1\le i\neq j\le \kappa$. 
Their main role for us is due to the formula:
\begin{align}\label{de-links}
&\de=\sum_{i=1}^\kappa \de_i+\sum_{i<j}I_{i,j}, \text{ where }
\de_i=\text{\,dim\,}_{\F}(\o_i/\r_i).
\end{align}
The $2${\small nd} sum here is 
dim\,$_{\F}\bigl(\oplus_{i=1}^\kappa \r_i\bigr)/\r$, which is for
the natural embedding $
\r\ni f\mapsto (f_1,\ldots,f_\kappa) \in 
\oplus_{i=1}^\kappa \r_i$.

The calculation of intersection numbers can be reduced to
proper valuations.
 Namely, 
 $I(f,g)=\nu_f(g)$, where $f$ is assumed 
irreducible over the algebraic closure $\overline{\F}$ of $\F$ 
and $g$ is not divisible by $f$. Here $\nu_f(g)=ord_z(g)$ for the 
uniformizing parameter $z$ for $f$ and the embedding
$\r_f=\F[[x,y]]/(f(x,y))\subset \o_f=\F[[z]]$, where  
the latter ring is the integral closure
of $\r_f$ in the field of its fractions; here $g$ is identified
with its image if $\o_f$.  To justify this formula, we calculate 
dim\,$_\F(\o_f/g\r_f)$ via 
$g\r_f\subset\r_f\subset \o_f$ and via
$g\r_f\subset g\o_f\subset \o_f$. Then
$I(f,g)=$\, dim\,$_\F(\r_f/g\r_f)$
and $I(f,g)\,\de_f=\de_f$\,dim\,$_\F(\o_f/g\o_f)=\de_f\nu_f(g)$,
where  $\de_f=$\,dim\,$_\F(\o_f/\r_f)=$
\,dim\,$_\F(g\o_f/g\r_f)$.
\vskip 0.2cm

{\sf Valuations.}
We have $\kappa$  valuations $\nu_i(f)=ord_{z_i}(f_i)$ 
in $\o\ni f=(f_1,\ldots,f_\kappa)$  and the corresponding semigroups $\Ga_i$; here
$\nu_i(0)=+\infty$. Then $I_{i,j}\equal I(F_i,F_j)=\nu_i(F_j)$.
Recall that $f_i=fe_i$ for 
$1\le i\le \kappa$. As an example,
let $F(u,v)=(u^a+v^b)(u^m-v^n)$ provided
gcd$(a,b)=1=$\,gcd$(m,n)$. Then,
$I_{1,2}=\nu_1(z_1^{na}+z_1^{mb})=
\min\{mb,na\}.$ 
If $m=a$ and $n=b$ in this formula,
 i.e. for $T(2a,2b)$, then $I_{1,2}=ab$.
The characteristic of $\F$ must be $p\neq 2$ here
to ensure that $F$ is square-free in the latter formula.
To avoid this restriction, take
$(u^m+\al v^n)$ instead of $(u^m-v^n)$ for $\al\neq 0,1$.

Let us use (\ref{de-links}) to re-calculate
$\de$ for $T(\kappa a,\kappa b)$, where gcd$(a,b)=1$ and
$N=\kappa a b$. Recall that $\de_{a,b}=(a-1)(b-1)/2$ 
for $T(a,b)$. One has:
{\small
\begin{align}\label{de-kappa}
&\de_{\kappa a,\kappa b}=\kappa\frac{(a-1)(b-1)}{2}+
\frac{\kappa(\kappa-1)}{2}ab=
\kappa \frac{(\kappa a b-a -b +1)}{2}.
\end{align}
}

It coincides with  
$\frac{(N-a)(N-b)}{2ab}+\frac{\kappa -1}{2}=
 \frac{(\kappa b-1)(\kappa a-1)+\kappa -1}{2}$ (see above).
\vskip 0.2cm 

Let $\vec{\nu}=(\nu_1,\ldots, \nu_\kappa)$ and
$\vec{\Ga}=\vec{\nu}(\r)\cap \Z_+^\kappa$. 
In the case of $\kappa=2$, $\vec{\Ga}$ can be used to calculate
$I_{1,2}$. Namely, let  $\e$ be the set of 
{\em extremal points} $(m,n)$ of $\vec{\Ga}$ defined as follows:
either  $m=\max\{m'\mid (m',n)\in \vec{\Ga}\}$ 
 or (equivalently) $n=\max\{n'\mid (m,n')\in \vec{\Ga}\}$. 
Then $|\e|=I_{1,2}$. See   
\cite{Bay} for an explicit
description of $\vec{\Ga}$ for $\kappa=2$ in terms of $\Ga_1,\Ga_2$ 
and $I_{1,2}$.  

We will not use $\vec{\Ga}$ in this paper. Instead, let
$\Ga\equal\{\nu(f), \r\ni f\neq 0\}$, where 
$\nu(e_i)=\nu^i$ and $\nu(f)\equal\min_{i=1}^\kappa
\{\nu_i(f_i)+\up^i\}$ for any $f\in \o$.  We assume that 
$\up^1=0<\up^2<\cdots<\up^\kappa<1$. Actually, they can be any 
numbers from $\R_+$ such that $\up^i-\up^j\not\in \Z$ for $i\neq j$.
Then, $\de=|\nu(\o)\setminus \Ga|.$ Note that
$|\nu(\mathfrak{m}_\o)\setminus \nu(\mathfrak{m}_\r)|=\de-\kappa+1$, 
where
$\mathfrak{m}_\o\equal\oplus_{i=1}^\kappa \mathfrak{m}_{\o_i}$ 
for $\mathfrak{m}_{\o_i}=z_i\o_i\subset \o_i$ and $\mathfrak{m}_\r\equal
\r\cap \mathfrak{m}_\o$. 
\vskip 0.2cm

Finally, we need the 
{\em conductor} $\mathfrak{c}$:\, the greatest 
$\o$-ideal that belongs to $\r$. It is $\o$-generated by 
$z_i^{n_i}\in \r$ for 
minimal $n_i>0$ such that $z_i^{n_i+j}\subset \r$ for all $j\ge 0$.
 Equivalently, the latter numbers are minimal such that  
$n_i+\up_i+\Z_+\subset \Ga$. Recall that 
dim\,$_\F(\o/\mathfrak{c})=2\de$ for plane curve
singularities (the defining property of Gorenstein rings).

\subsection{\bf Superpolynomials for links}\label{sec:sup-links}
For $\F=\C$, the (square-free) equation $F(u,v)=0$
gives the corresponding singularity (with $\kappa$ branches).
The corresponding link is 
$\{F(u,v)=0\}\cap S^3_\ep$ in $\C^2$
with the coordinates $u,v$; it has $\kappa$ components.
Its isotopy type gives the topological type of the singularity.
It is fully determined by the connected components and their
pairwise linking numbers, which coincide with the 
intersection numbers defined above. Thus, the topological
type is given by the sets $\{\Ga_i,1\le i\le \kappa\},
\{I_{i,j}, 1\le i<j\le \kappa\}.$

The passage from $\C$ to $\F_q$ is basically the same as in the 
unibranch case. Namely, we pick $x,y\in \Z[[z]]$ within
a given topological type and then switch to 
$\F_q$ for $q=p^m$ provided that $p$ is a prime of {\em
good reduction}. By definition, good $p$ are such that 
the corresponding $F_i$ remain irreducible and 
pairwise non-proportional over $\F_q$. The  semigroups 
$\Ga_i$ for $\r_i$
and the pairwise intersection numbers 
$\{I_{i,j}, 1\le i<j\le \kappa\}$ must remain unchanged.
Due to (\ref{de-links}), it suffices to assume that $\{\de_i\}$
and the total $\de=\de_\r$ remain unchanged, since the 
intersection numbers
are ``semi-continuous": do not decrease upon specializations. 

The {\em good reduction} here will be necessary for the
conjectural coincidence of motivic superpolynomials
with the DAHA superpolynomials and topology; the corresponding
$p$ must be good. This will be assumed in the Coincidence Conjecture
below.
\vskip 0.2cm

{\sf Standard modules.} Let $\F=\F_q$ and $\r\subset \o$ are 
as above. The {\em standard modules} $M\subset \o$ are $\r$-invariant
such that $M_i=M e_i$ are standard in $\o_i$. Equivalently,
$\o M=\o$. This means that for any 
$1\le i\le \kappa$, there exists $\phi\in M$ such that 
$\phi_i-e_i\in z_i\o$. Such $\phi$ can be different
for different $e_i$, but it is simple to check by induction 
that there always 
exists $\phi\in \o^*\cap M$ for any standard $M$,
where $\o^*=\prod_{i=1}^\kappa \o_i^*$ (the group of
invertibles in $\o$).

A special case is when  $\r_i$ are all isomorphic to $\r_\circ$.
Then $\r\cong \r_\circ$ and our $\o$ (in this case) and
$\Om_\circ=\o_{\circ}^{\oplus \kappa}$ from Section 
\ref{sec:ranks} become isomorphic as $\r$-modules. 
The only difference is
that the standard modules for $\Om_\circ$  are
such that 
$\phi_i-e_i\in z\Om$, a subset of those defined
in this section.

\vskip 0.2cm

The {\em motivic
superpolynomials} for uncolored links are
$$
\h^{mot}\equal\sum_{M} t^{\hbox{\tiny dim}(\o/M)}
\prod_{j=1}^{r\!k_q(M)-1}(1+aq^j)
\text{ summed over standard } M.
$$
Here $r\!k_q(M)=$\,dim\,$_{\F_q} M/\mathfrak{m}_\r M$ and
\,dim\,$(\o/M)=$\,dim\,$_{\F_q}(\o/M)$, which we generally
denote by $deg(M)$. One has:\, 
$deg(M)=
|\nu(\o)\setminus \De(M)|$, where $\De(M)\equal\{\nu(f) \mid
f\in M\}$ contains $\Ga$ for standard $M$.

Note that the action of $\Ga$ on  $\De(M)$ is given
by the ``twisted" formulas $(\nu^i+\ga_i)+(\nu^i+x_i)=
(\nu^i+\ga_i+x_i)$, which is due to our usage of $\{\nu^i\}$.
 One has $\G+\De=\De$ in this (twisted) sense
for standard $\De$. As above, standard $M$ contain 
the conductor $\mathfrak{c}$. Indeed,
 $M=\r M\supset \mathfrak{c}M\supset
\mathfrak{c}\,(1+z(\cdot))=\mathfrak{c}.$ This readily gives
that  there are finitely many standard modules for $\F=\F_q$. 

Extending the corresponding definition in the unibranch case, 
let $\j(\De)\equal \{ M \mid \De(M)-\De\}$ for standard 
$\Ga$-modules $\De$.  They are quasi-projective schemes over 
$\F_q$. We conjecture that they are {\em affine configuration 
spaces}, which is based on extensive calculations.
\vskip 0.2cm

Let us emphasize  that evaluations and gaps heavily depend on
$\{\nu^i,1\le i\le \kappa\}$. See Lemma 2.6 in \cite{ChP2} 
concerning the passage from one such set to another. It is
not always convenient to make $\nu_i<1$. For instance, 
let $n_1\ge n_2\ge \cdots\ge n_\kappa$ for the conductor 
$\mathfrak{c}=\oplus_{i=1}^\kappa z^{n_i}\o$, 
and $\nu^i=n_1-n_i-\ep_i$ for small sufficiently general
$\ep_i>0$, where $i\ge 2$ and $\nu^1=0$ as above. 
Then $\nu(\mathfrak{c})=\{v\in \nu(\o) \mid v>n_1-1\}$.
\vskip 0.2cm

The {\em invertible modules} are important here;
they are standard of $q$-rank $1$, equivalently,
those with one $\r$-generator, which is 
$\phi=1+\sum_{i=1}^{\kappa-1}\al_i e_i+\sum_{g} a_g g$
for $g\in G\equal\nu(\mathfrak{m}_\o)\setminus 
\nu(\mathfrak{m}_\r)$. The parameters $\al_i\in \F_q^*$
and $a_g\in \F_q$
are free and such modules constitute a group scheme
 $\mathbb{G}_m^{\kappa-1}\times
\mathbb{G}_a^{\de-\kappa+1}$ over $\F$. 
Their contribution to $\h^{mot}$ is
$(q-1)^{\kappa-1} q^{\de-\kappa+1} t^{\de}$. 

The greatest $q$-rank is that 
for $\o$,
which is $m_1+\cdots+m_\kappa$ for 
$m_i=\min\bigl\{\Ga_i-\nu^i\setminus \{0\}\bigr\}$,
where $\Ga_i$ are valuation
semigroups for $\r_i$. It is, generally,
not a unique standard module of maximal $q$-rank. 
\vskip 0.2cm

{\sf Coincidence Conjecture.} Its colored version will be
stated later. In the variant for uncolored algebraic {\em links},
we conjecture that 
\begin{align}\label{coin-links}
\h^{daha}_{\l,'\!\l}(q,t,a)=\h^{mot}(q,t,a) \text{ for the
corresponding } 
\r,
\end{align}
where the pair of trees $\{\l,'\!\l\}$
in $\h^{daha}=\hat{\h}$ from (\ref{h-polynoms-hat}) 
and $\r\subset \o$ over $\F_q$
are those associated with the 
corresponding algebraic link.

\subsection {\bf Reduction formula}
The following formula reduces $\h^{mot}$ to considering 
$M$ containing $\r$; such $M$ are standard automatically.
Let $\tilde{M}\equal (M+\mathfrak{m}_\o)/\mathfrak{m}_\o$
for standard $M$; it is a vector spaces over $\F_q$ of dimension
$\upsilon_M=|\De(M+\mathfrak{m}_\o)\setminus \De(\mathfrak{m}_\o)|$.
One has: $\upsilon_\o=\kappa, \upsilon_\r=1$. 
Also, we set $\tilde{U}=\o^*/(\F_q^*+\mathfrak{m}_\o)\cong
(\F_q^*)^{\kappa-1}$. 
The {\em generalized Jacobian} of $\r$ is 
$U\equal \o^*/\r^*$, which is isomorphic
to $\tilde{U}\times  \F_q^{\de -\kappa+1}$ as a group.

Let $I_M=\{\phi\in M\cap \o^*\} \mod \r^*$ and
$I_{\tilde{M}}= \{\tilde{M}\cap \tilde{\o}^*\} \mod \F_q^*$.
Since $M$ are standard, $I_M$ and $I_{\tilde{M}}$
are never empty. The groups
$U$ and $\tilde{U}$ act in $I_M$ and $\tilde{I}_M$ respectively. 
One has $I_{\tilde{M}}\cong (\F_q^*)^{\upsilon_M-1}$ for 
sufficiently general $M$, but the set $I_{\tilde{M}}$ can be 
smaller than this (see the example below). 
Generally, finding  such sets 
leads to interesting linear algebra problems; see e.g., \cite{Fr}.
Finally, let $i_{\tilde{M}}\equal |I_{\tilde{M}}|$, which is 
$(q-1)^{\upsilon_M-1}$ for generic $M$.

We will use the variety $\j^{ext}$ of
pairs $(M, \vph)$, where $M$ is standard and 
$\vph\in I_{\tilde{M}}$. There are $2$ natural surjective 
projections 
from $\j^{ext}$.
The $1${\small st}  is\, $\j^{ext}\to \{M'\supset \r\}$,  
sending $(M,\vph)\mapsto M'=M\vph^{-1}\supset \r$, which is
with the fibers of size 
$(q-1)^{\kappa-1}q^{\de-\kappa+1}$.
The $2${\small nd} one is the forgetting
map $(M,\vph)\mapsto M$ with the fibers of size 
$q^{dev(M)-\up_M+1}i_{\tilde{M}}$.
Recall that $dev(M)=\de-deg(M)$; we use a counterpart of
 Lemma \ref{lem:inv} for links.
Combining these maps, we
reduce the calculation of $\h^{mot}$ to the following weighted
sum over the image of the
first projection:
\begin{align}\label{M-one}
\h^{mot}\!=\!\sum_{M\supset \r}(qt)^{deg(M)}
\Bigl(\frac{\,(q-1)^{\kappa-1}}{q^{\kappa-\up_M}\,i_{\tilde{M}}}\Bigr)\,
(1\!+\!aq)\cdots(1\!+\!aq^{r\!k_q(M)-1}).
\end{align}
This formula is helpful practically, especially for $\kappa=1$,
when the weights $w_M$ of modules $M\supset \r$ are  $(qt)^{deg(M)}.$
Theoretically, when combined with the {\em reciprocity map}, it
allows reducing the calculation of $\h^{mot}$ to
counting certain {\em ideals} in $\r$. The weight $w_M$ is the
coefficient before the product of $(1+a q^i)$. 
See \cite{ChG} and Section \ref{sec:unibr} around
formula (\ref{Dedual}) (a comment on the reduction formula).
\vfil

{\sf Hopf 3-link.} Let us demonstrate how formula (\ref{M-one})
works for the Hopf $3$-link, which is a union of $\kappa=3$ unknots
with pairwise linking numbers $+1$, $T(3,3)$ in
our standard notations. It can be represented by 
$\r=\F_q[[e_1+e_2+e_3, z_1+z_2, z_1+z_3]]
\subset \o=\F_q[[e_i,z_j]]$. Here $\r$ contains all $z^2_i$
and $\de=$\, dim\,$_{\F_q}(\o/\r)=3$. The families of modules $M$
containing $\r$, their\, $r\!k_q,\, deg,\, \up$, the corresponding
weights $w=
q^{deg(M)-\kappa+\up_M}
(q-1)^{\kappa-1}\, t^{deg(M)}/i_{\tilde{M}}$, and 
the sizes $\#$ of these families are:\,
\vfil

(0) $\{M=\o\}$, which is with $r\!k_q=3$, $deg=0, 
\up=3, i=(q-1)^{2}, w=1$, and $\#=1$
(the total contribution to $\h^{mot}$ is $(1+qa)(1+q^2a)$);
\vskip 0.1cm

(1) $\{M=\r\}$ (a unique invertible module containing $\r$), where
$r\!k_q=1, deg=3, \up=1$, $i_{\tilde{M}}=1,$ 
$w=(q-1)^2 q t^3$, and $\#=1$ (only $\r$);
\vskip 0.1cm

(2) $\{M_\al=\lan e_1+e_2+e_3, e_2+ \al e_3\ran\}$ for
$\al\neq 0,1$, where
$r\!k_q=2, deg=1, i_{\tilde{M}}=q-2$, $\up=2$,
$w=(q-1)^2 t/(q-2)$, and of size $\#=q-2$;
\vskip 0.1cm

(2$'$) $\{M_i=\lan e_1+e_2+e_3, e_i\ran\}$ for $i=1,2,3$
of size $\#=3$ with
the same $r\!k_q=2, deg=1, \up=2$, but now for $i_{\tilde{M}}=q-1$ and
$w=(q-1)t$;
\vskip 0.1cm

(3)  $\{M=\lan e_1+e_2+e_3, z_1\ran\}$ (containing $z_{2}$ and $z_3$),
where $r\!k_q=2$, $deg=2, \up=1, i_{\tilde{M}}=1$,
$w=(q-1)^2t^2$, and $\#=1$ (only $1$ module).
\vskip 0.2cm

A direct calculation of these families (without the
usage of the reduction formula) is not difficult.
\vskip 0.1cm

Family (1) becomes $\{M=\lan e_1+ a_2 e_2 + a_3 e_3+b_1 z_1\ran\}$
for $a_2,a_3\in \F_q^*$ and $b_1\in \F_q$, which is of size
$(q-1)^2 q$.
\vskip 0.1cm

The modules of family (2) correspond to $2$-dimensional
subspaces in $\F_q e_1+\F_1 e_2+\F_q e_3$ that are not
$\lan e_i,e_j\ran$ for $i\neq j$. Their number is $(q^2+q+1)-3=
(q-1)(q+2)$. 
\vskip 0.1cm

Family (3) becomes  now $\{M=\lan e_1+ a_2 e_2 + a_3 e_3, z_1\ran\}$
such that  $a_2,a_3\in \F_q^*$, which is of size  $(q-1)^2.$

\vskip 0.2cm

Combining this families with the corresponding products of $(1+a q^i)$,
we obtain the following formula for the motivic superpolynomial
of $T(3,3)$ (Hopf $3$-plus-link):
\begin{align} \label{3-Hopf}
&\ \ \ \ \ \ \ \h^{mot}_{3,3}=\!\sum_{M\supset \r}
w_M(1\!+\!aq)\cdots(1\!+\!aq^{r\!k_q(M)-1})=\\
&\text{\small $(1\!+\!qa)(1\!+\!q^2a)
\!+\! (q\!-\!1)^2 q t^3 \!+\!
(1\!+\!qa)\Bigl((q\!-\!1)(q\!+\!2)t\!+\! (q\!-\!1)^2t^2\Bigr)$}.\notag
\end{align}
\vskip 0.2cm

\subsection{\bf Double trefoil} \label{sec:dtref}
One has: $\o=\F[[e_1,e_2,z_1,z_2]]$ and
$\r=\F[[1=e_1+e_2, x=z_1^2-z_2^2, y=z_1^3+z_2^3]]$. The corresponding
equation will be $F=(x^3+y^2)(x^3-y^2)=x^6-y^4=0$. The Milnor number is 
$\mu=$\,dim\,$_{\F}\bigl(\F[x,y]/(F_x,F_y)\bigr)=
\frac{(N-a)(N-b)}{ab}=15$,
where $a=2,b=3,N=12$; the characteristic $p$ of $\F$ must not be
 $2,3$ for this formula.
Using the Milnor number, $\de=(\mu+\kappa-1)/2=16/2=8$, which is for
any characteristic $p\neq 2$. This restriction can be avoided
by taking $x=z_1^2+\ep z_2^2$ and with the same $y$, where  $\ep\neq 0$
is such that $\ep^3\neq 1$. The formula $\de=8$ gives that  
the conductor is $\mathfrak{c}=\o z_1^8+\o z_2^8$.
 
Concerning the intersection number, it is 
$I_{1,2}=$\,dim\,$\F[x,y]/(x^3, y^2)=
6$, where we use directly the definition. This is an alternative
way to obtain $\de$; namely, 
$\de=\de_1+\de_2+I_{1,2}=2+6=8$.

Using now the valuation:\,
$I_{1,2}=\nu_1(x^3+y^2)=\nu_1(2z_1^6)=6$, where the images of
$x,y$ in $\r_1$ are $z_1^2,z_2^3$, which we already know.
 Alternatively, let us use
$\vec{\Ga}=$ {\small
$\{(0,0),(2,2),(3,3),(4,4),(5,5),(6,6+u),(6+u,6),
(7,7),(8,8+u),(8+u,8),\ldots\}$}, where $u$ is any number
from $\Z_+$. There are $6$ {\em maximal elements}, 
{\small $ \{(0,0),(2,2),(3,3),(4,4),(5,5),(7,7)\}$}, which gives 
$I_{1,2}=6$.
 
It is not accidental that $\de\!=\!8\!=\!\de'$, where $\de'$ is 
for $\r'=\F[[x=z^4,y=z^6+z^7]]$. Topologically, $T(6,4)$ is a 
degeneration
of $C\!ab(13,2)T(3,2)$. Algebraically, this gives  that $\Ga'$
coincides with $2\Ga$ for $\up^1=0,\up^2=1/2$, where $\Ga'=
\{0,4,6,8,10,12,13,14,16, 17,\ldots\}$.  Indeed,
\begin{align}\label{Ga-6-4}
\text{\small
$\Ga=\{0,2,3,\ldots\}\cup\{\frac{1}{2}+6,\,\frac{1}{2}+8,
\,\frac{1}{2}+10,\,\frac{1}{2}+11,\ldots\} \text{\, for\, } T(6,4)$}.
\end{align}
\vskip 0.2cm

{\sf DAHA superpolynomial for $T(6,4)$.}
This one 
was provided in\cite{ChD2}, Section 6.2.
In the notation there,
\begin{align}\label{T6-4}
&T(6,4):\ 
\l=\l_{(\{3,2\},\{3,2\})}^{\,\circ\rightrightarrows,\, 
(\square,\square)},\ \ 
\h^{daha}_{6,4}=\hat{\h}{}^{min}_{\l}\,(q,t,a)=
\end{align}

\renewcommand{\baselinestretch}{0.5} 
\noindent
{\small
\(
1-t+q t+q^2 t+q^3 t-q t^2+2 q^4 t^2-q^2 t^3-q^4 t^3
+2 q^5 t^3-q^3 t^4-q^5 t^4+2 q^6 t^4-q^4 t^5+q^7 t^5-q^5 t^6
+q^7 t^6-q^6 t^7+q^7 t^7-q^7 t^8+q^8 t^8+
a^3 \bigl(q^6-q^6 t+q^7 t-q^7 t^2+q^8 t^2\bigr)+
a^2 \bigl(q^3+q^4+q^5-q^3 t+q^5 t+2 q^6 t-q^4 t^2-q^5 t^2
+2 q^7 t^2
-q^5 t^3-q^6 t^3+q^7 t^3+q^8 t^3-q^6 t^4+q^8 t^4-q^7 t^5
+q^8 t^5\bigr)+a \bigl(q+q^2+q^3-q t+q^3 t+3 q^4 t+q^5 t-q^2 t^2
-q^3 t^2-q^4 t^2+3 q^5 t^2+q^6 t^2-q^3 t^3-q^4 t^3-2 q^5 t^3
+3 q^6 t^3+q^7 t^3-q^4 t^4-q^5 t^4-q^6 t^4+3 q^7 t^4-q^5 t^5
-q^6 t^5
+q^7 t^5+q^8 t^5-q^6 t^6+q^8 t^6-q^7 t^7+q^8 t^7\bigr).
\)
}
\renewcommand{\baselinestretch}{1.2} 
\smallskip

The $a$\~degree of this $\hat{\h}{}^{daha}_{\l}\,(q,t,a)$
is $3$, which is a particular case
of the general formula deg$_a=\ss(2|\la|)-|\la|=3$
from (\ref{deg-a-jj}). 

This polynomial is non-positive, which creates some problems
in Khovanov-Rozansky theory. The positivity for 
algebraic links
is part of Conjecture 5.3 in \cite{ChD2}; it is only
upon the 
division by $(1-t)^{\kappa-1}$ (in the uncolored case). 
In this example, 
$\h_{6,4}^{daha}/(1-t)$ is positive.

Our superpolynomial matches that {\em conjectured} in 
\cite{DMMSS} in Section 2.8. The main factor of
$-P_1^{T[4,6]}$ there coincides with our one upon the
following substitution:
$ a\mapsto A^2, q\mapsto q^2, t\mapsto t^2$. I.e.
their {\em non-bold\,} $A,q,t$ are essentially the
DAHA parameters.

Also, it coincides with the reduced (uncolored)
Khovanov-Rozansky polynomial defined 
via {\em Soergel modules} posted  in Example 1.3 of
\cite{HM}.
Their expression 
is $\h^{daha}_{6,4}(q\mapsto 1/t, t\mapsto q, a)$ 
multiplied by $\frac{1+a}{(1-q)^2}$. 
\vskip 0.2cm

{\sf Motivic superpolynomial for $T(6,4)$.}
Let us provide the $q$-rank decomposition of
the corresponding motivic superpolynomial,
which coincides with $\h^{daha}_{6,4}$. 
One has: $\h^{mot}_{6,4}=$

\renewcommand{\baselinestretch}{0.5} 
\noindent
{\small
\(
q^7(q-1)t^8+(1\!+\!a q) (1\!+\!a q^2) \bigl(1\!+\!a q^3)(1+(q-1) t+ 
q(q-1) t^2\bigr)+
(1\!+\!a q) (1\!+\!a q^2) \bigl(q^2 (1+q) t+q^2 (q^2-1) t^2+q^2 
(q+1)(q^2-1)t^3
+q^3 (q^2-1) t^4+q^4 (q-1) t^5\bigr)+(1\!+\!a q)\bigl(q^4 t^2+
q^3 (q-1)^2 t^3
+2q^5 (q-1) t^4+q^5 (q^2-1) t^5+q^5 (q^2-1) t^6+q^6 (q-1) t^7\bigr).
\)
}
\renewcommand{\baselinestretch}{1.2} 
\vskip 0.2cm

{\sf Rank $1$}. The generators of invertible modules are
$\phi=e_1+\om e_2+a_0 z_1+ a_1 z_2+a_2 z_2^2+a_3 z_2^3+a_4 z_2^4+
a_5 z_2^5+a_7 z_2^7$ for 
arbitrary $\om\in \F_q^*$ and $\{a_i\in \F_q\}$. They contribute
$q^7(q-1) t^8$ to $\h^{mot}_{6,4}$. Here and below the monomials
are ordered with respect to the valuation $\nu$. 

{\sf Rank $4$}. The standard modules of top $q$-rank (which is $4$) are 
$\o$, and the following $2$ families:\, 
$\{\lan e_1+\om e_2, z_1,z_2,z_1^2\ran$ with $\om\in \F_q^*\}$,\ 
and  \, $\{\lan e_1+\om e_2+\al z_2, z_1-\om z_2, 
z_1^2, z_1^3\ran\}$,
where $\om\in \F_q^* $, $\al\in \F_q$. The coefficient
$-\om$ in $z_1-\om z_2$ is necessary to ensure that $r\!k_q=4$. 
Here and below the span $\lan\cdot\ran$ is understood 
over $\r$ unless stated
otherwise. 
The total contribution of these families matches the formula
above. Namely:   
\begin{align}\label{4st6-4}
\h_{6,4}^{r\!k\!=\!4}\!=\!\text{\small
$\bigl(1+(q\!-\!1)t+q(q\!-\!1)t^2\bigr)
(1\!+a q)(1\!+\!a q^2)(1\!+\!a q^3)$}.
\end{align}


{\sf Rank $2$}. The coefficient of $(1\!+\!a q)$ 
is due to standard modules $M$ of $q$-rank $2$, which
we will now describe.
 See also Proposition \ref{thm:h-iter};
families $(1-6)$ correspond to standard modules for $\r'$, where the
total contribution of the latter must be multiplied by $(q-1)/q$. 
Family $(0)$ below comes from invertible modules for $T(3,2)$
of $\o$-rank $2$ (for $c=2$). 

We will use $\De',D'$ for $\r'=\F_q[[z^4,z^6+z^7]]$
using the identification above: 
 $\De\ni v\mapsto v'=2 v\in \De'$, where $\nu^1=0,
\nu^2=1/2$. Recall that
 $D'=\{1,2,3,5,7,9,11,15\}\cap \De'$. 

\vskip 0.2cm

(0) The modules from the family 
$\{\lan e_1+a z_1+b z_2, e_2+ c z_1+ d z_2\ran\}$ for any $a,b,c,d\in \F_q$.
contain $z_1^u,z_2^u$ for $u\ge 2$. This family  
contributes $q^4 t^2$ (times $(1+aq)$ here and below). 
The corresponding $\De'$
is generated by $0$ and $1$; accordingly, $D'=[1,5,7,9,11,15].$

\vskip 0.2cm
The remaining $6$ families of $r\!k_q=2$ are such that $1/2\not\in \De$,
equivalently, $1\not\in \De'$. 
Families $(1,2,3,4)$ are such and when 
$\nu(z_1)\!=\!1\! \not\in \De(M)$, equivalently, 
$2 \not\in \De'$. They are:

(1) $\{\lan e_1+\om e_2+a_0 z_1+ a_1 z_2+a_2z_2^2+a_3z_2^3+a_4z_2^4+
a_5z_2^5,\, z_2^7\ran \}$, where the corresponding contribution
 to $\h^{mot}$  is $q^6(q-1)t^7$ (multiplied by
 $(1+a q)$ here and below);

(2) $\{\lan e_1+\om e_2+a_0 z_1+ a_1 z_2+a_2z_2^2+a_3z_2^3+a_4z_2^4,\,
z_2^5\ran $ and $\lan e_1+\om e_2+a_0 z_1+ a_1 z_2+a_2z_2^2+a_3z_2^3+
a_5z_2^5, \,
z_2^4+b_5 z_2^5\ran \}$,
contributing $q^5(q^2-1)t^6$;

(3) $\{\lan e_1+\om e_2+a_0 z_1+ a_1 z_2+a_2z_2^2+a_4z_2^4,\,
z_2^3+b_4 z_2^4\ran \}$ and 
$\{\lan e_1+\om e_2+a_0 z_1+ a_1 z_2+a_3z_2^3+a_4z_2^4,\,
z_2^2+b_3 z_2^3+b_4 z_2^4\ran \}$,
contributing $q^5(q^2-1)t^5$;

(4) $\{\lan e_1+\om e_2+a_0 z_1+ a_1 z_2+a_3z_2^3,\,
z_2^2+b_3 z_2^3\ran\} $, 
contributing $q^4(q-1)t^4$.
\vskip 0.2cm

The considerations become more ramified if $1\in \De(M)\not\ni 1/2$. 
Equivalently,
$2\in \De'\not\ni 1$, where we switch from $\De$ to 
$\De'$ as above. 
\vskip 0.2cm

(5) For $D'=[2,7,11,15]$ the family is 
$\{\lan v_1=e_1+\om e_2+ a_1 z_2+a_2z_2^2+a_4z_2^4,\,
v_2=z_1+b_1z_2+b_2 z_2^2+b_4 z_2^4\ran\} $, where
the relation $u+b_1\neq 0$ is
necessary for the  occurrence of $z_2^3+\ldots$ in $M$. Indeed,
$y v_1-xv_2=z_1^3+uz_2^3-z_1^3+b_1 z_2^3+\ldots=(u+b_1)z_2^3+\ldots.$
Similarly,
$x^2 v_1-y v_2=z_1^4+u z_2^4-z_1^4-b_1 z_2^4+\ldots=(u-b_1)z_2^4+
\ldots$, which results in $u=b_1$ to prevent the occurrence of $9$
in $D'$. The characteristic must be $p\neq 2$. We obtain that 
the contribution of this family is $q^5(q-1) t^4$.

(6) Finally, $D'=[2,7,9,11,15]$ and this family is 
$\{\lan v_1=e_1+\om e_2+ a_1 z_2+a_2z_2^2,\,
v_2=z_1+b_1z_2+b_2 z_2^2\ran \}$ for
$b_1\neq \pm u$, contributing 
$q^3(q-1)^2 t^3$. 

\vskip 0.2cm
This is the end of the list, since $q$-rank is $\ge 3$ if $D'=
[2,9,11,15]$. The following is the sum of all contributions of modules
of $q$-rank $2$ (times $(1+a q)$ in $\h^{mot}$), 
where the corresponding family is shown:

{\small
\begin{align*}
&q^4 t^2\, [(0)]+q^6 (q-1) t^7\, [(1)]+
q^5 (q^2-1) t^6\, [(2)]+ q^5 (q^2-1) t^5\, [(3)]\\
&+q^5 (q-1) t^4\, [(4)]+ q^5 (q-1) t^4\, [(5)]+
q^3 (q-1)^2 t^3\, [(6)].
\end{align*}
}

\comment{
(-q^7+q^8) t^8+(1+a q) (1+a q^2) (1+a q^3)(1+(-1+q) t+(-q+q^2) t^2)+
(1+a q)(1+a q^2)(q^2 (1+q) t+q^2 (-1+q^2) t^2+q^2 (-1-q+q^2+q^3) t^3
+q^2 (-q+q^3) t^4+q^2 (-q^2+q^3) t^5)+(1+a q)(q^4 t^2
+q^3 (1-2 q+q^2) t^3+q^3 (-2 q^2+2 q^3) t^4+q^3 (-q^2+q^4) t^5
+q^3 (-q^2+q^4) t^6+q^3 (-q^3+q^4) t^7)
}

\vskip 0.2cm
{\sf Reduction to $M$ containing $\r$.}
It is instructional to employ formula (\ref{M-one}) to 
the modules of $q$-rank $2$. 
Then $u_M=1$ for types $(1-6)$ and $u_M=2$ for type $(0)$.
The parameters of the corresponding $M\supset \r$ are
those from the second generator; the number of such $M$
is after $(qt)^{\cdots}$. One has:
\begin{align*}
&\sum_{M\supset \r}^{r\!k_q(M)=2}\,(qt)^{deg(M)}(\frac{(q-1)}{q})^
{\kappa-u_M}=
(qt)^2 q^2+\frac{q-1}{q}\Bigr((qt)^7 \\
&+(qt)^6 (q+1)+ (qt)^5 q(q+1)
+2(qt)^4 q^2+ 
(qt)^3q(q-1)\Bigl).
\end{align*}

\subsection{\bf Iteration formula}
The above connection between cables and
links is a particular case of a 
general relations between iterated knots $C\!ab(\up \rr  \ss+1,\up)
T(\rr,\ss)$, torus links $T(\up \rr,\up \ss)$ and 
torus knots $T(\rr,\ss)$ colored be $m\om_1$,
 where\, gcd$(\rr,\ss)=1$. 

Topologically, a connection of invariants of
$C\!ab(\up \rr  \ss+1,\up)T(\rr,\ss)$ and
$T(\up \rr,\up \ss)$ can be expected, but this
is involved, especially for $\up>2$. This is not known 
for $KhR$ polynomials. 
Motivically, the corresponding $\de$
coincide:\, compare (\ref{Ga-cab}) and (\ref{de-kappa}).
Moreover, let $\Ga'$ be for $\r'$ for 
the cable above, which is 
$\r'=\F[[z^{\up\rr},z^{\up\ss}+z^{\up\ss+1}]]$. 
Then $\Ga'=\up \Ga$ for $\Ga$ for $T(\up \rr,\up \ss)$, where
we take $\nu^i=\nu(e_i)=\frac{(i-1)}{\up}$; cf. (\ref{Ga-6-4}).
For $\up=2$, this is
part of the Proposition \ref{thm:h-iter}, which is the case $m=0$ of the 
DAHA-relation from Proposition 4.3 and
formula (4.23) in \cite{ChW}. See \cite{ChP1,GMO} about the
corresponding Piontkowski cells.
The simplicity of the DAHA
formula is greatly
clarified in the motivic approach.

\begin{proposition}\label{thm:h-iter}
{\sf (a)} Let $\h^{mot}_{\rr,\ss;2}$ be the uncolored motivic 
superpolynomial for 
$\r'=\F_q[[x'\!=\!z^{2\rr}, y'\!=\!z^{2\ss}\! +
u\,\!z^{2\ss+1}]]$, where $u\neq 0$, 
gcd$(\rr,\ss)=1$ and the knot is $C\!ab(2\rr\ss+1,2)
T(\rr,\ss)$. {\sf (b)} Let $\h^{mot}_{2\rr,2\ss}$ be that for 
$\tilde{\r}=\F_q[[\tilde{x}=uz_1^\rr+z_2^\rr,\,
\tilde{y}=z_1^\ss+z_2^\ss]]
\subset \tilde{\o}=\F_q[[e_1,e_2,z_1,z_2]]$ as above, where 
$u^{\ss}\neq 0,1$. 
This is for $T(2\rr,2\ss)$.
{\sf (c)} Let $\h^{mot}_{\rr,\ss}(2)$ be
the motivic superpolynomial 
for $T(\rr,\ss)$ colored
by $2\om_1$, i.e. for 
 $\r_\circ=\F_q[[x=z^\rr,y=z^\ss]]$ acting diagonally in
$\Om=\F_q[[z]]\vep_1\oplus \F_q[[z]]\vep_2$. The following 
counterpart of the corresponding $\h^{daha}$-identity (for $m=0$) 
from \cite{ChW} holds:
\begin{align}\label{h-dtor}
(1+a q)\h^{mot}_{\rr,\ss}(2)+(q-1)\h^{mot}_{\rr,\ss;2}=
q\h^{mot}_{2\rr,2\ss}.
\end{align}
\end{proposition}
{\it Sketch of the proof.} We check this identity in this form:
\begin{align}\label{h-dtor1}
(q\!-\!1)\Bigl(\h^{mot}_{\rr,\ss;2}\!-\!(1\!+\!a q)\h^{mot}_{\rr,\ss}(2)\Bigr)\!=\!
q\Bigl(\h^{mot}_{2\rr,2\ss}\!-\!(1\!+\!a q)\h^{mot}_{\rr,\ss}(2)\Bigr).
\end{align}

The $2${\small nd} difference in (\ref{h-dtor1}) 
corresponds to the following
embedding $(1\!+\!a q)\h^{mot}_{\rr,\ss}(2)\subset 
\h^{mot}_{2\rr,2\ss}$. Generally, let 
$x_\circ=\tilde{x}(u=1)=\ze_1^{\rr}+\ze_2^{\rr}$,  
$y_\circ=y$ in {\sf (b)}. Then  $x_\circ^{\ss}=y_\circ^{\rr}$ and
$\tilde{\o}$ considered as a module over $\tilde{\r}_\circ=
\F_q[[x_\circ,y_\circ]]$  becomes $\Om=\o^{\oplus 2}$ upon the
identification  $\vep_i=e_i$.

We claim that the
standard modules in $\Om$ are exactly the specializations
at $u=1$ of $\tilde{M}\subset \tilde{\o}$ 
for  $\tilde{\r}$ under the following condition.
Modules $\tilde{M}$
must containing $e_1 \text{ mod } z\tilde{\o}$ and 
$e_2 \text{ mod } z\tilde{\o}$; equivalently, $\De(\tilde{M})\ni 0,1/2$.
 Moreover, this specialization identifies
the corresponding Piontkowski cells and preserves $r\!k_q$.
See Proposition
\ref{prop:link-rank} below for a more general setting and a
connection to \cite{ChP2}.

The remaining modules $\tilde{M}$ form the  $2${\small nd} 
difference. Here (and below) the factor $(1+a q)$ occurs 
because the products of $(1+a q^i)$ in $\h^{mot}(c)$ begin
with $(1+ aq^c)$. Let us consider the $1${\small st} difference.
\vskip 0.2cm

The interpretation of the embedding 
$(1\!+\!a q)\h^{mot}_{\rr,\ss}(2)\subset
\h^{mot}_{\rr,\ss;2}$
is basically as follows.
We use that 
the ring $\F_q[\sqrt{x'\,},\sqrt{y'\,}\in\o'=\F_q[[z]]$\, 
for  $x',y'\in \r'$ from {\sf (a)}  is isomorphic to 
$\r_\circ$ above for generic $u$.
 This gives that $\r'$-modules $M'\subset \o'$ 
such that $\De(M')\ni 0,1$ 
for $\r'_u=\F_q[[y=z^{2\ss}+u z^{2\ss+1}]]$ for generic $u$
can be identified
with standard modules in $\Om$ for  $\r'_{u=0}=\r_\circ$. 
This identification
preserves the  corresponding Piontkowski
cell and $r\!k_q$. Such embeddings seem a 
general feature of the iterated construction. In this case, 
we use that the generators of $M'$ with $\De(M')\ni 0,1$ can be
chosen uniformly in terms of $u$ (considered as
a formal variable), including $u=0$. 
\vskip 0.2cm

Finally, let us connect $\r'$ and $\tilde{\r}$. The key is 
the identification above of 
the $\De$-sets of the standard modules of $\tilde{\r}$
with the 
$\De'$-sets for $\r'$. Disregarding the modules 
of $\o$-rang $2$ coming from  $\r_\circ$ in the embeddings
above, i.e. those with $0,1\in \De'$ and $0,1/2\in \De$, 
this identification works at the level of
the parameters of the corresponding Piontkowski cells
with the following adjustment.
The parameter $\al$ of $z$ in the generator 
$\vph=1+\al z+\ldots\in M'$ for $\r'$ is arbitrary
from $\F_q$. However, the corresponding parameter $\al$ in
$e_1+\al e_2+\ldots \in M$ for $\r_{2\rr,2\ss}$ 
must be from $\F_q^*$ (arbitrary such).
This gives  $(q-1)$ and
$q$ in (\ref{h-dtor1}). \sq
\smallskip

The simplest version of this formula is the following 
identity
$$(q-1)\Bigl(\h^{mot}_{1,1;2}-(1+a q)\h^{mot}_{1,1}(2)\Bigr)=
q\Bigl(\h^{mot}_{2,2}-(1+a q)\h^{mot}_{1,1}(2)\Bigr).
$$
Here the iterated knot $C\!ab(3,2)T(1,1)$ is isotopic to
$T(3,2)$. Also, $\h^{mot}_{1,1}(2)=1$ since $T(1,1)$ is unknot, 
The identity becomes
$$(q-1)\Bigl(\h^{mot}_{3,2}-(1+aq)\Bigr)=
q\Bigl(\h^{mot}_{2,2}-(1+aq)\Bigr).$$

It clarifies very well the similarity of
$\h^{mot}_{3,2}=qt+(1+aq)$
and $\h^{mot}_{2,2}=(q-1)t+(1+aq)$, not immediate in
the DAHA approach. Recall that  $\h^{daha}_{1,1;2}=\h^{daha}_{3,2}$ 
corresponding to $C\!ab(3,2)T(1,1)=T(3,2)$ follows from the DAHA
identity 
$\{\tau_-\tau_+\tau_-(P_\la)\}_{ev}=
\{\tau_-\bigl(\tau_+\tau_-(P_\la)(1)\bigr)\}_{ev}$.

\vskip 0.2cm
{\sf Example}.
The following is the list of all $5^2$ cells such that
$0,1\in \De'$ for $\r'=\F_q[[z^6,z^8+z^{11}]]$ and their
contributions as $a=0$. Here $\de=18$,
the cable is $C\!ab(25,2)T(4,3)$, and $5$
is the number of standard $\De$ for $T(4,3)$. We 
provide only $D\setminus [1,7,9,13,15,17,19,21,23,27,29,35]$:
\vskip 0.2cm

\(
\begin{array}{cc}
 \{\} & q^{12} t^6 \\
 \{10\} & q^{11} t^5 \\
 \{11\} & q^{10} t^5 \\
 \{2,10\} & q^{11} t^4 \\
 \{3,11\} & q^{10} t^4 \\
 \{4,10\} & q^{10} t^4 \\
 \{5,11\} & q^9 t^4 \\  
 \{10,11\} & q^8 t^4 \\
\end{array}
\ \ \begin{array}{cc}
 \{2,4,10\} & q^9 t^3 \\
 \{2,10,11\} & q^9 t^3 \\
 \{3,5,11\} & q^8 t^3 \\  
 \{3,10,11\} & q^8 t^3 \\
 \{4,10,11\} & q^7 t^3 \\
 \{5,10,11\} & q^6 t^3 \\   
 \{2,3,10,11\} & q^8 t^2 \\
 \{2,4,10,11\} & q^7 t^2 \\
\end{array}
\ \ \begin{array}{cc}
 \{2,5,10,11\} & q^6 t^2 \\
 \{3,4,10,11\} & q^6 t^2 \\
 \{3,5,10,11\} & q^5 t^2 \\
 \{4,5,10,11\} & q^4 t^2 \\
 \{2,3,4,10,11\} & q^5 t \\
 \{2,3,5,10,11\} & q^4 t \\
 \{2,4,5,10,11\} & q^3 t \\
 \{3,4,5,10,11\} & q^2 t \\
 \{2,3,4,5,10,11\} & 1 \,.\\
\end{array}
\)

The last entry is for $M=\o$. The total  sum 
equals $\h_{4,3}(2)(a=0)$.

\subsection{\bf Trefoil linked with unknot}
We will end this section with $\h^{mot}$ for
$\r=\F[[1=e_1+e_2, x=z_1^2+z_2, y=z_1^3+z_2]]\subset \o=\F[[e_1, e_2,
z_1,z_2]]$. The corresponding equation is $(x^3-y^2)(x-y)=0$.
Topologically, this is $T(3,2)$ linked with
unknot with $lk=2$. 

One has:\, $\de=\de_1+\de_2+I_{1,2}=|G|$, where  
$G=\nu(\o)\setminus \Ga$ is the set
of potential gaps for $\nu^1=0,\nu^2=\frac{1}{2}$,
the valuations of $e_2,z_1,z_1^2$,
$G=\{\frac{1}{2},1,3\}$ and $\de=3$ in this case.
The conductor is 
$\mathfrak{c}=\o z_1^4+\o z_2^2$, which can be seen
directly:\,
$u=y^2-x^3=z_2^2-z_2^3$ implies
 $z_2^2=u+ux+ux^2+\ldots\in \r$, 
and $x^2=(z_1^2+z_2)^2=z_1^4+z_2^2$ implies $z_1^4\in \r$. 
Alternatively, \,dim\,$\o/\mathfrak{c}=2\de=6$  because
$\r$ is Gorenstein. We will use that all standard modules $M$ contain
$\mathfrak{c}$.

One has: $\h^{mot}=$
{\small
$q^2(q-1) t^3+(1\!+\!a q) (1\!+\!a q^2)
+(1\!+\!a q)\bigl(q^2t+(q-1) t+q(q-1) t^2\bigr).$
}
Recall that, generally,
$deg_a=m_1+m_2-1$ for $2$-branch singularities, where $m_i$ are 
the multiplicities of singularities.

The simplest families of standard modules are 
invertible $M$, contributing  $q^{\de-1}(q-1)t^\de=q^2(q-1) t^3$
to $\h^{mot}$, and $M=\o$, which is a unique standard module 
of $r\!k_q=3$ in this case (not true generally). 
The remaining (families of) modules are of $r\!k_q=2$; they are as
follows:
\vskip 0.2cm

(0) $\lan v_1=e_1+ a_1 z_1, v_2=e_2+ b_1 z_1 \ran$
for $a_1,b_1\in \F_q$,  containing 
$z_2^2$, $y v_2=z_2$\, and $ x^2 v_1=z_1^2$   
modulo $\mathfrak{c}$, which gives that  $deg=3-2=1$;

(1) $\lan v_1=e_1+ \om e_2, v_2=z_1 \ran$ for $\om\in F_q^*$,
which modules 
contain $z_2$ 
due to $\om\neq 0$ (consider $yv_1$) and are of the same degree $1$;
 
(2) $\lan v_1=e_1+ \om e_2+ a_1z_1, v_2=z_1^2\ran$ for $a_1\in \F_q$,
$\om\in \F_q^*$, containing
$z_2$ (consider $y v_1$) which are of degree $2$.
\vskip 0.2cm

Similarly, one can calculate $\tilde{\h}^{mot}$ for
$T(3,2)$ and the unknot with $lk=3$,
which is for $\tilde{\r}=\F[[e_1+e_2,x=z_1^2+z_2, y=z_1^3+z_2^2]]$.
 The corresponding equation is $(x^3-y^2)(x^2-y)=0$. One has:
$\tilde{\h}^{mot}=$
{\small
$q^3(q-1) t^4+(1\!+\!a q) (1\!+\!a q^2) 
\bigl(1+(q-1) t\bigr)+(1\!+\!a q) 
\bigl(q^2 t+q (q^2-1)t^2+q^2 (q-1) t^3\bigr).$
}

Now $\de=4$, $\mathfrak{c}=\o z_1^5+\o z_2^3$, and
$G=\nu(\o)\setminus \Ga=$ 
$\{1/2,1,2,4\}$, the valuations of $e_2,z_1, z_1^2, z_1^4$, where
we take $\nu^1=0, \nu^2=1/2$ as above.

Note that there are $2$ families of $M$ of maximal
$r\!k_q=3$ now:
 $\o$ and $\lan e_1+\om e_2+a_1 z_1, z_1^2,z_1^3\ran$ for
$\om\in \F_q^*$ and $a_1\in \F_q$.

\section{\sc Unibranch {\sl L}-functions}
Let us begin with the unibranch case: when 
$\r\subset \o=\F_q[[z]]$ is with $2$ generators and
the localization $\F_q((z))$. 
Recall that $\de=$dim\,${}_{\F_q}\o/\r$
$=|\,\Z_+\setminus \Ga\,|$.
V.M. Galkin (in 1973) and St\"ohr (in 1998)
studied zeta-functions and $L$-functions
for arbitrary Gorenstein rings in dimension one. See
\cite{Sto} for a discussion of prior works. 
They considered non-irreducible
singularities and allowed different extensions of $\F_q$ for 
different component. We assume that the same $\F_q$ serves
all irreducible components; the case of different extensions 
of $F_q$ for different components seems beyond algebraic links
and the DAHA superpolynomials. 
 
Adding $a$ and the 
case of arbitrary $\o$-ranks closely follow what we did for the 
superpolynomials. However, we will not consider 
{\em colored} $Z$ and $L$
in this paper.                                          

\subsection{\bf Flagged {\bf\em L}-functions}
The {\em admissible  flags} of ideals in unibranch 
$\r=\F_q[[x,y]]\subset \F_q[[z]$ are 
$\vect{M}\!=\!\{M_0\!\subset\! M_1\!\subset\!
\cdots\!\subset\!  M_\ell\!\subset\! \r\}$ such that 
$\{z^{-m_0}M_i\subset \o\}$ 
for $m_0=\text{Min}(\De_0)$ are {\em standard flags} in
$\o$ from Section \ref{sec:stand}.
Then, the  {\em flagged zeta function} is defined as follows:
$Z(q,t,a)\!\!\equal\!
\sum_{\vect{M}} 
a^\ell t^{\hbox{\tiny \,dim}(\r /M_{\ell})}$,
where the summation is over all admissible flags of ideals
$\vect{M}\subset \r$,\,and we set dim$=$dim$_{\F_q}$.
The corresponding {\em flagged $L$-function}  is 
$L(q,t,a)\equal(1-t)Z(q,t,a)$,
which is a polynomial in terms of $t$ (use the conductor of 
$\r$ to see this).

Similar to Theorem \ref{thm:a-formula},
we obtain the following decomposition:
\begin{align}\label{Z-exp}
&Z(q,t,a)=\!\!\sum_{M\subset \r} 
t^{\hbox{\tiny \,dim}(\r/M)}(1\!+\!aq)\cdots
(1\!+\!aq^{r\!k_q(M)-1}),
\end{align} 
summed over all ideals $ M\subset \r$,
where $r\!k_q(M)=$\,dim\,$M/\mathfrak{m}_\r M.$ 

There is an immediate application of this formula:
\begin{align}\label{z-zuniga}
Z(q,t,a\!\!=-1/q)\!=\! Z_{\hbox{\tiny prncpl}}
(q,t)=\sum_{(M=f\, \r\subset \r} 
t^{\hbox{\tiny \,dim}(\r /M)}.
\end{align}
The latter 
is the {\em Z\'u{\oldt{n}}iga 
zeta function}: the summation is only over {\em principle} 
ideals $f\, \r\subset\r$; see \cite{Sto,Zu}.
 Recall that \,dim\,$_{\F_q}(\r/(f\,\r)=
\nu_z(f).$ 
When $q\to 1$ (for the ``field" with $1$ element):\,
$Z_{\hbox{\tiny prncpl}}=\sum_{\nu\in \Ga}t^\nu$. Then
$L_{\hbox{\tiny prncpl}}=(1-t)Z_{\hbox{\tiny prncpl}}$
is essentially the {\em  Alexander polynomial}\,:
$$
\lim_{q\to 1}L_{\hbox{\tiny prncpl}}=(1-t)\bigl(\sum_{i=1}^{\de}
t^{g_i}\bigr)+t^{2\de} \text{ for } \{g_i\}=\Ga \setminus 
(2\de + \Z_+).
$$
We obtain that $\bigl(L_{\hbox{\tiny prncpl}}-t^{2\de}\bigr)/(1-t)$
becomes $\de$ when $q\to 1$ and $t=1$. 
\vskip 0.2cm

Let $\r(i)=\r\cap z^{\ga_i} \o$,
where $\Ga=\{\ga_0=0<\ga_1<\ga_2<\cdots\}$. Then \, 
dim\,$\r(i)/\r(i+1)=1$.
Following the proof of Theorem \ref{thm:2ndexp}$(ii)$, 
there is another decomposition:
\begin{align}\label{2ndZ-exp}
Z(q,t,a)\!=\! \sum_{M\subset \r} t^{deg(M)}\,
\Bigl(1\!+\!\frac{a}{t}\Bigr) \Bigl(1\!+\!\frac{a}{t}q\Bigr)\cdots 
\Bigl(1\!+\!\frac{a}{t}q^{\,d^\Box(M)-1}\Bigr),
\end{align}
where the summation is over all ideals $M\subset \r$ and 
we set $d^\Box(M)=$\, dim\,$_{\F_q}(M^\Box/M)$ 
for  $M^\Box\equal \bigl\{f\in \r(m) 
\mid f\mathfrak{m}_{\r}\subset M\bigr\}$, where
$m=v_M=\min\{\De(M)\}$ 
for $\De(M)=\nu(M)$.

This is a version of Theorem \ref{thm:2ndexp}
where $\o$ is replaced by $\r$;
the justification remains basically the same.
Recall, that $M=\o$ was a unique standard module such that
$M=M^\Diamond$  (and $d^\Diamond=0$). Similarly,  
$M=M^\Box$ and $d^\Box=0$ for $M\subset \r$
if and only if $M=\r(i)$ for some 
$i$.

As an application of (\ref{2ndZ-exp}), the {\em second 
decomposition}, we obtain that $Z(q,t,a=-t)=(1-t)^{-1}$ 
and $L(q,t,a=-t)=1$, which holds for
arbitrary $\r$, not only Gorenstein. 
\vskip 0.2cm

{\sf Functional equation.}
This is one of the key properties of $L(q,t,a)$:\, 
\begin{align}\label{functLuni}
t^{-\de}L(q,t,a) \text{ is invariant 
under } t\mapsto 1/(qt), q\mapsto q, a\mapsto a.
\end{align}
The transformation  $t\mapsto 1/(qt)$ is the one
for the {\em Hasse-Weil zeta-functions}.
Note that the functional equation does not hold for $Z(q,t,a)$.

When $a=0$ and $a=-1/q$, it suffices
to assume that $\r$ is {\em Gorenstein}.
Recall that it is Gorenstein 
if the {\em combinatorial reciprocity map} 
$g\mapsto g'=2\de\!-\!1\!-\!g$ is surjective; then this map 
identifies $\{g\in \Z_+\setminus \Ga\}$ (the set of ``gaps") with
$\{g'\in \Ga\setminus \{2\de+\Z_+\}\}$. For instance, the last
gap, which is $2\de-1$, maps to $g'=0$. 

The functional equation
for $L$ is directly related to the {\em reciprocity} for
$\r$-modules. The proof from \cite{Sto} is for $a=0$; 
the key  is Theorem 3.1 there. For any $a$ 
we use that $r\!k_q$ is
preserved by the reciprocity map, which requires plane curve 
singularities; see Theorem \ref{thm-feq}.


\subsection{\bf {\bf\em L}-functions and superpolynomials}
Geometrically, the functional equation for $L$
is somewhat surprising because there is no Poincar\'e duality for 
singular varieties, unless the intersection cohomology is used or so.
However, Tate's $p$-adic proof works well 
for curve singularities
(at least for Gorenstein rings).
St\"ohr found a short and entirely combinatorial proof of the
functional equation for $L$,  a significant
simplification of the approach from John Tate's thesis. 
\vskip 0.2cm

Let us state the {\em Coincidence Conjecture}  $\H\!=\!L$
in this
case, to be generalized later to non-unibranch 
plane curve singularities. We set
$\H^{mot}(q,t,a)\!\equal\!\h^{mot}(qt,t,a)$
for motivic $\h^{mot}$ above. The substitution 
$q_{new}=q/t$ was already used for {\em RH} for superpolynomials.

We conjecture that $\H^{mot}(q,t,a)\!=\!L(q,t,a)$ for 
any plane curve singularities. In particular, 
$\H^{mot}(q,t,a\!\!=-1/q)\!=\! L_{\hbox{\tiny prncpl}}
(q,t)$, which presumably holds for any sufficiently
general Gorenstein $\r\subset \o$.
Recall that formula (\ref{z-zuniga}) provides that
 $L(q,t,a=-1/q)=$
$L_{\hbox{\tiny prncpl}}
(q,t)$.

Composing the (conjectural) coincidences
 $\h^{daha}=\h^{mot}$ and 
$\H^{mot}=L$, we 
identify the superduality $q\leftrightarrow t^{-1}$
for $\h$ with the
functional equation $t\mapsto 1/(qt)$ for $L$, a remarkable
link from  DAHA  to  the Hasse-Weil-Deligne theory. 
See \cite{ChS} for further discussion, including 
(potential) interpretation
of this passage as an instance of the fundamental relation
between {\em SCFT}, super-conformal field theory (the DAHA side), 
and Landau-Ginzburg sigma model, {\em LGSM} 
(the motivic side).
\vskip 0.2cm

In contrast to the Hasse-Weil zetas for smooth projective 
curves, the Riemann Hypothesis 
for $L(q,t,a=0)$ holds only for sufficiently small $q$, 
which we proved above for $\H^{mot}(q,t,a=0)$ upon the assumption
that the dependence on $q$ is polynomial. 
Recall that the conjectural upper $RH$-bound for $\H(q,t,a\!=\!0)$
is $q\le 1/2$ in the unibranch uncolored case:\,
 far from ``arithmetic" values
$q=p^m$ for $L(q,t,0)$. 
\vskip 0.2cm

The coincidence $\H^{mot}(q,t,a)$ and $L(q,t,a)$ 
for $t=1$ is as follows (for any rings $\r\subset \o$, 
not only Gorenstein).
Any admissible flag of ideals $\vect{M}'\subset \r$
is $\,z^{m'}\, \vect{M}\,$ for standard $\vect{M}\subset \o$, where
$m',\vect{M}$ are uniquely determined  by $\vect{M}'$. 
Vice versa, given a standard flag $\vect{M}$, we set:\,
$$\{m \mid z^{m}\vect{M}\subset \r\}
=\{0\le m_1<m_2<\cdots<m_k<2\de\}\cup \{2\de+\Z_+\}$$ 
for proper $k=k_M$ and $\{m_i\}$.  
The contribution of $z^m \vect{M}\subset \r$ for $m$ in this list
to $L(q,t,a)=(1-t)Z(q,t,a)$ is 
$(1-t)t^{m-dev(M_\ell)}$. Recall:
$dev(M)\!=\!\de\! -\!\text{dim}\,\o/M$.
Given $\vect{M}$, the total contribution becomes
$(1-t)t^{-dev}
\bigl(\sum_{i=1}^{k_M} t^{m_i}+t^{2\de}/(1-t)\bigr)$, which 
approaches $1$ 
as $t\to 1$ (for any $q)$. Thus, we need to count standard
$\vect{M}\subset \o$ and arrive at the formula $L(q,t=1,a)=
\sum_{\ell=0}^{2\de-1} |\j_{\ell}(\F_q)|a^\ell$, which is
$\H^{mot}(q,t=1,a)$.
\vskip 0.2cm

{\sf Deforming $\de$ and $\mu$}. 
The reduction to $L_{prncpl}$ has many interesting features,
including the following deformation of $\de$ and the Milnor
number $\mu$ (the Witten index in the corresponding Landau-Ginzburg
theory). 

Let 
$G=\Z_+\!\!\setminus\! \Ga=$
$\cup_{i=1}^{\varpi}\, \{g_i\le x \le g'_i\}$, where
$g'_i+1\in \Ga$,
and $m_i\!=\!g'_i\!-\!g_i\!+\!1$;\, $\varpi$ counts the number
of $\ga\in \Ga$ such that $\ga+1\not\in \Ga$. 

Assuming ``\,$\H\!=\!L$" at $a=-1/q$ for
Gorenstein $\r\subset \C[[z]]$, let 
$\de(q,t)\!\equal\!$
{\large $\frac{\h^{mot}(q,t,a=-t/q)-(qt)^\de}{1-t}
=\frac{L_{pncpl}(q/t,t)-(qt)^\de}{1-t}$}. Then
\begin{align}\label{de-qt}
\de(q,t)\!=\! \text{\small $\frac{1\!-\!t^{g_1}}{1\!-\!t}+
\sum_{i=1}^{\varpi-1}\frac{t^{g_i'+1}-t^{g_{i+1}}}{1-t}
(\frac{q}{t})^{m_1\!+\cdots+\! m_i}$},\ \de(1,1)=\de.
\end{align}
The parameters $q,t$ from $\h^{daha}$ and $\h^{mot}$,
are used here, 
not $q_{new}=q/t$ from  the definition of $\H$ and $L$.
Let us deform $\mu$.


\vskip 0.2cm

One has:\,
$\mu(q,t)\equal\de(q,t)+(qt)^{\de-1} \de(t^{-1},q^{-1})=$
$\sum_{x=0}^{2\de-1} t^{v(x)-1}q^{g(x)}$, where
$v(x)=|\{\nu\!\in\! \Ga \mid 0\!\le\! \nu\!\le\! x\}|$ and 
$g(x)=|\{g\!\in\! G \mid 0\!\le\! g\! <\! x\}|$ for $G$ as above.
Then  $\mu(1,1)=\mu$.  The superduality of $\mu$ is granted by
construction:\, 
 $(qt)^{\de-1}\mu(1/t,1/q)=\mu(q,t)$.

All monomials are monic in $\de_{q,t}$, but not 
all of them are monic in $\mu(q,t)$. Exactly 
$\varpi$ of them are with coefficient $2$:\,
for $\ga\in \Ga \not\ni \ga+1$.
So this is not a ``full split" of $\mu$.

For example,
$\mu(q,t)\!=${\small 
$2 + q + q^2 + 2 q^3 t + 2 q^4 t^2 + 2 q^5 t^3 + 2 q^6 t^4 + q^7 t^5 + 
 q^7 t^6 + 2 q^7 t^7$}
for $\r=\C[[z^4,z^6\!+\!z^7]].$ The monomials with coefficients
$2$ correspond to $\ga=0,4,6,8,10,14$. 
Deformations of the {\em Witten invariant} and the 
{\em BPS-invariant} are
of interest to physicists.  

\subsection{\bf Quasi-rho for cables}\label{sec:rhocab}
Using the superduality, the following variant of the
formula for $\de(q,t)$ is quite reasonable algebraically:\,
$\varrho(q,t)\equal \frac{\h(q,t,a=-t/q)-q^{\de}}{(1-t)(1-q)}$.
The divisibility by $(1-t)$ is granted as above; it holds 
when any $q^\de t^m$ instead of $q^{\de}$
is subtracted. Taking $q^\de$ here
ensures the invariance of the numerator 
with respect to the superduality
$q\leftrightarrow 1/t$ and, accordingly, its divisibility
by $(1-q)$. Following \cite{ChS},
\vskip 0.2cm

\centerline{
$\varrho(q,t)=\sum_{x\in G}q^{g(x)}\frac{1\!-\!t^{v(x)}}{1-t}\!=\!
\!\sum_{G\ni x>y\in \Ga}q^{g(x)}t^{v(y)-1}=(q t)^{\de-1}
\varrho(t^{-1},q^{-1}).
$
}

\vskip 0.2cm
Presumably, 
$\varrho(1,1)$ is a version  of
the classical {\em $\rho_{ab}-$invariant} for cables $K$. The latter is
the {\em von Neumann rho-invariant} defined for the abelianization
representation $ab: \pi_1(S^3\setminus K)\to H_1(S^3\setminus K)=\Z$. 
This is,  basically, 
the {\em $\eta$-invariant}  of the $3$-fold $S^3\setminus K$ 
(it must be odd-dimensional) minus that for 
the line bundle over this manifold 
corresponding to the homomorphism $ab$. Also, one has that  
 $\rho_{ab}=\int_0^1 \si_K(e^{2\pi i x})dx$ for
the {\em Levine-Tristram signature} $\si_K$;  see e.g., \cite{Li}.
If a $K$-theoretical of the 
{\em Seifert pairing} can be defined (used here), 
it may be related to
our $q,t$-formulas.

For $r,s>0$ 
such that gcd$(r,s)=1$,
one has:
$
\rho_{r,s}=\varrho_{r,s}(1,1)=\frac{(r^2-1)(s^2-1)}{24}.
$
We note that this is the maximal size of the 
{\em $(r,s)$-core partitions}.                
The classical $\rho_{ab}$ is 
$-\frac{1}{3}\frac{(r^2-1)(s^2-1)}{rs}$.
Thus, we basically obtain the same formula up to some
renormalization.   A challenge is to interpret
our  $\rho(q,t)$ geometrically.

\vskip 0.2cm

The classical $\rho_{ab}$ 
for cables is known to be 
additive. For instance let $K=C\!ab(m,n)C\!ab(s,r)$.
Then 
$\rho_{ab}=-\frac{1}{3}\bigl(\frac{(m^2-1)(n^2-1)}{mn}
+\frac{(r^2-1)(s^2-1)}{rs}\bigr)$. Our $\varrho(1,1)$
is additive with some weights. This is as follows.


\vskip 0.2cm

Let $\r=\F[[z^{\upsilon r}, z^{\upsilon s} +z^{\upsilon s+p}]]$,
where gcd$(r,s)=1$ as above, $\upsilon>1$ and gcd$(\upsilon,p)=1$ for
$p\ge 1$. Then $\Ga=\lan \upsilon r,\,
\upsilon s,\,\upsilon r s\!+\!p\ran$,\ 
$2\de=\upsilon^2 r s\!-\!\upsilon (r\!+\!s)+(\upsilon\!-\!1)p+1$ and 
$K=C\!ab(m\!=\!\upsilon r s\!+\!p,\,
n\!=\!\upsilon)\,
C\!ab(s,r)$. One obtains:  $\varrho_K(1,1)=$
{\small $\frac{1}{24}\bigl((m^2-1)(n^2-1)
+\upsilon^2(r^2-1)(s^2-1)\bigr)$}. 
\vskip 0.1cm

More generally,
$\rho_k=\varrho_K(1,1)=$
{\small $\frac{1}{24}\sum_{i=1}^k\upsilon_i^2 (a_i^2-1)(r_i^2-1)$}
for the cable $K=C\!ab(a_k,r_k)\cdots C\!ab(a_2,r_2)C\!ab(a_1,r_1)$,
where  $1\le i\le k$, $\upsilon_i=r_k\cdots r_{i+1}$ and
$\upsilon_{k}=1$. We will post the
details elsewhere. Here
$\upsilon_i\!=$\,gcd$(u_1,\cdots,u_{i+1})$ for 
$\Ga=\lan u_1,\cdots,u_{k+1}\ran$, where  $u_i<u_{i+1}$ and
$\upsilon_{i+1}|\upsilon_i$; we note that 
$\de=\de(1,1)=$
{\small $\frac{1}{2}\sum_{i=1}^k\upsilon_i (a_i-1)(r_i-1)$}.

Importantly, the values of our  $\rho(1,1)$ are 
natural numbers and $\varrho(q,t)$ is a sum of
monic $q,t$-monomials. This is perfect for  
{\em categorification}! The $\rho_{ab}$\~invariant and
similar spectral ones  are not such generally.
 
\vskip 0.1cm

The invariants $\de(q,t), \mu(q,t)$
and $\varrho(q,t)$ have many symmetries.
 For instance, consider algebraic knots
$K$ and $K'=C\!ab(a,r)K$.
If $\de$ is that for the ring $\r$ of $K$, then 
 $\rho_{K}(q,t)$ is the sum of monomials $q^it^j$ in
$\rho_{K'}(q,t)$ such that $j<2\de-1$. This is just to
provide an example.

\subsection{\bf Super rho-invariants}
Let us now add $a$ to  $\varrho(q,t)$, which will be 
a ``triply-graded" deformation of $\varrho(1,1)$:
\begin{align*}
&R_K(q,t,a)\equal
\bigl(\h_K(qt,t,a)-t^{\de}\h_K(qt,1,a)\bigr)/
\bigl((1-q)(1-t)\bigr),\\
&\varrho_K(q,t)\equal R_K(q,t,a\!=\!-\!t/q),\ \text{where}\,\ 
\h_K(q,1,a\!=\!-1/q)=q^{\de}. 
\end{align*}

The superduality for $R$ reads:
$R(q,t,a)=(qt)^{\de-1}R(t^{-1},q^{-1},a)$.
The substitution $a\mapsto -\frac{t}{q}$ is compatible with
the superduality; $t/q$ is its invariant.  Generally, such 
specialization
is related to the {\em Heegaard-Floer homology} and 
Alexander polynomials (for $q\!=\!1,a\!=\!-1$).

\vskip 0.2cm

{\sf The case of C\!ab(13,2)C\!ab(2,3)}.
Here $\r\!=\!\C[[z^4,z^6\!+\!z^7]]$,
 $r\!=\!3,s\!=\!2,\upsilon\!=\!2$ (see above), $\de\!=\!8$. Then
$\rho=\varrho(1,1)=
\frac{(13^2-1)(2^2-1)^2}{24}+2^2\frac{(3^2-1)(2^2-1)}{24}=
25$ and  its refined
version (its $q,t$-split) is $\varrho(q,t)=$ 
\renewcommand{\baselinestretch}{1.2} 
{\small
\(
1+q+q^2+q^3+q^4+q^5+q^6+q^7+q^3 t+q^4 t+q^5 t+q^6 t+q^7 t+
q^4 t^2+q^5 t^2+q^6 t^2+q^7 t^2+q^5 t^3+q^6 t^3+q^7 t^3+
q^6 t^4+
q^7 t^4+q^7 t^5+q^7 t^6+q^7 t^7.
\)
}
\renewcommand{\baselinestretch}{1.2}

\vskip 0.1cm
We note that $RH$ holds for 
$\varrho(qt,t)$ when  $q\!<\!q_{sup}\!\approx\! 0.802$. 
Presumably, $\lim_{n\to\infty}q_{sup}\!=\!1$ 
for $C\!ab(13\!+\!2n,2)C\!ab(2,3)$;
for instance, $q_{sup}\!\approx\! 0.996$ for $n\!=\!2000$.

\vskip 0.2cm


The corresponding 
 $R(q,t,a)$ for  $\C[[z^4,z^6\!+\!z^7]]$ is as follows:
\renewcommand{\baselinestretch}{1.} 
{\footnotesize
\( 
1+q+q^2+q^3+q^4+q^5+q^6+q^7+t+2 q t+3 q^2 t+4 q^3 t+4 q^4 t+
4 q^5 t+4 q^6 t+q^7 t+t^2+2 q t^2+4 q^2 t^2+6 q^3 t^2+8 q^4 t^2+
8 q^5 t^2+4 q^6 t^2+q^7 t^2+t^3+2 q t^3+4 q^2 t^3+7 q^3 t^3+
10 q^4 t^3+8 q^5 t^3+4 q^6 t^3+q^7 t^3+t^4+2 q t^4+4 q^2 t^4+
7 q^3 t^4+7 q^4 t^4+6 q^5 t^4+4 q^6 t^4+q^7 t^4+t^5+2 q t^5+
4 q^2 t^5+4 q^3 t^5+4 q^4 t^5+4 q^5 t^5+3 q^6 t^5+q^7 t^5+t^6+
2 q t^6+2 q^2 t^6+2 q^3 t^6+2 q^4 t^6+2 q^5 t^6+2 q^6 t^6+q^7 t^6+
t^7+q t^7+q^2 t^7+q^3 t^7+q^4 t^7+q^5 t^7+q^6 t^7+q^7 t^7+
a^3 \Bigl(q^6+q^7+q^6 t+q^7 t\Bigr)+
a^2 \Bigl(q^3+2 q^4+3 q^5+3 q^6+3 q^7+q^3 t+
3 q^4 t+6 q^5 t+8 q^6 t+3 q^7 t+q^3 t^2+3 q^4 t^2+7 q^5 t^2+
6 q^6 t^2+3 q^7 t^2+q^3 t^3+3 q^4 t^3+3 q^5 t^3+3 q^6 t^3+
2 q^7 t^3+q^3 t^4+q^4 t^4+q^5 t^4+q^6 t^4+q^7 t^4\Bigr)+
a \Bigl(q+2 q^2+3 q^3+3 q^4+3 q^5+3 q^6+3 q^7+q t+3 q^2 t+6 q^3 t+
9 q^4 t+10 q^5 t+10 q^6 t+3 q^7 t+q t^2+3 q^2 t^2+7 q^3 t^2+
12 q^4 t^2+17 q^5 t^2+10 q^6 t^2+3 q^7 t^2+q t^3+3 q^2 t^3+
7 q^3 t^3+13 q^4 t^3+12 q^5 t^3+9 q^6 t^3+3 q^7 t^3+q t^4+
3 q^2 t^4+7 q^3 t^4+7 q^4 t^4+7 q^5 t^4+6 q^6 t^4+3 q^7 t^4+
q t^5+3 q^2 t^5+3 q^3 t^5+3 q^4 t^5+3 q^5 t^5+3 q^6 t^5+
2 q^7 t^5+q t^6+q^2 t^6+q^3 t^6+q^4 t^6+q^5 t^6+q^6 t^6+q^7 t^6\Bigr).
\)
}
\renewcommand{\baselinestretch}{1.2}

\subsection{\bf Using Hilbert schemes}
For a {\em rational} projective curve $C\subset \mathbb P^2$,
 the following formula was a starting point of quite a few developments
in modern enumerative geometry. It is due to
Gopakumar-Vafa and Pandharipande-Thomas \cite{PaT}:
$$
 \sum_{n\!\ge\! 0} q^{n\!+\!1\!-\!\de} e(C^{[n]})\!=\!
\sum_{0\le i\le \de} n_{C}{(i)}\bigl(\!\frac{q}
{(1\!-\!q)^2}\!\bigr)
^{i+1\!-\!\de},  
$$
which is for the Euler numbers of {\em Hilbert schemes} 
 $C^{[n]}$. The points of the latter are
zero-cycles of $C$, collections of  ideals at any points, of the
(total) colength $n$. Here $\de$ is the arithmetic genus of $C$,
$n_{C}(i)$ are some numbers. The passage from a series 
to a polynomial is far from obvious
even in this relatively simple case.
It is significantly more subtle 
to prove that $n_{C}{(i)}\!\in\! \Z_+$ 
(G\"ottsche, and then Shende for all $i$); 
versal deformations of 
singularities occurred in the Shende's proof \cite{Sh}.


Switching to local rings $\r$ of singularities, the 
following conjecture is 
for {\em nested\,}  Hilbert schemes
$H\!ilb^{\,[l\le l+m]}$ formed by pairs $\{I,I'\}$ of ideals in $\r$
of colengths $l,l+m$ 
such that $\mathfrak{m}I\!\subset\!
I'\!\subset\! I$
 for the maximal ideal $\mathfrak{m}\subset \r$.
One needs the {\em weight $t$-polynomial} 
$\, \mathfrak{w}(H\!ilb^{\,[l\le l+m]})$ defined for the weight
filtration of $H\!ilb^{\,[l\le l+m]}$
due to  Serre and Deligne.
 The
{\em Oblomkov-Rasmussen-Shende conjecture} from \cite{ORS}
 states that 
$$
\sum_{l,m\ge 0}q^{2l}\,a^{2m}\,
t^{m^2}\,\mathfrak{w}(H\!ilb^{\,[l\le l+m]})$$
is proportional to the Poincar\'e  series of 
the HOMFLY-PT triply graded
homology  of the
corresponding link. The connection with the perverse
filtration of $\j_0$  is due to Maulik-Yun 
and Migliorini-Shende \cite{MY,MS}. 
 The {\em ORS series} is a geometric
variant of $Z(q,t,a)$. See also \cite{KT}; the conjecture $\H\!=\!L$
is proven there for (uncolored) $T(3, 3n\pm 1)$.  

The ORS-conjecture adds $t$ to
the {\em Oblomkov-Shende conjecture}, which was
extended by adding  colors $\la$ and 
proved by Maulik \cite{Ma}.
\vskip 0.2cm

Our motivic superpolynomials presumably 
coincide with
the {\em reduced\,} KhR-polynomials.
The  Cherednik-Danilenko conjecture was that the latter coincide with
the uncolored 
$\h^{daha}(q,t,a)$. Recall that the polynomial
dependence on $q,t,a$   is manifest
for the DAHA superpolynomials (with some work to do for $a$).
Motivically, it is manifest for $t,a$; the
polynomial dependence on $q$  of $\h^{mot}$
remains a conjecture. 
\vskip 0.2cm


\section{\sc Unification: colored links}
We extend in this section 
the definition of motivic superpolynomials
to the most general ``characteristic polynomials"
$\prod_{i=1}^\kappa F_i(u,v)^{c_i}$ for 
pairwise non-proportional (absolutely) irreducible $F_i$ over $\F_q$.
The singularity remains the same:\, for 
$F(x,y)=\prod_{i=1}^\kappa F_i(u,v)$.
The DAHA superpolynomials for algebraic links colored
by the sequences of weights 
$\si=(c_i\om_1, 1\le i\le \kappa)$ are their counterparts
in the Coincidence Conjecture $\h^{daha}=\h^{mot}$.
The sequence $\si=(c_i)$ will be assumed ordered: 
$c_1\ge c_2\ge \cdots c_\kappa>0$, which
can be achieved by permuting $\{F_i\}$. 
\vskip 0.2cm

It will be convenient
to extend the idempotents $\{e_1,\ldots, e_\kappa\}$
to idempotents $\{\ep_j, 1\le j\le \tau\}$ setting
$e_1=\ep_1+\cdots+\ep_{c_1}$,
$e_2=\ep_{c_1+1}+\cdots+\ep_{c_1+c_2}$, and so on. 
The valuation
$\nu$ is now defined for the sequence 
$\nu^1=0<\nu^2<\cdots<\nu^\tau<1$,
where $\nu(\ze_j)=\nu_j$ as above. As a matter of fact,
any $\nu^j\in \R_+$ can be taken here
provided that $\nu^j-\nu^k\not\in \Z$ for
$j\neq k$. 

Accordingly,    
$z_1=\ze_1+\cdots+\ze_{c_1}, \,\ldots,\,
z_\kappa=\ze_{\tau-c_{\kappa}+1}+\cdots+
\ze_{\tau}$ for $\ze_j=\ze \ep_j$.  
This gives a natural
diagonal embedding of 
the prior $\r=\F_q[[x,y]]$ and $\o$ into
$\Om\!\equal\!\F[[\ep_1,\ldots, \ep_\tau,
\ze_1,\ldots, \ze_{\tau}]]=\sum_{i=1}^{\tau} \Om_i$, where
$\Om_j=\ep_j\Om$ for $1\le j\le \tau$, and
$\o_i=e_i \Om$ for $1\le i\le \kappa$.
Recall that  $\r\subset \o=\sum_{i=1}^\kappa \o_i$.

Let $x=\sum_{i=1}^\kappa x_i=
\sum_{j=1}^\tau \chi_j$ and   $y=\sum_{i=1}^\kappa y_i=
\sum_{j=1}^\tau \xi_j$, where we use
 $\chi_j\equal x \ep_j$ and  $\xi_j\equal y\ep_j$.

Actually, the fact that $\Om$ is a ring is not needed in the
definition of standard modules and the corresponding motivic 
superpolynomials.
The action of $\o$ is the key:\, $z^k_i\ze^m_j=\ze_j^{k+m}$ for the
indices $j$
from the segment corresponding to $e_i$
and $0$ otherwise. However, this is
not only for the sake of better readability 
as we will see in Proposition \ref{prop:link-rank}.

\vskip 0.2cm

\subsection{\bf General superpolynomials}
Similar to the prior definitions, 
{\em standard modules} $M$ are 
$\r$-invariant $\F_q$-subspaces $M\subset \Om$
such that $\o M=\Om$. 
As above, $r\!k_q(M)\!\equal$\,dim${}_{\F_q}
M/\mathfrak{m}_{\r}\,M$, where  
$\mathfrak{m}_{\r}=\r\cap z\o$.
The minimal $q$-rank is then $r\!k_{min}=c_1$ and
the maximal $q$-rank, $r\!k_{max}$, is that for $M=\Om$.
It equals $\sum_{i=1}^\kappa m_i c_i$, where $m_i$ are the
multiplicities
of singularities associated with
$F_i(x,y)=0$, which are $m_i=\min\{\Ga_i\setminus\{0\}\}$
for the corresponding $\Ga_i$. 

The 
{\em motivic superpolynomial} of $\r,\si$ is defined as follows:
$$
\h^{mot}_{\si}=\sum_{M} t^{\hbox{\tiny dim}(\Om/M)}
\prod_{j=r\!k_{min}}^{r\!k_q(M)-1}(1+aq^j)
\text{ summed over standard } M.
$$
The product $\Pi$ for $r\!k_q(M)=r\!k_{min}=c_1$ is $1$;
the dimension dim$_{\F_q}(\Om/M)$ will be  
denoted by $deg(M)$ (as above).
We note that
formula (\ref{2ndHfomula}),
the $(1\!+\!q^i a/t)$-decomposition,
does not work as such now.

The $a$-degree of
$\h^{mot}_{\si}$ is $\sum_{i=1}^\kappa m_i c_i-c_1$, where $-c_1$
is because $r\!k_q(M)$ begins with $c_1$. The 
standard modules of rank $c_1$ are not now 
``invertible" in any sense. For instance,
they can be unions of families of different degrees. 
The maximal 
$deg(M)$ among all standard $M$
is that for $M^\circ=\lan v^\circ_i\ran$, the span of the 
{\em sums} of the $i${\tiny th} $\ep$-vectors (correspondingly) in the groups
$\{\ep_1,\ldots, \ep_{c_1}\}, \{\ep_{c_1+1},\cdots, \ep_{c_1+c_2}\}$
and so on. For instance,  
$v^\circ_1=\ep_1+\ep_{1+c_1}+\ep_{1+c_1+c_2}+\cdots$\,\ldots, 
$v^\circ_{c_1}$ is the sum of the last $\ep$-vectors in the
groups with $c_1=c_2=\cdots$. 
\vskip 0.2cm

Let us state the Coincidence Conjecture in full generality
(within this work). Recall that its unibranch uncolored
version is Conjecture \ref{con:motknot}. The uncolored version
for links is Conjecture \ref{coin-links}.

\begin{conjecture}\label{conj:maincon}
The superpolynomial $\h^{mot}_\si$  depends polynomially on $q$
and coincides with the DAHA superpolynomial
$\h^{{\la}}(q,t,a)$ for the corresponding link colored
by the sequences $\{\la\}=\{c_1\om_1,c_2\om_1,\ldots,
c_\kappa \om_1\}$ (only pure rows). In particular, 
$\h^{mot}_\si$  are topological
invariants of the corresponding links and satisfy the
superduality in the uncolored case (when $\si=(1,\cdots,1)$):\,
$(qt)^\de \h^{mot}(1/t,1/q,a)=\h^{mot}(q,t,a)$.
 \sq
\end{conjecture}

{\sf The embedding $\h^{mot}_{\rr,\ss}(c)\subset
\h^{mot}_{c\,\rr,c\,\ss}$}\,. Its special version was part of 
Proposition \ref{thm:h-iter} and (\ref{h-dtor}). 
The justification in Proposition \ref{prop:link-rank}
utilizes some special features of $\r$ for torus
knots, but it can be true for any unibranch singularities $\r$.

\comment{
Let us fix $i\in \{1,\ldots,\kappa\}$..
One can assume that $z_i\in \Om_i\cong \o_i^{\oplus c_i}$ is  a 
{\em regular} diagonal matrix $\tilde{z}_i=diag(u_j)z_i$ acting there. 
Then $\Om_i$ becomes  
$\o[[\ze_i]]/\bigl(\prod_{j=1}^{c_i}(\ze_i-u_j z_i)\bigr)
\cong \oplus \o[[z_i]]/(z_i-u_j)$ as an algebra.
Accordingly, $\tilde{x}_i$ and $\tilde{y}_i$, will be given by the same
series in terms of $z_i$, where $z_i$ is replaced by 
 become the corresponding diagonal
matrices and the elements of this
algebra.  The equation $F_i(x_i,y_i)=0$ in $\r_i$ becomes
equation $\Phi_i(x_i,y_i)=0$ can be assumed
{\em square-free} for sufficiently general $u=(u_j, 1\le j\le c_1)$.
Extending this construction to all $i$, we obtain an algebra
$\o^u$ depending on $u=(u_j, 1\le j\le \tau)$, where 
$\Phi(x,y)=0$ for a square-free polynomial $\Phi$ for sufficiently
general $u$. 

The corresponding idempotents of
the components 
(depending on $u$) will be denoted $\vep_1,\ldots,\vep_\tau$.
We set $\ze=\sum_{j=1}^{\tau}\ze_i$
and, correspondingly,
extend  $x_i$ and $y_i$ to $\xi$ and $\chi$, acting
Assuming that $u=(u_i, 1\le i\le \tau)$ is sufficiently general, 
$\Phi(\xi,\chi)=0$ for a {\em square-free}
equation of degree $\tau$. 
We arrive at
the following proposition. 
}

Generally, let $\r=\F_q[[x,y]]$ be a unibranch ring  
diagonally embedded in $\Om=\o^{\oplus c}=
\oplus_{j=1}^c \o \ep_j$.
For a sequence $\{\pi_j\in \F_q, 1\le j\le c\}$, we set  
$\tilde{x}=\pi_1\chi_1+\cdots+\pi_{c}\chi_{c}$
for $\chi_j=x \ep_j$.  There is no change of $y$; we set 
$\tilde{y}=y=\sum_{j=1}^c \xi_j$,  where $\xi_j=y\ep_j$. 
Accordingly, 
Let $\tilde{\r}\equal\F_q[[\tilde{x},\tilde{y}]]\subset \Om$.
If $F(x,y)=0$ in $\r$, then $\tilde{F}(\tilde{x},\tilde{y})=0$
for a square-free $\tilde{F}(u,v)$, where $\{\pi_j\}$ 
are generic. The corresponding  link is with $c$ 
components isomorphic to the knot for  $F(u,v)=0$. 

For torus knots $T(\rr,\ss)$:\, 
$\r_{\rr,\ss}=\F_q[[z^{\rr},z^{\ss}]]$ and
 $F(u,v)=u^{\ss}-v^{\rr}$, where  $gcd(\rr,\ss)=1$. Accordingly,
$\tilde{F}(u,v)=\prod_{j=1}^c
(u^{\ss}-\pi_j^\ss v^{\rr})$ serves $T(c\,\rr,c\,\ss)$, where
$\pi_j^{\ss}\neq 0$ are pairwise distinct
($\tilde{F}$ must be square-free).

\begin{proposition}\label{prop:link-rank}
For $\r_{\rr,\ss}$, let $\tilde{\r}$ be as above 
for pairwise distinct $\pi_j^{\ss}\neq 0$ considered 
as free (formal) parameters. Consider 
standard modules $\tilde{M}\subset \Om$
for $\tilde{\r}$  containing $\ep_1,\ldots,\ep_c$
modulo $\ze\Om$. Then the specialization
map $\pi_j\mapsto 1 (1\le j\le c)$ sends $\tilde{M}$ to
$c$-colored standard $M$ over $\r_{\rr,\ss}$, and
any such $M\subset \Om$
can be obtained by this construction. Moreover, 
$deg(M)=\deg(\tilde{M})$, where 
$deg(\cdots)=$\,dim\,$_{\F_q}(\Om/\cdots)$,
and this map preserves $\De$ and $r\!k_q$.   
 Accordingly, 
$(1+aq)\cdots(1+aq^{c-1})\h_{\rr,\ss}^{mot}(c)$ is naturally
embedded into (uncolored) $\tilde{\h}^{mot}$, that for $\tilde{\r}$. 
\end{proposition}
{\em Sketch of the proof.}
Let $\pi_j(u)\in \F_q[[u]]$ for a parameter $u$
such that $\pi_j(0)=1$ and $\pi_j(1)=\pi_j$. 
Consider $\tilde{\r}$-modules
$\tilde{M}=M(u)\subset \Om$ with the coefficients in $\F_q[[u]]$.  
Special {\em fibers} $M=M(u=0)$ are those generated
by any linear combinations of elements
of $\tilde{M}$ 
divided by the greatest possible powers of $u$. This
definition gives that $deg(M)=deg(\tilde{M})$, but
the corresponding $\De$ can possibly change at $u=0$
from that at $u=1$. This procedure is explicit (at least)
for torus knots.

We follow the justification
of Theorem 3.1 from \cite{ChP2} for colored torus knots.
The coordinates of the corresponding Piontkowski cells
are the coefficients of the generators of the corresponding
standard $\r$-modules
$M\in \Om$. These generators are defined
in terms the valuation in $\Om$, which is the same for $\tilde{\r}$
and $\r_\circ$, subject to a certain {\em elimination
procedure}. Given $\De^c$,  the resulting coefficients of 
the generators of 
$M\subset \Om$ with  $\De^c(M)=\De^c$ 
serve $\tilde{M}$ if $\De(\tilde{M})\ni 0,1/2$. 
This is the required embedding.   
The term $(1+aq)\cdots(1+aq^{c-1})$ is because 
such product factors begin with $(1+aq^{c-1})$ in the 
$c$-colored case. 

There is another approach based on taking $\pi_j$ 
for $1\le j\le c$ such that\, 
$\pi_j^{\ss}$\, are pairwise distinct roots of unity of order $c$. Then
$\tilde{F}(u,v)=u^{\,c\,\ss}- v^{\,c\,\rr}$, and
the corresponding singularity is quasi-homogeneous. This
significantly simplifies the considerations. For instance,
pure monomials in terms of $\ze_j$ can be taken now as generators
for any given $\De$; for instance, this gives that the corresponding cells
are never empty. They are generally not affine spaces, but the 
{\em elimination procedure} works well in this case for any $\De$. 
The corresponding cells are affine spaces  for $\De\ni 0,1/2$
and can be identified with those in \cite{ChP2}. \sq

\comment{
Note that $c$ can be $1$
if $m>1$. For instance, one can take $\tilde{\r}\subset \F_q[[z]]$
 $\r=\F_q[[x=z^{m\rr},y z^{m\ss}\phi]]$, where
$\phi\in 1+z\o$ and $gcd(\rr,\ss)=1$. The setting needed in 
Proposition \ref{thm:h-iter} is when we begin with any unibranch
$\r$ and construct $\tilde{\r}$ in $\F[[\ze]]$ for $\ze^c=z$,
i.e. when $c=m$. 
This gives the embeddings generalizing those in
(\ref{h-dtor}) and, more subtle, 
relations between different $\tilde{\r}$ that can be obtained
when we use different rings $\r$ corresponding to the same
singularity, i.e. rings
equivalent formally analytically (or topologically). 
Here 
the relation between $z$ and $\ze$ can be $P(\ze)=0$ for
monic polynomials $P$ with coefficients in $z\o$ and 
$z$ as the free term.

The corresponding connection formulas, generalizing those 
in Proposition \ref{thm:h-iter}, can be complicated.
Recall that 
it corresponds to the usage of Macdonald polynomials for 
$\c\om_1$ in the DAHA theory, which become involved for $c>2$, 
The passage to 
$\tilde{\h}^{mot}-(1+aq)\cdots(1+aq^{c-1})\h^{mot}(c)$, 
the first step, is quite parallel.
However, such differences must be interpreted using 
further embedding of this kind in terms of  $c'<c$. This
especially concerns relation between different $\tilde{\h}^{mot}$ 
corresponding to equivalent $\r$ (presenting the same singularity),
which can be subtle.
}
\vskip 0.2cm

\subsection{\bf Two colored examples}
As an example of the deformation procedure above,
let us consider  
the case of $2$-colored trefoil in full detail.
Generally, see (\ref{trefa-t}) for arbitrary $c$ and a more
uniform approach to counting the dimensions. We will do this
here {\em directly}.
\vskip 0.2cm

{\sf Trefoil of rank 2.}
The case $c=2$ is for $\r=\F_q[[x=\ze_1^2+\ze_2^2, y=
\ze_1^3+\ze_2^3]]\subset \Om=\F_q[[\ep_1,\ep_2,\ze_1,\ze_2]]$.
The standard modules are:

$\{M_0=\Om\}$, the only one of $r\!k_q=4$, which contributes $1$ to 
$\h^{mot}$,
multiplied here and below by 
$\prod_{i=1}^{r\!k_q-1}(1+q^i a)$;

$\{M_1=\lan \ep_1+ a_1\ze_1+a_2\ze_2, \ep_2+b_1\ze_1+b_2\ze_2\ran\}$
for $a_i,b_i\in \F_q$, ``invertible" ones of $r\!k_q=2, deg=2$,
contributing $q^4 t^2$;

$\{M_2=\lan \ep_1+ a_2\ze_2, \ep_2+b_2\ze_2, \ze_1+c_2 \ze_2\ran\}$
for $a_2,b_2,c_2\in \F_q$, which are of $r\!k_q=3, deg=1$,
contributing $q^3 t$;

$\{M_3=\lan \ep_1+ a_1\ze_1, \ep_2+b_1\ze_1, \ze_2\ran\}$
for $a_1,b_1\in \F_q$, which are of $r\!k_q=3, deg=2$,
contributing $q^2 t$.

Thus, $\h^{mot}_{3,2}(c=2)= (1\!+\!a q^2)(1\!+\!a q^3)+q^4t^2+
(1\!+\!a q^2) (q^3 t+q^2 t)$.
\vskip 0.2cm

Equivalently, one can consider $\r=\F_q[[x=z^2,y=z^3]]\subset
\F_q[[z]]\subset \Om=\F_q[[\ze]]$, where we set
$z=\ze^2=z$, accordingly, $\xi=\ze^2$ and $\chi=\ze^3$. 
 Then the
standard modules $M$ will contain $1$ vector 
from $1+z\Om$ and $1$ from $\ze+z\Om$. This was the setting of
Section 3 of \cite{ChP2}.  
\vfil

The $\pi$-deformation above is 
 $\tilde{\r}=\F_q[[x=\pi_1\ze_1^2+\pi_2 \ze_2^2, y=
\ze_1^3+\ze_2^3]]\subset \Om=\F_q[[\ep_1,\ep_2,\ze_1,\ze_2]]$.
The corresponding equation for $x,y$ will be then
$(x^3-\pi_1^3 y^2)(x^3-\pi_2^3 y^2)=0$, which is square-free
for generic $\pi_1,\pi_2$.
So this is {\em double trefoil}. The generators in the families above
can be used as such for the action of $\tilde{\r}$. We obtain
{\em all} standard $\tilde{M}$ such that $\De(\tilde{M})\ni 0, 1/2$.
Arbitrary (all) standard $\tilde{\r}$-modules
$\tilde{M}$ are such that $\De(\tilde{M})\ni 0$, which is 
a significantly wider class. 
\vskip 0.2cm

{\sf Colored Hopf $2$-link.} Firstly, the simplest uncolored 
Hopf link is for $\r=\F_q[[1=e_1+e_2,x=z_1,y=z_2]]\subset \o=
\F_q[[e_1,e_2,
z_1, z_2]]$. Then standard $M$ can be $\o$ of $q$-rank $2$ or
$(e_1+\om e_2)\r$ for $\om\in \F_q^*$. This gives
 $\mathcal{H}^{mot}_\sigma=
(1+aq)+(q-1)t$. Here $(1+aq)$ is the contribution of the unknot,
the simplest instance  of Proposition \ref{prop:link-rank}.

The simplest colored Hopf $2$-link  is for 
$\sigma=(2,1)$. In this case,
$\mathcal{R}=\mathbb{F}_q[[x=\zeta_1+\zeta_2,y=\zeta_3]]$ and
$\tilde{\Omega}=\mathbb{F}_q[[\zeta_1,\zeta_2,\zeta_3]].$
One has: $\mathcal{H}^{mot}(2,1)=
(1+a q^2 )+(q^2-1)t$. The justification is as follows.

The standard modules
$M$ are $\Om=\ep_1\r+\ep_2\r+\zeta_1\r$ of $q$-rank $3$ and 
``invertible" ones: 
$M=\lan v_1=\ep_1+\om_1 \ep_3, v_2=\ep_2+\om_2 ep_3\ran$ for 
$\om_1,\om_2\in \F_q$ such that at least one of them is nonzero.
Since $x v_1=\ze_1, xv_2=\ze_2, yv_1=\om_1 \ze_3,  
yv_2=\om_2 \ze_3$, all these modules must contain 
$\ze_1,\ze_2,\ze_3$, so the latter can be omitted in the 
expressions for $v_1$ and $v_2$. 

Generalizing, the superpolynomials is
$\mathcal{H}^{mot}(m,1)\!\!=\!
(1\!+\!a q^m)\!+\!(q^m\!-\!1)t$ for 
the Hopf $2$-links with  $\si=(m,1)$.
The calculation is similar.
\vskip 0.2cm

Note that colored algebraic {\em links} were not  
considered in Proposition \ref{prop:link-rank}.
Its extension to $\r$ in this section is
the following conjecture.

\begin{conjecture}
For arbitrary $\r=\F_q[[x,y]]\subset \Om=
\F_q[[\ep_j,\ze_j]]$, where $1\le j\le \tau$ as in the beginning of this
section, let $\tilde{\r}=\F_q[[\tilde{x},\tilde{y}]]\subset \Om$
for $\tilde{x}=
\sum_{j=1}^\tau \pi_j\chi_j, \tilde{y}=y$ for sufficiently general
$\pi_j\in \F_q^*$. Then the standard modules $M\subset \Om$
can be identified with $\tilde{M}\subset \Om$ for $\tilde{\r}$
(with $\tau$ irreducible components)
such that $\o_i \tilde{M}=\Om_i$ for $1\le i \le \kappa$. 
As above, $c_1\ge c_2\ge \cdots\ge c_\kappa$.
Accordingly, $(1+a q^{1})\cdots (1+a q^{c_1-1})\,
\h^{mot}_\si$ for $\si=(c_1,\ldots, c_\kappa)$ is naturally
embedded in $\tilde{\h}$ for $\tilde{\r}$. \sq
\end{conjecture}

\vskip 0.2cm

For instance, $(1+aq)\mathcal{H}^{mot}(2,1)$ 
is the following portion of {\em uncolored}
superpolynomial from (\ref{3-Hopf}) for the Hopf $3$-link, 
$\h^{mot}_{3,3}\!=\!$
{\small 
$(1\!+\!a q)(1\!+\!a q^2)+ (q-1)^2 q t^3 +
(1\!+\!a q)\bigl((q-1)(q+2)t+ (q-1)^2t^2\bigr)$.
}
For $r\!k_q=2$ in $\h^{mot}_{3,3}$, the image of this embedding
is due to family $(2)$ there,
and $2$ modules from $(2')$. There are $3$ (not $2$) 
modules in $(2')$: $M_i=\lan e_1+e_2+e_3, e_i\ran$ for $i=1,2,3$. 
However,
$M_3$ does not satisfies the conditions from the conjecture.
Such modules contribute 
$(q-1)^2t+2(q-1)t=(q^2-1)t$, which is, indeed, part
of $\mathcal{H}^{mot}(2,1)$, that for $r\!k_q=2$.

\subsection{\bf Colored Hopf 3-link}
Finding  $\h^{mot}_\si$ for
Hopf links can be quite challenging. Let us
calculate $\h^{mot}_{(2,1,1)}$ 
 for the Hopf $3$-link with the linking numbers $\{+1,+1,+1\}$ and the
1{\small st} component colored
by $2\om_1=\yng(2)\,$. We take $\r=\F_q[[1,x,y]]$
for $ x=\ze_1+\ze_2+\ze_3, y=\ze_1+\ze_2+\ze_4$, 
satisfying the equation  $F(x,y)=(x-y)xy=0$. 
It is naturally embedded in 
$
\o=\F_q[[\ep_1+\ep_2, \ep_3,\ep_4,
\ze_1+\ze_2,\ze_3,\ze_4]]\subset \Om=\F_q[[\ep_i,\ze_j, 1\le i,j\le 4]].
$
One has:\, $\h^{mot}_{(2,1,1)}=$
\renewcommand{\baselinestretch}{1.2} 
{\small
\( 
1+a^2 q^5-2 t+q^2 t+q^3 t+t^2-2 q^2 t^2+q^4 t^2+q^2 t^3
-q^3 t^3-q^4 t^3+q^5 t^3 +a \bigl(q^2+q^3-2 q^2 t+q^4 t+
q^5 t+q^2 t^2-q^3 t^2-q^4 t^2+q^5 t^2\bigr)\, =\, 
(q-1)^2 (1+q) qt^2 (1+q t)\, +\, 
\bigl(\,(q-1) t (2+2 q+q^2+(q^2-1)t\,\bigr)
(1+aq^2)\, +\,  (1+aq^2)(1+aq^3).
\)
}
\renewcommand{\baselinestretch}{1.2} 

The justification is as follows.
The corresponding families of standard modules are
of $q$-ranks $2,3$ and $4$ (since $c_1=2$). Recall
that $deg_a\h^{mot}=\sum_{i=1}^\kappa c_im_i-c_1$, which is $2+1+1-2=2$,
for the multiplicities $m_i$ of unibranch
singularities $\{F_i(x_i,y_i)=0\}$. 

The ring $\r$ 
contains $xy=\ze_1^2+\ze_2^2$,
$x^2-xy=\ze_3^2, y^2-xy=\ze_4^2$. All 
$\ze_i^n$ for $n\ge 2$ belong to standard modules due to
the definition of the latter. Indeed, standard modules must contain
$v_1=\ep_1+a_3\ep_3+a_4\ep_4 \mod J$, where 
$J=\sum_{i=1}^4 \ze_i\Om$, 
and  $v_2=\ep_2+b_3 \ep_3+b_4\ep_4 \mod J$.
Here either $a_3$
or $b_3$ must be nonzero, and either $a_4$ or $b_4$ must be
nonzero; the projections of
$M$ onto $\o_3$ and $\o_4$ must be standard. Using the multiplication
by $xy, x^2,y^2$, we obtain that $\ze_i^2\in M$ for any $i$.
\vskip 0.2cm

{\sf Full description.}
Let us describe these families, the $q$-ranks,
the degrees $deg(M)$, and the cardinalities\, $\#$.  
\vskip 0.2cm

(0) The module $M=\Om$ is the only standard one of $r\!k_q=4$;
the elements $\ep_1,\ep_2,\ep_3, \ze_4$ are its $\r$-generators.
Its degree is $0$.

(1a) Here $r\!k_q=2$, $deg=2$; so 
$M=\{v_1,v_2\}$ for $v_i$ as above. Since
$xv_1-y v_1=a_3\ze_3-a_4\ze_4\in M$ and
$xv_2-y v_2=b_3\ze_3-b_4\ze_4\in M$, we obtain 
that $z_i\in M$
for all $i$ if
 $\det\begin{pmatrix}
a_3 &a_4 \\ b_3 &b_4\end{pmatrix}\neq 0$, which we impose. This
gives that ``nonzero conditions" for $a_i,b_i$ are satisfied, 
and  that the elements from $J$ can be disregarded
in the expansions of $v_i$. Thus,
there are  $\#=(q^2-1)(q^2-q)=q(q-1)(q^2-1)$ such modules.

(1b) Continuing, let us now assume that
$\det=0$. Then the corresponding $M$
are $\lan \ep_1+u(a_3\ep_3+a_4\ep_4) +b_3\ze_3,
 \ep_2+w(a_3\ep_3+a_4\ep_4)+c_3\ze_3\ran$. Here 
$a_3a_4\neq 0$ because
$\ep_3$ and $\ep_4$ must be present;
thus, $a_3,a_4\in \F_q^*$ and one can take $\{u,w\}$ to be either
$\{u\in \F_q,1\}$ or $\{0,1\}$. Since 
$b_3,c_3\in \F_q$ are arbitrary, 
$\#=q^3(q-1)^2(q+1)=q^2(q-1)(q^2-1)$,
$r\!k_q=2$, and $deg=3$ (due to $\det=0$ in this family).
\vskip 0.2cm

(2) Adding all $z_i$ to the modules
from (1b), we obtain the family with $r\!k_q=3,
deg=2$ and of cardinality  $\#=(q^2-1)(q-1)$.

(3) The remaining families are for $M_1=\lan \ep_1+a_4\ep_4, \ep_2+b_4\ep_4,
\ep_3+c_4 \ep_4$, where $c_4\neq 0$ and $\#_1=q^2(q-1)$,
$M_2= \lan \ep_1+a_4\ep_4, \ep_2+b_4\ep_4,
\ep_3\ran$, where the case $a_4=0=b_4$ is excluded,
with  $\#_2=(q^2-1)$,
and $M_3=\lan \ep_1+a_4\ep_4, \ep_2+b_4\ep_4, \ep_4\ran$ 
of (the same) size $\#_3=(q^2-1)$. Their total size
is $\#=q^3+q^2-2=(q-1)(2+2q+q^2)$. Here $r\!k_q=3$ and $deg=1$.
\vskip 0.2cm

Combining (0,1,2,3), we obtain:\, $\h^{mot}_{2,1,1}\!=\!$
\renewcommand{\baselinestretch}{1.2} 
{\small
\( 
(1\!+\!aq^2)(1\!+\!aq^3)+q(q-1)(q^2-1)t^2+
q^2(q-1)(q^2-1)t^3+(q^2-1)(q-1)t^2(1+aq^2)+(q-1)(2+2q+q^2))t(1\!+\!aq^2)
\, =\, (q-1)^2 (1+q) qt^2 (1+q t)\, +\, 
(q-1)t\,\bigl(\,(q^2-1)t+2+2q+q^2\,\bigr)(1\!+\!aq^2)+
(1\!+\!aq^2)(1\!+\!aq^3),
\)
}
\renewcommand{\baselinestretch}{1.2} 
as we claimed above.
\vskip 0.2cm

More uniformly, we need to describe 
linear subspaces $V$ of codimension $deg=1,2$ 
in $\Om/\mathfrak{m}_{\o}\Om=
\F_q^2\oplus \F_q\oplus \F_q$ with surjective projections
onto these $3$  components. The families above are
when  we count  the number of
possibilities for 
$\mathfrak{m}_\r V$ modulo $\{\ze_i^2\}$ for such $V$.

\subsection{\bf Multibranch {\bf\em L}-functions} \label{sec:multL}
Let us define $L(q,t,a)$ for 
{\em uncolored, possibly non-unibranch}
plane curve singularities. Concerning $c>1$, 
see \cite{ChS} for some examples.
For $a=0$, this $L$ is from \cite{Sto};
see the references there and the proof of
the functional equation 
(for any Gorenstein $\r$).
Recall that $\r\!=\!\F_q[[x,y]]/\bigl(F(x,y)\bigr)$, where
$F(x,y)\!=\!\prod_{i=1}^\kappa F_i(x,y)$ and 
$F_i(x,y)$ are 
pairwise non-proportional irreducible polynomials 
over $\F_q$. We assume them
absolutely irreducible, 
which is the case $r=1$ in the notation from \cite{Sto}.
Namely, $\r_i/\mathfrak{m}_i=\F_q$ for the maximal ideals
$\mathfrak{m}_i\subset \r_i$, where 
$\r_i\simeq\F_q[[x,y]]/\bigl(F_i(x,y)\bigr)$.

As above, $\mathcal{H}^{mot}$ and
 $\mathbf{H}^{mot}(q,t,a)\equal 
\mathcal{H}^{mot}(qt, t, a)$
are defined for 
the ring $\mathcal{R}=\mathbb{F}_q[[x,y]]$ embedded in $\o\!=\!
\mathbb{F}_q[[z,e_1,\ldots,e_\kappa]]$ for the idempotents $\{e_i\}$.
Let
$z_i=z e_i (1\le i\le \kappa)$, and 
$\mathcal{K}^\times=
\oplus_{i=1}^\kappa \mathbb{F}_q((z_i))^\times$, the corresponding
group of
(invertible) rational functions. We set $M'\sim M$ for $\r$-modules
$M,M'\in \o$ if $M'=fM$ for
$f\in \mathcal{K}^\times$. 

The {\em $a$-extended $L$-function} is then  
$L(q,t,a)\equal (1\!-\!t)^\kappa Z(q,t,a)$ for  
\begin{align}\label{Z-links}
&Z(q,t,a)\!=\!\sum_{M} 
t^{deg_\r(M)}\!\!\prod_{j=1}^{r\!k_q(M)-1}
\!(1+aq^j), 
\end{align}
where $M\subset\r$ are ideals 
such that $deg_\r(M)\equal$\,dim${}_{\mathbb{F}_q}(\r/M)<\infty$.
For any $a$, this definition (as such) 
is actually useful only for plane curve
singularities, which will be assumed in the next theorem.

\begin{theorem}\label{thm-feq}
The functional equation holds for  $L(q,t,a)$,
which is the invariance of $t^{-\de}L(q,t,a)$ with respect to
the Hasse-Weil transformation $t\mapsto 1/(qt)$, where 
$a\mapsto a$. Equivalently,
$q^{\de}t^{2\de}L(q,\frac{1}{qt},a)=L(q,t,a)$. This is for $c=1$ and
plane curve singularities, possibly multibranch.
\end{theorem}
{\it Proof.} The cases  $a=0$ and $a=-1/q$ are due to Section (3.4)
from \cite{Sto}, which is for any Gorenstein rings $\r$. Given standard
$M\subset \o$,
the contribution of the
set $\{M'\subset \r \mid M'\sim M \}$ 
to $L$ was explicitly identified there
with that for $M^\ast=\r: M$. We need only to check that
$r\!k_q(M)=r\!k_q(M^\ast)$, which does require 
plane curve singularities; see the example
of $\mathbb{F}_q[[z^4,z^5,z^6]]$ below. The 
{\em Auslander-Buchsbaum theorem} is used as in Lemma 10 in \cite{ObS}. 

Another approach is the $a$-extension 
of the motivic proof of the functional equation 
for Kapranov's motivic zeta-function due to Kapranov
and M.Artin; see \cite{Kap}. This follows Lemma 9 in \cite{ORS}.\sq

\begin{corollary}
Let $\r=\F[[z^{\rr},z^{\ss}]]$ for  gcd$(\rr,\ss)=1$.
Recall that the  $\Ga$-rank $r\!k_\Ga(\De)$ of a 
standard $Ga$-module $\De$
is the (minimal) number
of generators. Setting $\De^\bullet=\De^\ast-\min\{\De^\ast\}$ for
$\De^\ast$ from (\ref{Dedual}),
 $r\!k_\Ga(\De)=r\!k_\Ga(\De^\bullet)$. Also, if $r\!k_q(M)$
is not constant in the corresponding Piontkowski cell $\j_0(\De)\ni M$,
then the same holds for $\De^\bullet$. 
\end{corollary}
{\it Proof.} Any  $\j_0(\De)$ for standard $\De$ contains a monomial
module $M_\circ=\sum_{\de\in \De} \F z^{\de}$. Its $q$-rank $r\!k_q(M)$
coincides with $r\!k_\Ga(\De)$. Then the module 
$M^\bullet_\circ=z^{-m} M_\circ^\ast$ 
is monomial and standard for $m=\min\{\De^\ast\}$ in $\j_0(\De^\bullet)$
and its $q$-rank coincides with
$r\!k_q(M)$ due to the end of the proof of the theorem. 
This gives that $r\!k_\Ga(\De^\bullet)=r\!k_q(M^\bullet_\circ)=
r\!k_q(M_\circ)=r\!k_\Ga(\De)$. 

There is a pure combinatorial
proof of this coincidence based on  Theorem 2.12 from \cite{GM}.
The identification of standard $\De$ with Dyck paths above the 
diagonal in the
rectangle ``$\rr\times \ss$" from this paper is used; 
see also Section 7.5.3 of the book
\cite{ChB}. The passage to $\De^\bullet$ corresponds to the switch
$\rr\leftrightarrow \ss$ and the transposition the
Young diagram representing the corresponding 
path. The $\Ga$-rank is the number of the 
outer corners of the latter, which is preserved.

If there is $M\in \j_0(\De)$ with $r\!k_q(M)<r\!k_\Ga(\De)$
($r\!k_q$  can be only smaller), then 
$r\!k_q(M^\bullet)<r\!k_\Ga(\De^\bullet)$,
which gives the $2${\small nd} claim.

To give an example, 
the only such cell in the case of $T(7,3)$ from \ref{Table7-3}
is with $D=[4,8,11]$ and $\De=\Ga\cup D$, where 
$\Ga=\{0,3,6,7,9,10,12,13,\ldots\}$ and $G=\{1,2,3,5,8,11\}$.
Accordingly, $\De^\ast=\{6,9,10,12,13,\ldots\}$ and
$\De^\bullet=\De^\ast-6=\{0,3,4,6,7,\ldots\}=\De$, as
it has to be. \sq
\vskip 0.2cm

Let us state the conjecture $\,\H\!=\!L$ 
for any $\kappa$. See Conjecture 4.5 in \cite{ChW} (unibranch),
Section 5 in \cite{ChS} and, also, formula (\ref{conjKG}) below.   

\begin{conjecture}\label{conj:HLlinks}
For uncolored algebraic links and the corresponding
rings $\r$ of plane curve singularities over $\F_q$, 
$\mathbf{H}(q,t,a)=L(q,t,a)$, where we put $\h$ and $\H(q,t,a)$
omitting {\em\small mot} or {\em\small daha} due to
the conjectural coincidence 
$\h^{daha}=\h^{mot}$. Also, this coincidence
``$\mathbf{H}=L$\," is expected to hold for Gorenstein $\r$,
possibly sufficiently general. 
The Hasse-Weil transformation $t\mapsto 1/(qt)$ for $L$ and $\H$ 
corresponds 
to the transformation  $q\leftrightarrow t^{-1}$ of $\h(q,t,a)$
for ``\,{\em\small mot}" or ``\,{\em\small daha}". \sq
\end{conjecture}

Recall that the conjectural relation $\h^{mot}_\si=\h^{daha}_\si$ 
is for any colors $\si=(c_1,\cdots, c_\kappa)$. The corresponding
DAHA superduality requires the transposition of Young diagrams,
$m\om_1\mapsto \om_m$ in our setting. However, $\h^{mot}$ 
and the corresponding $L$-functions
are unknown for arbitrary Young diagrams by now. Moreover,
there is no general 
definition of $L$ for 
algebraic {\em  links} colored by arbitrary 
$\si$ at the moment, the one satisfying (conjecturally) 
``$\H=L$\,". 
\vskip 0.2cm

The {\em Weak Riemann Hypothesis} for $a=0$ from \cite{ChW} 
for algebraic uncolored {\em links}
is that $\H(q,t,a=0)$ has only 
$\kappa-1$ different pairs of conjugated $t$-zeros 
that do not satisfy $RH$, i.e. such that 
$|t|\neq q^{-1/2}$. Presumably, the proof of
Theorem  \ref{thm:q-rank} can be adjusted
to verify this claim. Let us mention here, that
 we conjectured that
the monomials in the
series $\h^{daha}(q,t,a)/(1-t)^{\kappa-1}$ for algebraic links
and rectangle Young
diagrams have only positive
coefficients, which conjecture remains open even in the case of
uncolored algebraic knots. 
\vskip 0.2cm

{\sf Quote schemes.}
Let us briefly
discuss the {\em Hilb-vs-Quot Conjecture} from \cite{KT};
see there for details and references. Their Conjecture 1 is our
``$\,\H=L$" for $a=0$ (from \cite{ChW}) in the setting with virtual 
weight polynomials instead
of counting $\mathbb{F}_q$-points (for a different choice of 
motivic measure). 
The polynomials $Q_S(q,t)$ there are essentially  
our uncolored motivic superpolynomials; their $S$ is our $\o$. 
Conjecture 2 there is, indeed, beyond \cite{ChW} and  
the current paper.

The passage to arbitrary $a$ requires the schemes 
$\mathcal{Q}uot_{r-len}(\r)$ and 
$\mathcal{Q}uot_{r-len}(\o)$, where 
$r\!k_q(M)=r$ is imposed. These schemes
are not discussed much in \cite{KT}, especially the latter.
They and their
virtual weight polynomials can be complicated even
for relatively simple $\r$.

Conjecture 4 
is a special case of our ``$\,\h^{mot}=$
{\em KhR-reduced}\,", which is 
for algebraic links colored by rows. We agree 
that the passage from $\mathcal{Q}uote(\r)$ from the original 
{\em ORS Conjecture}  to 
$\mathcal{Q}uote(\o)$, our
superpolynomials, has many advantages.
Recall that the strongest conjecture here is that the DAHA 
superpolynomials
$\h^{daha}$ 
coincide with reduced (colored) $KhR$-polynomials for iterated
torus links \cite{ChD2}. 

Our multibranch definition of $\h^{mot}(q,t,a)$ in this paper 
and in \cite{ChS} 
is actually quite similar to that from the original 
\cite{ChP1,ChP2} in the
unibranch case. However, it took time to incorporate $a$ properly 
when $F(x,y)$ is allowed to be reducible and not square-free.     
Topologically, this is the case of
arbitrary algebraic links colored by ``rows". 

The conjecture $\,\h^{daha}=\h^{mot}\,$, the main in this paper,
is based
on many examples and considerations, including the colored case.
Numerical confirmations of our (uncolored) 
``$\h^{mot}(qt,t,a)=L(q,t,a)\,$" are limited; practically, because 
$L$ are more involved than $\h^{mot}$.
 
The definition of  $L(q,t,a)$ for arbitrary $F(x,y)$, possibly not
square-free, 
remains a challenge by now. We note that
the Galkin-St\"ohr $L$-functions 
(our uncolored $L$ for $a=0$) are defined for any number 
of branches. Moreover, 
the residue fields can be different extensions
of $\F_q$ for different
branches in \cite{Sto}, which seems ``non-topological". 
\vfil

{\sf Parahoric affine Springer fibers.}
We refer the readers to
Chapter 7 of \cite{KT}; see also Section \ref{sec:asf-nil} above.
The definition is basically as follows.

Let $\ga\in \mathfrak{g}(\mathbb{F}(\!(x)\!))$
for a reductive finite-dimensional Lie algebra $\mathfrak{g}$, and
$\x_{\ga}^\circ=\{g\in G(\mathbb{F}((x)))/G(\mathbb{F}[[x)]])\mid 
Ad_g^{-1}(\ga)\in \mathfrak{g}(\mathbb{F}[[x]])\}$. Here
Lie$(G)=\mathfrak{g}$, $G$ is simply-connected,
 and $\mathbb{F}=\mathbb{F}_q$.

For $GL_n$,  the
scheme $\x_\ga^\circ$ can be identified with that 
formed by $\r$-invariant
submodules $M\subset \mathcal{K}=\prod_{i=1}^\kappa
\mathbb{F}((z_i))$, {\em lattices}, such that
$M \mathcal{K}=\mathcal{K}$, where 
$\r=\mathbb{F}[[x,y]]/(F_\ga(x,y))$.
$F_\ga(x,y)= $det$(\ga-y\mathbf{1})$. We assume that $F_\ga(x,y)$ 
is square-free  with $\kappa$ (absolutely) irreducible factors.
Mostly, only such $ASF$ are considered in the literature.

The corresponding {\em affine Springer fiber} is 
$\x_\ga=\x_{\ga}^\circ/\Lambda$ for the discrete 
group $\La$ due
to Kazhdan and Lusztig \cite{KL}. For $GL_n$,
this group acts via
$M\mapsto \prod_{i=1}^\kappa z_i^{n_i}M$ for $n_i\in \Z$
in the {\em lattice presentation} above.  
The division by $\La=\Z^{\kappa}$ reduces the consideration
to $M\subset \o$ of degree $0$ and results in our $\j_0$. 
Thus, the latter variety is a {\em fundamental domain} for $\La$; 
see e.g., \cite{KT}. 
Generally, this  $\La$  comes
from {\em principle adeles} for a rational planar complete curve
such that it has a unique singularity at $x=0=y$ 
with $\r$ as its local ring. 
 
Thus, we do extended $p$-adic orbital integrals
in this paper. Adding $t$ is, basically,
the passage to {\em Igusa integrals}; see \cite{Sto,Yun}.
Adding $a$ is the 
consideration of $q$-ranks. Adding colors $c$ is an 
extension to not square-free $F(x,y)$. Recall that
this is done directly in terms of  $\r$. 

\comment{
Let $\g$ be the group scheme over 
a smooth projective curve $E$ with the
Lie algebra $\mathfrak{G}$  such that 
$H^0(E,\g)$ and $H^1(E,\g)$ are trivial. The starting point 
is a subscheme $\t\subset \g$, which is assumed a maximal
subtorus at the generic point of $E$. 
Since $H^1(E,\g)=\{0\}$, any cocycle 
$\phi$ in the {\em generalized Jacobian}, which is
 $H^1(E,\t)$, becomes
the {\em boundary} $\{\phi_i\phi_j^{-1}\}$ for an open cover
$E=\cup_i U_i$ and $\phi_i\in H^0(U_i,\g)$. Then  
$\t_\phi=\phi_i^{-1}\t \phi_i\subset \g$ is another
toric subscheme with 
the characteristic polynomial coinciding with that of  $\t$. 
The $ASF$ is a compactification of the variety of such
subschemes. 

Generally, any $\t$ can be represented
as $\t=\mathbb G_m (C)$ for a projective curve 
$C$ covering $E$, {\sf\em possibly singular}. 
Thus $ASF$ is a natural compactification
the generalized Jacobian  
of $C$. The {\em Jacobian factors} will be the contributions 
of singular points of $C$ to $\overline{Jac}(C)$.

Now, we take $E=\mathbb P^1$ and consider a planar 
rational curve  $C$ with only one singularity, which is
with our $r$ as its local ring.
Then we arrive at our $\j_0$. 
}

\subsection{\bf Example: double trefoil}
The calculations of $Z(q,t,a)$ and $L(q,t,a)$ can be complicated even
when $\h(q,t,a)$ are relatively straightforward.
Let us discuss the coincidence of
{\em uncolored} $\H(q,t,a)$ and $L(q,t,a)$
for $T(6,4)$. Using (\ref{T6-4}), 
 $\H_{6,4}=$
\renewcommand{\baselinestretch}{1.2} 
{\small
\(
1-2 t+t^2+q^3 t^4-2 q^3 t^5+q^3 t^6+q^4 t^6-2 q^4 t^7+q^4 t^8
+q^5 t^8-2 q^5 t^9+q^5 t^{10}+q^6 t^{10}-2 q^6 t^{11}+q^7 t^{12}
+q^6 t^{14}-2 q^7 t^{15}+q^8 t^{16}+ 
\bigl(
t-t^2+q t^2-t^3-q t^3+q^2 t^3+t^4-q t^4-2 q^2 t^5+2 q^3 t^5+2 q t^6
-q^3 t^6+q^4 t^6-q t^7-3 q^3 t^7+q^4 t^7+2 q^2 t^8+q^3 t^8-q^4 t^8
+q^5 t^8-q^2 t^9-3 q^4 t^9+q^5 t^9+2 q^3 t^{10}-q^5 t^{10}+
q^6 t^{10}-2 q^5 t^{11}+2 q^6 t^{11}+q^4 t^{12}-q^5 t^{12}-q^5 t^{13}
-q^6 t^{13}+q^7 t^{13}-q^6 t^{14}+q^7 t^{14}+q^7 t^{15}
\bigr) 
(1+a q)+
\bigl(
t^3-t^4+q t^4+q^2 t^5-t^6-2 q t^6+t^7+q t^7+2 q^3 t^7-q t^8
-2 q^2 t^8-q^3 t^8+q t^9+q^2 t^9+2 q^4 t^9
-q^2 t^{10}-2 q^3 t^{10}+q^5 t^{11}-q^4 t^{12}+q^5 t^{12}+q^5 t^{13}
\bigr) 
(1+a q) (1+a q^2)+t^6 
\bigl(
t^6-t^7+q t^8-q t^9+q^2 t^{10}
\bigr) 
(1+a q) (1+a q^2) (1+a q^3).
\)
}
\renewcommand{\baselinestretch}{1.2} 

It coincides with $L_{6,4}(q,t,a)$, which is an
involved calculation (even in such a basic case) 
beyond $a=0$ and $a=-1/q$. 
For instance, $\H_{6,4}(q,t,-1/q)$
coincides with the {\em Z\'u{\oldt{n}}iga $L$-function} due to
Theorem 4.3 from \cite{Sto} for the latter.
Let us demonstrate  that the contribution  of 
the ideals $M\subset \r$ of maximal $r\!k_q=4$ to $Z(q,t,a)$
coincides with $t^6 \bigl(
t^6-t^7+q t^8-q t^9+q^2 t^{10}
\bigr) 
(1+a q) (1+a q^2) (1+a q^3)$ from 
 $\H_{6,4}$.

Recall that $\mathcal{R}=\mathbb{F}_q[[x=z_1^2-z_2^2,y=z_1^3+z_2^3]]
\subset \mathcal{O}=\mathbb{F}_q[[e_1,e_2,z_1,z_2]]$; all our
rings contain $1=\sum_i e_i$. Its 
{\em conductor} $\mathfrak{c}=z^8\,\mathcal{O}$ 
has $4$ generators, $z_{1,2}^8$ and $z_{1,2}^9$,
over $\mathcal{R}$ and $deg_\r(\mathfrak{c})=$
\,dim$_{\mathbb{F}_q}\mathcal{R}/\mathfrak{c}=8$.

The images of 
$\{e_1,z_1,z_2,z_1^2,z_1^3,z_1^4,z_1^5,z_1^7\}$
form a natural basis in $\mathcal{O}/\mathcal{R}$.
Note that 
$z_{1,2}^6\in \mathcal{R}$ due to different
signs in $x$ and $y$; so are all $z_{1,2}^m$ for $m\ge 8$. 
The natural $\F_q$-generators of $\r$ are
{\small
$1,\, z_1^2-z_2^2,\, z_1^3+z_2^3,\,
z_1^2+z_2^4,\, z_1^5-z_2^5,\,
z_i^6 (i=1,2),\, z_1^7+z_2^7;\, \mathfrak{c}=\{z_i^8,z_i^9,
\ldots$\}} for $i=1,2$.
\vfil

Let us begin with  {\em exceptional} ideals
$M\subset \r$ of $r\!k_q=4$\,:\, those not entirely in $\mathfrak{c}$.
They are
$M_1=\lan z_1^6-z_2^6, z_1^7+z_2^7, z_1^8,z_2^8\ran$,\, 
$M_2=\lan z_1^7+z_2^7, z_1^8-z_2^7, z_1^9,z_2^9\ran$,\, and 
$M_3=\lan z_1^7+z_2^7+\al z_1^8, z_1^8-z_2^8, z_2^9,z_2^{10}\ran$
for any $\al\in \F_q$.  The degrees $deg_\r$ are $6,7,8$.
Accordingly, their contribution to $Z$ will be
$\bigl(t^4+t^7+q t^8\bigr)(1+q a)(1+q^2 a)(1+q^3 a)$.
Potentially, one can take here ideals 
$\lan z_1^m\pm z_2^m, z_1^{m+1}\mp z_2^{m+1}, z_1^k,z_2^k\ran$ 
for $m=2,4,6$ and $k\ge 8$ but 
$r\!k_q<4$ for such ideals.

The {\em non-exceptional} ideals $M$ are those obtained by
the multiplication by $z_1^m+z_2^n$ for $m,n\ge 8$ of
{\em standard} $M_{st}$ of $r\!k_q=4$ described  
in Section \ref{sec:dtref}. Their contribution to $\h^{mot}_{6,4}(a=0)$
was calculated in (\ref{4st6-4}): 
$\h_{6,4}^{r\!k\!=\!4}(a=0)\!=\!\text{\small
$\bigl(1+(q\!-\!1)t+q(q\!-\!1)t^2\bigr)$}$. Accordingly, the 
contribution of the ideals corresponding to {\em non-exceptional}  
$M$ to $Z(q,t,a)$ will be $(1+a q) (1+a q^2) (1+a q^3)$ times
$\bigl(1+(q\!-\!1)t+q(q\!-\!1)t^2\bigr)\frac{t^8}{(1-t)^2}$; the
latter factor is due to multiplications by all
$z_1^m+z_2^n$.

Accordingly, the total contribution of ideals of $r\!k_q=4$ to
$Z(q,t,a)=\frac{L(q,t,a)}{(1-t)^2}$ is 
$\Bigl(\bigr(t^6+t^7+q t^8\bigl)+
\frac{t^8\bigl(1+(q\!-\!1)t+q(q\!-\!1)t^2\bigr)}{(1-t)^2}\Bigr)
(1+aq)(1+aq^2 )(1+aq^3)$, which, indeed, results in
$\bigl(
t^6-t^7+q t^8-q t^9+q^2 t^{10}\bigr)$ from $\H_{6,2}$.

The usage of some portions of
$\h^{mot}$ in the calculation of $Z(q,t,a)$ can be expected.
What is surprising  is that the Coincidence Conjecture
$\h(qt,t,a)\!=\!L(q,t,a)$ requires the substitution $q\mapsto qt$.

\vskip 0.2cm

The calculations with $Z,L$ become simpler when $a=0$.
Some such formulas are provided in \cite{Sto}. 
In the notation there:\, $L(\o)$
is our $L(q,t,a=0)$ and $L(\o,\o)$ is 
$L_{prncpl}$, which is $L(q,t,a=-1/q)$. 
Examples 6.3 and 6.2 are for the Hopf $2$-link and $3$-link.
They are, indeed,
our $\H(q,t,a)$ for $a=0$ and $a=-1/q$.   Theorems 4.1 and 4.3 
in \cite{Sto} give explicit formulas for  
$L_{prncpl}$ for $\kappa=1,2$ in terms of $\Ga,\vec{\Ga}$.
Using the latter theorem:
$\H_{6,4}(a\!=\!-1/q)\!=\!L_{6,4}(a\!=\!-1/q)\!=\!L_{6,4}(\o,\o)$.

Theorem 4.1 is basically our formula (\ref{de-qt}).
Using the latter for $\r'=\F_q[[z^4,z^6+z^7]]$:
$\h'(q,t,a\!=\!-t/q)=
1 - t + q^3 t - q^3 t^2 + q^4 t^2 - q^4 t^3 + q^5 t^3 - q^5 t^4 + 
 q^6 t^4 - q^6 t^5 + q^7 t^5 - q^7 t^8 + q^8 t^8=
L'(q/t,t,a=-t/q)=L'_{prncpl}(q/t,t)$,
where $\h'$ and $L'$ are the corresponding superpolynomial and
$L$-function. We use in this paragraph $q$
from $\h'$ (not from $L'$). 

\comment{
We note that such formulas can be reduced to the unibranch 
case using (\ref{h-dtor1}):
$
\h_{\rr,\ss}=(1\!+\!a q)\h_{\rr,\ss}(2)+
\frac{q-1}{q}\Bigl(\h_{\rr,\ss;2}-(1\!+\!a q)\h_{\rr,\ss}(2)
\Bigr).
$
It holds for $L(q,t,a)$ upon the
substitution $q\mapsto qt$; the definition of
colored unibranch $L(q,t,a)$ is parallel to that in
the uncolored case.
Letting $a\!=\!-t/q$, one has:\,
$\h_{3,2;2}(a\!=\!-t/q)=
1 - t + q^3 t - q^3 t^2 + q^4 t^2 - q^4 t^3 + q^5 t^3 - q^5 t^4 + 
 q^6 t^4 - q^6 t^5 + q^7 t^5 - q^7 t^8 + q^8 t^8$ 
and $\h_{3,2}(2)(a\!=\!-t/q)=1 - q t + q^3 t$.
The formula for $L_{3,2;2}(a\!=\!-1/q)=L_{prncpl}$ for 
$\r_{3,2;2}=\F_q[[z^4,z^6+z^7]]$ 
can be readily
extracted from (\ref{de-qt}) (or directly from \cite{Sto}); 
it coincides with 
$\H^{mot}_{3,2;2}(a\!=\!-1/q)$. 
}


\comment{
Let $\sigma=(2,1)$ for the Hopf $2$-link. Namely, 
$\mathcal{R}=\mathbb{F}_q[[x=\zeta_1+\zeta_2,y=\zeta_3]]$ and
$M\subset
\tilde{\Omega}=\mathbb{F}_q[[\zeta_1,\zeta_2,\zeta_3]]\cup 
\mathcal{R}(\epsilon_2+\epsilon_3).$
It is easy to calculate that $\mathcal{H}^{mot}_\sigma=
(1+aq^2)+(q^2-1)t$; the corresponding $L$ is as follows.

The modules $M$
of $r\!k_q(M)=2$ are those generated over $\mathcal{R}$
by $\{1,\zeta_1^u, u\ge 1\}, 
\{\epsilon_2+\epsilon_3,\ze_1^u,u\ge 1\}$, 
$\{1,\epsilon_2+\epsilon_3\}$, and 
$\{\zeta_1^u+\alpha\zeta_3^w,
\zeta_2^v+\beta\zeta_3^w\}$,
where  $u,v,w\ge 1$,\, $\alpha,\beta\in \mathbb{F}_q$ 
and the pair $\alpha=0=\beta$ is excluded. 

Their contribution to $Z_\sigma(q,t,a)$ will be
$\frac{(q^2-1)t^3}{(1-t)^3}+
\frac{1+t}{(1-t)}$. The submodules of $q$-rank $3$ are
those generated by 
$\{\zeta_1^u, \zeta_2^v,\zeta_3^w\}$,  which contribute
$(1+a q^2)\frac{t^2}{(1-t)^3}$. Combining,  $L_\sigma(q,t,a)=
(q^2-1)t^3+(1+t)(1-t)^2+t^2(1+aq^2)=q^2t^3+1-t+t^2 a q^2=(1+aq^2t^2)
+(q^2t^2-1)t.$ 
\vskip 0.2cm
}

\subsection{\bf Non-planar examples}
Generally, our main results and conjectures hold
only for plane curve singularities, with a reservation
about $\H(q, t, a=-1/q)=L(q,t,a=-1/q)$, which is
expected to hold for (certain?) 
Gorenstein rings $\r$. See Conjecture 4.5
and Section 4.2.3 in \cite{ChW}.  There are
also examples below of $L(q/t,t,a=-t/q)=\h^{daha}_K(q,t,a=-t/q)$,
for {\em non-algebraic} iterated torus links associated with
{\em non-planar} singularities. Recall that $L(q,t,a=-1/q)$
is the  {\em Z\'u{\oldt{n}}iga 
zeta function}: the summation is only over {\em principle} 
ideals $f\, \r\subset\r$ for $f\in \r$. Topologically,
the correspondence $\r\rightsquigarrow K$
 goes through the
lens spaces and their generalizations
(only for some families!). In this case below,  it goes via
$L(2,1)$ corresponding to $uv=w^2$ for $u=z^4, v=z^6, w=z^5$.
This is preliminary, and $a=-t/q$ is required 
for  non-algebraic cables. 

Let $\r=\F[[z^4,z^5,z^6]]\subset \o=\F[[z]]$, which is a
Gorenstein ring with $\de=4, \cc=z^8\o$. Then $L(q,t,a=-1/q)$ coincides with
(uncolored) $\h^{daha}_K(q, t, a=-t/q)$ upon the  substitution $q\mapsto q t$
for the cable $K=C\!ab(5,2)T(3,2)$. Recall that the smallest algebraic
knot in this family is $C\!ab(13,2)T(3,2)$. The general formula is due
to \cite{ChD1} and Proposition 4.2 in \cite{ChW}:\, $\h_K=$
\renewcommand{\baselinestretch}{0.5} 
{\small
\(
-\frac{a^3 q^5}{t}+a^2 (-q^5+q^4 t-\frac{q^4}{t}+q^3)+
a (q^4 t^3+q^4 t^2-q^4+q^3 t^2+2 q^3 t+q^2 t+q^2+q)+
q^4 t^4+q^3 t^3+q^3 t^2+q^2 t^2+q^2 t+q t+1.
\)
}
\renewcommand{\baselinestretch}{1.2} 

One has: $\h^{daha}_K(q,t,a=-t/q)=
1-t+q^3 t-q^3 t^4+q^4 t^4$, which coincides with
$L(q/t,t,a=-t/q)$ and with
$\h^{mot}(q,t,a=-t/q)$. 

The next cable in this family is $K'=C\!ab(7,2)T(3,2)$, which
corresponds to $\r'=\F_q[[z^4,z^6,z^7]]$. One has:
$\h^{daha}_{K'}(q,t)=1-t+q^3 t-q^3 t^2+q^4 t^2-q^4 t^5+q^5 t^5=
L(q/t,t)=\h^{mot}(q,t)$ for $a=-t/q$,
\vskip 0.2cm

Piontkowski's cells are affine spaces $\mathbb{A}^{dim}$
for quasi-homogeneous singularities. 
The corresponding $r\!k_q$ are constant in the cases above,
which holds only for relatively simple $\r$.
 Table \ref{Table456} contains
the list of standard $D=\De\setminus \Ga$ for $\r$, where
$\Ga=\{0,4,5,6,8,9,\cdots\}$, and the corresponding
$q^{dim}t^{deg}$, where the generators are marked ($0$ is
not shown). Recall that
deg$(M)=\dim_\F (\o/M)$ in superpolynomials.

{\footnotesize
\begin{table*}[ht!]
\[
\centering
\begin{tabular}{|l|l|}
\hline    
$D$-sets, $r\!k_q=1$ & $q,t$-terms  \\
\hline 
$\emptyset$ & $q^4t^4$ \\       
\hline                   
$D$-sets, $r\!k_q=2$ & $q,t$-terms \\
\hline 
$\widetilde{7}$   & $q^3t^3$ \\  
$\widetilde{1},7$ & $q^4 t^2$\\ 
$\widetilde{2},7$ & $q^3 t^2$\\
$\widetilde{3},7$ & $q^2 t^2$ \\

\hline
\end{tabular}
\hspace{0.1cm}
\begin{tabular}{|l|l|}
\hline    
$D$-sets, $r\!k_q=4$ & $q,t$-terms  \\
\hline 
$\widetilde{1},\widetilde{2},\widetilde{3},7$ & $1$ \\
\hline 
$D$-sets, $r\!k_q=3$ & $q,t$-terms  \\
\hline  
& \\ 
$\widetilde{1},\widetilde{2},7$ & $q^3 t$\\
$\widetilde{1},\widetilde{3},7$ & $q^2 t$ \\
$\widetilde{2},\widetilde{3},7$ & $q t$ \\
\hline   
\end{tabular}
\]
\caption{Piontkowski cells for $\F[[z^4,z^5,z^6]]$}
\label{Table456}
\end{table*}
}
By analogy with
the planar case, the contributions $q^{dim}t^{deg}$ of 
modules of $q$-rank $r$ will be multiplied by
$(1+a)\cdots (1+q^{r-1}a)$; though this definition 
seems of limited importance for 
non-planar singularities unless for some special $a$.  
Then $\h^{mot}=$
\renewcommand{\baselinestretch}{0.5} 
{\small
\(
q^4 t^4+(1+a q) (q^2 t^2+q^3 t^2+q^4 t^2+q^3 t^3)+
(1+a q) (1+a q^2) (q t+q^2 t+q^3 t)
+(1+a q) (1+a q^2) (1+a q^3).
\)
}
\renewcommand{\baselinestretch}{1.2}

In this example and generally,
 $\h^{mot}$ does not satisfy
the superduality (unless for plane curve singularities).
Namely,  
the sums of contributions from
\ref{Table456} for fixed $r\!k_q$ 
 are far from being invariant with respect to $q\leftrightarrow t^{-1}$.
The same is applicable to Table \ref{TableL456}. The reason 
for this is that the $q$-ranks
can be different for $M$ and $M^\ast$. Let us provide 
an example. 

The $q$-rank is $3$ for the monomial
$M_\circ=\lan 1,z,z^2,z^4,z^5,\ldots\ran$ with $\De=\{0,1,2,4,5,\cdots\}$
and $D=[1,2,7]$. However, $r\!k_q(M_\circ^\ast)=2$ for $M_\circ^\ast=
\lan z^4,z^8,z^9,\ldots \ran$ with 
$\De^\ast-4=\{0,4,5,\cdots\}$ for the $D$-set $[7]$.
\vskip 0.2cm

{\sf L-functions.}
We define $Z(q,t,a)$ and $L(q,t,a)=(1-t)Z(q,t,a)$ as for
plane curve singularities. To calculate them we need
to trace which $fM$ occur inside $\r$ for $f\in \o$ for standard 
$M\subset \o$. See Table \ref{TableL456}.
Generally, this results in restrictions for the (free) parameters
of the corresponding standard $M\subset \o$ 
until the stable range is reached, which occurs
when $\nu_z(f)$ become sufficiently large.
{\footnotesize
\begin{table*}[ht!]
\[
\centering
\begin{tabular}{|l|l|}
\hline    
$D$-sets, $r\!k_q=1$ & $q,t$-terms  \\
\hline 
$\emptyset$ & $1 - t + q^3 t^4 - q^3 t^7 + q^4 t^8$ \\       
\hline                   
$D$-sets, $r\!k_q=2$ & $q,t$-terms \\
\hline 
$\widetilde{7}$   & $q^2t^3-q^2t^6+q^3t^7$ \\  
$\widetilde{1},7$ & $q^2t^2-q^2t^4+q^4t^6$\\ 
$\widetilde{2},7$ & $q t^2-qt^3+q^2t^4-q^2t^5+q^3t^6$\\
$\widetilde{3},7$ & $qt^3-qt^5+q^2t^6$ \\
\hline    
$D$-sets, $r\!k_q=3$ & $q,t$-terms  \\
\hline  
$\widetilde{1},\widetilde{2},7$ & $t-t^2+q^3t^5$\\
$\widetilde{1},\widetilde{3},7$ & $t^2-t^3+q^2t^5$ \\
$\widetilde{2},\widetilde{3},7$ & $t^3-t^4+qt^5$ \\
\hline    
$D$-sets, $r\!k_q=4$ & $q,t$-terms  \\
\hline 
$\widetilde{1},\widetilde{2},\widetilde{3},7$ & $t^4$ \\  
\hline    
\end{tabular}
\]
\caption{$L$-function for $\F[[z^4,z^5,z^6]]$.}
\label{TableL456}
\end{table*}
}
For instance, let us provide the corresponding calculation 
for $r\!k_q=1$.
The standard {\em invertible} modules are $M=\phi\r\subset \o$ for
$\phi=1+az+bz^2+cz^3+dz^7$ with free $a,b,c,d\in \F$.
 The corresponding $fM\subset \r$
are $\r$, $(z^4+az^5+bz^6+dz^{11})\r$,  
$(z^5+az^6+cz^8+dz^{12})\r$,  $(z^6+bz^8+cz^9+dz^{13})\r$,
and then $z^{8+i}\phi M$ for any $i\ge 0$. Their contributions
to $Z$ and, respectively,  to $L$ are 
$1+q^3(t^4+t^5+t^6)+q^4 t^8/(1-t)$ and 
$1-t+q^3 t^4-q^3 t^7+q^4 t^8$.

Generally, this procedure can be complicated, but such sums for 
$a=0$ and $a=-1/q$ (for $q$ from $L$) 
can be explicitly found; see \cite{Sto}. 
The usage of the motivic measure $X\mapsto |X(\F_q)|$  
makes the calculations relatively elementary.
The other usual choices of motivic measures, are via {\em virtual weight 
polynomials}
and  the {\em Borel-Moore homology}.

Thus, we have $3$ different formulas. Their specializations
at $a=0$ are all different (but have many terms in common):
\begin{align*}
\h^{mot}(q t,t,a\!=\!0)=
&\text{\small
$1 + q t^2 + q^2 t^3 + (q^2+q^3) t^4 + q^3 t^5+ (q^3+q^4)t^6 + q^4 t^8$,} 
\\
L(q,t,a\!=\!0)=
&\text{\small
$1+(q+q^2)t^2+ q^2 t^3+q^3 t^4 +q^3 t^5+(q^3+q^4)t^6+q^4 t^8$,}
\\
\h^{daha}_K(q t, t, a\!=\!0)=
&\text{\small
$1 + q t^2 + q^2 t^3 + q^2 t^4 + q^3 t^5 + q^3 t^6 + q^4 t^8$.}
\end{align*}
Recall that all $3$ coincide for $a=-1/q$. 
In this form, it is easy to see that $q^4 t^8 L(q,1/(qt),a\!=\!0)=
L(q,t,a\!=\!0)$, the functional equation.
 The same holds for  $\h^{daha}_K(q t, t, a\!=\!0)$ and
Z\'u{\oldt{n}}iga's 
$L(q,t,a\!=\!-1/q)$, which is the 
contribution of modules of $r\!k_q=1$ in 
Table \ref{TableL456}. However, such invariance
under $t\to 1/(qt)$ fails for
$\h^{mot}(qt,t,a\!=\!0)$. 

\vskip 0.2cm

{\sf Non-algebraic cables.}
The relation between Gorenstein $\r$ 
and non-algebraic $K$ seems  of sufficiently
general nature, but we mostly deal with examples by now. 
For instance, let
$\r=\F_q[[z^{\upsilon \rr}, z^{\upsilon \ss}, z^{\tt}]]$, where 
$\tt=\upsilon\rr\ss+p$ for $p\in \Z$, 
$\rr,\ss,\tt, \upsilon>0$, gcd$(\rr,\ss)=1$, gcd$(|p|,\nu)=1$,
and the ring is assumed Gorenstein. Then $\Ga=\lan \rr,\ss,\tt\ran$,
$G=\Z_+\setminus \Ga$; compare with (\ref{Ga-cab}).
The corresponding
cable will be $K=C\!ab(\tt,\upsilon)T(\rr,\ss)$. We will use
$q$ from $\h$ in this section.

Given $\rr,\ss$, we expect that there exists $p_0\ge 0$ such that
 for $p\ge p_0$, 
\begin{align}\label{conjKG}
\h^{daha}_K(q,t,a\!=\!-\frac{t}{q})\!=\!L(q/t,t,a\!=\!-\frac{t}{q})\!=\!
(1\!-\!q)(1\!-\!t)\varrho(q,t)+ q^\de
\end{align}
for $\varrho(q,t)\!=\!
\!\sum_{G\ni x>y\in \Ga}q^{g(x)}t^{v(y)-1}=(q t)^{\de-1}
\varrho(t^{-1},q^{-1})$
from Section \ref{sec:rhocab}, where 
$v(x)=|\{\nu\!\in\! \Ga \mid 0\!\le\! \nu\!\le\! x\}|$ and 
$g(x)=|\{g\!\in\! G \mid 0\!\le\! g\! <\! x\}|$. 
Let us provide a couple of examples; we mostly follow \cite{ChW, ChS}.
\vskip 0.2cm

For $\rr=3,\ss=2,\upsilon=2,p=-4$, $K=C\!ab(9,2)T(3,2)$, and
$\h_K^{daha}=$
{\small $1+a^3 q^6+q t+q^2 t+q^3 t+q^2 t^2
+q^3 t^2+2 q^4 t^2+q^3 t^3+q^4 t^3+q^5 t^3
+q^4 t^4+q^5 t^4+q^5 t^5+q^6 t^6
+a^2 (q^3+q^4+q^5+q^4 t+2 q^5 t+q^6 t+q^5 t^2+q^6 t^2+q^6 t^3)
+a (q+q^2+q^3+q^2 t+2 q^3 t+3 q^4 t+q^5 t+q^3 t^2+2 q^4 t^2
+3 q^5 t^2+q^4 t^3+2 q^5 t^3+q^6 t^3+q^5 t^4+q^6 t^4+q^6 t^5).$
}

Then $\varrho=$ {\small $1+q+q^2+q^3+q^4+q^5+q^3 t+q^4 t+q^5 t+q^4 t^2
+q^5 t^2+q^5 t^3+q^5 t^4+q^5 t^5$} in this case, and
 $\h_K^{daha}(q,t,a=-t/q)=$
{\small $1-t+q^3 t-q^3 t^2+q^4 t^2-q^4 t^3+q^5 t^3-q^5 t^6+q^6 t^6$},
which 
coincides with $L(q/t,t,a=-t/q)$. 
\vskip 0.2cm

For $\rr=4,\ss=3,\upsilon=2,p=-13$, the ring $\r=\F_q[[z^6,z^8,z^{11}]]$
is still Gorenstein, $K=C\!ab(11,2)T(4,3)$, and the coincidence holds.
One has:
$\varrho=$ {\small 
$1+q+q^2+q^3+q^4+q^5+q^6+q^7+q^8
+q^9+q^10+q^5 t+q^6 t+q^7 t+q^8 t+q^9 t+q^10 t+q^6 t^2+q^7 t^2+q^8 t^2
+q^9 t^2+q^{10} t^2+q^8 t^3+q^9 t^3+q^{10} t^3+q^8 t^4+q^9 t^4
+q^{10} t^4
+q^9 t^5+q^{10} t^5+q^{10} t^6+q^{10} t^7+q^{10} t^8+q^{10} t^9
+q^{10} t^{10},$
}
and $L(q/t,t,a\!=\!-t/q)=$
{\small $1-t+q^5 t-q^5 t^2+q^6 t^2-q^6 t^3+q^8 t^3-q^8 t^5+q^9 t^5
-q^9 t^6+q^{10} t^6-q^10 t^11+q^{11} t^{11}$.}
\vfil

Let us provide the DAHA superpolynomial for the  $K$ above:\, 
$\h_K^{daha}=$
\renewcommand{\baselinestretch}{1.} 
{\footnotesize
\( 
t^{11} q^{11}-t^5 q^{11}+t^{10} q^{10}+t^9 q^{10}+t^8 q^{10}+
t^7 q^{10}+t^6 q^{10}-2 t^4 q^{10}-t^3 q^{10}+t^9 q^9+t^8 q^9
+2 t^7 q^9+2 t^6 q^9+3 t^5 q^9-3 t^3 q^9-t^2 q^9+t^8 q^8+t^7 q^8
+2 t^6 q^8+3 t^5 q^8+4 t^4 q^8+t^3 q^8-3 t^2 q^8-t q^8+t^7 q^7
+t^6 q^7+2 t^5 q^7+3 t^4 q^7+4 t^3 q^7-2 t q^7+t^6 q^6+t^5 q^6
+2 t^4 q^6+3 t^3 q^6+3 t^2 q^6-q^6+t^5 q^5+t^4 q^5+2 t^3 q^5
+2 t^2 q^5+t q^5+t^4 q^4+t^3 q^4+2 t^2 q^4+t q^4+t^3 q^3+t^2 q^3
+t q^3+t^2 q^2+t q^2+t q+a^5 \bigl(-\frac{q^{14}}{t}
-\frac{q^{13}}{t^2}-\frac{q^{12}}{t^3}\bigr)+a^4 \bigl(-t q^{14}
-q^{14}-t q^{13}-\frac{2 q^{13}}{t}-2 q^{13}-\frac{3 q^{12}}{t}
-\frac{2 q^{12}}{t^2}-2 q^{12}+t q^{11}-\frac{2 q^{11}}{t}
-\frac{2 q^{11}}{t^2}-\frac{q^{11}}{t^3}-\frac{q^{10}}{t^2}
-\frac{q^{10}}{t^3}+q^{10}\bigr)+a^3 \bigl(-t^2 q^{14}-t^3 q^{13}
-2 t^2 q^{13}-3 t q^{13}-q^{13}-3 t^2 q^{12}-5 t q^{12}
-\frac{q^{12}}{t}-5 q^{12}+t^5 q^{11}+t^4 q^{11}+2 t^3 q^{11}
+t^2 q^{11}-5 t q^{11}-\frac{5 q^{11}}{t}-\frac{q^{11}}{t^2}
-7 q^{11}+t^4 q^{10}+2 t^3 q^{10}+4 t^2 q^{10}+2 t q^{10}
-\frac{5 q^{10}}{t}-\frac{3 q^{10}}{t^2}-5 q^{10}+t^3 q^9+2 t^2 q^9
+4 t q^9-\frac{3 q^9}{t}-\frac{2 q^9}{t^2}-\frac{q^9}{t^3}+q^9
+t^2 q^8+2 t q^8-\frac{q^8}{t^2}+2 q^8+t q^7+q^7+q^6\bigr)
+a^2 \bigl(-t^4 q^{13}-t^3 q^{13}-t^2 q^{13}-t^4 q^{12}-4 t^3 q^{12}
-4 t^2 q^{12}-2 t q^{12}+t^8 q^{11}+t^7 q^{11}+2 t^6 q^{11}
+2 t^5 q^{11}+2 t^4 q^{11}-2 t^3 q^{11}-9 t^2 q^{11}-7 t q^{11}
-3 q^{11}+t^7 q^{10}+2 t^6 q^{10}+4 t^5 q^{10}+6 t^4 q^{10}
+6 t^3 q^{10}-2 t^2 q^{10}-11 t q^{10}-\frac{2 q^{10}}{t}-7 q^{10}
+t^6 q^9+2 t^5 q^9+5 t^4 q^9+8 t^3 q^9+9 t^2 q^9-2 t q^9
-\frac{4 q^9}{t}-\frac{q^9}{t^2}-9 q^9+t^5 q^8+2 t^4 q^8+5 t^3 q^8
+8 t^2 q^8+6 t q^8-\frac{4 q^8}{t}-\frac{q^8}{t^2}-2 q^8+t^4 q^7
+2 t^3 q^7+5 t^2 q^7+6 t q^7-\frac{q^7}{t}-\frac{q^7}{t^2}+2 q^7
+t^3 q^6+2 t^2 q^6+4 t q^6+2 q^6+t^2 q^5+2 t q^5+2 q^5+t q^4+q^4
+q^3\bigr)+a \bigl(-t^5 q^{12}-t^4 q^{12}-t^3 q^{12}+t^{10} q^{11}
+t^9 q^{11}+t^8 q^{11}+t^7 q^{11}+t^6 q^{11}-3 t^4 q^{11}-4 t^3 q^{11}
-2 t^2 q^{11}+t^9 q^{10}+2 t^8 q^{10}+3 t^7 q^{10}+4 t^6 q^{10}
+5 t^5 q^{10}+2 t^4 q^{10}-5 t^3 q^{10}-7 t^2 q^{10}-3 t q^{10}
+t^8 q^9+2 t^7 q^9+4 t^6 q^9+6 t^5 q^9+9 t^4 q^9+5 t^3 q^9-6 t^2 q^9
-7 t q^9-2 q^9+t^7 q^8+2 t^6 q^8+4 t^5 q^8+7 t^4 q^8+10 t^3 q^8
+5 t^2 q^8-5 t q^8-\frac{q^8}{t}-4 q^8+t^6 q^7+2 t^5 q^7+4 t^4 q^7
+7 t^3 q^7+9 t^2 q^7+2 t q^7-\frac{q^7}{t}-3 q^7+t^5 q^6+2 t^4 q^6
+4 t^3 q^6+6 t^2 q^6+5 t q^6-\frac{q^6}{t}+t^4 q^5+2 t^3 q^5
+4 t^2 q^5+4 t q^5+q^5+t^3 q^4+2 t^2 q^4+3 t q^4+q^4+t^2 q^3+2 t q^3
+q^3+t q^2+q^2+q\bigr)+1.
\)
}
\renewcommand{\baselinestretch}{1.2}


\vskip 0.2cm
However, there is no relation between $\h_K^{daha}(q,t,a)$
for (non-algebraic) $K=C\!ab(7,2)T(4,3)$ 
and $L(q/t,t,a)$ evaluated at $a=-t/q$  for
the corresponding Gorenstein 
ring, which is $\r=\F_q[[z^6,z^7,z^{8}]]$
for $\rr=4,\ss=3,\upsilon=2, p=-17$. They are 
{\small $1-t+q^5 t+q^6 t-q^7 t-q^5 t^2+q^6 t^2+q^7 t^2-q^8 t^2-
2 q^6 t^3+q^7 t^3+q^8 t^3-q^7 t^4+q^8 t^4-q^8 t^9+q^9 t^9$},
and {\small $1 - t + q^5 t - q^5 t^4 + q^8 t^4 - q^8 t^9 + 
q^9 t^9$} correspondingly. Both satisfy the superduality
$(q t)^\de \h(1/t,1/q,a)=\h(q,t,a)$.
The consideration
of sufficiently 
general Gorenstein rings (if they are non-planar) is needed for 
such a coincidence. The same reservation
may be needed to ensure that  
$\h^{daha}(q,t,a)=L(qt,t,a)$ at $a=-t/q$, which is not clear by now.


\subsection{\bf Concluding remarks}
As we mentioned, there is no
final definition  of $L$-functions of plane curve singularities
by now for general 
 $\si=(c_1,c_2,\cdots,c_\kappa)$. This is needed to have 
the most general version of the conjecture
$\H\!=\!L$, i.e. that for arbitrary $ASF$ in type $A$.
This seems quite doable, but we do not have
sufficient number of examples.
Obtaining explicit formulas for $L(q,t,a)$ for arbitrary $a$ 
is, generally, significantly
more involved than that for $\h^{mot}(q,t,a)$.
The conjecture $\h^{mot}\!=\!\h^{daha}$
is much better tested, including various $\si$.   

\vskip 0.2cm
A natural challenge is to extend motivic superpolynomials
and the Connection
Conjectures from $\{c_i\om_1\}$
to any Young diagrams. The DAHA superpolynomials
are defined in this generality; many examples of $\h^{daha}$ 
are known and posted with sufficiently general colors
(not only rows or columns). 
See e.g., \cite{ChD2}, \cite{ChP2}, \cite{ChW};
supplementary 
materials to the latter paper contain many examples.

\vskip 0.2cm
Another challenge is to prepare 
the next stage of the theory of {\em affine Springer fibers}\,:\,
its extension from $t=1,a=0,\si=\{1,\ldots, 1\}$ to any $t,a,\si$. 
The rationale is that our motivic superpolynomials
are simple to define. This is in type $A$ so far,
but similar theory is expected for {\em suitable} $AFS$ of
other types. Some examples of superpolynomials of classical types
and for  $E$ are discussed in \cite{CJ, ChE, ChS}.

A related problem is to understand the
{\em Fundamental Lemma} in this generality:\, with $q,t,a, \si$. 
 Let me mention in this regard papers
\cite{Cha, Yun}, where certain (related) extensions of orbital integrals
and the Fundamental Lemma are considered.  
\vskip 0.2cm

The theory of motivic  superpolynomials is very rich
even for relatively
simple plane curve singularities. For instance, the case
of colored {\em Hopf
links} triggers a lot of interesting {\em Schubert calculus}.
Generally, the 
combinatorial data of {\em Piontkowski cells} can be viewed as
generalized Young diagrams. Recall that  $\j_0$ belongs  the 
Grassmannian Gr$(2\de,\de)$, where the
subspaces of the middle dimension in
 $\o/(z^{2\de})$ must be invariant with respect
to the matrices corresponding to $x$ and $y$, which provides a link
to the theory of {\em commuting varieties}.
 The connection to Young diagrams is even
more direct for the {\em Gr\"obner cells}; see
\cite{ChG}. Also, the {\em Dyck paths}, 
important for the superpolynomials of torus knots \cite{GM},
are closely related to Young diagrams.
\vskip 0.2cm

Another important aspect concerns 
recurrence formulas for iterated torus links.
Potentially, they can be used to justify ``$\h^{daha}=\h^{mot}$\,".
The class of $\si$-colored algebraic
links in this paper is a natural one closed with respect to such 
formulas. Recall that, classically, the Jones-type invariants of
algebraic knots are calculated via a lot of
intermediate considerations of non-algebraic links, though {\em Soergel
modules} provide now a more``algebraic" approach. 

The Rosso-Jones-type recurrence formulas
require only $\si$-colored {\em algebraic} links. The latter 
occur in such formulas even if we want to calculate the
superpolynomials of uncolored algebraic {\em knots}, but
(at least) the resulting links are {\em algebraic}. 
A natural challenge here is to make the colors
arbitrary Young
diagrams. This was done for DAHA superpolynomials, but
remains an open problem in the motivic theory. 
\vfil

We note that in the DAHA setting, the recurrence formulas
are directly related to {\em Pieri-type
formulas} in general theory of Macdonald polynomials, symmetric
or not, including the stable ones
($J$-polynomials).  See \cite{ChW,ChD2,C101} 
and Ian Macdonald's \cite{Mac1,Mac2,Mac3}. He made many
contributions to mathematics, but developing and promoting
theory of {\em orthogonal polynomials}  is his great
achievement. This required vision and strong character. 
He was creating his theory in the landscape dominated by Lie groups
and algebras, when orthogonal polynomials were not considered
``main stream". 

\vfil
The technique of {\em rank decompositions} of motivic
superpolynomials is an important feature of this paper. 
The main ones for knots are in terms
of \,$\prod_{i\ge c}(1+q^i a)$\, and 
$\prod_{i\ge 0}(1+q^ia/t)$. When $c=1$ (the uncolored case), there
are $2$ more decompositions due to the superduality 
$q\leftrightarrow t^{-1}$:\, via 
$\prod_{i\ge 1}(1+t^{-i} a)$ and $\prod_{i\ge 0}(1+t q^{-i}a)$.
 
There are $3$ more due to the (conjectural)
``$\H=L$". The one corresponding to
the expansion of $L(q,t,a)$ in terms of 
$\prod_{i\ge 1}(1+q^ia)$, a counterpart of the
expansion from \cite{ORS}, becomes the decomposition of $\h^{mot}$
via $\prod_{i\ge 1}(1+q^i t^{-i}a)$, which is self-dual for 
$q\leftrightarrow t^{-1}$. Also, this passage provides
the decompositions
of $\h^{mot}$ in terms of  
$\prod_{i\ge 0}(1+q^i t^{-i-1}a)$, and its dual:\, 
via $\prod_{i\ge 0}(1+q^{i+1} t^{-i}a)$. 

Recall that the decomposition of $\h^{mot}$ in terms of
$\prod_{i\ge 0}(1+q^ia/t)$ requires the usage of $d^\Diamond(M)$,
which is the $q$-rank of 
$M^\Diamond\equal\{f\in \o \mid \mathfrak{m}_\r f\subset M\}$
for standard $M\subset \o$. 
Accordingly, its counterpart for  $L$, used in the last two expansions
above, requires the $q$-rank of 
$M^\Box$ from (\ref{2ndZ-exp}), a counterpart of $M^\Diamond$
adjusted to {\em ideals} $M\subset \r$.

\vskip 0.2cm

All such decompositions are meaningful
geometrically. They are expected
to be associated with {\em differentials} in certain
hyper-homology theories, the ones with 
degenerate {\em spectral sequences}.  
We mostly use the $1${\small st} decomposition above in this
paper, for instance, to
justify the {\em Weak Riemann Hypothesis} for algebraic knots, which 
is upon the specializations $a=0$ and $a=-t/q$ (for $q$ from
$\h^{mot}(q,t,a)$). 
\vskip 0.2cm

Generally, DAHA have a lot of specializations and implications.
Not all of them are of interest for superpolynomials, but some are.
Concerning $\h^{mot}$ and $L$,
the basic specialization is $a=0$, when the $q$-ranks are disregarded. The 
substitution
$a=-t/q$ (in the parameters from $\h^{daha}$ and $\h^{mot}$) is
equally interesting. This is, generally, the case of 
{\em Heegaard-Floer homology} and the
{\em Alexander polynomials} (as $a=-1, t=q$). The specialization
 $t=q$
is the passage to quantum groups and HOMFLY-PT polynomials,
$q=1$ is when the central charge is trivial in DAHA and
the case of ``field with $1$ element" motivically, and so on.
The consideration of roots of unity $q$ and $t=q^{m/n}$ is
important for DAHA and for  superpolynomials. 
\vskip 0.2cm

Recall that the coincidence of  $L(q,t,a=-1/q)$ with 
the Z\'u{\oldt{n}}iga 
$L$-function from \cite{Zu, Sto} (the consideration
of principle ideals only)
is immediate from the decomposition of $L$ via
$\prod_{i\ge 1}(1+q^i a)$ (for $q$ from $L$).
This was used in \cite{ChS} to obtain  
the $q$-split of the
{\em Witten index}.
We expect $\H=L$ at $a=-1/q$ to hold for sufficiently good
Gorenstein rings, not only for
plane curve singularities.  
\vfil 

A related 
application is
to the $\rho_{ab}$-invariant of $S^3\setminus K$
for algebraic knots $K$; superpolynomials provide its certain
$q,t,a$-deformation.  
This link may be not very 
surprising because the Alexander polynomials are the key in the theory
of  $\rho_{ab}$ via the {\em Levine-Tristram signature} $\si_K$.
However, this is an interesting development:\,
 $\rho$  belongs to a different
branches of geometry-topology, the theory of spectral invariants.

Generally, we think that a highly ramified theory of 
Khovanov-Rozansky polynomials, physics superpolynomials,
DAHA superpolynomials, motivic superpolynomials and
$L$-functions of algebraic links reaches now
new levels of uniformity. This includes
various related aspects in low-dimensional
topology, enumeration geometry and arithmetic
geometry. We hope that our paper contributes to this process.
 
These are quite different directions, with strong classical
origins and non-trivial overlapping. 
For instance,  definitions of superpolynomials, DAHA and motivic,
are simpler vs. those in Khovanov-
Rozansky theory, but we can do only iterated torus links and
related plane curve singularities. 
As we try to demonstrate in
this paper, the motivic theory can be
well developed almost from scratch,
though many problems certainly require advanced
methods. 

Beyond what we discuss in this paper,  
DAHA and motivic approaches can be presumably extended to
construct invariants of isolated {\em surface singularities}
associated with Seifert $3$-folds and other plumbed manifolds.
This is quite a challenge in low-dimensional topology-geometry
and physics. A related theme concerns possible applications
of the corresponding {\em superseries} in number theory;
see \cite{ChW} and \cite{ChS}.

\vskip 0.2cm
{\bf Acknowledgements.}
I thank Ian Philipp, David Kazhdan, and the referee for  
comments. Also, 
I would like to thank Tomoyuki Arakawa and RIMS 
for the kind invitation (it was a productive summer).


\vskip -2cm
\bibliographystyle{unsrt}

\medskip
\end{document}